\documentclass[a4paper,12pt]{amsart}
\usepackage{epsf,amssymb,theo,xy,french}
\title[G\'eom\'etrie d'Arakelov des vari\'et\'es toriques]
{G\'eom\'etrie d'Arakelov des vari\'et\'es toriques \\ et fibr\'es
en droites int\'egrables}
\author{Vincent Maillot}
\def\M#1{\mathbb#1}	

\def\C#1{\mathcal#1}
\def\E#1{\scr#1}	
\def\R{\M{R}}
\def\Z{\M{Z}}
\def\N{\M{N}}
\def\Q{\M{Q}}

\def\demo{\noindent {\bf D\'emonstration.}}
\def\hooklongrightarrow{\lhook\joinrel\longrightarrow}

\def\vfi{\varphi}
\def\epsi{\varepsilon}
\def\si{\sigma}
\def\t{\tau}
\def\P{\M{P}(\Delta)}
\def\PN{\M{P}(\nabla)}
\def\co{c\^one}
\def\cos{c\^ones}
\def\D{\Delta}
\def\Na{\nabla}
\def\op#1{\operatorname{#1}}
\def\SP{\op{Spec}\M{Z}}
\def\ssi{si et seulement si}
\def\PP{\M{P}(\D)_{\M{C}}}
\def\PPN{\M{P}(\nabla)(\M{C})}
\def\PPP{\P_{\geqslant}}
\def\TT{T(\M{C})}
\def\TTT{T_{\geqslant}}
\def\CT{\C{S}_{N}}
\def\ov#1{\overline{#1}}

\def\psh{plu\-ri\-sous\-har\-mo\-ni\-que}
\def\pshs{plu\-ri\-sous\-har\-mo\-ni\-ques}
\def\dd{dd^{c}}
\def\AA{\ov{A}}
\def\AAA{\ov{\ov{A}}}
\def\XR{X_{\R}}
\def\XC{X(\M{C})}
\def\arg{\op{Arg}}
\def\CH{\widehat{CH}_{\op{int}}}
\def\ZZ{\widehat{Z}_{\op{gen}}}
\def\odiv{\op{div}}
\DeclareMathOperator*{\Sup}{Sup}

\DeclareFontFamily{OT1}{rsfs}{}
\DeclareFontShape{OT1}{rsfs}{m}{n}{%
   <5>rsfs5%
   <6>rsfs7%
   <7>rsfs7%
   <8>rsfs7%
   <9>rsfs7%
   <10>rsfs10%
   <10.95>rsfs10%
   <12>rsfs10%
   <14.4>rsfs10%
   <17.28>rsfs10%
   <20.74>rsfs10%
   <24.88>rsfs10}{}
\DeclareMathAlphabet{\scc}{OT1}{rsfs}{m}{n}
\def\scr#1{{\ifmmode\scc{#1}\else{\fontfamily{rsfs}\selectfont{#1}}\fi}}
	
\xyoption{all}

\begin{document}
\maketitle
\bigskip~
\bigskip
\thispagestyle{empty}
\tableofcontents
\newpage~ 
\section{Introduction}~

Le but de cet article est l'\'etude, dans le cadre de la g\'eom\'etrie
d'Arakelov, des vari\'et\'es toriques projectives et lisses. Chemin faisant,
nous \'etendons sur certains points la th\'eorie d\'evelopp\'ee dans \cite{13}.
\medskip

Gr\^ace \`a Demazure (cf. \cite{8}), on sait associer \`a tout \'eventail $\D$ de
$\M{Z}^{d}$ un sch\'ema $\pi: \M{P}(\D) \rightarrow \op{Spec}\M{Z}$ que l'on
appelle {\it vari\'et\'e torique\/} associ\'ee \`a $\D$.
On trouve dans la litt\'erature (cf. \cite{4}, \cite{20} et \cite{11} pour des
r\'ef\'erences pr\'ecises) une description explicite de la vari\'et\'e complexe
$\M{P}(\D)_{\M{C}}$ en terme des propri\'et\'es combinatoires de $\D$. 
Revenant au point de vue originel de Demazure, nous montrons que cette
description s'\'etend sans difficult\'e \`a la situation sur $\M{Z}$. En
particulier, on dispose lorsque $\M{P}(\D)$ est projective et lisse, d'une
description agr\'eable bas\'ee sur un th\'eor\`eme de Jurkiewicz et Danilov, de
l'anneau de Chow $CH^{\ast}(\M{P}(\D))$ en terme de g\'en\'erateurs et
relations (th\'eor\`eme \ref{anneau_chow}). Plus pr\'ecis\'ement, nous montrons que
$CH^{\ast}(\M{P}(\D))$ est engendr\'e en tant qu'anneau par la premi\`ere
classe de Chern $c_{1}(L) \in CH^{1}(\M{P}(\D))$ des fibr\'es en droites $L$
sur $\M{P}(\D)$.
\medskip

Afin de donner une description analogue de l'anneau de Chow arithm\'etique 
$\widehat{CH}^{\ast}(\M{P}(\D))$, il nous est n\'ecessaire de munir tout 
fibr\'e en droites $L$ sur $\M{P}(\D)$ d'une m\'etrique $\|.\|_{L,\infty}$ 
``canonique'' permettant entre autre chose le calcul explicite du
produit~:
\[
\hat{c}_{1}(L,\|.\|_{L,\infty})^{p} \in \widehat{CH}^{p}(\M{P}(\D)), 
\]
o\`u $\hat{c}_{1}(L,\|.\|_{L,\infty})$ d\'esigne la premi\`ere classe de Chern
arithm\'etique de $(L,\|.\|_{L,\infty})$ (comparer avec \cite{33} o\`u un point de vue
analogue est d\'evelopp\'e pour les grassmanniennes).

Nous donnons plusieurs constructions d'une m\'etrique $\|.\|_{L,\infty}$
qui est canonique dans le sens o\`u
l'application $L \mapsto \|.\|_{L,\infty}$ poss\`ede des 
propri\'et\'es fonctorielles (proposition \ref{BT_fonct}) qui permettent de la 
caract\'eriser enti\`erement.
Malheureusement, la m\'etrique ainsi
construite n'est pas $C^{\infty}$ en g\'en\'eral. La premi\`ere ``forme'' de
Chern $c_{1}(L,\|.\|_{L,\infty})$ est un courant r\'eel de bidegr\'e (1,1), et
donc le produit $c_{1}(L,\|.\|_{L,\infty})^{p}$, et {\it a fortiori\/}
le produit $\hat{c}_{1}(L,\|.\|_{L,\infty})^{p}$, ne sont pas d\'efinis.

Nous sommes donc amen\'es dans un premier temps \`a \'etendre sur certains
points la th\'eorie d\'evelopp\'ee dans \cite{13}, afin de pouvoir consid\'erer
en g\'eom\'etrie d'Arakelov
des fibr\'es en droites munis de m\'etriques non-n\'ecessairement $C^{\infty}$. 
Comme, pour autant qu'il nous
soit permis d'en juger, ces d\'eveloppements poss\`edent un int\'er\^et propre,
nous avons choisi d'en donner une exposition valable en toute g\'en\'eralit\'e. 
\medskip

Quittons donc pour quelques temps l'univers torique, et consid\'erons $X$ une
vari\'et\'e arithm\'etique quelconque de dimension absolue $d+1$.
Reprenant une terminologie introduite par Zhang \cite{21}, nous dirons
d'un couple $(L,\|.\|)$ form\'e d'un fibr\'e en droites $L$ sur $X$ et d'une
m\'etrique hermitienne continue $\|.\|$ sur $L(\M{C})$ qu'il est 
{\it admissible\/} si $L$ est engendr\'e par ses sections globales, $\|.\|$ est
positive et peut \^etre approch\'ee uniform\'ement sur $X(\M{C})$ 
par des m\'etriques positives $C^{\infty}$.
Plus g\'en\'eralement, un fibr\'e en droites sur $X$ sera dit {\it int\'egrable\/}
s'il est diff\'erence de deux fibr\'es en droites admissibles.
Tout fibr\'e en droites $L$ sur $\M{P}(\D)$ muni de sa m\'etrique canonique
$\|.\|_{L,\infty}$ est int\'egrable (exemple \ref{exemple_decomp1}).
\medskip

Nous d\'eveloppons alors sur $X(\M{C})$ un formalisme de {\it formes
diff\'erentielles g\'en\'eralis\'ees\/}; 
en particulier la premi\`ere ``forme'' de Chern $c_{1}(\ov{L})$ d'un fibr\'e en
droites int\'egrable $\ov{L}$ est une forme diff\'erentielle g\'en\'eralis\'ee
\`a notre sens. En nous appuyant sur une th\'eorie d\'evelopp\'ee par
Bedford-Taylor \cite{1} puis Demailly \cite{6}, nous montrons comment l'on peut
donner un sens au produit de deux telles formes (proposition \ref{produit_generalise}).

Nous construisons ensuite
un groupe gradu\'e $\widehat{CH}^{\ast}_{\op{int}}(X)$
contenant l'anneau de Chow arithm\'etique usuel $\widehat{CH}^{\ast}(X)$, de
telle sorte que pour tout fibr\'e en droites int\'egrable $\ov{L}$ sur $X$, on
ait $\hat{c}_{1}(\ov{L}) \in \widehat{CH}^{1}_{\op{int}}(X)$.
Nous \'etendons alors partiellement le formalisme d\'evelopp\'e 
dans \cite{13} au groupe $\widehat{CH}^{\ast}_{\op{int}}(X)$, ce qui nous
permet de retrouver certains r\'esultats de Zhang \cite{21} concernant les fibr\'es int\'egrables.
Notre principal
r\'esultat dans cette direction est le suivant (th\'eor\`eme \ref{produit_bien_defini}). Il
existe un accouplement~:
\[
\widehat{CH}^{p}_{\op{int}}(X)\otimes \widehat{CH}^{q}_{\op{int}}(X) 
\longrightarrow \widehat{CH}^{p+q}_{\op{int}}(X)_{\M{Q}}, 
\]
qui prolonge celui d\'efini par Gillet-Soul\'e sur $\widehat{CH}^{\ast}(X)$ et
qui munit $\widehat{CH}^{\ast}_{\op{int}}(X)$ d'une structure d'anneau commutatif,
associatif et unif\`ere. De plus, on dispose d'un morphisme
$\op{\widehat{\op{deg}}}: \widehat{CH}_{\op{int}}^{d+1}(X) \rightarrow \M{R}$
qui prolonge celui d\'efini dans \cite{13}, et tel que (th\'eor\`eme \ref{gdthm}) si
$\ov{L}_{1},\dots,\ov{L}_{d+1}$ sont des fibr\'es int\'egrables sur $X$ et
$h_{\ov{L}_{1}, \dots,\ov{L}_{d+1}}(X)$ est la hauteur de $X$ relativement \`a
$\ov{L}_{1},\dots,\ov{L}_{d+1}$ telle qu'elle est d\'efinie dans \cite{21}, alors~:
\[
h_{\ov{L}_{1},\dots,\ov{L}_{d+1}}(X) =
\widehat{\op{deg}}(\hat{c}_{1}(\ov{L}_{1}) \dotsm \hat{c}_{1}(\ov{L}_{d+1})).
\]

Revenons au cas particulier des
vari\'et\'es toriques~: Soient $\ov{L}_{1},\dots,\ov{L}_{q}$ des
fibr\'es en droites sur $\M{P}(\D)$ munis de leur m\'etrique canonique; nous
donnons une formule explicite pour le produit g\'en\'eralis\'e
$c_{1}(\ov{L}_{1})\dotsm c_{1}(\ov{L}_{q})$ en terme de la combinatoire de $\D$
(th\'eor\`eme \ref{calcul_prod}). Nous montrons que cette formule conduit \`a un
algorithme particuli\`erement simple pour le calcul du volume mixte de
polytopes convexes (remarque \ref{algo_efficace}).

Lorsque $q = d+1$, nous montrons que la hauteur canonique 
$h_{\ov{L}_{1}, \dots,
\ov{L}_{d+1}}(\M{P}(\D))$ est nulle (proposition \ref{annulation_hauteur}),
puis nous approfondissons ce r\'esultat en montrant 
(th\'eor\`eme \ref{section_chow}) que le sous-anneau de 
$\widehat{CH}^{\ast}_{\op{int}}(\M{P}(\D))$ engendr\'e par 
les premi\`eres classes de Chern arithm\'etiques
$\hat{c}_{1}(L,\|.\|_{L,\infty})$ des fibr\'es en droites $L$ sur $\M{P}(\D)$
est isomorphe canoniquement \`a $CH^{\ast}(\M{P}(\D))$.
Ce r\'esultat est
l'analogue arithm\'etique du th\'eor\`eme de Jurkiewicz et Danilov. 
Nous en d\'eduisons
que la hauteur canonique d'une hypersurface dans
$\M{P}(\D)$ (i.e. sa hauteur relativement \`a des fibr\'es en droites sur $\M{P}(\D)$
munis de leur m\'etrique canonique)
est donn\'ee essentiellement par la mesure de Mahler du polyn\^ome
qui la d\'efinit (proposition \ref{hauteur_hypersurfaces}). 

En nous appuyant sur les r\'esultats pr\'ec\'edents, nous concluons cet article par la
d\'emonstration d'un analogue arithm\'etique du th\'eor\`eme de
Bernstein-Koushniren\-ko (th\'eor\`eme \ref{B_K} et corollaire \ref{coro_BK})
donnant une majoration de la hauteur des points d'intersections de $d$
hypersurfaces de $\M{P}(\D)$ en fonction de leur hauteur canonique et de leur
polyh\`edre de Newton.

\medskip

Passons maintenant en revue l'organisation de cet article.
\medskip

Le chapitre 2 est consacr\'e aux vari\'et\'es toriques sur $\op{Spec}\M{Z}$.
Nous rappelons au 2.1 quelques propri\'et\'es simples des c\^ones et des
\'eventails, et au 2.2 la construction d'apr\`es Demazure des vari\'et\'es
toriques et de leurs morphismes canoniques. Nous introduisons au 2.3 la notion
de diviseur invariant et de fonction support associ\'ee. Nous rappelons en les
adaptant \`a la situation sur $\M{Z}$ certains crit\`eres pour l'amplitude d'un
diviseur invariant. Au 2.4 nous montrons comment l'on peut construire de
mani\`ere canonique une vari\'et\'e torique projective \`a partir d'un polytope
convexe entier. Enfin le 2.5 est consacr\'e \`a la structure de l'anneau de
Chow d'une vari\'et\'e torique projective lisse. En utilisant une
d\'ecomposition cellulaire sur $\M{Z}$ d'une telle vari\'et\'e, nous montrons
comment la plupart des r\'esultats classiques (sur $\M{C}$) s'\'etendent \`a la
situation sur $\M{Z}$; nous en profitons pour prouver un th\'eor\`eme
d'annulation des groupes $CH^{p,p-1}(\M{P}(\D))$. C'est l\`a le seul point
vraiment nouveau de ce chapitre.
\medskip

Au chapitre 3 nous nous int\'eressons plus particuli\`erement aux vari\'et\'es
toriques complexes. La vari\'et\'e torique \`a coin $\M{P}(\D)_{\geqslant}$ 
est d\'efinie au 3.1. Nous pr\'esentons au 3.2 un recouvrement canonique de
$\M{P}(\D)(\M{C})$ introduit par Batyrev et Tschinkel \cite{2} puis \'etudions
certaines de ses propri\'et\'es. Nous utilisons ce recouvrement pour construire
au 3.3 la m\'etrique canonique $\|.\|_{L,\infty}$. Nous montrons que cette
m\'etrique v\'erifie certaines propri\'et\'es de multiplicativit\'e et de
fonctorialit\'e (proposition \ref{BT_fonct}). Nous en donnons ensuite deux autres
constructions, l'une bas\'ee sur un th\'eor\`eme de Zhang (th\'eor\`eme
\ref{zhang}), l'autre par image inverse (proposition \ref{decomposition}). Nous d\'eduisons de
cette derni\`ere construction un th\'eor\`eme d'approximation globale pour
$\|.\|_{L,\infty}$ (proposition \ref{approximation}).
\medskip

Le chapitre 4 est centr\'e sur l'\'etude des formes diff\'erentielles
g\'en\'eralis\'ees sur une vari\'et\'e complexe arbitraire. Apr\`es avoir
rappel\'e au 4.2 la th\'eorie d\'evelopp\'ee par Bedford-Taylor \cite{1} puis
Demailly \cite{6}, nous d\'efinissons au 4.3 et au 4.4 plusieurs types
remarquables de courants qui, gr\^ace \`a cette th\'eorie, peuvent \^etre
multipli\'es. Nous introduisons aux sections 4.5 et 4.7 les notions de fibr\'es
en droites admissibles et de fibr\'es en droites int\'egrables, et montrons au
4.6 que tout fibr\'e en droites ample muni d'une m\'etrique positive sur une
vari\'et\'e projective est admissible.
\medskip

Nous abordons au chapitre 5 le c\oe ur de notre sujet, \`a savoir la
construction de l'anneau $\widehat{CH}_{\op{int}}^{\ast}(X)$ pour toute
vari\'et\'e arithm\'etique $X$. Apr\`es avoir rappel\'e au 5.1 la th\'eorie
classique de Gillet-Soul\'e et introduit au 5.2 le groupe de Picard
arithm\'etique g\'en\'eralis\'e $\widehat{\op{Pic}}_{\op{int}}(X)$, nous
d\'efinissons au 5.3 les groupes de Chow arithm\'etiques g\'en\'eralis\'es
$\widehat{CH}^{p}_{\op{int}}(X)$ (d\'efinition \ref{def_gch_gen}) puis nous construisons
au 5.4 l'accouplement $\widehat{CH}^{p}_{\op{int}}(X) \otimes
\widehat{CH}^{q}_{\op{int}}(X) \rightarrow
\widehat{CH}^{p+q}_{\op{int}}(X)_{\M{Q}}$.
Finalement nous relions au 5.5 cette construction \`a celle de Zhang \cite{21}
(th\'eor\`eme \ref{gdthm}), puis nous d\'emontrons un r\'esultat g\'en\'eral de
positivit\'e (proposition \ref{positivite_arithmetique}) qui nous servira au chapitre 8 lors de la
d\'emonstration d'un analogue arithm\'etique du th\'eor\`eme de
Bernstein-Koushnirenko.
\medskip

Nous retournons \`a partir du chapitre 6 aux vari\'et\'es toriques. 
Apr\`es avoir donn\'e au 6.2 une expression explicite, pour tout fibr\'e en
droites $\ov{L}$ sur $\M{P}(\D)$ muni de sa m\'etrique canonique, du courant
$c_{1}(\ov{L})$, nous g\'en\'eralisons ce r\'esultat en donnant au 6.3 une
formule explicite pour le produit g\'en\'eralis\'e de tels courants de Chern
(th\'eor\`eme \ref{calcul_prod}). Nous d\'emontrons enfin au 6.4 un r\'esultat
d'annulation pour un tel produit (corollaire \ref{trivialite2}).
\medskip

Le chapitre 7 est consacr\'e \`a l'\'etude de la g\'eom\'etrie d'Arakelov des
vari\'et\'es toriques. Nous d\'emontrons au 7.1 un th\'eor\`eme d'annulation
des multihauteurs canoniques de $\M{P}(\D)$ (proposition \ref{annulation_hauteur}) et
\'etablissons au 7.2 une formule reliant la hauteur canonique d'une
hypersurface dans $\M{P}(\D)$ \`a la hauteur de Mahler du polyn\^ome qui la
d\'efinit (proposition \ref{hauteur_hypersurfaces}). Le 7.3 pr\'esente un exemple
int\'eressant~: le
calcul de la mesure de Mahler d'une droite dans $\M{P}^{2}$. Nous exhibons au 7.4 une
section canonique du morphisme d'anneaux application cycle $\zeta:
\widehat{CH}^{\ast}_{\op{int}}(\M{P}(\D)) \rightarrow CH^{\ast}(\M{P}(\D))$.
\medskip

Au 8.1, nous associons de mani\`ere fonctorielle une constante r\'eelle
positive $L(\nabla)$ \`a tout polytope convexe entier $\nabla$ de $\M{R}^{d}$,
puis nous en donnons une majoration explicite dans le cas o\`u $\nabla$ est
absolument simple (proposition \ref{majoration_cas_AS}). Finalement nous d\'emontrons par
r\'ecurrence au 8.2 un analogue arithm\'etique du th\'eor\`eme de
Bernstein-Koushnirenko faisant intervenir dans son \'enonc\'e les constantes
$L(\nabla)$ pr\'ec\'edemment introduites (th\'eor\`eme \ref{B_K} et corollaire
\ref{coro_BK}).
\medskip

Une partie des r\'esultats pr\'esent\'es ici avaient \'et\'e annonc\'es dans
\cite{34}.
\medskip

L'auteur tient \`a exprimer sa profonde gratitude \`a J.-B. Bost sans l'aide
duquel ce travail n'aurait jamais vu le jour. Il lui est \'egalement agr\'eable
de remercier 
J. Cassaigne, A. Chambert-Loir, D. Harari, M.
Laurent, E. Leichtnam, J. Pfeiffer, C. Soul\'e et B. Teissier pour des
discussions int\'eressantes.
\bigskip

\newpage~
\section{Vari\'et\'es toriques sur $\SP$}~

\subsection{C\^ones et \'eventails}~

Pour plus de d\'etails sur les d\'efinitions et les d\'emonstrations des
propositions \'enonc\'ees ici, on peut consulter (\cite{20}, \S 1.1) et aussi 
(\cite{11}, \S 1.1 et 1.2).

Soit $N \simeq \M{Z}^{d}\label{lab1}$ un $\M{Z}$-module libre de rang $d$ dans lequel on a
choisi une base $e_{1}, \dots, e_{d}$. On note $M =
\operatorname{Hom}_{\M{Z}}(N,\M{Z})\label{lab2}$ son $\M{Z}$-module dual. On obtient un
accouplement non d\'eg\'en\'er\'e~:
\[
\label{lab3}
<~\, ,~>\;: M\times N \longrightarrow \M{Z}.
\]
On pose $N_{\M{R}} = N \otimes_{\M{Z}} \M{R}$ et $M_{\M{R}} =
M\otimes_{\M{Z}}\M{R}$. Ce sont des $\M{R}$-espaces vectoriels de dimension
$d$. L'accouplement ci-dessus s'\'etend en une forme bilin\'eaire
non-d\'eg\'en\'er\'ee~:
\[
<~\, ,~>\;: M_{\M{R}}\times N_{\M{R}} \longrightarrow \M{R}.
\]
\begin{defn}
On appelle {\it c\^one polyh\'edral rationnel dans $N$} ou plus simplement {\it
c\^one} tout ensemble $\sigma \subseteq N_{\M{R}}$ de la forme~:
\[
\sigma = \sum_{i \in I}\M{R}^{+}n_{i}, 
\]
o\`u $(n_{i})_{i \in I}$ est une famille finie d'\'el\'ements de $N$. La {\it
dimension} de $\sigma$ est d\'efinie comme la dimension de l'espace vectoriel
r\'eel engendr\'e par les points du c\^one $\sigma$~:
\[
\dim \sigma = \dim_{\M{R}}(\operatorname{Vect}(\sigma)) = \dim_{\M{R}}(\sigma +
(- \sigma)).
\]
\end{defn}
\begin{rem}
On d\'efinit de la m\^eme fa\c con les {\it c\^ones polyh\'edraux rationnels
dans $M$}.
\end{rem}
On d\'efinit le {\it dual} $\sigma^{\ast}$ (resp. {\it l'orthogonal}
$\sigma^{\perp}$) du c\^one $\sigma$ de la fa\c con suivante~:
\begin{alignat*}{3}
&\sigma^{\ast}& &=& &\{v \in M_{\M{R}}\, :\; <v,x> \;\geqslant 0, \quad \forall x \in
\sigma\} \subseteq M_{\M{R}}, \\
&\sigma^{\perp}& &=& &\{v \in M_{\M{R}}\, :\; <v,x> \;= 0, \quad \forall x \in \sigma\}
\subseteq M_{\M{R}}.
\end{alignat*}
\begin{defn}
On dit que $\tau \subseteq \sigma$ est une {\it face} de $\sigma$, et on note
$\tau < \sigma$, si l'on peut trouver $v \in \sigma^{\ast}$ tel que~:
\[
\tau = \sigma \cap \{v\}^{\perp}
\]
\end{defn}
\begin{rem}
Dans un tel cas, on peut toujours choisir $v \in \sigma^{\ast}\cap M$ (voir par
exemple \cite{20}, prop. 1.3 ou \cite{11}, p. 13).
\end{rem}
\begin{defn}
Un c\^one $\sigma$ est dit {\it strict} s'il ne contient aucune droite
r\'eelle.
\end{defn}
Les assertions suivantes d\'ecoulent ais\'ement des d\'efinitions~:
\begin{prop}~
\begin{itemize}
\item Toute face d'un c\^one est un c\^one.
\item Le dual d'un c\^one $\sigma$ est un c\^one, et de plus
${(\sigma^{\ast})}^{\ast} = \sigma$.
\item Pour tout c\^one strict $\sigma$, $\dim \sigma^{\ast} = n$.
\item On a $\dim \sigma + \dim \sigma^{\perp} = n$ pour tout c\^one $\sigma$.
\end{itemize}
\end{prop}
\begin{defn}
Un {\it \'eventail} de $N_{\M{R}}$ est une famille finie $\Delta = \{\sigma\}$ de
c\^ones stricts de $N_{\M{R}}$ tels que~:
\begin{itemize}
\item Si $\sigma \in \Delta$, alors toute face $\tau$ de $\sigma$ appartient
\`a $\Delta$.
\item Si $\sigma, \sigma' \in \Delta$, alors $\sigma \cap \sigma'$ est une face
\`a la fois de $\sigma$ et de $\sigma'$.
\end{itemize}
\end{defn}
La r\'eunion $|\Delta| = \bigcup_{\sigma \in \Delta}\sigma$ est appel\'ee
{\it support} de $\Delta$. On note~:
\[
\Delta(j) = \{\sigma\}_{
\substack{
\sigma \in \Delta \\ \dim \sigma = j}}
\]
le {\it $j$-squelette} de
$\Delta$. On notera \'egalement $\Delta_{\rm{max}} = \Delta(d)$.
\begin{defn}
Soit $\sigma$ un c\^one strict de $N_{\M{R}}$, on note $\E{S}_{\sigma} = M \cap
\sigma^{\ast}$.
\end{defn}
La proposition suivante donne une caract\'erisation alg\'ebrique des ensembles
$\E{S}_{\sigma}$.
\begin{prop}
Soit $\sigma$ un c\^one strict de $N_{\M{R}}$. L'ensemble $\E{S}_{\sigma}$
est un semi-groupe (additif), satur\'e, de type fini. De plus, $\E{S}_{\sigma}$
engendre $M$ en tant que groupe. R\'eciproquement, si $\E{S} \subset M$ est un
tel semi-groupe, alors il existe $\sigma$ un c\^one strict de $N_{\M{R}}$ tel
que $\E{S} = \E{S}_{\sigma}$.
\end{prop}
\medskip
\begin{center}
\fbox{
\begin{minipage}{12cm}
Dans toute la suite du texte, $N$ est un $\M{Z}$-module libre de rang $d$
fix\'e une fois pour toute.
\end{minipage}
}
\end{center}
\bigskip
\medskip

\subsection{Construction des vari\'et\'es toriques sur Spec $\M{Z}$}~

On suit ici essentiellement (\cite{8}, \S 4). Demazure s'int\'eresse aux
vari\'et\'es toriques lisses, mais les d\'emonstrations des propositions
donn\'ees ici s'\'etendent imm\'ediatement au cas g\'en\'eral. On peut
\'egalement consulter (\cite{20}, \S 1) et (\cite{11}, \S 1 et \S 2) pour les d\'emonstrations
de ces propositions dans le cas o\`u le corps de base est $\M{C}$.
\begin{defn} (Vari\'et\'es toriques affines associ\'ee \`a un c\^one $\sigma$).
Soit $\sigma$ un c\^one strict de $N_{\M{R}}$, on note $U_{\sigma}$ le
$\M{Z}$-sch\'ema~:
\[
U_{\sigma} = \operatorname{Spec}\left(\M{Z}\left[\E{S}_{\sigma}\right]\right)
\]
\end{defn}
\begin{prop}
Le $\M{Z}$-sch\'ema $U_{\sigma}$ est affine, normal, plat sur
$\operatorname{Spec}\M{Z}$, \`a fibres g\'eom\'etriquement int\`egres.
\end{prop}
\begin{prop}
Soit $\tau < \sigma$, l'inclusion $\E{S}_{\sigma} \subseteq \E{S}_{\tau}$
induit un morphisme canonique $U_{\tau} \rightarrow U_{\sigma}$; c'est une
immersion ouverte de $\M{Z}$-sch\'emas.
\end{prop}
\demo\ Voir (\cite{8}, \S 4, lemme 1), et aussi (\cite{20}, th.
1.4).
\begin{expl}
Prenons $\si = \{0\}$, on obtient $U_{\{0\} } = T$, le $\M{Z}$-tore dual de
$N$.
\end{expl}
Soient $\si$ et $\si'$ deux \cos\ d'un m\^eme \'eventail $\Delta$. On a le
diagramme suivant~:
\begin{center}
\mbox{
\xymatrix{
& U_{\si} \\
U_{\si\cap\si'} \ar@{^{(}->}[ur] \ar@{^{(}->}[dr]& \\
& U_{\si'}
}}
\end{center}
La proposition pr\'ec\'edente permet de recoller $U_{\si}$ et $U_{\si'}$ le long
de $U_{\si \cap \si'}$. Plus g\'en\'eralement on pose la d\'efinition
suivante~:
\begin{defn}
Soit $\D$ un \'eventail. On appelle {\it vari\'et\'e torique associ\'ee \`a
l'\'e\-ventail $\D$\/} et on note $\P$ le sch\'ema obtenu par recollement des
$U_{\si}$, $\si$ parcourant $\D$, \`a l'aide des immersions ouvertes 
$U_{\si\cap\si'} \hookrightarrow U_{\si}$ et $U_{\si\cap\si'} \hookrightarrow
U_{\si'}$, pour $\si$, $\si' \in \D$.
\end{defn}
\begin{prop}
Le $\M{Z}$-sch\'ema $\P$ est plat sur $\SP$, normal, s\'epar\'e, int\`egre, 
de dimension absolue $d+1$,
\`a fibres g\'eom\'etriquement int\`egres.
\end{prop}
\demo\ Voir (\cite{8}, \S 4, prop. 1), et aussi (\cite{20}, th. 1.4); on 
peut \'egalement consulter (\cite{11}, \S 1.4) pour certains aspects.
\medskip

La proposition suivante justifie le nom de {\it vari\'et\'e torique} pour un
$\M{Z}$-sch\'ema de la forme $\P$~:
\begin{prop}
Notons $T$ le $\M{Z}$-tore dual de $N$; l'action~:
\[
T \otimes_{\op{Spec}\M{Z}} U_{\{0\}} \simeq T \otimes_{\op{Spec}\M{Z}} T \longrightarrow T
\]
de $T$ sur $U_{\{0\}}$, d\'efinie par la structure de sch\'ema en groupe de
$T$, se prolonge en une action~:
\[
T \otimes_{\op{Spec}\M{Z}} \P \longrightarrow \P
\]
de $T$ sur $\P$.
\end{prop}
\demo\ Voir (\cite{8}, p. 559). On peut \'egalement consulter (\cite{20}, p.
9).
\medskip

Comme $T = U_{\{0\}}$ est un ouvert dense de $\P$, tout mon\^ome $x^{m}$ 
pour $m \in M$  s'\'etend en une fonction rationnelle sur $\P$. On note
$\chi^{m}$ cette fonction que l'on appelle {\it caract\`ere} associ\'e \`a $m$.
La proposition suivante est une cons\'equence
imm\'ediate de la d\'efinition de $\P$~:
\begin{prop}
Soit $\si \in \D$ et $m \in M \cap \si^{\ast} = \E{S}_{\si}$, alors $\chi^{m}$
est r\'eguli\`ere sur $U_{\si}$. De plus si $m'$, $m'' \in \E{S}_{\si}$ alors
$\chi^{m' + m''} = \chi^{m'}\chi^{m''}$ sur $U_{\si}$.
\end{prop}
\begin{defn}
On notera $\epsilon : \SP \rightarrow \P$ la section nulle de $T$ vue comme
section de $\P$~:
\begin{center}
\mbox{
\xymatrix{
\qquad \quad \; T \subseteq \P \\
\SP \ar@/^/[u]_{\epsilon}}
}
\end{center}
\end{defn}
Les deux propositions suivantes donnent des crit\`eres simples sur $\D$ pour
que $\P$ soit propre (resp. lisse)~:
\begin{prop}
La vari\'et\'e torique $\P$ est propre sur $\SP$ si et
seulement si $\D$ est {\it complet} (i.e. $|\D| = N_{\M{R}}$).
\end{prop}
\begin{prop}
\label{lissite}
La vari\'et\'e torique $\P$ est lisse sur $\SP$ si et seulement si tout \co\
$\si \in \D$ est engendr\'e par une partie d'une base de $N$. Dans ce cas, $\D$
est dit {\it r\'egulier}.
\end{prop}
\demo\ Voir (\cite{8}, \S 4, prop. 4) et (\cite{8}, \S 4, def. 1 et prop. 1).
On peut \'egalement consulter (\cite{20}, th. 1.10) et (\cite{11}, \S 2.1 et
2.4).
\begin{rem}
Soit $k$ un corps quelconque et notons $\M{P}(\D)_{k} = \P
\otimes_{\SP}\op{Spec}k$ la
vari\'et\'e torique sur $k$ associ\'ee \`a $\D$ comme dans \cite{4}.
Les trois assertions
suivantes sont \'equivalentes~:
\begin{itemize}
\item {Le sch\'ema $\P$ est lisse sur $\SP$.}
\item {Le sch\'ema $\P$ est r\'egulier.}
\item {Le sch\'ema $\M{P}(\D)_{k}$ est lisse sur $k$.}
\end{itemize}
\end{rem}
\begin{expl}
\label{projectif}
On prend $N = \Z^{2}$. On note $e_{0} = - e_{1} - e_{2}$ et on pose $\si_{1} =
\R^{+}e_{1} + \R^{+}e_{2}$, $\si_{2} = \R^{+}e_{0} + \R^{+}e_{2}$ et $\si_{3} =
\R^{+}e_{0} + \R^{+}e_{1}$. 
\bigskip
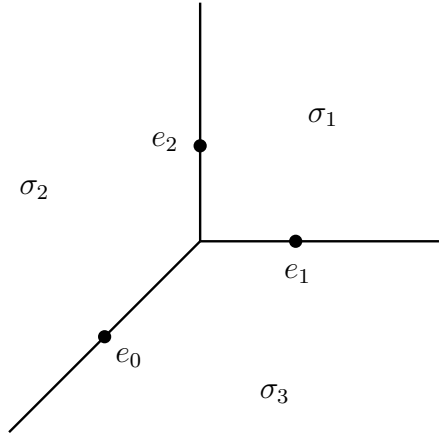
\begin{center}
\begin{picture}(0,0)%
\includegraphics{figure2.pstex}%
\end{picture}%
\setlength{\unitlength}{0.00083300in}%
\begingroup\makeatletter\ifx\SetFigFont\undefined
\def\x#1#2#3#4#5#6#7\relax{\def\x{#1#2#3#4#5#6}}%
\expandafter\x\fmtname xxxxxx\relax \def\y{splain}%
\ifx\x\y   
\gdef\SetFigFont#1#2#3{%
  \ifnum #1<17\tiny\else \ifnum #1<20\small\else
  \ifnum #1<24\normalsize\else \ifnum #1<29\large\else
  \ifnum #1<34\Large\else \ifnum #1<41\LARGE\else
     \huge\fi\fi\fi\fi\fi\fi
  \csname #3\endcsname}%
\else
\gdef\SetFigFont#1#2#3{\begingroup
  \count@#1\relax \ifnum 25<\count@\count@25\fi
  \def\x{\endgroup\@setsize\SetFigFont{#2pt}}%
  \expandafter\x
    \csname \romannumeral\the\count@ pt\expandafter\endcsname
    \csname @\romannumeral\the\count@ pt\endcsname
  \csname #3\endcsname}%
\fi
\fi\endgroup
\begin{picture}(2744,2744)(3579,-3983)
\put(4276,-3511){\makebox(0,0)[lb]{\smash{\SetFigFont{12}{14.4}{rm}$e_{0}$}}}
\put(5476,-2011){\makebox(0,0)[lb]{\smash{\SetFigFont{12}{14.4}{rm}$\si_{1}$}}}
\put(5326,-2986){\makebox(0,0)[lb]{\smash{\SetFigFont{12}{14.4}{rm}$e_{1}$}}}
\put(5176,-3736){\makebox(0,0)[lb]{\smash{\SetFigFont{12}{14.4}{rm}$\si_{3}$}}}
\put(3676,-2461){\makebox(0,0)[lb]{\smash{\SetFigFont{12}{14.4}{rm}$\si_{2}$}}}
\put(4501,-2161){\makebox(0,0)[lb]{\smash{\SetFigFont{12}{14.4}{rm}$e_{2}$}}}
\end{picture}

\end{center}
\bigskip
En prenant $\si_{1}$, $\si_{2}$, $\si_{3}$ ainsi que leurs faces $\R^{+}e_{0}$,
$\R^{+}e_{1}$, $\R^{+}e_{2}$ et $\{0\}$, on obtient un \'eventail complet et
r\'egulier $\D_{2}$. Dans $M$, les \cos\ duaux sont donn\'es par~:
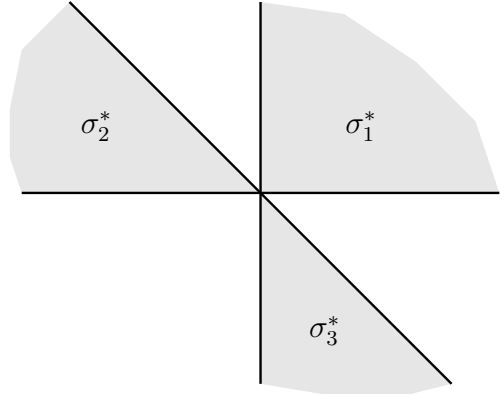
\bigskip
\begin{center}
\begin{picture}(0,0)%
\includegraphics{figure3.pstex}%
\end{picture}%
\setlength{\unitlength}{0.00083300in}%
\begingroup\makeatletter\ifx\SetFigFont\undefined
\def\x#1#2#3#4#5#6#7\relax{\def\x{#1#2#3#4#5#6}}%
\expandafter\x\fmtname xxxxxx\relax \def\y{splain}%
\ifx\x\y   
\gdef\SetFigFont#1#2#3{%
  \ifnum #1<17\tiny\else \ifnum #1<20\small\else
  \ifnum #1<24\normalsize\else \ifnum #1<29\large\else
  \ifnum #1<34\Large\else \ifnum #1<41\LARGE\else
     \huge\fi\fi\fi\fi\fi\fi
  \csname #3\endcsname}%
\else
\gdef\SetFigFont#1#2#3{\begingroup
  \count@#1\relax \ifnum 25<\count@\count@25\fi
  \def\x{\endgroup\@setsize\SetFigFont{#2pt}}%
  \expandafter\x
    \csname \romannumeral\the\count@ pt\expandafter\endcsname
    \csname @\romannumeral\the\count@ pt\endcsname
  \csname #3\endcsname}%
\fi
\fi\endgroup
\begin{picture}(3098,2498)(3225,-7337)
\put(5101,-6961){\makebox(0,0)[lb]{\smash{\SetFigFont{12}{14.4}{rm}$\si_{3}^{\ast}$}}}
\put(3676,-5686){\makebox(0,0)[lb]{\smash{\SetFigFont{12}{14.4}{rm}$\si_{2}^{\ast}$}}}
\put(5326,-5686){\makebox(0,0)[lb]{\smash{\SetFigFont{12}{14.4}{rm}$\si_{1}^{\ast}$}}}
\end{picture}

\end{center}
\bigskip
Les ouverts $U_{\si_{1}}$, $U_{\si_{2}}$ et $U_{\si_{3}}$ sont des plans 
affines que l'on recolle par les applications~:
\begin{alignat*}{9}
\Theta_{1,2}:\; & U_{\si_{1}}& &\rightarrow & &U_{\si_{2}}\qquad & 
\Theta_{1,3}:\; & U_{\si_{1}}& &\rightarrow & &U_{\si_{3}}\qquad &
\Theta_{2,3}:\; & U_{\si_{2}}& &\rightarrow & &U_{\si_{3}}\qquad \\
(&x,y) & &\mapsto & \bigg( & \frac{y}{x},\frac{1}{x}\bigg) & 
(&x,y) & &\mapsto & \bigg( & \frac{x}{y},\frac{1}{y}\bigg) & 
(&x,y) & &\mapsto & \bigg( & \frac{1}{x},\frac{y}{x}\bigg)
\end{alignat*}
On en d\'eduit que $\M{P}(\D_{2})$ s'identifie au plan projectif $\M{P}^{2}_{\Z}$.
\end{expl}
La proposition suivante donne une d\'ecomposition de $X$ sous forme d'une
r\'eu\-nion disjointe de tores~:
\begin{prop}
\label{decomposition1}
Soit $\D$ un \'eventail de $N_{\M{R}}$ et $\P$ la vari\'et\'e torique
associ\'ee. Pour tout $\si \in \D$ on consid\`ere le tore~:
\[
\E{O}(\si) = \operatorname{Spec}\left(\M{Z}\left[ M \cap \si^{\perp}
 \right]\right).
\]
Le tore $\E{O}(\si)$ se plonge de mani\`ere canonique dans l'ouvert $U_{\si}$ (et donc
dans $\P$) par le morphisme~:
\[
i_{\si}: \E{O}(\si) = \operatorname{Spec}\left(\M{Z}\left[ M \cap \si^{\perp}
 \right]\right) \hooklongrightarrow 
\operatorname{Spec}\left(\M{Z}\left[ M \cap \si^{\ast}
 \right]\right) = U_{\si},
\]
obtenu par prolongement par z\'ero (i.e. induit par le morphisme $i_{\si}:
\Z\left[M\cap \si^{\ast}\right] \rightarrow \Z \left[ M\cap\si^{\perp}
\right]$ d\'efini par $i_{\si}(\chi^{m}) = \chi^{m}$ si $m \in \si^{\perp}$ et
$i_{\si}(\chi^{m}) = 0$ sinon). De plus~:
\begin{enumerate}
\item {Soit $k$ un corps alg\'ebriquement clos, toute $T(k)$-orbite de $\M{P}(\D)(k)$ est de la
forme $\E{O}(\si)(k)$ avec $\si
\in \P$.}
\item {Le sch\'ema $\P$ est r\'eunion disjointe des sous-sch\'emas $\E{O}(\si)$
pour $\si$ parcourant $\D$ et cette d\'ecomposition est respect\'ee par
l'action de $T$.}
\item {On a~:
	\begin{itemize}
		\item{$\E{O}\left(\{O\}\right) = U_{\{O\}} = T.$}
		\item{$\t < \si \Leftrightarrow \E{O}(\si) \subset
\overline{\E{O}(\t)}.$}
		\item{$U_{\si} = \bigcup_{\t < \si}\E{O}(\t)$.}
	\end{itemize}
	}
\item {Pour tout $\si \in \D$, notons $V(\si) = \overline{\E{O}(\si)}$; c'est
une vari\'et\'e torique et l'on a~:
\[
V(\si) = \bigcup_{\si < \t}\E{O}(\t).
\]
	}
\end{enumerate}
\end{prop}
\demo\ Voir (\cite{8}, \S 4, prop. 2). On pourra \'egalement consulter
(\cite{20}, prop. 1.6) et (\cite{11}, \S 3.1).
\medskip

\begin{rem}
\label{rem_decomposition1}
Soient $\si \in \D$ et $N_{\si} \subset N$ le $\M{Z}$-module engendr\'e par
$\si \cap N$. On note $N(\si)$ le $\M{Z}$-module quotient $N /
N_{\si}$ et $M(\si) = \si^{\perp} \cap M$ son dual.

On appelle {\it \'etoile\/} du c\^one $\si$ l'ensemble des c\^ones $\t \in \D$
contenant $\si$. Soit $\t$ un tel c\^one; on note $\ov{\t}$ son image dans
$N(\si)_{\R}$, c'est-\`a-dire~:
\[
\ov{\t} = (\t + (N_{\si})_{\R})\left/ (N_{\si})_{\R} \subset N(\si)_{\R}\right.
.
\]
L'ensemble $\{\ov{\t}: \t \in \D, \; \si < \t\}$ forme un \'eventail de
$N(\si)$ que l'on note $\D(\si)$. On dispose de l'inclusion canonique
$\E{O}(\si) \subset \M{P}(\D(\si))$ correspondant au c\^one $\{0\} \in
\D(\si)$.

On peut montrer (voir \cite{11}, \S 3.1; ou \cite{20}, cor. 1.7)
que le morphisme $i_{\si}: \E{O}(\si) \hookrightarrow U_{\si}$ introduit \`a la
proposition (\ref{decomposition1}) s'\'etend de mani\`ere canonique en une immersion ferm\'ee
$i_{\si}: \M{P}(\D(\si)) \hookrightarrow \P$ qui a pour image $V(\si)$. Dans
toute la suite, on identifie $\M{P}(\D(\si))$ \`a $V(\si)$ par  
cet isomorphisme canonique.
\end{rem}

La d\'efinition et la proposition suivantes d\'ecrivent les morphismes
naturels entre vari\'et\'es toriques~:
\begin{defn}
Soient $(N,\D)$ et $(N',\D')$ deux \'eventails, avec $N \simeq \M{Z}^{d}$ et
$N' \simeq \M{Z}^{d'}$; un {\it morphisme d'\'eventails\/} $\varphi: (N',\D')
\rightarrow (N,\D)$ est un morphisme de $\M{Z}$-module $\varphi: N'\rightarrow N$
telle que l'application induite~: $\varphi_{\M{R}}: N'_{\M{R}} \rightarrow
N_{\M{R}}$, d\'efinie par extension des scalaires \`a partir de $\varphi$,
v\'erifie~: pour tout $\si'\in\D'$, il existe $\si\in\D$ tel que $\varphi(\si')
\subset \si$.
\end{defn}

Soit $\varphi: (N',\D') \rightarrow (N,\D)$ un tel morphisme d'\'eventails. On
construit \`a partir de $\varphi$ un morphisme \'equivariant $\varphi_{\ast}:
\M{P}(\D') \rightarrow \P$ de la fa\c con suivante~: On note $^{t}\varphi: M
\rightarrow M'$ la transpos\'ee de $\varphi$ et $^{t}\varphi_{\M{R}}:
M_{\M{R}} \rightarrow M'_{\M{R}}$ l'application d\'efinie par extension des
scalaires \`a partir de $^{t}\varphi$. Soient $\si' \in \D'$ et $\si \in \D$
tels que $\varphi_{\M{R}}(\si') \subset \si$. De l'inclusion
$^{t}\varphi_{\M{R}}(\si^{\ast}) \subset (\si')^{\ast}$, on tire que
$^{t}\varphi(\E{S}_{\si}) \subset \E{S}_{\si'}$. On dispose donc d'une
application $^{t}\varphi: \E{S}_{\si} \rightarrow \E{S}_{\si'}$ qui induit un
morphisme \'equivariant $\varphi_{\ast}: U_{\si'} \rightarrow U_{\si}$. En
particulier, si l'on prend $\si' = \{0'\}$ et $\si = \{0\}$, on obtient le
morphisme de tore~:
\[
\varphi_{\ast}: T' = \op{Spec}(\M{Z}[M']) \longrightarrow \op{Spec}(\M{Z}[M]) = T
\]
induit par l'application $^{t}\varphi : M \rightarrow M'$. La proposition 
suivante affirme qu'on peut recoller ces constructions locales pour obtenir un
morphisme global \'equivariant $\varphi_{\ast}$ et donne une condition
n\'ecessaire et suffisante sur $\varphi$ pour que $\varphi_{\ast}$ soit
propre~:

\begin{prop}
Soit un morphisme d'\'eventails $\varphi:(N',\D')\rightarrow (N,\D)$. Le
morphisme de tore alg\'ebrique~:
\[
\varphi_{\ast}:\,T' = \op{Spec}\left(\Z[M']\right) \longrightarrow
\op{Spec}\left(\Z [M]\right) = T,
\]
induit par l'application duale $^{t}\varphi : M \rightarrow M'$, se prolonge 
en un morphisme~:
\[
\varphi_{\ast}:\, \M{P}(\D') \longrightarrow \P.
\]
Le morphisme $\varphi_{\ast}$ est \'equivariant sous l'action de $T'$ et $T$.
De plus, $\varphi_{\ast}$ est propre \ssi~:
\[
\varphi_{\R}^{-1}\left(|\D|\right) = |\D'|.
\]
\end{prop}
\demo\ On peut consulter (\cite{20}, prop. 1.13 et 1.15) et aussi (\cite{11},
\S 1.4 et 2.4).
\begin{expl}
\label{definition_endo}
Soit $p$ un entier sup\'erieur ou \'egal \`a un et prenons $N' = N$ et $\D' =
\D$, et soit~: 
\begin{alignat*}{3}
[p]:\,&N& &\longrightarrow& &N \\
&n& &\longmapsto& &p\,n.
\end{alignat*}
On note encore $[p]$ l'endomorphisme de $\P$ induit par $[p]$. D'apr\`es la
proposition pr\'ec\'edente, le morphisme $[p]$ est propre. Sa restriction \`a
chacun des tores $\E{O}(\si)$ est le morphisme {\it puissance $p$-i\`eme\/}~:
\begin{alignat*}{3}
[p]:\, \E{O}&(\si) &&\longrightarrow &\E{O}&(\si) \\
&x &&\longmapsto &&x^{p}.
\end{alignat*}
On en d\'eduit que pour tout corps $k$, le morphisme $[p]_{k}: \M{P}(\D)_{k}
\rightarrow \M{P}(\D)_{k}$ obtenu par extension des scalaires est fini de
degr\'e $p^{d}$.
\end{expl}
\begin{defn}
\label{raffinement_1}
Soient $\D$ et $\D'$ deux \'eventails de $N$. On dit que $\D'$ est {\em plus
fin} que $\D$ ou encore que $\D'$ est un {\em raffinement} de $\D$ si pour tout
$\si'\in\D'$ il existe $\si \in \D$ tel que $\si' \subset \si$, et si de plus
$|\D'| = |\D|$. L'inclusion induit un morphisme \'equivariant propre canonique
$i_{\ast}: \M{P}(\D') \rightarrow \M{P}(\D)$.
\end{defn}
\bigskip

\subsection{Diviseurs invariants sur $\P$}~

On consid\`ere $\D$ un \'eventail complet, de sorte que la vari\'et\'e 
torique associ\'ee $\P$ est propre. On note $\tau_{1}, \dots, \tau_{r}$ les
\'el\'ements de $\D(1)$, c'est-\`a-dire les demi-droites de $\D$ et 
$u_{1}, \dots, u_{r}$ leur g\'en\'erateur dans $N$, c'est-\`a-dire les
\'el\'ements de $N$ tels que $\t_{i} \cap N = \M{N}u_{i}$. A tout $\t_{i}$
on a associ\'e pr\'ec\'edemment $V(\t_{i}) = \overline{\E{O}(\t_{i})}$ un
sch\'ema irr\'eductible de codimension $1$ invariant sous l'action de $T$.
\begin{defn}
On appelle {\it diviseur invariant \'el\'ementaire\/} ou plus simplement {\it 
diviseur \'el\'ementaire\/} et on note $D_{i}$ le cycle donn\'e par
$V(\t_{i})$.
\end{defn}
\begin{prop}
\label{inv_intro}
Tout diviseur de Weil $D$ sur $\P$ horizontal invariant par $T$ (i.e. dont
la restriction \`a la fibre g\'en\'erique est laiss\'ee invariante par
$T_{\M{Q}}$) est de la forme~:
\[
D = \sum_{i = 1}^{r}a_{i}D_{i} \qquad (a_{i} \in \Z).
\]
\end{prop}
\demo\ Cela d\'ecoule directement de la d\'ecomposition donn\'ee dans la
proposition (\ref{decomposition1}).
\medskip

\begin{prop}
\label{ordre}
Soit $m \in M$ et $\chi^{m}$ le caract\`ere associ\'e; l'ordre de $\chi^{m}$ en
$D_{i}$ est donn\'e par~:
\[
\operatorname{ord}_{D_{i}}(\chi^{m}) = <m,u_{i}>.
\]
\end{prop}
\demo\ Voir par exemple (\cite{11}, lemme p. 61).
\medskip

\begin{defn}
On dit qu'une fonction $\psi : N_{\R} \rightarrow \R$ est lin\'eaire par
morceaux sur $\D$ (ou plus simplement lin\'eaire par morceaux) si la 
restriction de $\psi$ \`a chacun des \cos\ de $\D$ est d\'efinie par
une forme lin\'eaire $m_{\psi,\si} \in M$.
\end{defn}
\begin{prop}
Un diviseur de Weil horizontal $T$-invariant $D = \sum_{i = 1}^{r}a_{i}D_{i}$ sur
$\P$ provient d'un diviseur de Cartier \ssi\ il existe une fonction $\psi :
N_{\R} \rightarrow \R$ continue et lin\'eaire par morceaux sur $\D$ telle que
$\psi(u_{i}) = -a_{i}$ pour $(1 \leqslant i \leqslant r)$. Si elle existe, une
telle fonction $\psi$ est unique.
\end{prop}
\demo\ Voir (\cite{20}, prop. 2.1) et (\cite{11}, p. 66).
\medskip

\begin{defn}
Soit $D$ un diviseur de Cartier sur $\P$ horizontal et $T$-inva\-riant, on
appelle {\it fonction support associ\'ee \`a $D$\/} et on notera $\psi_{D}$ la
fonction d\'efinie par la proposition ci-dessus. On notera $m_{D,\si}$
la forme lin\'eaire d\'efinissant $\psi_{D}$ sur $\si \in \D$.
\end{defn}
\begin{rem}
Lorsque $\P$ est lisse, tout diviseur de Weil est un diviseur de Cartier. On
remarquera que sous cette hypoth\`ese, l'existence d'une fonction support
d\'ecoule du crit\`ere de lissit\'e (\ref{lissite}).
\end{rem}
\begin{prop}
\label{intro_inverse}
Si $\varphi: (N',\D') \rightarrow (N,\D)$ est un morphisme d'\'eventails et $D$
un diviseur de Cartier $T$-invariant sur $\P$ de fonction support $\psi_{D}$, alors
le diviseur de Cartier $T$-invariant $(\varphi_{\ast})^{\ast}(D)$ 
sur $\M{P}(\D')$ admet $\psi_{d}\circ
\varphi$ comme fonction support.
\end{prop}
\demo\ Soient $\si' \in \D'$ et $\si \in \D$ tels que $\varphi(\si') \subset
\si$.
On pose $D' = (\varphi_{\ast})^{\ast}(D)$ et on note $m_{D,\si}$ (resp.
$m_{D',\si'}$) la forme lin\'eaire d\'efinissant $\psi_{D}$ sur $\si$ (resp.
$\psi_{D} \circ \varphi$ sur $\si'$). D'apr\`es la proposition (\ref{ordre}) le
diviseur $D$ est de la forme $\op{div}(\chi^{m_{D,\si}})$ sur $U_{\si}$, et
donc $D'$ est de la forme $\op{div}(\chi^{m_{D,\si}}\circ \varphi_{\ast})$ sur
$U_{\si'}$. Il suffit alors de remarquer que~: 
\[
\chi^{m_{D,\si}}\circ \varphi_{\ast} = \chi^{^t\varphi(m_{D,\si})} =
\chi^{m_{D',\si'}}, 
\]
et d'appliquer une nouvelle fois la proposition (\ref{ordre}) pour conclure.
\medskip

Le lemme suivant est une simple cons\'equence de la proposition (\ref{ordre}).
\begin{lem}
\label{locale}
Soit $D = \sum_{i=1}^{r}a_{i}D_{i}$ un diviseur de Cartier horizontal
$T$-invariant et $\C{O}(D)$ le faisceau inversible associ\'e \`a $D$. Pour
tout $\si \in \D$, on pose~:
\[
P_{D}(\si) = \{ v\in M_{\R}:\quad <v,u> \geqslant \psi_{D}(u), \quad \forall
u \in \si \} = \si^{\ast} + m_{D, \si}.
\]
Le $\Z$-module des sections de $\C{O}(D)$ sur $U_{\si}$ est donn\'e par~:
\[
\Gamma(U_{\si},\C{O}(D)) = \bigoplus_{m \in P_{D}(\si) \cap M}\Z \chi^{m}.
\]
\end{lem}
\medskip

La proposition suivante, qui d\'ecrit les sections globales d'un 
diviseur de Cartier $T$-invariant, d\'ecoule directement du lemme
pr\'ec\'edent~:
\begin{prop}
\label{sections_globales}
Soit $D = \sum_{i =1}^{r}a_{i}D_{i}$ un diviseur de Cartier horizontal
$T$-invariant et $\C{O}(D)$ le faisceau inversible associ\'e \`a $D$. On note
$K_{D}$ le polytope convexe de $M_{\R}$ d\'efini par les in\'equations
suivantes~:
\begin{align*}
K_{D} = &\{v \in M_{\R}, \quad <v,u_{i}> \geqslant -a_{i}, \quad 0 \leqslant i
\leqslant r\} \\
 = &\{v \in M_{\R}, \quad <v,u> \geqslant \psi_{D}(u),\quad \forall u \in N_{\R}\}.
\end{align*}
Le $\Z$-module des sections globales de $\C{O}(D)$ est donn\'e par~:
\[
\Gamma \left(\P,\C{O}(D)\right) = \bigoplus_{m \in K_{D} \cap M}\Z \chi^{m}.
\]
\end{prop}
Pour plus de d\'etails, on peut consulter (\cite{20},
lemme 2.3 et \cite{11}, p. 66), les arguments donn\'es s'\'etendant
imm\'ediatement \`a la situation sur $\SP$.
\medskip

\begin{defn}
Une fonction $\psi: N_{\R} \rightarrow \R$ est dite {\it concave} si~:
\[
\psi(tx + (1-t)y) \geqslant t\psi(x) + (1-t)\psi(y), \qquad \forall t \in [0,1]
\]
pour tout $x$, $y \in N_{\R}$.
\end{defn}
\begin{defn} {\bf (Minkowski).} Soit $K$ un compact convexe non vide de
$M_{\R}$. On appelle {\it fonction d'appui associ\'ee \`a $K$\/} et on note
$\psi_{K}$ la fonction $\psi_{K}: N_{\R} \rightarrow \R$ d\'efinie par~:
\[
\psi_{K}(u) = \op{inf}\{<v,u>, \quad v \in K\}
\]
pour tout $u \in N_{\R}$.
\end{defn}

On a alors le r\'esultat suivant (voir par exemple \cite{20}, th. A. 18)~:
\begin{thm}
\label{correspondance}
Soit $\C{C}(M_{\R})$ l'ensemble des compacts convexes non vides de $M_{|R}$
et notons $\C{S}(N_{\R})$ l'ensemble des fonctions $\psi: N_{\R} \rightarrow
\R$ positivement homog\`enes (i.e. telles que $\psi(cu) = c\psi(u)$ pour tout
$c \in \R^{+}$ et $u \in N_{\R}$) et concaves. Pour tout $\psi \in
\C{S}(N_{\R})$, on d\'efinit $K_{\psi}$ le {\it compact connexe associ\'e \`a
$\psi$\/} par les in\'egalit\'es suivantes~:
\[
K_{\psi} = \{v \in M_{\R}:\quad <v,u> \geqslant \psi(u), \quad \forall u \in
N_{\R}\}.
\]
On a alors~:
\begin{itemize}
\item {Les applications $\C{C}(M_{\R}) \rightarrow \C{S}(N_{\R})$ 
et $\C{S}(N_{\R}) \rightarrow \C{C}(M_{\R})$ qui envoient respectivement
$K$ sur $\psi_{K}$ et $\psi$ sur $K_{\psi}$ sont r\'eciproques l'une de
l'autre.}
\item {Par la correspondance biunivoque d\'ecrite ci-dessus, la somme de
Minkowski $K + K'$ de deux compacts convexes $K, K' \in \C{C}(M_{\R})$ (resp.
le dilat\'e $cK$, $c \in \R^{+}$) est associ\'ee \`a la somme $\psi_{K} +
\psi_{K'}$ (resp. \`a $c\psi_{K}$).}
\end{itemize}
\end{thm}
On peut alors \'enoncer~:
\begin{thm}
\label{polytope_et_vt}
Soit $D$ un diviseur de Cartier horizontal $T$-invariant sur $\P$. Le faisceau inversible
$\C{O}(D)$ est engendr\'e par ses sections globales \ssi\ $\psi_{D}$ la
fonction support de $D$ est concave. De plus, si pour tout \co\ $\si \in \D$
on note $m_{D, \si}$ 
l'\'el\'ement de $M$ d\'efinissant la
restriction de $\psi_{D}$ \`a $\si$, le polytope convexe
$K_{D}$ associ\'e \`a $D$ est alors l'enveloppe convexe des $m_{D, \si}$ dans
$M_{\R}$ pour $\si$ parcourant $\D_{\op{max}}$. En particulier, $K_{D}$ est \`a
sommets entiers. Enfin $K_{D}$ et $\psi_{D}$ sont images l'un de l'autre par la
correspondance d\'efinie \`a la proposition $($\ref{correspondance}$)$.
\end{thm}
\demo\ Voir (\cite{20}, th 2.7) ou (\cite{11}, p. 68).
\medskip

\begin{prop}
\label{additivite_1}
Soient $D$ et $D'$ deux diviseurs horizontaux $T$-invariants et engendr\'es par leurs
sections globales. Soit $D'' = D + D'$, on a~:
\begin{align*}
\psi_{D''} &= \psi_{D} + \psi_{D'} \\
K_{D''} &= K_{D} + K_{D'}.
\end{align*}
\end{prop}
\demo\ C'est une cons\'equence des \'enonc\'es pr\'ec\'edents. On peut
consulter (\cite{11}, p. 69) pour plus de d\'etails.
\medskip

\begin{rem}
On remarquera que l'on a pas n\'ecessairement~:
\[
K_{D''} \cap M = \left(K_{D} \cap M\right) + \left(K_{D'} \cap M\right).
\]
\end{rem}
\begin{defn}
Soit $\psi$ une fonction concave lin\'eaire par morceaux sur un \'eventail $\D$
complet. On dit que
$\psi$ est {\it strictement concave relativement \`a $\D$\/} (ou plus
simplement {\it strictement concave\/} lorsque aucune confusion n'est \`a
craindre)
\ssi\ $\D$ est l'\'eventail complet le plus grossier dans $N$ tel que
$\psi_{|\si}$ soit lin\'eaire pour tout $\si \in \D$.
\end{defn}
\begin{thm}
Soit $D$ un diviseur de Cartier horizontal $T$-invariant sur $\P$ et $\psi_{D}$
la fonction
support associ\'ee. Les trois assertions suivantes sont \'equivalentes~:
\begin{enumerate}
\item{Le diviseur $D$ est ample.}
\item{La fonction support $\psi_{D}$ est strictement concave relativement \`a
$\D$.}
\item{Le polytope $K_{D}$ est de dimension $d$, et si pour tout $\si \in \D$,
on note $m_{D, \si}$ l'\'el\'ement de $M$ donnant la restriction de $\psi_{D}$
\`a $\si$, alors les sommets de $K_{D}$ sont donn\'e par $\{m_{D,\si}, \si \in
\D_{\op{max}}\}$. De plus $m_{D,\si} \not= m_{D,\t}$ pour $\si$, $\t \in 
\D_{\op{max}}$ d\`es
que $\si \not= \t $.}
\end{enumerate}
\end{thm}
\demo\ Comme $\op{Spec}\M{Z}$ est un sch\'ema affine, le diviseur $D$ est ample
sur $\P$ si et seulement s'il est ample relativement \`a $\op{Spec}\M{Z}$ d'apr\`es
(\cite{27}, cor. 4.6.6). 

Par ailleurs, d'apr\`es (\cite{27}, cor. 4.6.4) et le crit\`ere d'amplitude
donn\'e dans (\cite{28}, th. 4.7.1), l'amplitude de $D$ sur $\P$ relativement
\`a $\op{Spec}\M{Z}$ est \'equivalente \`a l'amplitude de $D$ sur
$\P_{\M{F}_{p}}$ pour tout nombre premier $p$.

Enfin pour tout corps $k$, (\cite{20}, cor. 2.14) ou (\cite{11}, p. 70)
modifi\'es de fa\c con \'evidente montrent que les trois assertions du
th\'eor\`eme sont \'equivalentes sur $\P_{k}$.
\medskip

On remarque en particulier que sur une vari\'et\'e torique propre, tout
diviseur de Cartier horizontal $T$-invariant ample est engendr\'e par ses sections
globales.

Concernant les fibr\'es tr\`es amples, on a~:
\begin{thm}
Soit $D$ un diviseur de Cartier horizontal $T$-invariant sur $\P$ et $\psi_{D}$ sa
fonction support associ\'ee. Les trois assertions suivantes sont
\'equivalentes~:
\begin{enumerate}
\item{Le diviseur $D$ est tr\`es ample (relativement \`a $\op{Spec}\M{Z}$).}
\item{La fonction $\psi_{D}$ est strictement concave relativement \`a $\D$. De
plus, pour tout $\si \in \D_{\op{max}}$, l'ensemble $\left(M \cap K_{D}\right) -
m_{D,\si}$ engendre le semi-groupe $M\cap \si^{\ast} = \E{S}_{\si}$.}
\item{Le polytope $K_{D}$ est de dimension $d$, l'ensemble de ses sommets est
donn\'e par $\{m_{D,\si}, \si \in \D_{\op{max}}\}$. De plus, $\left( M\cap K_{D}\right)
- m_{D, \si}$ engendre le semi-groupe $M \cap \si^{\ast} = \E{S}_{\si}$ pour
tout $\si \in \D_{\op{max}}$.}
\end{enumerate}
\end{thm}
\demo\ Il suffit de modifier de fa\c con \'evidente la preuve de (\cite{20}, th 2.13).
\medskip

On a enfin le th\'eor\`eme suivant~:
\begin{thm}
\label{simplicite}
{\rm \bf (Demazure).} Soit $\P$ une vari\'et\'e torique propre et lisse
de dimension relative $d$. Soit
$D$ un diviseur horizontal $T$-invariant sur $\P$ et $\psi_{D}$ sa fonction support
associ\'ee. Les quatre assertions suivantes sont \'equivalentes~:
\begin{enumerate}
\item{Le diviseur $D$ est ample.}
\item{Le diviseur $D$ est tr\`es ample.}
\item{La fonction $\psi_{D}$ est strictement concave relativement \`a $\D$.}
\item{Le polytope $K_{D}$ est de dimension $d$, et si pour tout $\si \in \D$,
on note $m_{D,\si}$ l'\'el\'ement de $M$ donnant  la restriction de $\psi_{D}$
\`a $\si$, alors les sommets de $K_{D}$ sont donn\'es par $\{m_{D,\si}, \si \in
\D_{\op{max}}\}$. De plus $m_{D,\si} \not= m_{D,\t}$ pour $\si$, $\t \in
\D_{\op{max}}$ d\`es
que $\si \not= \t$.}
\end{enumerate}
Lorsque ces conditions sont satisfaites, le polytope $K_{D}$ est {\rm absolument
simple} dans le sens o\`u chaque sommet $m_{D,\si}$ rencontre exactement $d$
ar\^etes et o\`u,
si $m_{1,\si}, \dots, m_{d,\si}$ sont les points de $M$ les plus proches de $m_{D,\si}$
sur chacune de ces diff\'erentes ar\^etes, alors
$\{m_{1,\si}-m_{D,\si}, \dots, m_{d,\si}-m_{D,\si}\}$ est une
base de $M$.
\end{thm}
\demo\ Voir (\cite{8}, \S 4, th. 2 et cor. 1) et aussi (\cite{20}, cor. 2.15).
\bigskip

\subsection{Vari\'et\'e torique et fibr\'e en droites associ\'e \`a un
polytope}~

Dans ce paragraphe, on suit \cite{17}; on peut \'egalement consulter
(\cite{20}, th. 2.22, et \cite{11}, \S 5.5).

On a vu comment associer \`a tout diviseur de Cartier horizontal $T$-invariant $D$ sur une
vari\'et\'e torique propre $\P$ un polytope convexe $K_{D} \subset M_{\R}$.
Inversement le probl\`eme suivant se pose~: \'etant donn\'e un
polytope convexe \`a sommets entiers $K \subset M_{\R}$, peut-on construire une
vari\'et\'e torique propre $\M{P}(\D_{K})$ et un diviseur horizontal $T$-invariant $E$ sur 
$\M{P}(\D_{K})$ tels que le polytope $K_{E}$ soit \'egal \`a $K$ ? La r\'eponse
est donn\'ee par le th\'eor\`eme suivant~: 
\begin{thm}
\label{construction_inverse}
Soit $K \subset M_{\R}$ un polytope convexe d'int\'erieur non vide dont les sommets
sont dans $M$. Il existe un unique \'eventail complet $\D$ dans $N_{\R}$ 
et un unique diviseur de Cartier $E$
horizontal
$T$-invariant sur $\M{P}(\D)$ tels que~:
\begin{enumerate}
\item{$K_{E} = K$.}
\item{Le diviseur $E$ est ample.}
\end{enumerate}
L'\'eventail $\D$ est le plus petit \'eventail complet tel que la fonction
d'appui $\psi_{K}$ est lin\'eaire par morceau relativement \`a $\D$.
De plus, $\M{P}(\D)$ est lisse \ssi\ le polytope $K$ est absolument simple;
dans ce cas, le diviseur $E$ est tr\`es ample.
\end{thm}
\demo\ On suit ici (\cite{20}, th. 2.22) et (\cite{17}, \S 3, th. 1). 
On d\'efinit l'\'eventail $\D$ comme dans le dernier alin\'ea de l'\'enonc\'e. 
La fonction d'appui $\psi_{K}$ est continue et concave. Soient $u_{1}, \dots,
u_{r}$ les g\'en\'erateurs dans $N$ des demi-droites $\t_{1},\dots, \t_{r}$ de
$\D(1)$; le diviseur~:
\[
E = -\sum_{i=1}^{r}\psi_{K}(u_{i})V(\t_{i})
\]
est un diviseur de Cartier horizontal $T$-invariant engendr\'e par ses sections 
globales d'apr\`es (\ref{polytope_et_vt}).
De plus, par construction, $\psi_{K}$ est strictement concave relativement \`a
$\D$, donc $E$ est ample et $\D$ est l'unique \'eventail \`a
satisfaire cette condition. On a $K_{E} = K_{\psi_{K}} = K$. Enfin,
$\M{P}(\D)$ est lisse \ssi\ $K$ est absolument simple d'apr\`es
(\ref{lissite}) et (\ref{simplicite}).
\medskip
\begin{expl}
\label{exemple_intro_1}
Lorsque $d = 2$,
consid\'erons $K_{2} \subset M_{\R} = \R^{2}$ le polytope convexe absolument simple
d\'efini par les in\'equations~:
\[
x_{1} \geqslant 0, \quad x_{2} \geqslant 0 \quad \text{et} \quad x_{1} + x_{2}
\leqslant 1.
\]
\bigskip
\begin{center}
\begin{picture}(0,0)%
\includegraphics{figure4.pstex}%
\end{picture}%
\setlength{\unitlength}{0.00083300in}%
\begingroup\makeatletter\ifx\SetFigFont\undefined
\def\x#1#2#3#4#5#6#7\relax{\def\x{#1#2#3#4#5#6}}%
\expandafter\x\fmtname xxxxxx\relax \def\y{splain}%
\ifx\x\y   
\gdef\SetFigFont#1#2#3{%
  \ifnum #1<17\tiny\else \ifnum #1<20\small\else
  \ifnum #1<24\normalsize\else \ifnum #1<29\large\else
  \ifnum #1<34\Large\else \ifnum #1<41\LARGE\else
     \huge\fi\fi\fi\fi\fi\fi
  \csname #3\endcsname}%
\else
\gdef\SetFigFont#1#2#3{\begingroup
  \count@#1\relax \ifnum 25<\count@\count@25\fi
  \def\x{\endgroup\@setsize\SetFigFont{#2pt}}%
  \expandafter\x
    \csname \romannumeral\the\count@ pt\expandafter\endcsname
    \csname @\romannumeral\the\count@ pt\endcsname
  \csname #3\endcsname}%
\fi
\fi\endgroup
\begin{picture}(3024,2424)(3889,-6973)
\put(5101,-6061){\makebox(0,0)[lb]{\smash{\SetFigFont{12}{14.4}{rm}$K_{2}$}}}
\put(4501,-5236){\makebox(0,0)[lb]{\smash{\SetFigFont{12}{14.4}{rm}$e_{2}^{\ast}$}}}
\put(4576,-6586){\makebox(0,0)[lb]{\smash{\SetFigFont{12}{14.4}{rm}$0$}}}
\put(5926,-6586){\makebox(0,0)[lb]{\smash{\SetFigFont{12}{14.4}{rm}$e_{1}^{\ast}$}}}
\end{picture}

\end{center}
\bigskip
L'\'eventail complet de $N_{\R}$ associ\'e \`a $K_{2}$ est l'\'eventail
$\D_{2}$ d\'ej\`a consid\'er\'e \`a l'exemple (\ref{projectif}) dont la vari\'et\'e
torique associ\'ee est le plan projectif $\M{P}_{\Z}^{2}$. La fonction
$\psi_{K_{2}}$ est d\'efinie sur $\si_{1}$, $\si_{2}$ et $\si_{3}$ par
respectivement $m_{1} = (0,0)$, $m_{2} = (1,0)$ et $m_{3} = (0,1)$.
En particulier $\psi_{K_{2}}(e_{1}) = \psi_{K_{2}}(e_{2}) = 0$ et 
$\psi_{K_{2}}(e_{0}) = -1$. On en d\'eduit que le diviseur $E$ sur
$\M{P}_{\Z}^{2}$ associ\'e \`a $K_{2}$ par le th\'eor\`eme (\ref{construction_inverse}) 
est un
hyperplan coordonn\'e. On a notamment $\C{O}(E) \simeq \C{O}(1)$.

Plus g\'en\'eralement, consid\'erons $K_{d} \subset M_{\R} = \R^{d}$ le
simplexe standard d\'efini par les in\'equations~:
\[
x_{1} \geqslant 0, \dots, x_{d} \geqslant 0 \quad {\text{et}} \quad x_{1} +
\dots + x_{d} \leqslant 1.
\]
C'est un polytope convexe absolument simple. La vari\'et\'e torique
associ\'ee $\M{P}(K_{d})$ s'identifie avec l'espace projectif $\M{P}_{\Z}^{d}$.
Le diviseur $E$ sur $\M{P}_{\Z}^{d}$ associ\'e \`a $K_{d}$ par le th\'eor\`eme
(\ref{construction_inverse}) est un hyperplan coordonn\'e. On a donc $\C{O}(E)
\simeq \C{O}(1)$ (cf. \cite{20}, \S 2.4 et aussi \cite{11}, \S 1.4 et 1.5).
\end{expl}
Si l'on consid\`ere plusieurs polytopes convexes, on a le
r\'esultat suivant~:
\begin{thm}
\label{construction_inverse2}
Soient $K_{1},\dots,K_{m}$ des polytopes convexes de $M_{\R}$ \`a sommets dans
$M$. Posons $K = K_{1} + \dots + K_{m}$. On suppose que l'int\'erieur de $K$ est
non vide. Soient $\D$ l'\'eventail de $N_{\R}$ et $E$
le diviseur de Cartier horizontal $T$-invariant sur $\M{P}(\D)$
associ\'es \`a $K$.
Il existe des diviseurs de Cartier horizontaux $T$-invariants $E_{j}$ pour $(1 \leqslant
j \leqslant m)$ tels que $K_{E_{j}} = K_{j}$ et les faisceaux inversibles
$\C{O}(E_{j})$ soient engendr\'es par leurs sections globales. On a de plus $E
= E_{1} + \dots + E_{m}$.
\end{thm}
\demo\ Il suffit de remarquer que les fonctions $\psi_{K_{j}}$ sont lin\'eaires
par morceaux sur $\D$. On pose alors~:
\[
E_{j} = - \sum_{i = 1}^{r}\psi_{K_{j}}(u_{i})V(\t_{i}) \qquad (1 \leqslant j
\leqslant m).
\]
Comme les fonctions $\psi_{K_{j}}$ sont concaves, $E_{j}$ est un diviseur de
Cartier horizontal $T$-invariant et engendr\'e par ses sections globales; de plus,
$K_{\psi_{K_{j}}} = K_{j}$. Enfin, on a $\psi_{K} = \psi_{K_{1}} + \dots +
\psi_{K_{m}}$, et donc $E = E_{1} + \dots + E_{m}$ d'apr\`es (\ref{additivite_1}).
\medskip
\begin{rem}
\label{construction_inverse3}
Reprenons les hypoth\`eses et les notations du th\'eor\`eme
(\ref{construction_inverse2}). D'apr\`es la r\'esolution torique des
singularit\'es (cf. \cite{29}, th. 11, p. 94; voir aussi \cite{30}, th. 11, p.
273) il existe un raffinement $\D'$ de $\D$ tel que $\M{P}(\D')$ est projective et
lisse. Si l'on note $i_{\ast}: \M{P}(\D') \rightarrow \M{P}(\D)$ le
morphisme propre \'equivariant induit par l'inclusion $i: \D' \hookrightarrow \D$
comme \`a la d\'efinition (\ref{raffinement_1}) alors les diviseurs $E' =
(i_{\ast})^{\ast}(E)$, $E_{1}' = (i_{\ast})^{\ast}(E_{1}), \dots, E_{m}' =
(i_{\ast})^{\ast}(E_{m})$ sont tels que $K_{E'} = K$, $K_{E'_{1}} = K_{1},
\dots, K_{E_{m}'} = K_{m}$ d'apr\`es (\ref{intro_inverse}), et les faisceaux
inversibles $\C{O}(E')$, $\C{O}(E_{1}'), \dots, \C{O}(E_{m}')$ sont engendr\'es
par leurs sections globales (mais $E'$ n'est pas n\'ecessairement ample). De
plus $E' = E_{1}' + \dots + E_{m}'$.
\end{rem}
\bigskip

\subsection{Groupe de Picard et anneau de Chow d'une vari\'et\'e torique
projective
lisse}~

\subsubsection{Pr\'eliminaires}~

On d\'emontre une l\'eg\`ere g\'en\'eralisation d'un th\'eor\`eme d\^u \`a Gillet
et Soul\'e (\cite{14}, prop. 3.1.4).
\begin{thm} 
\label{iso_cellulaire}
Soit $X$ un sch\'ema
quasi-projectif sur $\SP$ admettant  une d\'ecomposition cellulaire,
c'est-\`a-dire tel qu'il existe une suite~: 
\[ X = X_{n} \supset X_{n-1}
\supset \dots \supset X_{0} \supset X_{-1} = \emptyset 
\]
de sous-sch\'ema
ferm\'es tels que $(X_{i} - X_{i-1})$ soit r\'eunion finie disjointe d'ouverts
affines $U_{i,j}$ isomorphes \`a $\M{A}^{i}_{\M{Z}}$. 
\begin{enumerate}
\item{Pour tout $l$ entier positif, on a des isomorphismes de groupes~: 
\[ CH^{l}(X) \stackrel{b}{\simeq} CH^{l}(X_{\Q})
\stackrel{b'}{\simeq} CH^{l}(X_{\M{C}}) \stackrel{cl}{\simeq}
H_{2n - 2l}(X(\M{C}),\Z), 
\] 
o\`u les isomorphismes $b$ et $b'$ sont d\'eduits des
morphismes de changement  de base $X_{\M{Q}} \rightarrow X$ et $X_{\M{C}}
\rightarrow X_{\M{Q}}$ et o\`u $cl$ est l'application cycle. De plus,
$CH^{l}(X)$ est un $\Z$-module libre de type fini,
et $H_{\op{impair}}(X(\M{C}),\M{Z}) = 0$.} 
\item{Pour tout couple $(r,s)$ d'entiers positifs tels que $r > s$, 
on a $CH^{r,s}(X) = 0$,
les groupes $CH^{r,s}(X)$ \'etant ceux d\'efinis dans (\cite{12},\S 8).}
\end{enumerate} 
\end{thm} 
\demo\
On commence par d\'emontrer l'assertion 2.
Du fait de l'invariance des groupes $CH^{r,s}$ par homotopie (cf. \cite{12}, th.
8.3), on a pour tout $r > s \geqslant 0$~:
\[
CH^{r,s}(\M{A}_{\Z}^{i}) = CH^{r,s}(\M{A}_{\Z}^{0}) = 0,
\]
et donc pour tout $i \in \{1,\dots,n\}$, 
\begin{equation}
\label{eq_intro_1}
CH^{r,s}(X_{i} - X_{i-1}) = \bigoplus_{j}CH^{r,s}(U_{i,j}) = 0.
\end{equation}
Par ailleurs, d'apr\`es la
longue suite exacte d'excision (\cite{12}, th. 8.1), on a~:
\begin{multline*}
\dots \longrightarrow CH^{r+1,s}(X_{i+1} - X_{i})
\stackrel{\partial}{\longrightarrow} CH^{r,s}(X_{i}) \longrightarrow \\
CH^{r,s}(X_{i+1}) \longrightarrow CH^{r,s}(X_{i+1} - X_{i}) 
\stackrel{\partial}{\longrightarrow} CH^{r-1,s}(X_{i}) \longrightarrow \dotsm
\end{multline*}
On tire de (\ref{eq_intro_1}) que pour tout $r > s \geqslant 0$,
\[
CH^{r,s}(X_{i+1}) \simeq CH^{r,s}(X_{i}).
\]
Par r\'ecurrence, on est donc ramen\'e \`a montrer le r\'esultat pour $X_{0}$,
ce qui est imm\'ediat car $X_{0}$ est r\'eunion finie disjointe de
$\M{A}_{\Z}^{0} = \op{Spec}\M{Z}$. 

On montre maintenant que $b$ est un
isomorphisme. Pour cela, on suit (\cite{14}, prop. 3.1.4) et on effectue une
r\'ecurrence sur la dimension. On remarque que le m\^eme raisonnement que
pr\'ec\'edemment appliqu\'e \`a $X_{\M{Q}}$ et \`a sa d\'ecomposition
cellulaire montre que $CH^{l+1,l}((X_{n} - X_{n-1})_{\M{Q}}) = 0$
pour tout entier positif $l$.
On a le diagramme commutatif suivant, form\'e de deux
suites exactes~: 
\medskip

\mbox{
\xymatrix{
  & CH(X_{n-1}) \ar[d]\ar[r] & CH(X_{n}) \ar[d] \ar[r]& CH(X_{n} - X_{n-1}) \ar[r]\ar[d] & 0 \\
(\ast) \quad 0 \ar[r] & CH((X_{n-1})_{\Q})\ar[r] & CH((X_{n})_{\Q})\ar[r] & CH((X_{n} - X_{n-1})_{\Q})\ar[r] & 0
}
}
\medskip

\noindent
o\`u les fl\`eches verticales sont induites par le morphisme $X_{\M{Q}}
\rightarrow X$ et ses restrictions et o\`u
la ligne $(\ast)$ est une suite exacte du fait 
de la longue suite exacte d'excision (\cite{12}, th. 8.1) et de 
la nullit\'e du groupe
$CH^{l+1,l}((X_{n} - X_{n-1})_{\Q})$. Puisque $X_{n-1}$ est encore un 
sch\'ema sur $\SP$ admettant une
d\'ecomposition cellulaire, la fl\`eche $CH(X_{n-1}) \rightarrow
CH((X_{n-1})_{\Q})$ est un isomorphisme par hypoth\`ese de r\'ecurrence. Comme
les groupes $CH$ sont invariants par homotopie (\cite{12}, th. 8.1), la
fl\`eche $CH(X_{n} - X_{n-1}) \rightarrow CH((X_{n} - X_{n-1})_{\Q})$ est
\'egalement un isomorphisme, et donc $CH(X_{n}) \rightarrow CH((X_{n})_{\Q})$
est un isomorphisme. De la suite exacte $(\ast)$, on tire imm\'ediatement que
$CH(X_{\Q})$, et donc $CH(X)$, est un $\Z$-module libre de type fini. Enfin
$b'$ et $cl$ sont des isomorphismes d'apr\`es 
(\cite{10}, \S 1.9.1 et \S 19.1.11), ce r\'esultat pouvant par
ailleurs se montrer par le m\^eme argument que pr\'ec\'edemment.
\medskip

Le r\'esultat suivant (voir \cite{9}, p. 144, ou encore \cite{11}, \S
5.2, la
d\'emonstration sur $\M{C}$ s'\'etendant imm\'ediatement \`a la situation sur
$\SP$) nous permet d'appliquer ce qui pr\'ec\`ede aux vari\'et\'es toriques
projectives lisses~:
\begin{thm}
Toute vari\'et\'e torique projective lisse $\P$ sur $\op{Spec}\M{Z}$ admet une
d\'ecomposition cellulaire.
\end{thm}
\medskip

On d\'eduit imm\'ediatement des deux th\'eor\`emes pr\'ec\'edents le corollaire
suivant~:
\begin{cor}
\label{absolu}
Soit $\P$ une vari\'et\'e torique projective lisse. 
\begin{enumerate}
\item{Pour tout entier $l$ positif, $CH^{l}(\P)$ est un $\Z$-module libre de type fini
et on a les isomorphismes de groupes
suivants~:
\[
CH^{l}(\P) \stackrel{b}{\simeq} CH^{l}(\M{P}(\D)_{\M{Q}})
\stackrel{b'}{\simeq} CH^{l}(\M{P}(\D)_{\M{C}}) \stackrel{cl}{\simeq}
H^{2l}(\M{P}(\D)(\M{C}), \Z),
\]
o\`u $b$ et $b'$ sont donn\'es par changement
de base de $\M{Z}$ \`a $\M{Q}$ et de $\M{Q}$ \`a $\M{C}$, et 
o\`u $cl$ est l'application cycle.}
\item{Pour tout $r$ entier strictement positif, on a~:
\[
CH^{r,r-1}(\P) = 0.
\]
}
\end{enumerate}
\end{cor}
\begin{rem}
Tous les isomorphismes de groupes introduits au corollaire (\ref{absolu}) 
sont compatibles avec l'intersection;
ce sont donc des isomorphismes pour les anneaux gradu\'es associ\'es aux
groupes consid\'er\'es.
\end{rem}

\subsubsection{Groupe de Picard}~

La proposition suivante caract\'erise les diviseurs de Cartier horizontaux
$T$-invariants principaux sur
une vari\'et\'e torique compl\`ete~:
\begin{prop}
\label{picard_0}
Soit $\P$ une vari\'et\'e torique compl\`ete et soit $D$ un diviseur de Cartier
horizontal $T$-invariant sur $\P$. Le diviseur $D$ est principal \ssi\ $\psi_{D}$ est
lin\'eaire sur tout $N_{\R}$.
\end{prop}
\demo\ En reprenant les notations de (\ref{inv_intro}), on \'ecrit $D$ sous la forme $D
= \sum_{i=1}^{r}a_{i}D_{i}$. Si $D$ est principal sur $\P$, alors $D_{\M{C}} =
\sum_{i=1}^{r}a_{i}(D_{i})_{\M{C}}$ est principal sur $\PP$ et d'apr\`es 
(\cite{20}, prop. 2.4) il existe donc $m \in M$ tel que $m(u_{i}) = - a_{i}$
pour tout $1 \leqslant i \leqslant r$. On en d\'eduit que $D =
\op{div}(\chi^{m})$. La r\'eciproque est imm\'ediate.
\medskip

\begin{prop}
\label{picard}
Soit $\P$ une vari\'et\'e torique lisse projective. On note
$\t_{1},\dots,$ $\t_{r}$ les \'el\'ements de $\D(1)$ et $u_{1},\dots,u_{r}$ leur
g\'en\'erateur dans $N$. On a la suite exacte~:
\[
0 \longrightarrow M \stackrel{\iota}{\longrightarrow} \bigoplus_{i =1}^{r}\Z V(\t_{i})
\stackrel{s}{\longrightarrow} \op{Pic}(\P) \longrightarrow 0,
\]
la fl\`eche $\iota$ \'etant donn\'ee par $\iota: m \rightarrow
\sum_{i=1}^{r}<m,u_{i}> V(\t_{i})$ et la fl\`eche $s$ d\'esignant la surjection
canonique sur $\op{Pic}(\P)$ vu comme groupe quotient par l'\'equivalence 
lin\'eaire. Donc $\op{Pic}(\P)$ est un $\Z$-module libre de rang $(\#\D(1) -
d)$.
\end{prop}
\demo\ C'est une cons\'equence directe de (\cite{20}, cor. 2.5) et et 
de l'isomorphisme $b: CH^{1}(\P) \simeq CH^{1}(\P_{\Q})$ donn\'e par le
corollaire (\ref{absolu}).
On peut \'egalement consulter (\cite{11}, \S 3.4).
\medskip

\subsubsection{Anneau de Chow}~

On d\'ecrit maintenant un th\'eor\`eme de Jurkiewicz et Danilov (\cite{4}, th.
10.8 et rem. 10.9) donnant la structure de l'anneau de Chow (et donc de
l'anneau d'homologie) d'une vari\'et\'e torique projective lisse~:
\begin{thm}
\label{anneau_chow}
Soient $\P$ une vari\'et\'e torique projective lisse,
$\t_{1},\dots,$ $\t_{r}$ les \'el\'ements de $\D(1)$ et $u_{1},\dots,u_{r}$
leurs
g\'en\'erateurs respectifs dans $N$. Consid\'erons l'anneau
de polyn\^omes en $r$ ind\'etermin\'ees~:
\[
\C{S} = \Z\left[t_{\t_{i}}\right]_{1 \leqslant i \leqslant r}.
\]
Soient $\C{I}$ l'id\'eal de $\C{S}$ engendr\'e par l'ensemble~:
\[
\{t_{\rho_{1}}t_{\rho_{2}}\dots t_{\rho_{s}}: \quad \rho_{1},\dots,\rho_{s} \in
\D(1)\text{ deux \`a deux distincts et }\rho_{1} + \dots + \rho_{s} \notin \D \}
\]
et $\C{J}$ l'id\'eal de $\C{S}$ engendr\'e par~:
\[
\left\{ 
\sum_{i=1}^{r}<m,u_{i}> t_{\t_{i}}, \quad m \in M 
\right\}.
\]
Soit $[~]: \C{S} \rightarrow CH^{\ast}(\P)$ le morphisme d'anneau d\'efini
par~:
\[
[t_{\t}] := [V(\t)] \in CH^{\ast}(\P),
\]
pour tout $\t \in \D(1)$. 
On a $\op{Ker}[~] = (\C{I} + \C{J})$ et le morphisme~:
\[
\C{S}/(\C{I} + \C{J}) \stackrel{[~]}{\longrightarrow} CH^{\ast}(\P),
\]
est un isomorphisme d'anneaux gradu\'es.
\end{thm}
\demo\ C'est une cons\'equence de (\cite{4}, th. 10.8) ou encore 
(\cite{20}, p. 134). On peut \'egalement 
consulter (\cite{11}, p. 106).
\medskip

\section{Vari\'et\'es toriques complexes}~

\subsection{Vari\'et\'e \`a coin associ\'ee \`a une vari\'et\'e torique}~

Dans ce paragraphe, on suit essentiellement (\cite{11}, \S 4). On peut
\'egalement consulter (\cite{20}, prop. 1.8).

Soit $\D$ un \'eventail de $N$. Pour tout \co\ $\si$ de $\D$ on note~:
\[
(U_{\si})_{\geqslant} = \op{Hom}_{\op{sg}}(\si^{\ast}\cap M, \R^{+})
\]
l'ensemble des morphismes de semi-groupe avec \'el\'ement neutre 
de $\si^{\ast}\cap M$ vers $(\R^{+},\times)$.
Du fait de l'inclusion $\R^{+} \subset \M{C}$, l'ensemble
$(U_{\si})_{\geqslant}$ peut \^etre vu comme sous-ensemble ferm\'e de
$U_{\si}(\M{C})$. En effet, on a~:
\begin{multline*}
(U_{\si})_{\geqslant} = \op{Hom}_{\op{sg}}(\si^{\ast}\cap M, \R^{+})\\
\subset
\op{Hom}_{\op{sg}}(\si^{\ast}\cap M, \M{C}) =
\op{Hom}_{\text{$\M{C}$-alg\`ebre}} \left(\M{C}[\si^{\ast}\cap M],\M{C}\right) 
= U_{\si}(\M{C}).
\end{multline*}
L'application module $|~|:\M{C} \rightarrow \R^{+}$ induit par composition 
une r\'etraction~:
\[
U_{\si}(\M{C}) = \op{Hom}_{\op{sg}}(\si^{\ast}\cap M, \M{C})
\longrightarrow \op{Hom}_{\op{sg}}(\si^{\ast}\cap M, \R^{+}) =
(U_{\si})_{\geqslant} 
\]
que l'on note encore $|~|$. 
Si $\si$ et $\si'$ sont deux \'el\'ements de $\D$, alors $(U_{\si})_{\geqslant}
\cap (U_{\si'})_{\geqslant} = (U_{\si \cap \si'})_{\geqslant}$. Il existe 
donc une partie ferm\'ee $\PPP$ de $\P(\M{C})$ telle que pour tout $\si \in \D$
on ait
$\PPP \cap U_{\si}(\M{C}) = (U_{\si})_{\geqslant}$.

L'application $|~|$ s'\'etend en une r\'etraction continue~:
\[
|~|: \P(\M{C}) \longrightarrow \PPP
\]
\begin{prop}
Si $\D$ est r\'egulier, $\PPP$ est une sous-vari\'et\'e \`a coins 
analytique r\'eelle de $\P(\M{C})$, de
dimension r\'eelle $d$.
\end{prop}
\demo\ C'est une cons\'equence directe des d\'efinitions (voir par 
exemple \cite{20}, prop. 1.8).
\medskip
\begin{expl}
On note $\TTT = (U_{\{0\}})_{\geqslant} = \op{Hom}_{\op{sg}}(M,\R^{+})$. On
v\'erifie que $\TTT$ est dense dans $\PPP$.
\end{expl}

Soient $\TT$ le tore analytique complexe associ\'e \`a $T$
et $\CT$ le {\it tore compact\/}~:
\[
\CT = \op{Hom}(M,\C{S}_{1}) \subset \op{Hom}(M,\M{C}^{\ast}) = \TT,
\]
o\`u $\C{S}_{1} \subset \M{C}$ est le cercle unit\'e. 
C'est le sous-groupe compact maximal de $\TT$; en fait on a~:
\[
\CT = \op{Hom}(M,\C{S}_{1}) = N\otimes_{\Z}\C{S}_{1} \cong
{\C{S}_{1}}^{d} \subset \TT = \op{Hom}(M,\M{C}^{\ast}) = 
N\otimes_{\Z}\M{C}^{\ast} \cong (\M{C}^{\ast})^{d}.
\]
L'action naturelle de $\TT$ sur $\P(\M{C})$ induit une action de $\CT$ sur
$\P(\M{C})$.
La d\'ecomposition polaire $\M{C}^{\ast} = \C{S}_{1} \times
\R^{+\ast}$ et l'isomorphisme $\log: \R^{+} \rightarrow \R$ d\'eterminent
un isomorphisme $\CT$-\'equivariant~:
\[
\log : \TT = \CT \times \op{Hom}_{\op{sg}}(M,\R^{+}) \longrightarrow \CT \times
\op{Hom}(M,\R) = \CT \times N_{\R}.
\]
On note $\op{pr_{2}} : \CT \times N_{\R} \rightarrow N_{\R}$ la seconde
projection. Le diagramme suivant commute~:
\begin{center}
\mbox{
\xymatrix{ \TT \ar[rr]^(0.4){\log} \ar[d]_{|~|} & &\CT \times N_{\R} \ar[d]^{\op{pr_{2}}}
\\
\TTT \ar[rr]^{\log} & &N_{\R}
}}
\end{center}
et les fl\`eches horizontales sont des isomorphismes de groupes de Lie r\'eels.\\
On note $\op{orb}$ la surjection canonique~:
\[
\op{orb}: \P(\M{C}) \longrightarrow \P(\M{C})/\CT.
\]
On remarque qu'on peut identifier~:
\[
\TTT = \op{Hom}_{\op{sg}}(M,\R^{+}) = \TT/\CT.
\]
Plus g\'en\'eralement, on a le r\'esultat suivant~:
\begin{thm}
Il existe un unique hom\'eomorphisme $\kappa : \PPP \rightarrow \P(\M{C})/\CT$ tel 
que le diagramme suivant commute~:
\begin{center}
\mbox{
\xymatrix{
& \P(\M{C}) \ar[dl]_{|~|} \ar[dr]^{\op{orb}} & \\
\PPP \ar[rr]^{\kappa}& & \P(\M{C})/\CT}
}
\end{center}
De plus, la restriction de $\kappa$ \`a $\TTT$ co\"\i ncide avec
l'identification ci-dessus.
\end{thm}
\demo\ Voir (\cite{11}, p. 79) et (\cite{20}, prop. 1.8).
\bigskip

\subsection{Un recouvrement canonique de $\P(\M{C})$}~

\begin{defn}{\bf (Batyrev et Tschinkel).}
Pour tout $\si \in \D$, on d\'efinit $C_{\si} \subset \P(\M{C})$ 
de la fa\c con suivante~:
\[
C_{\si} = \{x \in \P(\M{C}): 
\forall m \in \E{S}_{\si}= \si^{\ast}\cap M, 
\quad \chi^{m} \text{ est r\'egulier en $x$ et }|\chi^{m}(x)|
\leqslant 1\}.
\]
\end{defn}
\begin{prop}
Pour tout $\si \in \D$, on a $C_{\si} \subset U_{\si}(\M{C})$ et
$C_{\si}$ est compact. De plus, si $\t, \t' \in \D$, alors~: 
\[
C_{\t} \cap C_{\t'} = C_{\t \cap \t'}
\]
\end{prop}
\demo\ Soit $x \in C_{\si}$, la fonction $\chi^{m}$ est d\'efinie en $x$ pour
tout $m \in \E{S}_{\si} = \si^{\ast}\cap M$, et donc $x \in U_{\si}(\M{C})$.

Soit $\{m_{1},\dots,m_{q}\}$ une famille g\'en\'eratrice du semi-groupe
$\E{S}_{\si}$. On dispose de l'immersion ferm\'ee~: 
\begin{alignat*}{3}
\varphi : \; U_{\si}&(\M{C}) & &\longrightarrow & \;&\M{C}^{q} \\
&x & &\longmapsto & &(\chi^{m_{1}}(x),\dots,\chi^{m_{q}}(x)).
\end{alignat*}
On remarque que $\varphi (C_{\si}) = \varphi(U_{\si}(\M{C})) \cap B(0,1)^{q}$. 
Comme $\varphi (U_{\si}(\M{C}))$ est un ferm\'e analytique, on conclut
que $\varphi(C_{\si})$, et donc $C_{\si}$, sont compacts.

Enfin la derni\`ere assertion r\'esulte de l'identit\'e $\t^{\ast} +
{\t'}^{\ast} = (\t \cap \t')^{\ast}$ (Voir \cite{20}, th. A.1).
\medskip
\begin{prop}
\label{recouvrement}
Si $\D$ est complet, alors
les compacts $C_{\si}$ forment un
recouvrement de $\P(\M{C})$ lorsque $\si$ parcoure $\D_{\op{max}}$. 
\end{prop}
\demo\ On suit ici \cite{2}. Puisque $\TT$ est dense dans $\P(\M{C})$, il suffit de 
d\'emontrer que les $C_{\si}$, lorsque $\si$ parcourt
$\D_{\text{max}}$, recouvrent $\TT$. Soit $x \in \TT$; puisque $\D$ est
complet, il existe $\si \in \D_{\text{max}}$ tel que $-\log|x| \in \si$.
Pour tout $m \in \E{S}_{\si}$, on a~:
\[
<m,-\log|x|> = -\log|\chi^{m}(x)| \geqslant 0.
\]
On en conclut que $|\chi^{m}(x)| \leqslant 1$ pour tout $m \in \E{S}_{\si}$, 
c'est-\`a-dire que $x \in C_{\si}$.
\medskip

Les compacts de la forme $C_{\si}$ sont globalement invariants
sous l'action de $\CT$. Pour mieux comprendre la structure des $C_{\si}$, 
on peut \'etudier l'image de $C_{\si} \cap \TT$ dans $N_{\R}$ par 
l'application $\log |.|$. C'est 
l'objet de la proposition suivante~:
\begin{prop}
\label{application_log}
Pour tout $\si \in \D$, on note $\stackrel{\circ}{C}_{\si} = C_{\si} \cap \TT$.
On a alors~: 
\[
-\log |\stackrel{\circ}{C}_{\si}| = \si.
\]
\end{prop}
\demo\ Soit $x \in \TT$. Pour tout $m \in M$, on a $<m,-\log|x|> =
-\log|\chi^{m}(x)|$, et donc $x \in \stackrel{\circ}{C}_{\si}$ \ssi~: 
\[
<m,-\log|x|> \geqslant 0, \qquad \forall m \in \E{S}_{\si},
\]
ce qui est \'equivalent \`a~:
\[
- \log|x| \in (\si^{\ast})^{\ast} = \si .
\]
\medskip

L'\'etude des $C_{\si}$ est donc ramen\'ee \`a l'\'etude de $\D$ ``plong\'e''
dans $\PPP$.
\begin{expl}
On a repr\'esent\'e ici $\M{P}^{2}(\M{C})/\CT$~:
\bigskip
\begin{center}
\begin{picture}(0,0)%
\includegraphics{figure1.pstex}%
\end{picture}%
\setlength{\unitlength}{0.00083300in}%
\begingroup\makeatletter\ifx\SetFigFont\undefined
\def\x#1#2#3#4#5#6#7\relax{\def\x{#1#2#3#4#5#6}}%
\expandafter\x\fmtname xxxxxx\relax \def\y{splain}%
\ifx\x\y   
\gdef\SetFigFont#1#2#3{%
  \ifnum #1<17\tiny\else \ifnum #1<20\small\else
  \ifnum #1<24\normalsize\else \ifnum #1<29\large\else
  \ifnum #1<34\Large\else \ifnum #1<41\LARGE\else
     \huge\fi\fi\fi\fi\fi\fi
  \csname #3\endcsname}%
\else
\gdef\SetFigFont#1#2#3{\begingroup
  \count@#1\relax \ifnum 25<\count@\count@25\fi
  \def\x{\endgroup\@setsize\SetFigFont{#2pt}}%
  \expandafter\x
    \csname \romannumeral\the\count@ pt\expandafter\endcsname
    \csname @\romannumeral\the\count@ pt\endcsname
  \csname #3\endcsname}%
\fi
\fi\endgroup
\begin{picture}(3345,3356)(7178,-8795)
\put(9826,-5911){\makebox(0,0)[lb]{\smash{\SetFigFont{12}{14.4}{rm}$C_{\si_{1}}$}}}
\put(8686,-6646){\makebox(0,0)[lb]{\smash{\SetFigFont{12}{14.4}{rm}$C_{\si_{2} \cap \si_{3}}$}}}
\put(9901,-7531){\makebox(0,0)[lb]{\smash{\SetFigFont{12}{14.4}{rm}$C_{\si_{3}}$}}}
\put(8086,-5911){\makebox(0,0)[lb]{\smash{\SetFigFont{12}{14.4}{rm}$C_{\si_{2}}$}}}
\put(8851,-6061){\makebox(0,0)[lb]{\smash{\SetFigFont{12}{14.4}{rm}$C_{\si_{1} \cap \si_{2}}$}}}
\put(9691,-6706){\makebox(0,0)[lb]{\smash{\SetFigFont{12}{14.4}{rm}$C_{\si_{1} \cap \si_{3}}$}}}
\end{picture}

\end{center}
\bigskip
\end{expl}
\begin{prop}
Soit $p$ un entier strictement positif. Le morphisme $[p]: \P(\M{C})\rightarrow
\P(\M{C})$ envoie $\PPP$ dans lui-m\^eme. 
Pour tout $\si \in D$, le compact $C_{\si}$ est laiss\'e stable par $[p]$.
Enfin le
diagramme suivant commute~:
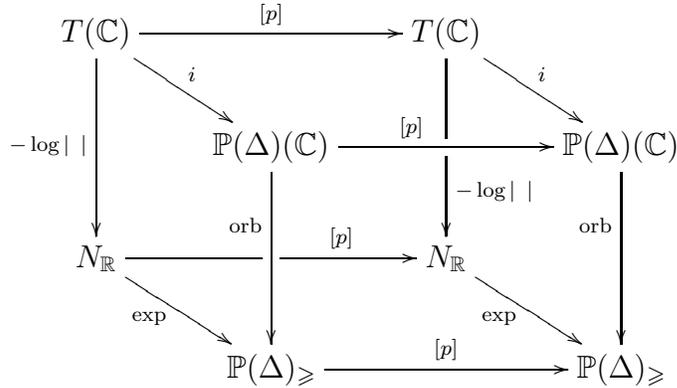
\begin{center}
\mbox{
\xymatrix{
\TT \ar[rr]^{[p]} \ar[dd]_{-\log|~|} \ar[dr]^{i} & & \TT \ar'[d]
[dd]^(0.4){-\log|~|}
\ar[dr]^{i} \\
& \P(\M{C}) \ar[dd]_(0.35){\op{orb}} \ar[rr]^(0.4){[p]} & & \P(\M{C}) \ar[dd]_(0.35){\op{orb}} \\
N_{\R} \ar'[r] [rr]^(0.4){[p]} \ar[rd]_(0.4){\exp} & & N_{\R} \ar[dr]_(0.4){\exp} \\
& \PPP \ar[rr]^{[p]} & & \PPP }
}
\end{center}
\medskip
o\`u $\exp$ est l'inverse du morphisme $\log: \TTT \rightarrow N_{\R}$, 
et $[p]: N_{\R} \rightarrow N_{\R}$ d\'esigne la multiplication par $p$.
\end{prop}
\demo\ La premi\`ere partie de la proposition est une cons\'equence de
l'\'egalit\'e $[p](x) = x^{p}$ pour tout $x \in \TT$ (on peut \'egalement
consulter \cite{11}, p. 80). \\
Soient $\si \in \D$ et $x \in C_{\si}$, pour tout $m \in \E{S}_{\si}$, on a 
$|\chi^{m}([p](x))| = |\chi^{m}(x)|^{p} \leqslant 1$, c'est-\`a-dire $[p](x) \in
C_{\si}$; et donc $C_{\si}$ est laiss\'e stable par $[p]$. \\
Enfin la commutativit\'e du diagramme r\'esulte des propositions \'enonc\'ees 
au paragraphe pr\'ec\'edent.
\medskip

\begin{defn}
Pour tout $\si \in \D$, on pose~:
\[
C_{\si}^{\op{int}} = \{ x \in C_{\si}\,: \quad |\chi^{m}(x)| < 1, \quad \forall
m \in (\si^{\ast} - \si^{\perp}) \cap M \}.
\]
\end{defn}

La proposition suivante rassemble diverses propri\'et\'es des ensembles
$C_{\si}$ et $C_{\si}^{\op{int}}$
\begin{prop}~
\begin{enumerate}
\item{Pour tout $\si \in \D$, on a $- \log | C_{\si}^{\op{int}} \cap \TT | = 
\stackrel{\circ}{\si}$,
o\`u $\stackrel{\circ}{\si}$ est l'int\'erieur relatif du c\^one $\si$.}
\item{Les ensembles $C_{\t}^{\op{int}}$ sont deux \`a deux disjoints. Pour tout
$\si \in \D$, on a~:
\begin{equation}
\label{an_eq_1}
C_{\si} = \bigcup_{\t < \si}C_{\si}^{\op{int}}.
\end{equation}
}
\item{Si l'\'eventail $\D$ est complet, alors les $C_{\t}^{\op{int}}$ pour $\t$
parcourant $\D$ forment une partition de $\P(\M{C})$.}
\end{enumerate}
\end{prop}
\demo\ D\'emontrons tout d'abord l'assertion (1). Soit $x \in \TT \cap
C_{\si}^{\op{int}}$. Pour tout $m \in (\si^{\ast} - \si^{\perp}) \cap M$, on a
$<m,-\log|x|> \; > 0$, ce qui \'equivaut \`a \'ecrire que $- \log |x| \in
\stackrel{\circ}{\si}$. On en d\'eduit que $- \log |C_{\si}^{\op{int}} \cap
\TT| = \stackrel{\circ}{\si}$.

On s'int\'eresse maintenant \`a l'assertion (2). On reprend ici les notations de la
remarque (\ref{rem_decomposition1}). Soient $\t$ et $\t'$ deux c\^ones
distincts de $\D$. D'apr\`es l'assertion (1), on a $- \log
|C_{\t}^{\op{int}}\cap
C_{\t'}^{\op{int}}\cap \TT | = \stackrel{\circ}{\t} \cap
\stackrel{\circ}{\t}^{\lower7pt\hbox{$\scriptstyle '$}}
= \emptyset$, ce dont on d\'eduit que les ensembles $C_{\t}^{\op{int}}$ et $C_{\t'}^{\op{int}}$ sont
disjoints sur $\TT$. En remarquant que pour tout $\si \in \D$, on a~:
\[
C_{\t}^{\op{int}} \cap V(\si)(\M{C}) = C_{\ov{\t}}^{\op{int}}, \qquad 
(\text{resp.} \quad C_{\t'}^{\op{int}} \cap V(\si)(\M{C}) = C_{\ov{\t}'}^{\op{int}})
\]
o\`u $\ov{\t} \in \D(\si)$ (resp. $\ov{\t}' \in \D(\si)$) et o\`u $C_{\ov{\t}}^{\op{int}}
\subset \P(\D(\si))(\M{C}) = V(\si)(\M{C})$ 
(resp. $C_{\ov{\t}'}^{\op{int}} \subset V(\si)(\M{C})$), 
on obtient la relation~:
\[
C_{\t}^{\op{int}} \cap C_{\t'}^{\op{int}} \cap V(\si)(\M{C}) = C_{\ov{\t}}^{\op{int}}
\cap C_{\ov{\t}'}^{\op{int}}.
\]
Pa r\'ecurrence, on d\'eduit de cela, de la d\'ecomposition en tores disjoints 
donn\'ee \`a la proposition (\ref{decomposition1}) et de la
remarque (\ref{rem_decomposition1}) que $C_{\t}^{\op{int}}$ et
$C_{\t'}^{\op{int}}$ sont disjoints sur $\P(\M{C})$. La d\'ecomposition
(\ref{an_eq_1}) se montre par r\'ecurrence de fa\c con similaire en utilisant
la relation $\si = \bigcup_{\t < \si}\stackrel{\circ}{\t}$ et la proposition
(\ref{application_log}).

Il reste \`a prouver l'assertion (3). D'apr\`es ce qui pr\'ec\`ede,
il suffit de montrer
que les $C_{\t}^{\op{int}}$ pour $\t \in \D$ forment un recouvrement de
$\P(\M{C})$; or ceci se d\'eduit directement de la d\'ecomposition
(\ref{an_eq_1}) et de la proposition (\ref{recouvrement}).
\medskip

Dans le cas o\`u $\D$ est complet et r\'egulier, on peut pr\'eciser la structure des 
ensembles $C_{\si}$ et $C_{\si}^{\op{int}}$~:
\begin{prop}
\label{dissection}
Soient $\D$ un \'eventail complet et r\'egulier et $\si$ un \'el\'ement de $\D$.
\begin{enumerate}
\item{
L'ensemble $C_{\si}^{\op{int}}$ est une sous-vari\'et\'e
analytique r\'eelle lisse de dimension r\'eelle $d + \op{dim} \si$ de
$\P(\M{C})$.}
\item{le compact $C_{\si}$ est une sous-vari\'et\'e r\'eelle \`a coins
de dimension r\'eelle $d + \op{dim} \si$ de $\P(\M{C})$, et son bord 
$\partial C_{\si}$ est une 
vari\'et\'e \`a coins v\'erifiant la relation~:
\[
\partial C_{\si} = C_{\si} - C_{\si}^{\op{int}} =
\bigcup_{\substack{
\t < \si \\ \t \not= \si}}
C_{\si}^{\op{int}}\; .
\]
}
\end{enumerate}
\end{prop}
\demo\ On pose $q = \dim \si$. Soit $\t \in \D_{\op{max}}$ tel que $\si < \t$.
D'apr\`es la proposition (\ref{lissite}), on peut trouver $\{m_{1}, \dots, m_{d}\}$ une base de $M$
telle que $\{m_{1},\dots,m_{q}\}$ (resp. $\{m_{1}, \dots, m_{d}\}$) soit 
une famille g\'en\'eratrice du semi-groupe $\E{S}_{\si}$ (resp. $\E{S}_{\t}$).

Dans la carte affine $\varphi: U_{\t}(\M{C}) \rightarrow \M{C}^{d}$ donn\'ee
par $\varphi(x) = (\chi^{m_{1}}(x),\dots,\chi^{m_{d}}(x))$, les ensembles
$C_{\si}$ et $C_{\si}^{\op{int}}$ sont d\'efinis par les conditions~:
\begin{align*}
C_{\si} &= \{x \in \M{C}^{d}: \quad |x_{1}| \leqslant 1, \dots, |x_{d-q}|
\leqslant 1, 
|x_{d-q+1}| = 1, \dots, |x_{d}| = 1\} \\
\intertext{et}
C_{\si}^{\op{int}} &= \{x \in \M{C}^{d}: \quad |x_{1}| < 1, \dots, |x_{d-q}|
< 1, |x_{d-q+1}| =
1, \dots, |x_{d}| = 1\}.
\end{align*}
On en d\'eduit directement les assertions \'enonc\'ees.
\bigskip

\subsection{M\'etriques canoniques sur les faisceaux inversibles \'equivariants
au-dessus de $\P(\M{C})$}~

Dans toute cette section, $\D$ d\'esigne un \'eventail 
complet de $N$. 
Pour tout diviseur de Cartier $D$ horizontal et $T$-invariant sur $\P$, on 
construit de mani\`ere canonique une m\'etrique sur le faisceau inversible
$\C{O}(D)(\M{C})$.

Plusieurs constructions \'equivalentes sont indiqu\'ees.

\subsubsection{Construction de Batyrev et Tschinkel}~

On pr\'esente ici, sous une forme
l\'eg\`erement diff\'erente de 
(\cite{2}, \S 2.1), une construction due \`a Batyrev et Tschinkel.
\begin{prop_defn}
\label{BT_construction}
Soit $s$ une section holomorphe de $\C{O}(D)(\M{C})$ au-dessus d'un ouvert
$\Omega \subset \P(\M{C})$. Pour tout $x \in \Omega$, soit $\si \in \D$ tel que
$x \in C_{\si} \subset U_{\si}(\M{C})$,
et $m_{D,\si} \in M$ la restriction de $\psi_{D}$ \`a $\si$. Le quotient~:
\begin{equation}
\label{BT_definition}
\| s(x) \|_{\op{BT}} = \left| \frac{s}{\chi^{m_{D,\si}}}(x) \right|
\end{equation}
est bien d\'efini et est appel\'e {\it norme de $s$ au point $x$ au sens de
Batyrev-Tschinkel\/}. Cette norme d\'efinit une m\'etrique continue $\CT$-invariante
sur $\C{O}(D)(\M{C})$, que l'on note $\| . \|_{\op{BT}}$.
\end{prop_defn}
\demo\ L'existence du quotient est une cons\'equence directe du lemme (\ref{locale}).
Enfin, le second membre de (\ref{BT_definition}) est ind\'ependant du choix de
$\sigma$ car~:
\[
|\chi^{m_{D,\si}}(x)| = | \chi^{m_{D,\si'}}(x)| 
\]
pour tout $x \in C_{\si} \cap C_{\si'}$ du fait de la continuit\'e de la
fonction support $\psi_{D}$ sur $N_{\R}$, et le m\^eme argument montre que 
$\| . \|_{\op{BT}}$ est bien continue sur $\P(\M{C})$.
\medskip

On donne dans la proposition suivante quelques propri\'et\'es de la
construction de Batyrev-Tschinkel~:
\begin{prop}~
\label{BT_fonct}
\begin{enumerate}
\item{
{\rm (Multiplicativit\'e).}
Soient $\D$ un \'eventail complet dans $N_{\R}$ et $D_{1}$, $D_{2}$ deux diviseurs de
Cartier horizontaux et $T$-invariants sur $\P$. L'isomorphisme~:
\[
\C{O}(D_{1}) \otimes \C{O}(D_{2}) \simeq \C{O}(D_{1} + D_{2}), 
\]
est compatible aux m\'etriques de Batyrev-Tschinkel.}
\item{
{\rm (Fonctorialit\'e).}
Soient $\varphi: \D_{1} \rightarrow \D_{2}$ un morphisme d'\'eventails complets
et
$\varphi_{\ast}: \M{P}(\D_{1}) \rightarrow \M{P}(\D_{2})$ le morphisme de
vari\'et\'es toriques associ\'e. Soient \'egalement $D_{2}$ un diviseur de Cartier
horizontal $T$-invariant sur $\M{P}(\D_{2})$ et $\psi_{2}$ sa fonction
support, et notons $D_{1} = (\varphi_{\ast})^{\ast}D_{2}$ le diviseur de Cartier
$T$-invariant sur $\M{P}(\D_{1})$ dont la fonction support est donn\'ee 
par $\psi_{1} = \psi_{2}\circ \varphi$.
L'isomorphisme~:
\[
\C{O}(D_{1}) \simeq (\varphi_{\ast})^{\ast}\C{O}(D_{2}),
\]
est une isom\'etrie lorsque $\C{O}(D_{1})$ et $\C{O}(D_{2})$ sont munis de leur
m\'etrique de Batyrev-Tschinkel.}
\end{enumerate}
\end{prop}
\demo\
On d\'emontre tout d'abord l'assertion (1). Soient $s_{1}$ et $s_{2}$ des
sections holomorphes, sur un ouvert $\Omega \subset \P(\M{C})$, des faisceaux
$\C{O}(D_{1})$ et $\C{O}(D_{2})$ respectivement.
Pour tout $x \in \Omega$, soit $\si \in \D$ tel que $x \in C_{\si}$.
On a~:
\[
\|s_{1}(x)\|_{\op{BT}}\cdot \|s_{2}(x)\|_{\op{BT}} = 
\left|
\frac{s_{1}\otimes s_{2}(x)}{\chi^{m_{D_{1},\si} + m_{D_{2},\si}}(x)}\right|
= 
\left|
\frac{s_{1}\otimes s_{2}(x)}{\chi^{m_{D_{1}+D_{2}},\si}(x)}\right|
= \|s_{1} \otimes s_{2}(x)\|_{\op{BT}}, 
\]
ce qui \'etablit l'\'enonc\'e recherch\'e.

On s'int\'eresse maintenant \`a l'assertion (2). Soit $s_{2}$ une section
holomorphe du faisceau $\C{O}(D_{2})$ sur un ouvert $\Omega \subset 
\M{P}(\D_{2})(\M{C})$
et notons $s_{1} = s_{2} \circ \varphi_{\ast}$ la section holomorphe de
$\C{O}(D_{1})$ au-dessus de $(\varphi_{\ast})^{-1}(\Omega)$ obtenue par image
r\'eciproque de $s_{2}$ par $(\varphi_{\ast})^{\ast}$.

Pour tout $x \in (\varphi_{\ast})^{-1}(\Omega)$, soit $\si_{1} \in \D_{1}$ tel
que $x \in C_{\si_{1}}$ et choisissons $\si_{2} \in \D_{2}$ tel que
$\varphi(\si_{1}) \subset \si_{2}$. D'apr\`es la proposition
(\ref{intro_inverse}), on sait que $\psi_{D_{1}} = \psi_{D_{2}}\circ \varphi$,
et donc que $m_{D_{1},\si} = {}^{t}\varphi(m_{D_{2},\si})$.

On d\'eduit de la proposition (\ref{application_log}) et du fait que
$\varphi_{\ast}$ est $T$-\'equivariant, l'inclusion~:
\[
\varphi_{\ast}(T_{1}(\M{C}) \cap C_{\si_{1}}) 
\subset T_{2}(\M{C}) \cap C_{\si_{2}},
\]
ce qui, comme $\varphi_{\ast}$ est continue, montre que~:
\[
\varphi_{\ast}(C_{\si_{1}}) \subset C_{\si_{2}}.
\]
On peut donc \'ecrire~:
\[
\|s_{1}(x)\|_{\op{BT}} = 
\left|\frac{s_{1}(x)}{\chi^{m_{D_{1},\si}}(x)}\right|
= 
\left|
\frac{s_{2}\circ \varphi_{\ast}(x)}{\chi^{{}^t\varphi(m_{D_{2},\si})}(x)}\right|
= 
\left|
\frac{s_{2}\circ \varphi_{\ast}(x)}{\chi^{m_{D_{2},\si}}\circ
\varphi_{\ast}(x)}\right|
= 
\|s_{2}(\varphi_{\ast}(x))\|_{\op{BT}},
\]
ce qui termine la d\'emonstration.
\medskip

\subsubsection{Construction d'apr\`es Zhang}~

On suit dans ce paragraphe (\cite{21}, th. 2.2). Soient $p$ un entier 
sup\'erieur
ou \'egal \`a $2$ et $D$ un diviseur de Cartier horizontal 
$T$-invariant sur $\P$. Comme $[p]$ est l'endomorphisme
\'equivariant de $\P$ associ\'e \`a  l'endomorphisme de $\D$ 
d\'efini par la multiplication par $p$, on a d'apr\`es (\ref{intro_inverse}) 
l'\'egalit\'e
des diviseurs de Cartier~:
\[
[p]^{\ast}D = p\,D,
\]
et donc un isomorphisme de faisceaux sur $\P$~:
\[
\Phi_{D,p}: \C{O}(D)^{\otimes p} \simeq [p]^{\ast}\C{O}(D).
\]

Comme $\P(\M{C})$ est compacte, on peut
munir $\C{O}(D)(\M{C})$ d'une m\'etrique continue que l'on notera $\|.\|_{0}$.
On d\'efinit alors par r\'ecurrence une suite de m\'etriques $(\|.\|_{n})_{n
\in \M{N}}$ sur $\C{O}(D)(\M{C})$ de la fa\c con suivante $(n \geqslant 1)$~:
\[
\| .\|_{n} = \left(\Phi^{\ast}_{D,p}\,[p]^{\ast}\,\|.\|_{n-1}\right)^{1/p} .
\]
On a alors le th\'eor\`eme suivant~:
\begin{thm}~
\label{zhang}
\begin{enumerate}
\item{Les m\'etriques $\| .\|_{n}$ convergent uniform\'ement vers une
m\'etrique $\|.\|_{\op{Zh},p}$ sur $\P(\M{C})$ (i.e $\log
\frac{\|.\|_{n}}{\|.\|_{0}}$ converge uniform\'ement sur $\P(\M{C})$ vers $\log
\frac{\|.\|_{\op{Zh},p}}{\|.\|_{0}}$).}
\item{La m\'etrique $\|.\|_{\op{Zh},p}$ est l'unique m\'etrique continue sur
$\C{O}(D)(\M{C})$ telle que~:
\[
\|.\|_{\op{Zh},p} = \left(\Phi_{D,p}^{\ast}\,[p]^{\ast}\,\|.\|_{\op{Zh},p}
\right)^{1/p}.
\]
}
\end{enumerate}
\end{thm}
\demo\ C'est une cons\'equence directe de (\cite{21}, th. 2.2).
\medskip
\begin{rem}
Zhang raisonne pour des vari\'et\'es projectives, mais son argument ne
n\'ecessite en fait que la compacit\'e.
\end{rem}
Le th\'eor\`eme suivant nous permet d'identifier les deux m\'etriques
introduites pr\'ec\'edemment~:
\begin{thm}
Soit $D$ un diviseur de Cartier horizontal et $T$-invariant 
sur $\P$. Pour tout entier $p \geqslant 2$, on a l'\'egalit\'e des m\'etriques~:
\[
\|.\|_{\op{BT}} = \|.\|_{\op{Zh},p}
\]
sur le faisceau inversible $\C{O}(D)(\M{C})$.
\end{thm}
\demo\ Il suffit de v\'erifier que $\| . \|_{\op{BT}}$ satisfait \`a la
condition donn\'ee au (2) du th\'eor\`eme (\ref{zhang}). Cela d\'ecoule des
propri\'et\'es de multiplicativit\'e et de fonctorialit\'e \'enonc\'ees \`a la
proposition (\ref{BT_fonct}).
\medskip

\subsubsection{Construction par image inverse}~
\label{chap_image_inverse}

Dans ce paragraphe, $\D$ d\'esigne un \'eventail poss\'edant une fonction
support strictement concave relativement \`a $\D$. Par cons\'equent la
vari\'et\'e torique $\P$ est {\em projective}.

Soit $D$ un diviseur de Cartier $T$-invariant sur $\P$ et $\C{O}(D)$ le faisceau
inversible associ\'e. On suppose dans un premier temps que $\C{O}(D)$ est engendr\'e
par ses sections globales. Le choix d'un ordre sur 
les \'el\'ements de $K_{D} \cap M$ permet de d\'efinir 
un morphisme $T$-\'equivariant associ\'e 
\`a $D$, que l'on note $\phi_{D}$, de la fa\c con suivante~: 
\begin{alignat*}{3}
\phi_{D}:\,\M{P} & (\D) & &\longrightarrow & & \;\M{P}_{\M{Z}}^{k_{D}} \\
&x & &\longrightarrow & &(\chi^{m}(x))_{m \in K_{D} \cap M}
\end{alignat*}
o\`u $k_{D}$ est un entier positif d\'efini par $k_{D} = \#(K_{D}\cap M) - 1$.
On note $\C{O}(1)$ le fibr\'e de Serre sur $\M{P}_{\M{Z}}^{k_{D}}$ 
et on le munit de la m\'etrique d\'efinie pour toute section m\'eromorphe de
$\C{O}(1)(\M{C})$ par~:
\begin{equation}
\label{image_inverse}
\| s(x) \|_{\infty} = \frac{|s(x)|}{\sup_{1 \leqslant i \leqslant k_{D} + 1}|x_{i}|}\; .
\end{equation}
Cette m\'etrique est la m\'etrique de Batyrev-Tschinkel ou de Zhang pour le
faisceau $\C{O}(1)$ sur $\M{P}_{\M{Z}}^{k_{D}}$ consid\'er\'ee comme
vari\'et\'e torique comme dans l'exemple (\ref{exemple_intro_1}). 
On note $\overline{\C{O}(1)}_{\infty}$ le faisceau
$\C{O}(1)$ muni de cette m\'etrique.
Comme $\|.\|_{\infty}$ est invariante si l'on permute les \'el\'ements de
$K_{D} \cap M$, la m\'etrique sur $\C{O}(D)$ d\'efinie par~:
\[
\|.\|_{D,\infty} = \phi_{D}^{\ast}\,\|.\|_{\infty}
\]
est ind\'ependante du choix effectu\'e pour d\'efinir $\phi_{D}$. On note
$\overline{\C{O}(D)}_{\infty}$ le faisceau $\C{O}(D)$ muni de la m\'etrique
$\|.\|_{D,\infty}$. On a alors la proposition suivante~:
\begin{prop}
\label{decomposition_0}
Soient $D$ et $E$ deux diviseurs $T$-invariants sur $\P$ dont les faisceaux
associ\'es sont engendr\'es par leurs sections globales, on a~:
\[
\overline{\C{O}(D)}_{\infty} \otimes \ov{\C{O}(E)}_{\infty} =
\ov{\C{O}(D + E)}_{\infty}.
\]
\end{prop}
\demo\ On d\'emontre tout d'abord le lemme suivant~:
\begin{lem}
\label{metrique}
Soit $D$ un diviseur $T$-invariant sur $\P$. On notera $K_{D}(0)$ les points
entiers extr\'emaux du polytope $K_{D}$.
Pour tout $x \in \P(\M{C})$, on a~:
\[
\sup_{m \in K_{D} \cap M}|\chi^{m}(x)| = \sup_{m \in K_{D}(0)}|\chi^{m}(x)|. 
\]
\end{lem}
\demo\ Soient $m_{0},
\dots, m_{q}$ les \'el\'ements de $K_{D}(0)$ et $m$ un \'el\'ement quelconque
de $K_{D} \cap M$. Comme $K_{D}$ est convexe, on peut trouver des r\'eels
positifs $\alpha_{0}, \dots, \alpha_{q}$ tels que~:
\[
\left\{
\begin{array}{l}
m = \alpha_{0}m_{0} + \dots + \alpha_{q}m_{q} \\
1 = \alpha_{0} + \dots + \alpha_{q}.
\end{array}
\right.
\]
On en d\'eduit l'in\'egalit\'e~:
\[
|\chi^{m}(x)| \leqslant \sup_{0 \leqslant i \leqslant q}|\chi^{m_{i}}(x)|,
\]
ce qui joint \`a l'inclusion $K_{D}(0) \subset K_{D} \cap M$ donne le
r\'esultat annonc\'e.
\medskip

On passe maintenant \`a la d\'emonstration de la proposition.
Soit $s$ (resp. $r$) une section holomorphe de $\C{O}(D)(\M{C})$ 
(resp. de $\C{O}(E)(\M{C})$) au-dessus d'un ouvert $\Omega$ de $\P(\M{C})$. 
Pour tout $x \in \Omega$, on a~:
\[
\|s(x)\|_{D,\infty} = \frac{|s(x)|}{\sup_{m \in K_{D} \cap M}|\chi^{m}(x)|}
\quad \left(\text{resp.} \quad \|r(x)\|_{E,\infty} = \frac{|r(x)|}{\sup_{m \in
K_{E} \cap M}|\chi^{m}(x)|} \right).
\]
En utilisant le lemme (\ref{metrique}) et le fait que~:
\[
K_{D + E}(0) = (K_{D} + K_{E})(0) \subseteq K_{D}(0) + K_{E}(0) \subseteq K_{D}
+ K_{E} = K_{D + E},
\]
on obtient l'\'egalit\'e~:
\[
\sup_{m \in K_{D + E}}|\chi^{m}(x)| = \left(\sup_{m \in
K_{D}}|\chi^{m}(x)|\right)
\left(\sup_{m \in K_{E}}|\chi^{m}(x)|\right)
\]
et la proposition est d\'emontr\'ee.
\medskip

\begin{prop_defn}
\label{decomposition}
Soit $D$ un diviseur de Cartier horizontal et $T$-invariant sur $\P$. 
Il existe $E$ et $F$ des diviseurs de Cartier horizontaux et
$T$-invariants sur $\P$, dont le faisceau associ\'e est engendr\'e par
ses sections globales, et tels que~:
\[
D = E - F.
\]
La m\'etrique $\| .\|_{E,\infty}\, . \,\| .\|_{F,\infty}^{-1}$ induite sur
$\C{O}(D) = \C{O}(E) \otimes \C{O}(F)^{-1}$ 
par les m\'etriques canoniques $\| .\|_{E,\infty}$ et
$\| .\|_{F,\infty}$ est ind\'ependante de $E$ et de $F$. On note
$\| .\|_{D,\infty}$ cette m\'etrique et on note $\ov{\C{O}(D)}_{\infty}$ le
faisceau $\C{O}(D)$ muni de la m\'etrique $\| .\|_{D,\infty}$.
\end{prop_defn}
\demo\ 
La premi\`ere partie de l'\'enonc\'e est classique (Serre) et est tr\`es facile
\`a g\'en\'eraliser dans le cadre \'equivariant (prendre $H$ diviseur de
Cartier $T$-invariant
tel que $\C{O}(H)$ est engendr\'e par ses sections globales
et consid\'erer $E = nH$ et $F = - D + nH$ pour $n$ assez grand). Soient
$E$, $F$, $E'$ et $F'$ des diviseurs de Cartier $T$-invariants dont les faisceaux
associ\'es sont engendr\'es par leurs sections globales et tels que~:
\[
D = E - F = E' - F'.
\]
On tire $E + F' = E' + F$, et donc des isomorphismes isom\'etriques~:
\[
\ov{\C{O}(E)}_{\infty}\otimes \ov{\C{O}(F')}_{\infty} = 
\ov{\C{O}(E) \otimes \C{O}(F')}_{\infty} = 
\ov{\C{O}(E') \otimes \C{O}(F)}_{\infty} = 
\ov{\C{O}(E')}_{\infty} \otimes \ov{\C{O}(F)}_{\infty},
\]
on conclut que~:
\[
\ov{\C{O}(E)}_{\infty} \otimes (\ov{\C{O}(F)}_{\infty})^{-1} \simeq
\ov{\C{O}(E')}_{\infty} \otimes (\ov{\C{O}(F')}_{\infty})^{-1}, 
\]
ce qui termine la d\'emonstration.
\medskip

\begin{prop}
Soient $E$ et $F$ deux diviseurs de Cartier horizontaux et $T$-invariants 
sur $\P$. On a un isomorphisme isom\'etrique~:
\[
\ov{\C{O}(E)}_{\infty} \otimes \ov{\C{O}(F)}_{\infty} \simeq 
\ov{\C{O}(E + F)}_{\infty}.
\]
\end{prop}
\demo\ La proposition est vraie dans le cas o\`u $\C{O}(E)$ et $\C{O}(F)$ 
sont engendr\'es par leurs sections globales.
Le cas g\'en\'eral s'obtient par diff\'erence \`a partir des propositions
(\ref{decomposition_0}) et (\ref{decomposition}).
\medskip

La construction qui vient d'\^etre donn\'ee co\"\i ncide avec les deux autres
pr\'esent\'ees pr\'ec\'edemment~:
\begin{thm}
Pour tout diviseur de Cartier horizontal $T$-invariant $D$ sur $\P$, les
m\'etriques $\|.\|_{\op{BT}}$, $\|.\|_{\op{Zh}}$ et $\|.\|_{D,\infty}$ sur
$\C{O}(D)$ co\"\i ncident.
\end{thm}
\demo\ 
Cela d\'ecoule directement des propri\'et\'es de multiplicativit\'e et de
fonctorialit\'e de la m\'etrique $\|.\|_{\op{BT}}$ \'enonc\'ees \`a la 
proposition (\ref{BT_fonct}), ajout\'e au fait que $\|.\|_{\op{BT}}$ et
$\|.\|_{D,\infty}$ co\"\i ncident par construction lorsque $\P$ est l'espace
projectif $\M{P}^{n}_{\Z}$ et que $\C{O}(D) = \C{O}(1)$.
\medskip

Dans toute la suite de ce paragraphe, $\P$ d\'esigne une vari\'et\'e 
torique {\em projective} et {\em lisse}.

On constate que les m\'etriques canoniques d\'efinies pr\'ec\'edemment
ne sont pas $C^{\infty}$ en g\'en\'eral (consid\'erer par exemple 
$\P = \M{P}^{n}_{\Z}$ et
$\ov{\C{O}(D)}_{\infty} = \ov{\C{O}(1)}_{\infty}$). 
On a n\'eanmoins le r\'esultat suivant~:
\begin{prop}
\label{psh}
Soient $\P$ une vari\'et\'e torique projective et lisse, et 
$D$ un diviseur de Cartier horizontal $T$-invariant sur $\P$
tel que $\C{O}(D)$ soit engendr\'e par ses sections globales.
Pour tout ouvert $\Omega \subset \P(\M{C})$ et pour toute 
section holomorphe $s$ de $\C{O}(D)(\M{C})$ sur $\Omega$ ne s'annulant pas, 
la fonction
$x \mapsto - \log \|s(x)\|_{D,\infty}^{2}$ est continue et \psh\ sur $\Omega$.
\end{prop}
\demo\ La continuit\'e d\'ecoule des d\'efinitions. Comme la condition de
plurisousharmonicit\'e est pr\'eserv\'ee par changement de variable holomorphe
(voir par exemple \cite{7}, th. 1.5.9), il suffit de d\'emontrer le r\'esultat
pour le faisceau $\ov{\C{O}(1)}_{\infty}$ sur $\M{P}^{n}_{\Z}$; et
d'apr\`es la formule (\ref{image_inverse}) cela r\'esulte
de (\cite{7}, th. 1.5.6 et
exemple 1.5.10). On peut \'egalement consulter (\cite{5}, \S 3).
\medskip

On montre maintenant un th\'eor\`eme d'approximation des m\'etriques canoniques
par des m\'etriques $C^{\infty}$ sur $\P(\M{C})$~:
\begin{prop}
\label{approximation}
Soient $\P$ une vari\'et\'e torique projective et lisse, et 
$D$ un diviseur de Cartier horizontal $T$-invariant sur $\P$
tel que $\C{O}(D)$ soit engendr\'e par ses sections globales.
Il existe une suite de m\'etriques
$\left( \|.\|_{n}\right)_{n \in \M{N}}$ sur $\C{O}(D)(\M{C})$ convergeant
uniform\'ement vers $\|.\|_{D,\infty}$ sur $\P(\M{C})$ et v\'erifiant les conditions
suivantes~:
\begin{itemize}
\item{Les m\'etriques $\|.\|_{n}$ sont $C^{\infty}$.}
\item{Si $\Omega \subset \P(\M{C})$ est un ouvert et $s$ une section holomorphe
de $\C{O}(D)(\M{C})$
sur $\Omega$ ne s'annulant pas, la fonction $x \mapsto \log \|s(x)\|_{n}^{2}$ est 
continue et \psh\ sur $\Omega$; en
particulier le courant~:
\[
c_{1}(\C{O}(D)(\M{C}),\|.\|_{n}) = - d d^{c}\,\log \|s(x)\|^{2}_{n},
\]
est positif.}
\item{Pour tout $x \in \Omega$, la suite $( - \log \|s(x)\|^{2}_{n})_{n \in
\M{N}}$ est
d\'ecroissante et converge vers $- \log \|s(x)\|^{2}_{D,\infty}$.}
\end{itemize}
\end{prop}
\demo\ La condition de plurisousharmonicit\'e \'etant pr\'eserv\'ee par changement
de variable holomorphe, il suffit d'apr\`es la construction par image inverse
de d\'emontrer la proposition pour $\ov{\C{O}(1)}_{\infty}$ sur
$\M{P}^{m}_{\Z}$. On consid\`ere alors la famille de m\'etriques $\|.\|_{n}$ 
d\'efinies
par~:
\[
\|s(x)\|_{n} = \frac{|s(x)|}{\left( \sum_{i = 0}^{m}|x_{i}|^{n}\right)^{1/n}}\;
,
\]
pour toute section locale holomorphe $s$ de $\C{O}(1)(\M{C})$. Les m\'etriques
$\|.\|_{n}$ sont $C^{\infty}$ sur $\P(\M{C})$. De plus, la fonction~:
\[
(x_{0}, \dots, x_{m}) \longmapsto \log (e^{x_{0}} + \dots + e^{x_{m}})
\]
\'etant convexe et croissante en chacun des $x_{i}$ et la fonction $t \mapsto \log
|t|$
\'etant sousharmonique, la fonction  $\log \left(
\sum_{i=0}^{m}|x_{i}|^{n}\right)^{1/n}$ est \psh\ (voir par exemple 
\cite{7}, th.
1.5.6 et \cite{5}, \S 3).
Enfin la suite $\left(\sum_{i
=0}^{m}|x_{i}|^{n}\right)^{1/n}$, $(n \in \M{N})$, est d\'ecroissante et converge
vers $\sup_{0 \leqslant i \leqslant m}|x_{i}|$, et la positivit\'e de
$c_{1}(\C{O}(1)(\M{C}),\|.\|_{n})$ d\'ecoule directement de (\cite{7}, exemple 3.1.18).
\bigskip

\subsection{M\'etriques canoniques sur les fibr\'es en droites sur $\P$.}~

Dans cette section, $\D$ d\'esigne un \'eventail complet et r\'egulier de $N$
poss\'edant une fonction support strictement concave. En d'autres termes, on
suppose que $\P$ est une vari\'et\'e torique projective lisse.

Pour tout fibr\'e en droites $L$ sur $\P$, on construit de mani\`ere canonique
une m\'etrique sur $L(\M{C})$.
\begin{prop_defn}
\label{metrique_ind}
Soit $L$ un fibr\'e en droites sur $\P$. Il existe un diviseur horizontal
$T$-invariant $D$ sur $\P$ et un isomorphisme~:
\[
\Phi: L \longrightarrow \C{O}(D).
\]
La m\'etrique $\Phi^{\ast}\|.\|_{D,\infty}$ sur $L$ est ind\'ependante des
choix de $D$ et $\Phi$. On l'appelle {\it m\'etrique canonique\/} sur $L$ et
on la note $\|.\|_{L,\infty}$, ou plus simplement $\|.\|_{\infty}$
lorsqu'aucune
confusion n'est \`a craindre. On note $\ov{L}_{\infty} =
(L,\|.\|_{L,\infty})$ le fibr\'e $L$
muni de sa m\'etrique canonique.
\end{prop_defn}
\demo\
L'existence de $D$ et $\Phi$ est une reformulation de la surjectivit\'e de $s$
dans la proposition (\ref{picard}). Soit maintenant $D'$ un second diviseur
horizontal $T$-invariant de $\P$ tel qu'il existe un isomorphisme
$\Phi': L \simeq \C{O}(D')$. Le diviseur $D-D'$ est principal, et il existe un
unique \'el\'ement $m \in M$ tel que $D-D' = \op{div}\chi^{m}$ d'apr\`es la
proposition (\ref{picard}).

Comme les unit\'es globales sur $\P$ sont $\{1,-1\}$, les isomorphismes $\Phi$
et $\Phi'$ sont uniques au signe pr\`es. Pour montrer que~:
\[
\Phi^{\ast}\|.\|_{D,\infty} = {\Phi'}^{\ast}\|.\|_{D',\infty}
\]
il suffit donc de v\'erifier que l'isomorphisme~:
\[
\C{O}(D') \simeq \C{O}(D)
\]
d\'efini par la multiplication par $\chi^{m}$ transporte $\|.\|_{D',\infty}$
sur $\|.\|_{D,\infty}$. Cela d\'ecoule de l'expression (\ref{BT_definition}) 
de Batyrev et Tschinkel pour ces m\'etriques.
\medskip

La proposition suivante est une cons\'equence imm\'ediate du th\'eor\`eme
(\ref{zhang}).
\begin{prop}
\label{relation_zhang}
Soit $p$ un entier sup\'erieur ou \'egal \`a $2$. Pour tout fibr\'e en droites
$L$ sur $\P$, on a un isomorphisme isom\'etrique~:
\[
[p]^{\ast}(\ov{L}_{\infty}) \simeq (\ov{L}_{\infty})^{\otimes p},
\]
et $\|.\|_{L,\infty}$ est l'unique m\'etrique continue telle que
$\ov{L}_{\infty} = (L,\|.\|_{L,\infty})$
v\'erifie cette propri\'et\'e.
\end{prop}
\bigskip

\section{Produits de courants}~

\subsection{Motivation}~

Soient $\P$ une vari\'et\'e torique projective lisse et 
$\ov{L}_{\infty}$ un fibr\'e en droites engendr\'e par ses
sections globales sur $\P$ et muni de sa
m\'etrique canonique. Soit $s$ une section holomorphe de $L$ sur un ouvert $U
\subset \P(\M{C})$. En g\'en\'eral la fonction $x \mapsto -\log
\|s(x)\|^{2}_{L,\infty}$ n'est pas de classe $C^{\infty}$ mais seulement \psh\
sur $U$ (cf. prop. (\ref{psh})); la premi\`ere ``forme'' de Chern 
$c_{1}(\ov{L}_{\infty}) = -\dd \log
\|s\|^{2}_{L,\infty}$ n'est alors d\'efinie qu'au sens des distributions (c'est
un courant de bidegr\'e (1,1)). On ne peut esp\'erer d\'efinir en toute 
g\'en\'eralit\'e le produit de deux courants; ici n\'eanmoins,
$c_{1}(\ov{L}_{\infty})$ est un courant positif
(voir \cite{7}, prop. 3.1.18 et 3.1.14). Une construction
due \`a Bedford et Taylor (cf. \cite{1}) permet d'associer 
\`a des fibr\'es en droites $\ov{L}_{1,\infty}, \dots,
\ov{L}_{p,\infty}$ engendr\'es par leurs
sections globales sur $\P$ et munis de leur
m\'etrique canonique, un produit~:
\[
c_{1}(\ov{L}_{1,\infty}) \dotsm c_{1}(\ov{L}_{p,\infty})
\]
poss\'edant des propri\'et\'es satisfaisantes.

La construction pr\'esent\'ee
dans ce qui suit \'etant en fait tout \`a fait g\'en\'erale, on abandonne dans
cette partie et la suivante le point de vue particulier des vari\'et\'es
toriques.
\bigskip

\subsection{Th\'eorie de Bedford-Taylor-Demailly}~

On suit l'expos\'e donn\'e dans (\cite{5}, \S 3; \cite{6}, \S 1 et \S 2; et
\cite{7}, \S 3.3). Tous les r\'esultats pr\'esent\'es ici sont dus \`a Bedford
et Taylor (cf. \cite{1}) et Demailly (cf. \cite{5}, \cite{6} et \cite{7}). 

On rappelle tout d'abord quelques d\'efinitions et notations.
\medskip

Dans toute cette partie, $X$ d\'esigne une vari\'et\'e analytique complexe de dimension
complexe $d$. On note $A^{p,q}(X)$ (resp. $D^{p,q}(X)$) l'espace vectoriel
des formes diff\'erentiables $C^{\infty}$ complexes (resp. l'espace vectoriel
des courants complexes) de type $(p,q)$ sur $X$. Soit $Y \subset X$ un
cycle analytique irr\'eductible de codimension $p$, on note
$\delta_{Y} \in D^{p,p}(X)$ le courant d'int\'egration sur $Y$. Pour
tout \'el\'ement $T \in D^{p,p}(X)$, on note $\op{Supp} T \subset X$ le
support de $T$. 
Enfin pour tout ouvert $U \subset X$, on notera $\op{Psh}(U)$ 
l'ensemble des fonctions
\pshs\ de $U$ vers $[-\infty,+\infty[$ semi-continues sup\'erieurement. 

\begin{defn} {\bf (Lelong).}~
Soit $p$ un entier positif et $U$ un ouvert de $X$. Un courant $T \in
D^{p,p}(U)$ est dit {\it positif\/} (ou {\it faiblement positif\/}) et on note
$T \geqslant 0$, \ssi\ pour
tout choix de $(1,0)$ formes $\alpha_{1}, \dots, \alpha_{d-p}$ de classe
$C^{\infty}$ \`a support compact sur $U$, la distribution~:
\[
T \wedge (i\alpha_{1}\wedge\ov{\alpha}_{1}) \wedge \dots \wedge
(i\alpha_{d-p}\wedge\ov{\alpha}_{d-p})
\]
est une mesure positive sur $U$. On note
$D_{+}^{p,p}(U) \subset D^{p,p}(U)$ l'ensemble des courants positifs de type
$(p,p)$ sur $U$.
\end{defn}
\medskip

On pose alors~:
\begin{defn} {\bf (Bedford-Taylor).}~
Soit $T \in D_{+}^{p,p}(U)$ un courant positif {\em fer\-m\'e} de type $(p,p)$ et
$u \in \op{Psh}(U)$ une fonction \psh\ localement born\'ee sur $U$ un ouvert de
$X$. Le produit $(\dd u) \wedge T$ est d\'efini par la formule~:
\[
(\dd u) \wedge T = \dd (uT)
\]
\end{defn}
\begin{rem}
Le produit $uT$ est bien d\'efini puisque $u$ est localement born\'ee et $T$
est d'ordre $0$, i.e. \`a coefficients mesures (pour une d\'emonstration de ce
fait, voir par exemple \cite{7}, 3.1.14).
\end{rem}
\begin{prop}
Le produit $(\dd u) \wedge T$ d\'efini ci-dessus est un courant positif ferm\'e
de bidegr\'e $(p+1, p+1)$. Il \'etend la d\'efinition usuelle du produit $(\dd
u) \wedge T$ dans le cas o\`u $u$ est une fonction $C^{\infty}$ sur $U$.
\end{prop}
\demo\ Voir par exemple (\cite{5}, prop. 5.1; \cite{6}, prop. 3.3.2 ou \cite{7},
prop. 1.2).
\medskip

De mani\`ere plus g\'en\'erale, \'etant donn\'ees $T \in D_{+}^{p,p}(U)$ et 
$u_{1}, \dots,
u_{q}$ des fonctions \pshs\ localement born\'ees sur $U$, on peut d\'efinir par
r\'ecurrence sur $q$ le courant~:
\[
(\dd u_{1}) \wedge \dots \wedge (\dd u_{q}) \wedge T = \dd (u_{1}(\dd u_{2})
\wedge \dots \wedge (\dd u_{q}) \wedge T).
\]
C'est encore un courant positif ferm\'e, de bidegr\'e $(p+q,p+q)$.
\medskip

Le produit ainsi d\'efini a un comportement agr\'eable vis-\`a-vis de la limite
uniforme~:
\begin{thm}
\label{lim_uniforme}
Soient $u_{1}, \dots, u_{q}$ des fonctions \pshs\ continues sur $U$. Soient
$(u_{1}^{(k)})_{k \in \M{N}}, \dots, (u_{q}^{(k)})_{k \in \M{N}}$, $q$ suites de
fonctions \pshs\ localement born\'ees sur $U$ convergeant uniform\'ement
sur tout compact de $U$ vers $u_{1}, \dots, u_{q}$ respectivement, et
$(T_{k})_{k \in \M{N}}$ une suite de courants positifs ferm\'es convergeant
faiblement vers $T$ sur $U$. Alors~:
\begin{alignat*}{3}
u_{1}^{(k)}(\dd u_{2}^{(k)}) \wedge \dots \wedge (\dd u_{q}^{(k)}) \wedge
&T_{k} & &\text{\quad tend vers\quad}& &u_{1}(\dd u_{2}) \wedge \dots \wedge
(\dd u_{q}) \wedge T \\
\text{et \quad} (\dd u_{1}^{(k)}) \wedge \dots \wedge (\dd u_{q}^{(k)})\wedge
&T_{k}& &\text{\quad vers\quad} & &(\dd u_{1}) \wedge \dots \wedge (\dd
u_{q}) \wedge T, 
\end{alignat*}
au sens de la convergence faible des courants.
\end{thm}
\demo\ Voir (\cite{6}, cor. 1.6) ou (\cite{7}, cor. 3.3.6).
\medskip

Au vue de la prop. (\ref{approximation}), on souhaite affaiblir les
hypoth\`eses du th\'eor\`eme pr\'ec\'edent et remplacer la convergence uniforme
par la convergence simple d\'ecroissante; c'est l'objet du th\'eor\`eme
suivant~:
\begin{thm}{\rm \bf (Bedford-Taylor).}~
\label{BeT}
Soient $u_{1}, \dots, u_{q}$ appartenant \`a $\op{Psh}(U)$ et localement 
born\'ees sur $U$, 
et $u_{1}^{(k)}, \dots, u_{q}^{(k)}$, $q$ suites d\'ecroissantes de fonctions
dans $\op{Psh}(U)$
localement born\'ees sur $U$
et convergeant simplement sur $U$ vers
$u_{1}, \dots, u_{q}$ respectivement. On a~:
\begin{alignat*}{3}
u_{1}^{(k)}(\dd u_{2}^{(k)}) \wedge \dots \wedge (\dd u_{q}^{(k)}) \wedge
&T & &\text{\quad tend vers\quad}& &u_{1}(\dd u_{2}) \wedge \dots \wedge
(\dd u_{q}) \wedge T \\
\text{et \quad} (\dd u_{1}^{(k)}) \wedge \dots \wedge (\dd u_{q}^{(k)})\wedge
&T & &\text{\quad vers\quad} & &(\dd u_{1}) \wedge \dots \wedge (\dd
u_{q}) \wedge T, 
\end{alignat*}
au sens de la convergence faible des courants.
\end{thm}
\demo\ Voir par exemple (\cite{6}, \S 1.7 ou \cite{7}, th. 3.3.7).
\medskip

Par r\'egularisation (voir par exemple \cite{7}, th. 1.5.5), on d\'eduit
imm\'ediatement de ce th\'eor\`eme le corollaire~:
\begin{cor}
\label{BeT_commutativite}
Soient $u_{1}, \dots, u_{q}$ des fonctions \pshs\ localement born\'ees sur $U$
et $T \in D_{+}^{p,p}(U)$ un courant positif ferm\'e; le produit $u_{1}(\dd u_{2})
\wedge \dots \wedge (\dd u_{q}) \wedge T$ (resp. le produit $(\dd u_{1}) \wedge
\dots \wedge (\dd u_{q}) \wedge T$) est ind\'ependant de l'ordre des $u_{2},
\dots, u_{q}$ (resp. de l'ordre des $u_{1}, \dots, u_{q}$).
\end{cor}
\begin{defn}
Soit $U$ un ouvert de $X$ et $u \in \op{Psh}(U)$. On appelle {\it lieu non
born\'e\/} de $u$ dans $U$ et on note $L(u)$ l'ensemble des points $x \in U$
tels que $u$ n'est born\'ee sur aucun voisinage de $x$ dans $U$.
\end{defn}
\medskip

Le r\'esultat suivant, d\^u \`a Demailly, montre que sous certaines conditions,
on peut abandonner l'hypoth\`ese selon laquelle les
fonctions \pshs\ $u_{1}, \dots, u_{q}$ dans les \'enonc\'es (\ref{BeT}) et
(\ref{BeT_commutativite}) doivent \^etre choisies {\it localement born\'ees\/}.
\begin{thm} {\rm \bf (Demailly).}~
\label{demailly}
Soit $U$ un ouvert de $X$ et soient $T \in D_{+}^{p,p}(U)$ ferm\'e et $u_{1},
\dots, u_{q}$ appartenant \`a $\op{Psh}(U)$. On suppose que $q \leqslant d - p$
et que pour tout choix d'indices $j_{1} < \dots < j_{m}$ dans $\{1,\dots,q\}$
l'intersection $L(u_{j_{1}}) \cap \dots \cap L(u_{j_{m}}) \cap \op{Supp}T$ est
contenue dans un ensemble analytique de dimension complexe inf\'erieure ou
\'egale \`a $(d - p - m)$. On peut construire 
des courants $u_{1}(\dd u_{2}) \wedge \dots \wedge
(\dd u_{q}) \wedge T$ et $(\dd u_{1}) \wedge \dots \wedge (\dd u_{q}) \wedge T$
de masse localement finie sur $U$ et caract\'eris\'es de mani\`ere unique
par la propri\'et\'e suivante~: Pour toutes suites d\'ecroissantes
$(u_{1}^{(k)})_{k \in \M{N}}, \dots, (u_{q}^{(k)})_{k \in \M{N}}$
de fonctions \pshs\ convergeant simplement vers $u_{1}, \dots,
u_{q}$ respectivement, on a~:
\begin{alignat*}{3}
u_{1}^{(k)}(\dd u_{2}^{(k)}) \wedge \dots \wedge (\dd u_{q}^{(k)}) \wedge
&T & &\text{\quad tend vers\quad}& &u_{1}(\dd u_{2}) \wedge \dots \wedge
(\dd u_{q}) \wedge T \\
\text{et \quad} (\dd u_{1}^{(k)}) \wedge \dots \wedge (\dd u_{q}^{(k)})\wedge
&T & &\text{\quad vers\quad} & &(\dd u_{1}) \wedge \dots \wedge (\dd
u_{q}) \wedge T, 
\end{alignat*}
au sens de la convergence faible des courants sur $U$.
\end{thm}
\demo\ Voir (\cite{6}, th. 2.5 et prop. 2.9 et \cite{7}, th. 3.4.5 et prop.
3.4.9). On peut aussi consulter (\cite{5}, th. 5.4).
\medskip

On d\'eduit imm\'ediatement du th\'eor\`eme pr\'ec\'edent les corollaires
suivants~:
\begin{cor}
Soient $u_{1}, \dots, u_{q}$ et $T$ comme au th\'eor\`eme (\ref{demailly}); le produit
$u_{1}(\dd u_{2})\wedge \dots \wedge (\dd u_{q}) \wedge T$ (resp. le produit
$(\dd u_{1}) \wedge \dots \wedge (\dd u_{q}) \wedge T$) ne d\'epend pas de
l'ordre de $u_{2}, \dots, u_{q}$ (resp. de l'ordre de
$u_{1}, \dots u_{q}$).
\end{cor}
\begin{cor}
\label{produit}
Soient $u_{1}, \dots, u_{q}$ \pshs\ sur $U$ et telles que pour tout $1\leqslant
i \leqslant q$, $L(u_{i})$ est contenu dans un ensemble analytique $A_{i}
\subset U$. Si pour tout choix d'indices $j_{1} < \dots j_{m}$ dans $\{1,\dots,
q\}$, on a~:
\[
\op{codim}A_{j_{1}} \cap \dots \cap A_{j_{m}} \geqslant m, 
\]
alors les courants $u_{1}(\dd u_{2}) \wedge \dots \wedge (\dd u_{q})$ et $(\dd
u_{1}) \wedge \dots \wedge (\dd u_{q})$ sont bien d\'efinis et v\'erifient les
propri\'et\'es d'approximation \'enonc\'ees au th\'eor\`eme (\ref{demailly}).
\end{cor}
On expose enfin une construction due \`a Gillet-Soul\'e (cf. \cite{13}, \S
2.1.5)~:
\\
Soient $X$ une vari\'et\'e projective complexe, $Z$ un cycle de codimension $q$
de $X$ et $g_{Z}$ un courant de Green localement $L^{1}$ pour $Z$ (i.e. un
\'el\'ement de $D^{p-1,p-1}(X)$, localement $L^{1}$ sur $X$, $C^{\infty}$ sur
$X - |Z|$ et tel que $\dd g_{Z} + \delta_{Z} = \omega_{Z}$ est $C^{\infty}$ sur
$X$). Choisissons une m\'etrique sur $X$ et notons pour tout $\varepsilon \in
\R^{+\ast}$~:
\[
N_{\epsi}(Z) = \{x \in X: \quad d(x,Z) < \epsi\}.
\]
Pour tout $\epsi \in \R^{+\ast}$, soit $\rho_{\epsi}: X \rightarrow \R$ une
fonction $C^{\infty}$ telle que~:
\begin{itemize}
\item{$0 \leqslant \rho_{\epsi} \leqslant 1$ sur $X$,}
\item{$\rho_{\epsi} = 1$ en dehors de $N_{\epsi}(Z)$,}
\item{$\rho_{\epsi} = 0$ sur $N_{\epsi/2}(Z)$,}
\end{itemize}
et posons $g_{Z}^{(\epsi)} = \rho_{\epsi}g_{Z}$.
\begin{prop}
\label{construction_GS_lissage}
Pour tout $\epsi \in \R^{+\ast}$, les assertions suivantes sont v\'erifi\'ees~:
\begin{itemize}
\item{$g_{Z}^{(\epsi)}$ est une forme $C^{\infty}$ sur $X$.}
\item{On a~: $\dd g_{Z}^{(\epsi)} + \delta_{Z}^{(\epsi)} = \omega_{Z}$, o\`u
$\delta_{Z}^{(\epsi)}$ est une forme ferm\'ee $C^{\infty}$ dont le support est
contenu dans $\ov{N_{\epsi}(Z)}$.}
\end{itemize}
De plus, on a les limites suivantes~:
\begin{itemize}
\item{$\lim_{\epsi \rightarrow 0}g_{Z}^{(\epsi)} = g_{Z}$,}
\item{$\lim_{\epsi \rightarrow 0}\delta_{Z}^{(\epsi)} = \delta_{Z}$,}
\end{itemize}
au sens de la topologie faible des courants.
\end{prop}
\demo\ Voir \cite{13}, \S 2.1.5.
\bigskip

\subsection{Formes diff\'erentielles g\'en\'eralis\'ees}~

Les r\'esultats \'enonc\'es dans la section pr\'ec\'edente motivent la
d\'efinition suivante~:
\begin{defn}
Soit $U$ un ouvert de $X$ et $u_{1}, \dots, u_{q}$ des fonctions de $U
\rightarrow [-\infty, +\infty [$. Le $q$-uplet $(u_{1}, \dots, u_{q})$ est dit
{\it admissible sur $U$\/} \ssi~:
\begin{enumerate}
\item{Les fonctions $u_{1}, \dots, u_{q}$ sont des \'el\'ements de
$\op{Psh}(U)$.}
\item{Pour tout $1 \leqslant i \leqslant q$, l'ensemble $L(u_{i})$ est contenu
dans un ensemble analytique $A_{i} \subset U$.}
\item{Pour tout choix d'indices $j_{1} < \dots < j_{m}$ dans $\{1, \dots, q\}$,
on a~:
\[
\op{codim}A_{j_{1}} \cap \dots \cap A_{j_{m}} \geqslant m.
\]
}
\end{enumerate}
\end{defn}
\medskip

Pour tout $q$-uplet $(u_{1}, \dots, u_{q})$ admissible sur un ouvert $U \subset
X$, le corollaire (\ref{produit}) nous permet de d\'efinir sur $U$ les courants
$u_{1}(\dd u_{2}) \wedge \dots \wedge (\dd u_{q})$ et $(\dd u_{1}) \wedge \dots
\wedge (\dd u_{q})$.

On pose~:
\begin{defn}
Soit $p$ un entier positif. On note $\ov{\ov{A}}^{p,p}(X) \subset D^{p,p}(X)$
(resp. $\ov{\ov{A}}^{p,p}_{\op{log}}(X) \subset D^{p,p}(X)$) l'espace vectoriel
complexe form\'e des \'el\'ements de $D^{p,p}(X)$ qui, sur tout ouvert $U$
d'un recouvrement suffisamment fin de $X$, peuvent s'\'ecrire sous la forme~:
\begin{align*}
& \sum_{i = 1}^{n}\omega_{i}(\dd u_{i,1})\wedge \dots \wedge (\dd u_{i,q_{i}})
\\
\big( \text{resp. \quad} &\sum_{i =1}^{n}\omega_{i}u_{i,1}(\dd u_{i,2}) \wedge
\dots \wedge (\dd u_{i, q_{i}}) \big) , 
\end{align*}
o\`u pour tout $1 \leqslant i \leqslant n$, on a $\omega_{i} \in A^{p-q_{i}, p -
q_{i}}(U)$ (resp. $\omega_{i} \in A^{p -q_{i}+1, p - q_{i} + 1}(U)$) et le
$q_{i}$-uplet $(u_{i,1}, \dots, u_{i,q_{i}})$ est admissible sur l'ouvert $U$
consid\'er\'e.
\end{defn}
\medskip

On pose aussi la d\'efinition suivante~:
\begin{defn}
\label{formes_generalisees}
Soit $p$ un entier positif. On dit que $\alpha \in \ov{\ov{A}}^{p,p}(X)$ est
une {\it forme diff\'erentielle g\'en\'eralis\'ee\/} de type $(p,p)$ \ssi\ sur
tout ouvert $U$ d'un recouvrement suffisamment fin de $X$, on peut \'ecrire la restriction de
$\alpha$ \`a $U$ sous la forme~:
\[
\alpha_{/U} = \sum_{i=1}^{n}\omega_{i}(\dd u_{i,1}) \wedge \dots \wedge (\dd
u_{i,q_{i}}),
\]
o\`u pour tout $1 \leqslant i \leqslant n$, on a $\omega_{i} \in A^{p-q_{i},p -
q_{i}}(U)$ et $u_{i,1}, \dots, u_{i,q_{i}}$ sont \pshs\ et {\it 
localement born\'ees\/}
sur $U$. On note $\ov{A}^{p,p}(X)$ l'espace vectoriel complexe des formes
diff\'erentielles g\'en\'eralis\'ees de type $(p,p)$ sur $X$.
\end{defn}
\medskip

On pose \'egalement~:
\[
\ov{A}^{\ast}(X) = \bigoplus_{p \geqslant
0}\ov{A}^{p,p}(X), \qquad \ov{\ov{A}}^{\ast}(X) = \bigoplus_{p \geqslant 0}
\ov{\ov{A}}^{p,p}(X), \text{\quad et \quad} \ov{\ov{A}}_{\op{log}}^{\ast}(X) = \bigoplus_{p \geqslant 0}
\ov{\ov{A}}_{\op{log}}^{p,p}(X).
\]
\begin{rem}
On a imm\'ediatement les inclusions suivantes~:
\[
\ov{A}^{\ast}(X) \subset \ov{\ov{A}}^{\ast}(X) \subset \ov{\ov{A}}_{\op{log}}^{\ast}(X)
\subset D^{\ast}(X).
\]
\end{rem}
\begin{rem}
L'op\'erateur $\dd: D^{p,p}(X) \rightarrow D^{p+1,p+1}(X)$ induit par
restriction une application $\dd: \AAA^{p,p}_{\log}(X) 
\rightarrow \AAA^{p+1,p+1}_{\log}(X)$.
\end{rem}
\medskip

Soient $x \in \AA^{p,p}(X)$ et $y \in \AAA^{q,q}_{\op{log}}(X)$. Sur tout
ouvert $U \subset X$ assez petit, on peut \'ecrire~:
\begin{align*}
x &= \sum_{i = 1}^{n}\omega_{i}(\dd u_{i,1}) \wedge \dots \wedge (\dd
u_{i,r_{i}}) \\
\text{et \quad} y &= \sum_{j = 1}^{m}\eta_{j}v_{j,1}(\dd v_{j,2}) \wedge \dots
\wedge (\dd v_{j, s_{j}}), 
\end{align*}
o\`u pour tout $1 \leqslant i \leqslant n$, on a $\omega_{i} \in A^{p-r_{i}, p
- r_{i}}(U)$ et $u_{i,1}, \dots, u_{i, r_{i}}$ sont des fonctions \pshs\
localement born\'ees sur $U$; et o\`u pour tout $1 \leqslant j \leqslant m$,
on a $\eta_{j} \in A^{q -s_{j}+1, q - s_{j}+1}(U)$ et le multiplet $(v_{j,1}, \dots,
v_{j, s_{j}})$ est admissible sur $U$.

D'apr\`es le th\'eor\`eme (\ref{demailly}), l'expression~:
\[
\sum_{i=1}^{n}\sum_{j =1}^{m} \omega_{i}\eta_{j} v_{j,1}(\dd u_{i,1}) \wedge
\dots \wedge (\dd u_{i,r_{i}}) \wedge (\dd v_{j,2}) \wedge \dots \wedge (\dd
v_{j, s_{j}}), 
\]
a bien un sens et d\'efinit sur $U$ un courant de type $(p+q,p+q)$ que l'on
note provisoirement $[x\cdot y](U)$.
\begin{prop}
\label{produit_formes}
Le courant $[x\cdot y](U)$ d\'efini ci-dessus ne d\'epend que de $x$ et de $y$.
\end{prop}
\demo\ Soient $x'$ et $y'$ deux courants sur $U$ tels que $x' = x_{/U}$ et $y' =
y_{/U}$, et tels qu'on puisse \'ecrire~:
\begin{align*}
x' &= \sum_{i = 1}^{n'}\omega_{i}' (\dd u'_{i,1}) \wedge \dots \wedge (\dd
u'_{i, r'_{i}}) \\
\text{et \quad} y' &= \sum_{j = 1}^{m'}\eta'_{j} v'_{j,1} (\dd v'_{j,2}) \wedge
\dots \wedge (\dd v'_{j, s'_{j}}), 
\end{align*}
o\`u pour tout $1 \leqslant i \leqslant n'$, on a $\omega'_{i} \in A^{p-r'_{i}, p
- r'_{i}}(U)$ et $u'_{i,1}, \dots, u'_{i, r'_{i}}$ sont des fonctions \pshs\
localement born\'ees sur $U$; et o\`u pour tout $1 \leqslant j \leqslant m'$,
on a $\eta'_{j} \in A^{q -s'_{j}, q - s'_{j}}(U)$ et le multiplet $(v'_{j,1}, \dots,
v'_{j, s'_{j}})$ est admissible sur $U$.

On d\'efinit les courants~:
\begin{align*}
&[x'\cdot y](U) = \sum_{i=1}^{n'}\sum_{j =1}^{m} \omega'_{i}\eta_{j} 
v_{j,1}(\dd u'_{i,1}) \wedge
\dots \wedge (\dd u'_{i,r'_{i}}) \wedge (\dd v_{j,2}) \wedge \dots \wedge (\dd
v_{j, s_{j}}) \\
&\text{et~}\\
&[x'\cdot y'](U) = \sum_{i=1}^{n'}\sum_{j =1}^{m'} \omega'_{i}\eta'_{j} 
v'_{j,1}(\dd u'_{i,1}) \wedge
\dots \wedge (\dd u'_{i,r'_{i}}) \wedge (\dd v'_{j,2}) \wedge \dots \wedge (\dd
v'_{j, s'_{j}}).
\end{align*}
On veut d\'emontrer que le courant $\delta$ d\'efini sur $U$ par
l'\'egalit\'e~:
\[
\delta = [x'\cdot y'](U) - [x\cdot y](U) = ([x'\cdot y'](U) - [x'\cdot y](U)) +
([x'\cdot y](U) - [x\cdot y](U)), 
\]
est nul. On va prouver pour cela que $[x'\cdot y](U) - [x\cdot y](U)= 0$. On
d\'emontrerait de m\^eme que $[x'\cdot y'](U) - [x'\cdot y](U) = 0$.

Par r\'egularisation, on peut trouver pour tout $1 \leqslant j \leqslant m$ et
$1 \leqslant k \leqslant s_{j}$, une suite d\'ecroissante $(v_{j,k}^{(n)})_{n \in \M{N}}$
d'\'el\'ements de $\op{Psh}(U) \cap C^{\infty}(U)$ convergeant simplement sur
$U$ vers $v_{j,k}$. Pour tout $n \in \M{N}$, on note $y^{(n)}$ la forme
diff\'erentielle $C^{\infty}$ de type $(q,q)$ sur $U$ d\'efinie par~:
\[
y^{(n)} = \sum_{j=1}^{m}\eta_{j}v_{j,1}^{(n)}(\dd v_{j,2}^{(n)})\wedge \dots
\wedge (\dd v_{j,s_{j}}^{(n)}). 
\]
De l'\'egalit\'e $x = x'$, on d\'eduit que $x\cdot y^{(n)} = x' \cdot y^{(n)}$
pour tout $n \in \M{N}$, le produit \'etant le produit habituel d'un courant par
une forme diff\'erentielle. Or, d'apr\`es (\ref{demailly}), on a~:
\begin{alignat*}{3}
&x\cdot y^{(n)}& &\text{\quad tend vers \quad}& &[x\cdot y](U) \\
\text{et \quad} &x'\cdot y^{(n)}& & \text{\quad tend vers \quad}& &[x'\cdot
y](U), 
\end{alignat*}
au sens de la topologie faible. On en conclut que $[x'\cdot y](U) - [x\cdot
y](U) = 0$ et la proposition est d\'emontr\'ee.
\medskip

L'assertion suivante est une cons\'equence directe de ce qui pr\'ec\`ede.
\begin{prop}
\label{produit_generalise}
Soient $x \in \AA^{p,p}(X)$ et $y \in \AAA^{q,q}_{\log}(X)$. Il existe un
(unique) \'el\'ement de $\AAA_{\log}^{p+q,p+q}(X)$ que l'on notera $x\cdot y$
et dont la restriction \`a tout ouvert $U \subset X$ assez petit co\"\i ncide
avec le courant $[x\cdot y](U)$ d\'efini \`a la proposition
(\ref{produit_formes}).
\end{prop}
\begin{defn}
Soient $x \in \AA^{p,p}(X)$ et $y \in \AAA_{\log}^{q,q}(X)$. On appelle {\it
produit de $x$ et de $y$\/} le courant $x\cdot y \in \AAA_{\log}^{p+q,p+q}(X)$
d\'efini \`a la proposition pr\'ec\'edente.
\end{defn}
\medskip

Plus g\'en\'eralement, soient $x \in \AA^{\ast}(X)$ et $y \in
\AAA^{\ast}_{\log}(X)$, on appelle {\it produit de $x$ et de $y$\/} et l'on note
encore $x\cdot y$ le produit gradu\'e de $x$ et de $y$.

Les propositions suivantes sont des cons\'equences imm\'ediates des
d\'efinitions~:
\begin{prop}
\label{produit2}
Soient $x_{1} \in \AA^{p,p}(X)$ et $x_{2} \in \AA^{q,q}(X) \subset
\AAA^{q,q}_{\log}(X)$; on a~: $x_{1}\cdot x_{2} \in \AA^{p+q,p+q}(X)$.
\end{prop}
\begin{prop}
Soient $x_{1}$ et $x_{2}$ des \'el\'ements de $\AA^{\ast}(X)$ et $y \in
\AAA_{\log}^{\ast}(X)$; on a les relations suivantes~:
\begin{enumerate}
\item{$x_{1}\cdot (x_{2} \cdot y) = (x_{1}\cdot x_{2}) \cdot y \quad
\text{\rm (associativit\'e)}.$}
\item{$x_{1}\cdot x_{2} = x_{2} \cdot x_{1} \qquad \qquad 
\text{\rm (commutativit\'e)}.$}
\end{enumerate}
\end{prop}
\begin{prop}
Le produit d\'efini \`a la proposition $($\ref{produit2}$)$ munit $\AA^{\ast}(X)$
d'u\-ne structure d'alg\`ebre associative commutative unif\`ere et gradu\'ee 
et $\AAA_{\log}^{\ast}(X)$ d'une structure de $\AA^{\ast}(X)$-module, qui
\'etend sa structure usuelle de $A^{\ast}(X)$-module.
La restriction \`a $A^{\ast}(X)$ de la structure
d'alg\`ebre sur $\AA^{\ast}(X)$ co\"\i ncide avec la structure d'alg\`ebre
usuelle sur $A^{\ast}(X)$.
\end{prop}

Soient $X$ et $Y$ deux vari\'et\'es complexes et $f: Y \rightarrow X$ un
morphisme lisse de vari\'et\'es analytiques (i.e. une submersion holomorphe). 
On consid\`ere $U \subset X$ et $V \subset Y$ deux ouverts
tels que $f(V) \subset U$. \\
Un r\'esultat classique (voir par exemple
\cite{7}, th. 1.5.9), affirme que si $ u \in \op{Psh}(U)$ alors 
$f^{\ast}(u) = u \circ f \in \op{Psh}(V)$. Plus
g\'en\'eralement, on a le r\'esultat suivant~:
\begin{prop}
Soit $(u_{1}, \dots, u_{q})$ un $q$-uplet admissible de fonctions r\'eelles sur
$U$; le $q$-uplet $(u_{1}\circ f, \dots, u_{q}\circ f)$ est admissible sur $V$.
\end{prop}
\demo\ D'apr\`es la remarque ci-dessus, on sait que les fonctions $u_{1}\circ
f, \dots, u_{q} \circ f$ sont \pshs\ sur $V$. Pour tout $1 \leqslant i
\leqslant q$, on a~: $L(u_{i}\circ f) = f^{-1}(L(u_{i})) \subset
f^{-1}(A_{i})$. Enfin, du fait de la lissit\'e de $f$,
les ensembles $f^{-1}(A_{i})$ sont analytiques 
et on a pour tout choix
d'indices $j_{1} < \dots < j_{m}$ dans $\{1, \dots,q\}$~:
\[
\op{codim} \, f^{-1}(A_{j_{1}}) \cap \dots \cap f^{-1}(A_{j_{m}}) =
\op{codim} A_{j_{1}} \cap \dots \cap A_{j_{m}}
\geqslant m ,
\]
d\`es que les intersections consid\'er\'ees sont non vides.
\medskip

La proposition suivante montre que l'image inverse du courant 
$u_{1}\dd u_{2} \wedge \dots
\wedge \dd u_{q}$ s'exprime simplement~:
\begin{prop}
\label{pullback_formes}
Soit $(u_{1},\dots,u_{q})$ un $q$-uplet admissible de fonctions r\'eelles sur
$U$; on a sur $V$ l'\'egalit\'e~:
\[
f^{\ast}(u_{1}\dd u_{2} \wedge \dots \wedge \dd u_{q})
= 
(u_{1}\circ f)(\dd u_{2}\circ f) \wedge \dots \wedge (\dd u_{q} \circ f).
\]
\end{prop}
\demo\ Par r\'egularisation, on peut trouver pour tout $1 \leqslant i \leqslant
q$ une suite d\'ecroissante $\left(u_{i}^{(n)}\right)_{n \in \M{N}}$
d'\'el\'ements de $\op{Psh}(U) \cap C^{\infty}(U)$ convergeant simplement sur
$U$ vers $u_{i}$. Pour tout $n \in \M{N}$ on note $x^{(n)}$ la forme
diff\'erentielle $C^{\infty}$ de type $(q-1,q-1)$ sur $U$ d\'efinie par~:
\[
x^{(n)} = u_{1}^{(n)}\dd u_{2}^{(n)} \wedge \dots \wedge \dd u_{q}^{(n)}.
\]
D'apr\`es (\ref{demailly}) la suite $x^{(n)}$ (resp. la suite
$f^{\ast}(x^{(n)})$) tend vers $u_{1}\dd u_{2} \wedge \dots \wedge \dd u_{q}$
(resp. vers $(u_{1}\circ f) (\dd u_{2}\circ f) \wedge \dots \wedge (\dd u_{q}
\circ f)$) au sens de la convergence faible des courants quand $n$ tend vers
$+\infty$. Le r\'esultat d\'ecoule alors de la continuit\'e faible de
$f^{\ast}$.
\medskip

On d\'eduit de ce r\'esultat que toute application analytique lisse 
$f : Y \rightarrow X$ induit un
morphisme de groupes $f^{\ast}: \AAA^{p,p}_{\log}(X) \rightarrow
\AAA_{\log}^{p,p}(Y)$, qui induit un morphisme gradu\'e
$f^{\ast}: \AAA_{\log}^{\ast}(X) \rightarrow
\AAA_{\log}^{\ast}(Y)$.\\
Les deux propositions suivantes sont des cons\'equences directes de
(\ref{pullback_formes})~:
\begin{prop}
On a~: $f^{\ast}\left(\AAA^{\ast}(X)\right) \subset \AAA^{\ast}(Y)$ et
$f^{\ast}\left(\AA^{\ast}(X)\right) \subset \AA^{\ast}(Y)$.
\end{prop}
\begin{prop}
\label{prop_formes1}
Soient $x \in \ov{A}^{\ast}(X)$ et $y \in \AAA^{\ast}_{\log}(X)$, on a~: 
\[
f^{\ast}(x\cdot y) = f^{\ast}(x)\cdot f^{\ast}(y).
\]
En particulier, le morphisme $f^{\ast}: \ov{A}^{\ast}(X) \rightarrow \ov{A}^{\ast}(Y)$
est un morphisme d'alg\`ebres.
\end{prop}
\medskip

Enfin les deux propositions suivantes sont valables en toute g\'en\'eralit\'e~:
\begin{prop}
\label{prop_formes2}
Soient $g: Z \rightarrow Y$ et $f: Y \rightarrow X$ deux morphismes lisses de
vari\'et\'es complexes, on a l'identit\'e~: $(f \circ g)^{\ast} = g^{\ast}
\circ f^{\ast}$.
\end{prop}

\begin{prop}
\label{prop_formes3}
Soient $\alpha \in A_{c}^{\ast}(Y)$ et $x \in \AAA_{\log}^{\ast}(X)$, on a~:
\[
f_{\ast}\left(\alpha \cdot f^{\ast}(x)\right) = f_{\ast}(\alpha)\cdot x, 
\]
l'image directe $f_{\ast}\left(\alpha \cdot f^{\ast}(x)\right)$ \'etant prise
au sens des courants.
\end{prop}
\bigskip

\subsection{Convergence au sens de Bedford-Taylor et image inverse}~

On introduit dans ce paragraphe de nouvelles classes remarquables de courants.

\begin{defn}
\label{def_formeB_1}
Soit $p$ un entier positif et $Z$ un cycle de codimension $q$ de $X$. On note
$B^{p,p}(X)$ (resp. $B_{\log}^{p,p}(X)$, resp. $B_{Z}^{p,p}(X)$, resp.
$B_{Z,\log}^{p,p}(X)$) l'espace vectoriel complexe form\'e des \'el\'ements de
$D^{p,p}(X)$, qui, sur tout ouvert $U$ d'un recouvrement suffisamment fin de
$X$, peuvent s'\'ecrire sous la forme~:
\begin{align*}
&\sum_{i=1}^{n}\omega_{i}(\dd u_{i,1}) \wedge \dotsm \wedge (\dd u_{i,q_{i}}), \\
\bigg(resp. \qquad
&\sum_{i=1}^{n}\omega_{i}u_{i,1}(\dd u_{i,2}) \wedge \dotsm \wedge (\dd u_{i,q_{i}})
\bigg), \\
\bigg(resp. \qquad
&\sum_{i=1}^{n}\omega_{i}(\dd u_{i,1}) \wedge \dotsm \wedge (\dd u_{i,q_{i}})
\wedge \delta_{Z}, 
\bigg), \\
\bigg(resp. \qquad
&\sum_{i=1}^{n}\omega_{i}u_{i,1}(\dd u_{i,2}) \wedge \dotsm \wedge (\dd u_{i,q_{i}})
\wedge \delta_{Z}, 
\bigg),
\end{align*}
o\`u pour tout $1 \leqslant i \leqslant n$, on a $\omega_{i} \in
A^{p-q_{i},p-q_{i}}(U)$ (resp. $\omega_{i} \in A^{p-q_{i}+1,p-q_{i}+1}(U)$, 
resp. $\omega_{i} \in A^{p-q-q_{i},p-q-q_{i}}(U)$,
resp. $\omega_{i} \in A^{p-q-q_{i}+1,p-q-q_{i}+1}(U)$) et $u_{i,1},\dots,
u_{i,q_{i}}$ sont des fonctions \pshs\ {\em continues} sur $U$.
\end{defn}
\begin{rem}
\label{rem_BT1}
On a imm\'ediatement les inclusions $B^{p,p}(X) \subset B_{\log}^{p,p}(X)$,\\
$B_{Z}^{p,p}(X) \subset B_{Z,\log}^{p,p}(X)$ et l'\'egalit\'e
$B_{\log}^{p,p}(X) = B_{X,\log}^{p,p}(X)$.
\end{rem}
En reprenant les hypoth\`eses et les notations de la d\'efinition
pr\'ec\'edente, on pose \'egalement~:
\begin{defn}
On note $B_{0}^{p,p}(X)$ (resp. $B_{\log,0}^{p,p}(X)$, 
resp. $B_{Z,0}^{p,p}(X)$, \\
resp. $B_{Z,\log,0}^{p,p}(X)$) le sous-espace
vectoriel complexe de $B^{p,p}(X)$ 
(resp. de $B_{\log}^{p,p}(X)$, 
resp. de $B_{Z}^{p,p}(X)$, resp. de $B_{Z,\log}^{p,p}(X)$) form\'e des
\'el\'ements qui, sur tout ouvert $U$ d'un recouvrement suffisamment fin de
$X$, peuvent s'\'ecrire sous la forme~:
\begin{align*}
&\sum_{i=1}^{n}\omega_{i}(\dd u_{i,1}) \wedge \dotsm \wedge (\dd u_{i,q_{i}}), \\
\bigg(resp. \qquad
&\sum_{i=1}^{n}\omega_{i}u_{i,1}(\dd u_{i,2}) \wedge \dotsm \wedge (\dd u_{i,q_{i}})
\bigg), \\
\bigg(resp. \qquad
&\sum_{i=1}^{n}\omega_{i}(\dd u_{i,1}) \wedge \dotsm \wedge (\dd u_{i,q_{i}})
\wedge \delta_{Z}, 
\bigg), \\
\bigg(resp. \qquad
&\sum_{i=1}^{n}\omega_{i}u_{i,1}(\dd u_{i,2}) \wedge \dotsm \wedge (\dd u_{i,q_{i}})
\wedge \delta_{Z}, 
\bigg),
\end{align*}
o\`u les $u_{i,j}$ et les $\omega_{i}$ sont comme \`a la d\'efinition
(\ref{def_formeB_1})
et o\`u de plus,
pour tout $1 \leqslant i \leqslant n$, la forme $\omega_{i}$ est {\em
ferm\'ee}.
\end{defn}
\begin{rem}
On a les inclusions $B_{0}^{p,p}(X) \subset B_{\log,0}^{p,p}(X) =
B_{X,\log,0}^{p,p}(X)$ et $B_{Z,0}^{p,p}(X) \subset B_{Z,\log,0}^{p,p}(X)$.
\end{rem}
\begin{defn}
\label{conv_BT1}
Soit $(x^{(k)})_{k \in \M{N}}$ une suite d'\'el\'ements de $B_{Z,\log}^{p,p}(X)$.
On dit que $(x^{(k)})_{k \in \M{N}}$ converge vers $x \in B_{Z,\log}^{p,p}(X)$ {\it
au sens de Bedford-Taylor\/} (en abr\'eg\'e {\it au sens BT\/}) 
si pour tout ouvert $U$ d'un recouvrement
suffisamment fin de $X$, on peut \'ecrire~:
\begin{align*}
x^{(k)} & = 
\sum_{i=1}^{n}\omega^{(k)}u_{i,1}^{(k)}(\dd u_{i,2}^{(k)}) \wedge \dotsm 
\wedge (\dd u_{i,q_{i}}^{(k)})\wedge \delta_{Z} \\
\text{et} \qquad 
x &= 
\sum_{i=1}^{n}\omega u_{i,1}(\dd u_{i,2}) \wedge \dotsm 
\wedge (\dd u_{i,q_{i}})\wedge \delta_{Z},
\end{align*}
o\`u $n$ est ind\'ependant de $k$ et o\`u, pour tout $1 \leqslant 
i \leqslant n$ et tout $k \in \M{N}$, on a
$\omega_{i}^{(k)} \in A^{p-q-q_{i}+1,p-q-q_{i}+1}(U)$, 
$\omega_{i} \in A^{p-q-q_{i}+1,p-q-q_{i}+1}(U)$ et $u_{i,1}^{(k)},
\dots,u_{i,q_{i}}^{(k)}, u_{i,1}, \dots,u_{i,q_{i}}$ sont des fonctions \pshs\
{\em continues} sur $U$; que $\omega_{i}^{(k)}$ tend vers
$\omega_{i}$ 
au sens des formes $C^{\infty}$ quand $k$ tend vers $+\infty$ et que de plus $(u_{i,1}^{(k)})_{k
\in \M{N}}, \dots,(u_{i,q_{i}}^{(k)})_{k\in \M{N}}$ convergent uniform\'ement sur $U$
vers
$u_{i,1}, \dots, u_{i,q_{i}}$ respectivement quand $k$ tend vers $+\infty$.
\end{defn}
\begin{defn}
\label{conv_BT2}
Soit $(x^{(k)})_{k \in \M{N}}$ une suite d'\'el\'ements de $B_{Z}^{p,p}(X)$.
On dit que $(x^{(k)})_{k \in \M{N}}$ converge vers $x \in B_{Z}^{p,p}(X)$ {\it
au sens BTR\/}, si pour tout ouvert $U$ d'un recouvrement
suffisamment fin de $X$, on peut \'ecrire~:
\begin{align*}
x^{(k)} & = 
\sum_{i=1}^{n}\omega^{(k)}(\dd u_{i,1}^{(k)}) \wedge \dotsm 
\wedge (\dd u_{i,q_{i}}^{(k)})\wedge \delta_{Z} \\
\text{et} \qquad 
x &= 
\sum_{i=1}^{n}\omega (\dd u_{i,1}) \wedge \dotsm 
\wedge (\dd u_{i,q_{i}})\wedge \delta_{Z},
\end{align*}
o\`u $n$ est ind\'ependant de $k$ et o\`u, pour tout $1 \leqslant 
i \leqslant n$ et tout $k \in \M{N}$, on a
$\omega_{i}^{(k)} \in A^{p-q-q_{i},p-q-q_{i}}(U)$, 
$\omega_{i} \in A^{p-q-q_{i},p-q-q_{i}}(U)$ et $u_{i,1}^{(k)},
\dots,u_{i,q_{i}}^{(k)}, u_{i,1}, \dots,u_{i,q_{i}}$ sont des fonctions \pshs\
{\em continues} sur $U$; que $\omega_{i}^{(k)}$ tend vers
$\omega_{i}$ 
au sens des formes $C^{\infty}$ quand $k$ tend vers $+\infty$ et que de plus $(u_{i,1}^{(k)})_{k
\in \M{N}}, \dots,(u_{i,q_{i}}^{(k)})_{k\in \M{N}}$ convergent uniform\'ement sur $U$
vers
$u_{i,1}, \dots, u_{i,q_{i}}$ respectivement quand $k$ tend vers $+\infty$.
\end{defn}
En reprenant les hypoth\`eses et les notations de la d\'efinition
pr\'ec\'edente,
on pose \'egalement~:
\begin{defn}
\label{conv_BT3}
Soit $(x^{(k)})_{k \in \M{N}}$ une suite d'\'el\'ements de $B_{Z,\log,0}^{p,p}(X)$.
On dit que $(x^{(k)})_{k \in \M{N}}$ converge vers $x \in B_{Z,\log,0}^{p,p}(X)$
{\it fortement au sens BT\/} si pour tout ouvert $U$ d'un recouvrement
suffisamment fin de $X$, il existe $n \in \M{N}^{\ast}$
tel que pour tout $k \in \M{N}$ on puisse \'ecrire~:
\begin{align*}
x^{(k)} & = 
\sum_{i=1}^{n}\omega^{(k)}u_{i,1}^{(k)}(\dd u_{i,2}^{(k)}) \wedge \dotsm 
\wedge (\dd u_{i,q_{i}}^{(k)})\wedge \delta_{Z} \\
\text{et} \qquad 
x &= 
\sum_{i=1}^{n}\omega u_{i,1}(\dd u_{i,2}) \wedge \dotsm 
\wedge (\dd u_{i,q_{i}})\wedge \delta_{Z},
\end{align*}
o\`u pour tout $1 \leqslant i \leqslant n$, $\omega_{i}^{(k)}$ tend vers
$\omega_{i}$ dans $A^{\ast}(U)$ quand $k$ tend vers $+\infty$ et
$u_{i,1}^{(k)},\dots,u_{i,q_{i}}^{(k)}$ converge uniform\'ement sur $U$ vers $u_{i,1},
\dots, u_{i,q_{i}}$ respectivement quand $k$ tend vers $+\infty$, et o\`u pour
tout $1 \leqslant i \leqslant n$ et tout $k \in \M{N}$ les formes
$\omega_{i}^{(k)}$ et $\omega_{i}$ sont {\em ferm\'ees}.
\end{defn}
\begin{defn}
\label{conv_BT4}
Soit $(x^{(k)})_{k \in \M{N}}$ une suite d'\'el\'ements de $B_{Z,0}^{p,p}(X)$.
On dit que $(x^{(k)})_{k \in \M{N}}$ converge vers $x \in B_{Z,0}^{p,p}(X)$
{\it fortement au sens BTR\/} si pour tout ouvert $U$ d'un recouvrement
suffisamment fin de $X$, il existe $n \in \M{N}^{\ast}$
tel que pour tout $k \in \M{N}$ on puisse \'ecrire~:
\begin{align*}
x^{(k)} & = 
\sum_{i=1}^{n}\omega^{(k)}(\dd u_{i,1}^{(k)}) \wedge \dotsm 
\wedge (\dd u_{i,q_{i}}^{(k)})\wedge \delta_{Z} \\
\text{et} \qquad 
x &= 
\sum_{i=1}^{n}\omega (\dd u_{i,1}) \wedge \dotsm 
\wedge (\dd u_{i,q_{i}})\wedge \delta_{Z},
\end{align*}
o\`u pour tout $1 \leqslant i \leqslant n$, $\omega_{i}^{(k)}$ tend vers
$\omega_{i}$ dans $A^{\ast}(U)$ quand $k$ tend vers $+\infty$ et
$u_{i,1}^{(k)},\dots,u_{i,q_{i}}^{(k)}$ converge uniform\'ement sur $U$ vers $u_{i,1},
\dots, u_{i,q_{i}}$ respectivement quand $k$ tend vers $+\infty$, et o\`u pour
tout $1 \leqslant i \leqslant n$ et tout $k \in \M{N}$ les formes
$\omega_{i}^{(k)}$ et $\omega_{i}$ sont {\em ferm\'ees}.
\end{defn}
\begin{rem}
D'apr\`es le th\'eor\`eme (\ref{lim_uniforme}) les quatre notions de convergence
introduites aux
d\'efinitions (\ref{conv_BT1}), (\ref{conv_BT2}), (\ref{conv_BT3}) et
(\ref{conv_BT4}) sont plus fortes que la convergence faible au sens des
courants.
\end{rem} 
\begin{rem}
L'application $\dd: B_{\log,0}^{p,p}(X) \rightarrow B_{0}^{p+1,p+1}(X)$
envoie l'ensemble des suites convergeant fortement au sens BT 
dans l'ensemble des suites convergeant fortement au sens BTR.
\end{rem}
\begin{rem}
L'application $\dd: B^{p,p}(X) \rightarrow B_{0}^{p+1,p+1}(X)$ envoie 
l'ensemble des suites convergeant au sens BTR 
dans l'ensemble des suites convergeant fortement au sens BTR.
\end{rem}
\begin{defn}
\label{formes_adherentes}
On note $C^{p,p}(X)$ (resp. $C^{p,p}_{\log}(X)$) l'ensemble 
des limites des suites de $A^{p,p}(X)$ convergeant 
dans $B^{p,p}(X)$ (resp. dans $B_{\log}^{p,p}(X)$) au sens BTR (resp. au sens
BT).

On note $C^{p,p}_{0}(X)$ (resp. $C^{p,p}_{\log,0}(X)$)l'ensemble 
des limites des suites de $A^{p,p}(X)$ convergeant fortement
dans $B_{0}^{p,p}(X)$ (resp. dans $B_{\log,0}^{p,p}(X)$) au sens BTR
(resp. au sens BT).
\end{defn}

Soient $x \in B^{p,p}(X)$ et $y \in B^{r,r}_{Z,\log}(X)$. 
Sur un ouvert $U \subset X$ assez petit, on peut \'ecrire~:
\begin{align*}
x &= \sum_{i=1}^{n}\omega_{i}(\dd u_{i,1}) \wedge \dotsm \wedge (\dd
u_{i,q_{i}}) \\
\text{et}\qquad y &= \sum_{j=1}^{m}\eta_{j}v_{j,1}(\dd v_{j,2}) \wedge \dotsm
\wedge (\dd v_{j,r_{j}})\wedge \delta_{Z}, 
\end{align*}
o\`u pour tout $1\leqslant i \leqslant n$, $\omega_{i} \in
A^{p-q_{i},p-q_{i}}(U)$ et $u_{i,1},\dots,u_{i,q_{i}}$ sont des fonctions
\pshs\ continues sur $U$; et o\`u pour tout $1 \leqslant j \leqslant n$, on a
$\eta_{j} \in A^{p-q-r_{j}+1,p-q-r_{j}+1}(U)$ et $v_{j,1}, \dots , v_{j,r_{j}}$
sont des fonctions \pshs\ continues sur $U$. A partir de ces donn\'ees, on
d\'efinit un courant de $B_{Z,\log}^{p+r,p+r}(U)$ que l'on note provisoirement $[x\cdot
y](U)$, par la formule~:
\begin{multline*}
[x\cdot y](U) = \\
\sum_{i=1}^{n}\sum_{j=1}^{m}
\omega_{i}\eta_{j} v_{j,1} (\dd u_{i,1}) \wedge \dotsm \wedge 
(\dd u_{i,q_{i}}) \wedge 
(\dd v_{j,2}) \wedge \dotsm \wedge (\dd v_{j,r_{j}})\wedge \delta_{Z}.
\end{multline*}
\begin{prop}
\label{produit_uniforme}
Le courant $[x\cdot y](U)$ ainsi d\'efini ne d\'epend que de $x$ et de $y$.
\end{prop}
\demo\ Par lin\'earit\'e, on se ram\`ene au cas o\`u $Z$ est effectif. Le
probl\`eme \'etant local, on se restreint \`a un ouvert $U' \subset U$ tel
qu'il existe une suite $(\delta_{Z}^{(n)})_{n \in \M{N}}$ de formes $C^{\infty}$
ferm\'ees et positives sur $U'$ convergeant faiblement au sens des courants
vers $\delta_{Z}$. On suit alors {\it mutatis mutandis\/} la d\'emonstration de
la proposition (\ref{produit_formes}) en utilisant le th\'eor\`eme (\ref{lim_uniforme}).
\medskip

On d\'eduit imm\'ediatement de ce qui pr\'ec\`ede~:
\begin{prop_defn}
\label{produit_uniforme2}
Soient $x \in B^{p,p}(X)$ et $y \in B_{Z,\log}^{r,r}(X)$. Il existe un unique
\'el\'ement de $B_{Z,\log}^{p+r,p+r}(X)$ que l'on note $x\cdot y$ et qu'on
appelle produit de $x$ et de $y$, dont la restriction \`a tout ouvert $U
\subset X$ assez petit co\"\i ncide avec le courant $[x\cdot y](U)$ d\'efini
\`a la proposition (\ref{produit_uniforme}).
\end{prop_defn}

Plus g\'en\'eralement, soient $x \in B^{\ast}(X)$ et $y \in
B_{Z,\log}^{\ast}(X)$, on appelle {\it produit de x et de y\/} et on note
encore $x\cdot y$ le produit gradu\'e de $x$ et de $y$.
\begin{rem}
On d\'efinit de fa\c con similaire un produit~:
\[
(~\cdot~): B_{\log}^{\ast}(X) \times B_{Z}^{\ast}(X) \longrightarrow
B_{Z,\log}^{\ast}(X).
\]
\end{rem}
\begin{rem}
Du fait de l'inclusion $B^{\ast}(X) \subset B_{\log}^{\ast}(X) =
B_{X,\log}^{\ast}(X)$, on dispose d'un produit $(~\cdot~): B^{\ast}(X) \times
B^{\ast}_{\log}(X) \rightarrow B_{\log}^{\ast}(X)$. 
Ce produit co\"\i ncide avec le produit d\'efini \`a la proposition
(\ref{produit_generalise}).
\end{rem}
La proposition suivante est une cons\'equence directe des d\'efinitions~:
\begin{prop}
Le produit d\'efini \`a la proposition (\ref{produit_uniforme2}) munit $B^{\ast}(X)$ d'une
structure d'alg\`ebre associative commutative unif\`ere et gradu\'ee et
$B_{Z,\log}^{\ast}(X)$ d'une structure de $B^{\ast}(X)$-module qui \'etend sa
structure usuelle de $A^{\ast}(X)$-module.
\end{prop}
La proposition suivante montre que le produit d\'efini \`a la proposition
(\ref{produit_uniforme2}) se comporte bien vis-\`a-vis des convergences
ordinaires ou fortes au sens BT et BTR~:
\begin{prop}
Le produit~:
\[
(~\cdot~): B^{\ast}(X) \times B_{Z,\log}^{\ast}(X) \longrightarrow 
B_{Z,\log}^{\ast}(X), 
\]
est compatible avec les convergences au sens BTR et BT sur $B^{\ast}(X)$
et $B_{Z,\log}^{\ast}(X)$ respectivement. Sa restriction~:
\[
(~\cdot~): B_{0}^{\ast}(X) \times B_{Z,\log,0}^{\ast}(X) \longrightarrow 
B_{Z,\log,0}^{\ast}(X),
\]
est compatible avec les convergences fortes au sens BTR et BT sur $B^{\ast}_{0}(X)$
et $B_{Z,\log,0}^{\ast}(X)$ respectivement.
\end{prop}
\demo\ C'est une cons\'equence imm\'ediate des d\'efinitions.
\medskip

\begin{rem}
Soient $\varphi \in C_{\log,0}^{p,p}(X)$, $\vfi' \in C_{\log,0}^{q,q}(X)$ et
$\vfi'' \in C_{0}^{r,r}(X)$, on a $\vfi \dd \vfi' \in
C_{\log,0}^{p+q+1,p+q+1}(X)$ et $\vfi \vfi'' \in C_{\log,0}^{p+r,p+r}(X)$.
\end{rem}

\begin{thm}
\label{th_image_inverse}
Soit $f: X \rightarrow Y$ un morphisme de vari\'et\'es projectives complexes.
Il existe un unique morphisme de groupes gradu\'es~:
\[
f^{\ast}: B^{\ast}_{\log}(Y) \longrightarrow B_{\log}^{\ast}(X), 
\]
tel que les assertions suivantes soient v\'erifi\'ees~:
\begin{enumerate}
\item{$f^{\ast}$ restreint \`a $A^{\ast}(Y)$ co\"\i ncide avec l'application
image inverse usuelle.}
\item{$f^{\ast}$ respecte la convergence au sens BT.}
\item{On a $f^{\ast}(B^{\ast}(Y)) \subset B^{\ast}(X)$,
$f^{\ast}(B_{0}^{\ast}(Y)) \subset B_{0}^{\ast}(X)$ et
$f^{\ast}(B_{\log,0}^{\ast}(Y)) \subset B_{\log,0}^{\ast}(Y)$.}
\item{$f^{\ast}$ restreint \`a $B^{\ast}(Y)$ (resp. \`a $B_{0}^{\ast}(Y)$,
resp. \`a $B_{\log,0}^{\ast}(Y)$) respecte la convergence au sens BTR (resp. la
convergence forte au sens BTR, resp. la convergence forte au sens BT).}
\item{Si $x \in B^{\ast}(Y)$ et $y \in B_{\log}^{\ast}(Y)$, alors on a~:
\[
f^{\ast}(x\cdot y) = f^{\ast}(x)\cdot f^{\ast}(y).
\]
}
\item{Soit $g : V \rightarrow X$ un autre morphisme de vari\'et\'es projectives
complexes. On a~: $(f \circ g)^{\ast} = g^{\ast}\circ f^{\ast}$.}
\end{enumerate}
\end{thm}
\demo\ 
Le morphisme $f: X \rightarrow Y$ peut se factoriser comme la composition d'une
immersion ferm\'ee $i: X \hookrightarrow \M{P}^{n}_{Y}$ pour $n$ assez grand,
et de la projection $\pi: \M{P}_{Y}^{n} \rightarrow Y$ qui est un morphisme
{\em lisse}.
On va d\'efinir $\pi^{\ast}$ et $i^{\ast}$, puis poser $f^{\ast} = i^{\ast}
\circ \pi^{\ast}$. 

Comme $\pi$ est lisse, $\pi^{\ast}: D^{\ast}(Y) \rightarrow D^{\ast}(X)$ est
d\'efini en toute g\'en\'eralit\'e et sa restriction \`a $B_{\log}^{\ast}(Y)$
v\'erifie les assertions $(1)$ \`a $(6)$ d'apr\`es les propositions (\ref{pullback_formes}), 
(\ref{prop_formes1}),
(\ref{prop_formes2}) et (\ref{prop_formes3}).

On donne \`a pr\'esent une construction de $i^{\ast}$. Soit $ x \in
B_{\log}^{\ast}(Y)$; sur tout ouvert $U$ d'un recouvrement suffisamment fin de
$Y$, on peut \'ecrire~:
\[
x = \sum_{j = 1}^{n}\omega_{j}u_{j,1}(\dd u_{j,2})\wedge \dotsm \wedge (\dd
u_{j,q_{j}}),
\]
o\`u pour tout $1 \leqslant j \leqslant n$, $\omega_{j} \in A^{\ast}(U)$ et
$u_{j,1}, \dots , u_{j,q_{j}}$ sont des fonctions \pshs\ {\em continues} sur
$U$. On pose~:
\[
i^{\ast}(x)(i^{-1}(U)) = 
\sum_{i=1}^{n}
i^{\ast}(\omega_{j})(u_{j,1}\circ i)(\dd u_{j,2}\circ i)\wedge \dotsm \wedge
(\dd u_{j,q_{j}}\circ i) \in B_{\log}^{\ast}(i^{-1}(U)).
\]
Le morphisme image directe $i_{\ast}: D^{\ast}(X) \rightarrow 
D^{\ast}(Y)$ \'etant injectif, le courant $i^{\ast}(x)$ est l'unique
courant tel que~:
\[
i_{\ast}(i^{\ast}(x)) = x\cdot \delta_{i(X)}, 
\]
ce qui montre que le
morphisme $i^{\ast}$ est bien d\'efini. On d\'eduit ais\'ement des
d\'efinitions que $i^{\ast}$ v\'erifie les assertions $(1)$ \`a $(6)$.

Les morphismes $\pi^{\ast}$ et $i^{\ast}$ \'etant d\'efinis, on pose $f^{\ast}
= i^{\ast}\circ \pi^{\ast}$. Par composition, $f^{\ast}$ v\'erifie les
assertions $(1)$ \`a $(5)$. Il faut montrer que $f^{\ast}$ ne d\'epend pas de la
d\'ecomposition $f = \pi \circ i $ utilis\'ee pour le d\'efinir. 
Pour cela, on montre que les assertions $(1)$ et $(2)$ caract\'erisent
$f^{\ast}$ de mani\`ere unique.

Soit $x \in B_{\log}^{\ast}(X)$. En utilisant une partition de l'unit\'e
associ\'ee \`a un recouvrement suffisamment fin de $Y$, on peut supposer que le
support de $x$ est contenu dans un ouvert $U \subset X$ tel que l'on puisse
\'ecrire~:
\[
x = \sum_{j=1}^{n}\omega_{j}u_{j,1}(\dd u_{j,2})\wedge \dotsm \wedge (\dd
u_{j,q_{j}}), 
\]
o\`u pour tout $1 \leqslant j \leqslant n$, $\omega_{j}$ est une forme
$C^{\infty}$ dont le support est contenu dans $U$ et
$u_{j,1},\dots,u_{j,q_{j}}$ sont des fonctions \pshs\ continues sur $U$.
Par r\'egularisation, on peut trouver pour tout $1 \leqslant j \leqslant n$ des
suites $(u_{j,1}^{(k)})_{k \in \M{N}}, \dots, $$(u_{j,q_{j}}^{(k)})_{k \in \M{N}}$ de fonctions \pshs\ $C^{\infty}$ convergeant
uniform\'ement vers $u_{j,1},\dots,$$u_{j,q_{j}}$ respectivement sur $U$.

Posons, pour tout $k \in \M{N}$, 
\[
x^{(k)} = \sum_{j=1}^{n} \omega_{j}u_{j,1}^{(k)}(\dd u_{j,2}^{(k)})\wedge
\dotsm \wedge (\dd u_{j,q_{j}}^{(k)}).
\]
On d\'eduit de l'assertion (2) que la suite $(f^{\ast}(x^{(k)}))_{k \in \M{N}}$
converge faiblement vers $f^{\ast}(x)$. Comme d'autre part les formes images
inverses $f^{\ast}(x^{(k)})$ 
pour $k \in \M{N}$ ne d\'ependent pas, d'apr\`es l'assertion $(1)$, 
de la d\'ecomposition choisie, il en est de m\^eme pour $f^{\ast}(x)$.

L'assertion $(6)$ se montre de fa\c con similaire.
\medskip

\begin{rem}
Les op\'erateurs $f^{\ast}$ et $\dd$ commutent.
\end{rem}
\begin{rem}
Si $\vfi \in B_{\log,0}^{0}(X)$, alors $f^{\ast}(\vfi) = \vfi \circ f$.
\end{rem}
\begin{rem}
Si $i: X \hookrightarrow Y$ est une immersion ferm\'ee et $\vfi \in
B_{\log}^{\ast}(Y)$, on d\'eduit de la construction donn\'ee au th\'eor\`eme
(\ref{th_image_inverse}) que~:
\[
i_{\ast}(i^{\ast}(\vfi)) = \vfi \cdot \delta_{i(X)}.
\]
\end{rem}
\bigskip

\subsection{M\'etriques admissibles}~

Dans cet article, on appellera
{\it fibr\'e en droites hermitien\/} sur $X$ tout
couple $\ov{L} = (L,\|.\|)$ form\'e d'un fibr\'e en droites holomorphe $L$ sur $X$ et d'une
m\'etrique hermitienne {\em continue} sur $L$. 

On consid\`ere $U \subset X$ un
ouvert et $s_{1}$ et $s_{2}$ deux sections holomorphes de $L$ ne s'annulant pas sur
$U$. Il existe alors une fonction $h$ holomorphe et ne s'annulant pas sur $U$
telle que $s_{2} = h\cdot s_{1}$. On en d\'eduit que~:
\[
\dd (- \log \|s_{2}\|^{2}) = - \dd \log |h|^{2} + \dd ( - \log \|s_{1}\|^{2}) =
\dd ( - \log \|s_{1}\|^{2}).
\]
Cette remarque justifie la d\'efinition suivante~:
\begin{defn}
\label{classe_chern}
Soit $\ov{L} = (L,\|.\|)$ un fibr\'e en droites hermitien sur $X$. On appelle
{\it premier courant de Chern\/} de $\ov{L}$ et on note $c_{1}(\ov{L}) \in
D^{1,1}(X)$ le courant d\'efini localement par l'\'egalit\'e~:
\[
c_{1}(\ov{L}) = \dd ( - \log \|s\|^{2}), 
\]
o\`u $s$ est une section locale holomorphe et ne s'annulant pas du fibr\'e
$L$.
\end{defn}
\begin{prop}
\label{additivite_chern}
Soient $\ov{L}_{1}$ et $\ov{L}_{2}$ des fibr\'es en droites hermitiens sur $X$.
On a la relation~:
\[
c_{1}(\ov{L}_{1}\otimes \ov{L}_{2}) = c_{1}(\ov{L}_{1}) + c_{1}(\ov{L}_{2}).
\]
\end{prop}
\demo\ C'est une cons\'equence imm\'ediate de la d\'efinition
(\ref{classe_chern}).
\medskip

L'\'enonc\'e suivant est une extension imm\'ediate de la formule de
Poincar\'e-Lelong classique au cas des fibr\'es hermitiens quelconques~:
\begin{thm}{\bf \rm (formule de Poincar\'e-Lelong
g\'en\'eralis\'ee).}~
\label{PL_generalisee}
Soit $\ov{L} = (L,\|.\|)$ un fibr\'e en droites hermitien sur
$X$ et $s$ une section m\'eromorphe de $L$ sur $X$ non identiquement nulle sur
chaque composante connexe de $X$. On a l'\'egalit\'e entre
courants~:
\[
\dd (- \log \|s\|^{2}) + \delta_{\op{div}s} = c_{1}(\ov{L}).
\]
\end{thm}

\begin{defn}
Soit $\ov{L} = (L,\|.\|)$ un fibr\'e en droites hermitien. La m\'etrique
$\|.\|$ est dite {\it positive\/} (resp. {\it strictement positive\/}) \ssi\
$c_{1}(L,\|.\|) \geqslant 0$ sur $X$ (resp. pour toute forme de K\"ahler
$\alpha$ sur $X$ et tout ouvert $U$ d'un recouvrement suffisamment fin de $X$, 
il existe $\epsi \in \R^{+\ast}$ tel que $c_{1}(L,\|.\|)
\geqslant \epsi\alpha$ sur $U$).
\end{defn}
\begin{defn}
\label{def_admissibilite1}
Soit $L$ un fibr\'e en droites holomorphe sur $X$. 
Une m\'etrique $\|.\|$ continue et positive sur $L$ 
est dite {\it admissible\/}
s'il existe une suite $\left(\|.\|_{n}\right)_{n \in \M{N}}$ de
m\'etriques positives de classe $C^{\infty}$ sur $L$, convergeant uniform\'ement
vers $\|.\|$ sur $X$.\\
On appelle {\it fibr\'e admissible sur $X$\/} un fibr\'e en
droites holomorphe muni d'une
m\'etrique admissible sur $X$.
\end{defn}
\begin{prop}
\label{exemple_produit_formes}
Soient $\ov{L}_{1} = (L_{1},\|.\|_{1}), \dots, \ov{L}_{q} = (L_{q}, \|.\|_{q})$
des fibr\'es en droites admissibles sur $X$ et $s_{1},\dots,s_{q}$ des sections
m\'eromorphes non identiquement nulles, sur chaque composante connexe de $X$,
de $L_{1},\dots,L_{q}$ respectivement, et telles que les cycles $\op{div}s_{1},
\dots, \op{div}s_{q}$ soient d'intersection propre (i.e. tels que
$\op{codim}(\op{div}s_{j_{1}} \cap \dotsm \cap \op{div}s_{j_{m}}) \geqslant m$
pour tout choix d'indices $j_{1} < \dots < j_{m}$ dans $\{1,\dots,q\}$).
Pour tout $1 \leqslant i \leqslant q$, le courant~:
\[
(- \log \|s_{i}\|_{i}^{2})c_{1}(\ov{L}_{1})\dotsm c_{1}(\ov{L}_{i-1})\cdot
\delta_{\op{div}s_{i+1}}\dotsm \delta_{\op{div}s_{q}}, 
\]
est bien d\'efini et est un \'el\'ement de $\AAA_{\log}^{q-1,q-1}(X)$.
\end{prop}
\demo\
Le probl\`eme \'etant local, on se ram\`ene par lin\'earit\'e au cas o\`u
$s_{1}, \dots,s_{q}$ sont holomorphes au-dessus d'un ouvert $U$. 
La proposition est alors une cons\'equence directe de la formule de
Poincar\'e-Lelong g\'en\'eralis\'ee (\ref{PL_generalisee}), du fait que pour tout $1 \leqslant
i \leqslant q$, $c_{1}(\ov{L}_{i}) \in \AA^{1,1}(X)$, et des d\'efinitions.
\medskip

\begin{prop}
\label{approx_uniforme1}
Soit $\ov{L} = (L,\|.\|)$ un fibr\'e admissible sur $X$ compacte. Il existe une suite
croissante de m\'etriques $C^{\infty}$ positives convergeant
uniform\'ement vers $\|.\|$.
\end{prop}
\demo\ 
Soit $(\|.\|_{n})_{n \in \M{N}}$ une suite de m\'etriques $C^{\infty}$ positives
sur $L$ convergeant uniform\'ement vers $\|.\|$.
Pour tout $n \in \M{N}$, on note $\varphi(n)$ le plus petit entier positif tel
que~:
\[
\left| 
\frac{\;\;\;\|.\|_{\varphi(n)}}{\|.\|\;\;} - 1 \right| < \frac{1}{2^{n+2}}.
\]
On choisit alors un r\'eel $\lambda_{n}$ dans $]0,1[$ tel que~:
\[
\left(1 - \frac{1}{2^{n}}\right) < 
\lambda_{n} \frac{\;\;\;\|.\|_{\varphi(n)}}{\|.\|\;\;} < \left(1 -
\frac{1}{2^{n+1}}\right).
\]
On a donc construit une
application croissante $\varphi: \M{N} \rightarrow \M{N}$ 
et une suite de
r\'eels $(\lambda_{n})_{n \in \M{N}}$ dans $]0,1[$ tendant vers $1$ tels 
que la suite des m\'etriques $\|.\|_{n}' =
\lambda_{n}\|.\|_{\varphi(n)}$ pour $n \in \M{N}$ soit croissante sur $X$. Comme
$c_{1}(L,\|.\|'_{n}) = c_{1}(L,\|.\|_{\varphi(n)}) \geqslant 0$ pour tout $n
\in \M{N}$ et que $\left(\|.\|'_{n}\right)_{n \in \M{N}}$ converge uniform\'ement
vers $\|.\|$ sur $X$, la proposition est d\'emontr\'ee.
\medskip

\begin{expl}~
\label{exemple_adm1}
\begin{itemize}
\item{Toute m\'etrique positive $C^{\infty}$ est
admissible (prendre $\|.\|_{n} = \|.\|$).}
\item{
Soient $\P$ une vari\'et\'e torique projective lisse et
$L$ un fibr\'e en droites engendr\'e par ses sections
globales au-dessus de $\P$; la m\'etrique canonique $\|.\|_{L,\infty}$ 
introduite \`a la proposition (\ref{metrique_ind}) est
admissible (voir prop. \ref{approximation}).}
\item{Plus g\'en\'eralement, soit $L$ un fibr\'e en droites holomorphe engendr\'e 
par ses
sections globales holomorphes au-dessus d'une vari\'et\'e $X$ que l'on 
suppose compacte. Soit $\varphi: X
\rightarrow X$ un morphisme surjectif tel que $\varphi^{\ast}(L)
\stackrel{\Phi_{L}}{\simeq} L^{k}$, avec $k$ un entier $> 1$; 
on munit $L$ d'une m\'etrique de Zhang $\|.\|_{\op{Zh}}$ pour $\varphi$ (i.e.
une m\'etrique continue telle que $\varphi^{\ast}\left(L,\|.\|_{\op{Zh}}\right)
\stackrel{\Phi_{L}}{\simeq} \left(L,\|.\|_{\op{Zh}}\right)^{k}$). Le fibr\'e
hermitien $\left(L,\|.\|_{\op{Zh}}\right)$ est admissible. (Prendre sur $L$ une
m\'etrique $\|.\|_{0}$ positive $C^{\infty}$, puis consid\'erer la suite de
m\'etriques $(\|.\|_{n})_{n \in \M{N}}$ d\'efinie par la r\'ecurrence~:
\[
\|.\|_{n} = \left(\Phi_{L}^{\ast}\,\varphi^{\ast}\, \|.\|_{n-1}\right)^{1/k}.
\]
Les m\'etriques $\|.\|_{n}$ sont positives et $C^{\infty}$, et la suite
$(\|.\|_{n})_{n \in \M{N}}$ converge uniform\'ement vers $\|.\|_{\op{Zh}}$ sur $X$
d'apr\`es (\cite{21}, th. 2.2)).}
\end{itemize}
\end{expl}
\bigskip

\subsection{Un th\'eor\`eme d'approximation globale}~

Le th\'eor\`eme suivant montre que pour un fibr\'e ample, les notions de
m\'etriques positives et de m\'etriques admissibles sont \'equivalentes~:
\begin{thm}
\label{admi_positivite}
Soit $X$ une vari\'et\'e complexe projective et $L$ un fibr\'e en droites
ample sur $X$. Toute m\'etrique positive sur $L$ est admissible.
\end{thm}
Pour d\'emontrer le th\'eor\`eme (\ref{admi_positivite}) on utilise une
m\'ethode de r\'egularisation et de recollement essentiellement due \`a
Richberg (cf. \cite{22}). On s'inspire ici de la pr\'esentation donn\'ee dans
(\cite{7}, \S 1.5.C).

Soit $\theta : \R \rightarrow \R^{+}$ une fonction $C^{\infty}$ dont le
support est inclus dans $[-1,1]$ et telle que $\int_{\R}\theta(t) \, dt = 1$ et
$\int_{\R}t\theta(t)\, dt = 0$.
\begin{lem}
\label{approx_richberg}
Pour tout $\eta = (\eta_{1}, \dots, \eta_{p}) \in ]0, + \infty[^{p}$,
la fonction~:
\[
M_{\eta}(t_{1}, \dots, t_{p}) = \int_{\R^{p}}\op{max}\{t_{1}+h_{1}, \dots,
t_{p} + h_{p}\} \prod_{1 \leqslant j \leqslant
p}\theta\left(h_{j}/\eta_{j}\right)\, dh_{1}\dots dh_{p}, 
\]
v\'erifie les propri\'et\'es suivantes~:
\begin{enumerate}
\item{$M_{\eta}(t_{1}, \dots, t_{p})$ est croissante en chacune des variables
et est convexe et $C^{\infty}$ sur $\R^{p}$.}
\item{On a~: $\op{max}\{t_{1}, \dots, t_{p}\} \leqslant M_{\eta}(t_{1}, \dots,
t_{p}) \leqslant \op{max}\{t_{1} + \eta_{1}, \dots, t_{p} + \eta_{p}\}$.}
\item{On a~: $M_{\eta}(t_{1}, \dots, t_{p}) = M_{(\eta_{1}, \dots,
\widehat{\eta_{j}}, \dots, t_{p})}(t_{1}, \dots, \widehat{t_{j}}, \dots, t_{p})$ d\`es
que $t_{j} + \eta_{j} \leqslant \max_{k \neq j}\{t_{k} - \eta_{k}\}$.}
\item{$M_{\eta}(t_{1}+a, \dots, t_{p} +a) = M_{\eta}(t_{1}, \dots, t_{p}) + a$
pour tout $a \in \R$.}
\item{Si $u_{1}, \dots, u_{p}$ sont \pshs\ sur $X$ et v\'erifient $\dd u_{j}
\geqslant \alpha$, o\`u $\alpha \in A^{1,1}(X)$, alors $u = M_{\eta}(u_{1},
\dots, u_{p})$ est \psh\ et satisfait \`a l'in\'egalit\'e $\dd u \geqslant
\alpha$.}
\end{enumerate}
\end{lem}
\demo\ Voir (\cite{7}, Lemme 1.5.16).
\medskip

On peut passer \`a la d\'emonstration du th\'eor\`eme~: \\
On note $\|.\|$ la m\'etrique positive consid\'er\'ee. 
Quitte \`a consid\'erer une puissance tensorielle assez grande de $L$, on peut
supposer $L$ tr\`es ample.\\
On suppose dans un
premier temps que la m\'etrique $\|.\|$ est strictement positive. On 
note $\alpha$ une forme de K\"ahler sur $X$ telle que $c_{1}(L,\|.\|) \geqslant
\alpha$.
On choisit $S = \{s_{1}, \dots, s_{N}\}$ une $\Z$-base des sections
holomorphes de $L$ au-dessus de $X$. Pour tout $i \in \{1,\dots, N\}$ on
note $u_{i} = - \log \|s_{i}\|^{2} \in \op{Psh}(X)$. Pour $i$ et $j$
\'el\'ements de $\{1,\dots, N\}$, on note $f_{i,j}$ la fonction m\'eromorphe
d\'efinie par $s_{i} = f_{i,j}s_{j}$. Pour tout couple $(\Omega, x)$ form\'e d'un point $x
\in X$ et d'un ouvert $\Omega \subset X$ le contenant, on dit que $\Omega$ est
une boule centr\'ee en $x$ de rayon $r$ s'il existe une carte $(V,\varphi)$
de $X$ contenant $\Omega$ telle que $\varphi(x) = 0$ et que
$\varphi(\Omega)$ soit une boule centr\'ee en $0$ de rayon $r$ dans
$\M{C}^{\dim X}$.

On choisit $\left(\Omega_{\alpha}\right)_{\alpha \in A}$ un recouvrement ouvert
fini de $X$ v\'erifiant les propri\'et\'es suivantes~:
\begin{enumerate}
\item{Pour tout $\alpha \in A$, l'ouvert $\Omega_{\alpha}$ est une boule
centr\'ee en un point $0_{\alpha} \in \Omega_{\alpha}$ de rayon $r_{\alpha}$
dans une carte $(V_{\alpha}, \varphi_{\alpha})$ de $X$.}
\item{Pour tout $\alpha \in A$, il existe $s_{i_{\alpha}} \in S$ qui ne
s'annule pas sur un
voisinage $U_{\alpha}$ de $\ov{\Omega_{\alpha}}$.}
\end{enumerate}
Comme $L$ est tr\`es ample et donc engendr\'e par ses sections globales, un tel
recouvrement existe toujours.

On choisit $\lambda$ un \'el\'ement de $]0,1[$. 
Pour tout $\alpha \in A$, on peut construire par r\'egularisation (voir par
exemple \cite{7}, th. 1.5.5) une famille
$\left(u_{\alpha}^{(\varepsilon)}\right)_{\varepsilon \in \R^{+\ast}}$
d'\'el\'ements de $\op{Psh}(\Omega_{\alpha}) \cap C^{\infty}(\Omega_{\alpha})$ telle que
$u_{\alpha}^{(\varepsilon)}$ soit une fonction croissante de $\varepsilon$ et
que~:
\[
\forall \varepsilon \in \R^{+\ast}, \qquad \|u_{i_{\alpha}} -
u_{\alpha}^{(\varepsilon)}\|_{\ov{\Omega}_{\alpha}, \infty} < \varepsilon.
\]
Pour tout $\alpha \in A$, on choisit $\Omega_{\alpha}'' \subset
\Omega_{\alpha}' \subset \Omega_{\alpha}$ des boules concentriques de rayon
respectif $r''_{\alpha} < r'_{\alpha} < r_{\alpha}$ dans la carte $V_{\alpha}$,
et telles que la famille $\left(\Omega''_{\alpha}\right)_{\alpha \in A}$ forme
encore un recouvrement de $X$.

Soient $\epsi_{\alpha}$ et $\gamma_{\alpha}$ deux nombres r\'eels strictement
positifs que l'on fixera par la suite. Pour tout $z \in \ov{\Omega}_{\alpha}$,
on pose (dans le syst\`eme de coordonn\'ees relatif \`a la carte
$(V_{\alpha},\varphi_{\alpha})$), 
\[
v_{\alpha}^{(\epsi_{\alpha},\gamma_{\alpha})}(z) =
u_{\alpha}^{(\epsi_{\alpha})}(z) + \gamma_{\alpha}({r'}_{\alpha}^{2} - |z|^{2}).
\]
Pour $\epsi_{\alpha} < \epsi_{\alpha,0}$ et $\gamma_{\alpha} <
\gamma_{\alpha,0}$ assez petits, on a
$v_{\alpha}^{(\epsi_{\alpha},\gamma_{\alpha})} \leqslant u_{i_{\alpha}} +
\lambda / 2$ et $\dd v_{\alpha}^{(\epsi_{\alpha},\gamma_{\alpha})} \geqslant
(1-\lambda)\alpha$ sur $\ov{\Omega}_{\alpha}$. On pose~:
\[
\eta_{\alpha} = \gamma_{\alpha} \min \left\{{r'}_{\alpha}^{2} - {r''}_{\alpha}^{2},
(r_{\alpha}^{2} - \left. {r'}_{\alpha}^{2})\right/ 2\right\}.
\]
On choisit tout d'abord $\gamma_{\alpha} < \gamma_{\alpha,0}$ tel que
$\eta_{\alpha} < \lambda / 2$, puis on choisit $\epsi_{\alpha} <
\epsi_{\alpha,0}$ suffisamment petit pour que l'on ait~:
\[
u_{i_{\alpha}} \leqslant u_{\alpha}^{(\epsi_{\alpha})} < u_{i_{\alpha}} +
\eta_{\alpha}
\]
sur $\ov{\Omega}_{\alpha}$. Comme $\gamma_{\alpha}\left({r'}_{\alpha}^{2} -
|z|^{2}\right)$ est inf\'erieur \`a $-2\eta_{\alpha}$ sur
$\partial\Omega_{\alpha}$ et strictement sup\'erieur \`a $\eta_{\alpha}$ sur
$\ov{\Omega}''_{\alpha}$, on a~: $v_{\alpha}^{(\epsi_{\alpha},\gamma_{\alpha})}
< u_{i_{\alpha}} - \eta_{\alpha}$ sur $\partial \Omega_{\alpha}$ et 
$v_{\alpha}^{(\epsi_{\alpha},\gamma_{\alpha})} > u_{i_{\alpha}} +
\eta_{\alpha}$ sur $\ov{\Omega}''_{\alpha}$.

Pour tout $i \in \{1,\dots,N\}$, on d\'efinit maintenant sur
$\ov{\Omega}_{\alpha}$ la fonction~:
\[
\widetilde{u}_{i,\alpha} = v_{\alpha}^{(\epsi_{\alpha},\gamma_{\alpha})} +
\left( - \log |f_{i,i_{\alpha}}|^{2} \right).
\]
D'apr\`es ce qui pr\'ec\`ede, les fonctions $\widetilde{u}_{i,\alpha}$
v\'erifient les assertions suivantes~:
\begin{enumerate}
\item{$\widetilde{u}_{i,\alpha} \leqslant u_{i} + \lambda / 2$ sur
$\ov{\Omega}_{\alpha}$.}
\item{$\widetilde{u}_{i,\alpha} < u_{i} - \eta_{\alpha}$ sur
$\partial\Omega_{\alpha}$.}
\item{$\widetilde{u}_{i,\alpha} > u_{i} + \eta_{\alpha}$ sur
$\ov{\Omega}''_{\alpha}$.}
\item{$\widetilde{u}_{i,\alpha}$ est $C^{\infty}$ sur
$\ov{\Omega}_{\alpha} - \op{div}s_{i}$, et on a~: $\dd 
\widetilde{u}_{i,\alpha} \geqslant (1-\lambda)\alpha$.}
\end{enumerate}
Gr\^ace aux $(2)$ et $(3)$ ci-dessus et au $(3)$ du lemme
(\ref{approx_richberg}), on peut d\'efinir pour tout $i \in \{1,\dots,N\}$ la
fonction~:
\begin{alignat*}{3}
\widetilde{u}_{i}: X - &\op{div} s_{i}& &\longrightarrow  & &\R \\
&z & &\longrightarrow &
&M_{(\eta_{\alpha})}\left(\widetilde{u}_{i,\alpha}(z)\right).
\end{alignat*}
La fonction $\widetilde{u}_{i}$ est $C^{\infty}$ d'apr\`es le $(1)$ du lemme
(\ref{approx_richberg}), et \psh\ d'apr\`es le $(5)$ de
(\ref{approx_richberg}). De plus, elle v\'erifie l'encadrement~:
\begin{equation}
\label{encadrement1}
u_{i} < u_{i} + \min_{\alpha \in A}\eta_{\alpha} < \widetilde{u}_{i} \leqslant
u_{i} + \lambda.
\end{equation}
Enfin, pour tout $i$ et $j$ \'el\'ements de $\{1,\dots,N\}$, on a~:
\begin{equation}
\label{egalite1}
\widetilde{u}_{j} = \widetilde{u}_{i} + \left(- \log | f_{j,i}|^{2}\right),
\end{equation}
d'apr\`es le $(4)$ du lemme (\ref{approx_richberg}).

Soit $U$ un ouvert de $X$ suffisamment petit pour qu'il existe $s_{i} \in S$ ne
s'annulant pas sur $U$. Soit $s$ une section r\'eguli\`ere de $L$ au-dessus de
$U$. On pose~:
\[
\|s\|_{\lambda} = \left| \frac{s}{s_{i}}\right| e^{-\widetilde{u}_{i}/2}.
\]
Gr\^ace \`a l'\'egalit\'e (\ref{egalite1}), la d\'efinition ci-dessus ne
d\'epend pas du choix de $s_{i}$ et d\'efinit par recollement une norme
$C^{\infty}$ positive sur $L$. De plus, on tire des d\'efinitions l'\'egalit\'e~:
\[
\frac{\;\,\|s\|_{\lambda}}{\|s\|} = e^{\frac{u_{i} - \widetilde{u}_{i}}{2}}, 
\]
ce qui donne l'encadrement~:
\[
e^{- \lambda /2} \leqslant \frac{\;\,\|s\|_{\lambda}}{\|s\|} \leqslant 1, 
\]
d'apr\`es l'encadrement (\ref{encadrement1}).

En faisant tendre $\lambda$ vers $0$,
on extrait de la famille $\left(\|.\|_{\lambda}\right)_{\lambda \in
]0,1[}$ une suite croissante de m\'etriques $C^{\infty}$ positives convergeant
uniform\'ement vers $\|.\|$ sur $X$. Le th\'eor\`eme est donc d\'emontr\'e dans
le cas o\`u $\|.\|$ est strictement positive sur $X$.

On s'int\'eresse d\'esormais au cas o\`u $\|.\|$ est suppos\'ee simplement
positive.
Comme $L$ est tr\`es ample, il
existe sur $L$ une m\'etrique $\|.\|'$ qui est $C^{\infty}$ et strictement
positive.
Pour tout $n \in \M{N}$, la
m\'etrique $\|.\|_{n} = (\|.\|)^{1-1/n}\cdot (\|.\|')^{1/n}$ est strictement positive.
D'apr\`es le r\'esultat que l'on vient de d\'emontrer, elle est donc
approchable uniform\'ement par des m\'etriques $C^{\infty}$ positives sur $X$. 
Comme $\|.\|_{n}$ tend uniform\'ement vers $\|.\|$ quand $n$ tend vers $+
\infty$, on en d\'eduit que $\|.\|$ est approchable uniform\'ement par des
m\'etriques $C^{\infty}$ positives.
\bigskip

\subsection{Fibr\'es en droites int\'egrables}~

Suivant Zhang, on introduit \`a pr\'esent une nouvelle classe de fibr\'es en droites plus 
g\'en\'erale que celle des fibr\'es
admissibles. On conserve ici la terminologie de Zhang (cf. \cite{21}, \S
1.5)~:
\begin{defn}
\label{decomposable}
Soit $\ov{L} = (L,\|.\|)$ un fibr\'e en droites hermitien sur $X$. Le fibr\'e
$\ov{L}$ est dit {\it int\'egrable sur $X$\/} ou plus simplement {\it
int\'egrable\/} si et seulement s'il existe $\ov{E}_{1}$ et
$\ov{E}_{2}$ deux fibr\'es admissibles sur $X$ tels que~:
\[
\ov{L} = \ov{E}_{1} \otimes (\ov{E}_{2})^{-1}.
\]
\end{defn}
\begin{expl}~
\label{exemple_decomp1}
\begin{itemize}
\item{Supposons $X$ projective; tout fibr\'e en droites $\ov{L} = (L,\|.\|)$
muni d'une
m\'etrique $C^{\infty}$ est int\'egrable (consid\'erer
$\ov{L}\otimes\ov{H}^{n}$, o\`u $\ov{H}$ est un fibr\'e ample muni d'une
m\'etrique $C^{\infty}$ strictement positive sur $X$).}
\item{Soient $\P$ une vari\'et\'e torique projective lisse et
$\ov{L}_{\infty} = (L,\|.\|_{L,\infty})$ un fibr\'e en droites
sur $\P$ muni de sa m\'etrique canonique, $(L(\M{C}),\|.\|_{L,\infty})$ est
int\'egrable sur $\P(\M{C})$ (cf. prop. (\ref{decomposition}) et
(\ref{approximation})).}
\end{itemize}
\end{expl}
\begin{prop}
\label{produit_decomposables}
Soient $\ov{E}$ et $\ov{F}$ deux fibr\'es admissibles (resp. int\'egrables)
sur $X$; leur produit tensoriel $\ov{E} \otimes \ov{F}$ est admissible (resp.
int\'egrable) sur $X$.
\end{prop}
\demo\ On suppose tout d'abord que $\ov{E} = (E,\|.\|_{E})$ 
et $\ov{F} = (F,\|.\|_{F})$ sont admissibles.
On a $c_{1}(\ov{E}\otimes\ov{F}) = c_{1}(\ov{E}) + c_{1}(\ov{F}) \geqslant
0$ d'apr\`es
(\ref{additivite_chern}). De plus, si $\left(\|.\|_{E}^{(n)}\right)_{n \in \M{N}}$
et $\left(\|.\|_{F}^{(n)}\right)_{n \in \M{N}}$ sont deux suites de
m\'etriques $C^{\infty}$ positives convergeant uniform\'ement sur $X$ vers $\|.\|_{E}$
et $\|.\|_{F}$ respectivement, alors la suite donn\'ee par
$\left(\|.\|_{E}^{(n)} \otimes \|.\|_{F}^{(n)}\right)_{n \in \M{N}}$ est une suite
de m\'etriques positives $C^{\infty}$ convergeant uniform\'ement sur $X$ vers $\|.\|_{E}
\otimes \|.\|_{F}$; on en d\'eduit que $\ov{E}\otimes\ov{F}$ est admissible sur
$X$.

On suppose maintenant que $\ov{E}$ et $\ov{F}$ sont int\'egrables sur $X$. 
On peut donc trouver $\ov{E}_{1}$, $\ov{E}_{2}$, $\ov{F}_{1}$ et $\ov{F}_{2}$
des fibr\'es admissibles sur $X$ tels que $\ov{E} = \ov{E}_{1} \otimes
(\ov{E}_{2})^{-1}$ et $\ov{F} = \ov{F}_{1} \otimes (\ov{F}_{2})^{-1}$. On
tire~:
\[
\ov{E} \otimes \ov{F} = \left(\ov{E}_{1} \otimes \ov{F}_{1}\right) \otimes
\left(\ov{E}_{2} \otimes \ov{F}_{2}\right)^{-1}, 
\]
et comme $\ov{E}_{1} \otimes \ov{F}_{1}$ et $\ov{E}_{2} \otimes \ov{F}_{2}$
sont admissibles d'apr\`es ce qui pr\'ec\`ede, on en d\'eduit que $\ov{E}
\otimes \ov{F}$ est int\'egrable.
\medskip

\begin{prop}
\label{image_reciproque}
Soient $X$ et $Y$ deux vari\'et\'es complexes et $f: Y \rightarrow X$ une
application holomorphe. Si $\ov{L}$ est un fibr\'e en droites admissible (resp.
int\'egrable) sur $X$, alors $f^{\ast}(\ov{L})$ est admissible (resp.
int\'egrable) sur $Y$.
\end{prop}
\demo\ C'est une cons\'equence directe des d\'efinitions.
\medskip

\begin{prop}
Soit $\ov{L}$ un fibr\'e en droites hermitien int\'egrable sur $X$; on a~:
\[
c_{1}(\ov{L}) \in C_{0}^{1,1}(X) \subset \AA^{1,1}(X).
\]
\end{prop}
\demo\ C'est une cons\'equence directe des d\'efinitions
(\ref{formes_generalisees}), (\ref{formes_adherentes}) et (\ref{decomposable}).
\medskip
\begin{thm} 
\label{coho_courant}
Supposons $X$ projective. Soient $\ov{L}_{1}, \dots, \ov{L}_{p}$
des fibr\'es en droi\-tes hermitiens int\'egrables sur $X$ et $\alpha \in
\AAA^{\ast}_{\log}(X)$ un courant ferm\'e; on a l'\'egalit\'e des classes~:
\[
[c_{1}(\ov{L}_{1}) \dotsm c_{1}(\ov{L}_{p})\cdot \alpha] = c_{1}(L_{1}) \dotsm
c_{1}(L_{p})\cdot [\alpha], 
\]
en cohomologie de de Rham des courants, o\`u l'on a not\'e $[c_{1}(\ov{L}_{1})
\dotsm c_{1}(\ov{L}_{p})\cdot \alpha]$ et $[\alpha ]$ les classes des courants
$c_{1}(\ov{L}_{1})\dotsm c_{1}(\ov{L}_{p})\cdot \alpha$ et $\alpha$, et o\`u
$c_{1}(L_{1}), \dots, c_{1}(L_{p})$ d\'esignent les premi\`eres classes de
Chern des fibr\'es $L_{1}, \dots, L_{q}$ respectivement.
\end{thm}
\demo\ Par polarisation, il suffit de d\'emontrer le r\'esultat pour
$(L_{1},\|.\|_{1}),$ $\dots, (L_{p},\|.\|_{p})$ des fibr\'es en droites
admissibles.
Soient $\left(\|.\|_{1}^{(n)}\right)_{n \in \M{N}}, \dots, 
\left(\|.\|_{p}^{(n)}\right)_{n \in \M{N}}$ des suites croissantes de m\'etriques
$C^{\infty}$ positives sur $L_{1}, \dots, L_{p}$ convergeant vers les m\'etriques 
$\|.\|_{1}, \dots,
\|.\|_{p}$ respectivement. On a~:
\[
c_{1}\left(L_{1},\|.\|_{1}^{(n)}\right)\dotsm 
c_{1}\left(L_{p},\|.\|_{p}^{(n)}\right)\cdot \alpha \text{\quad tend vers
\quad} c_{1}(\ov{L}_{1})\dotsm c_{1}(\ov{L}_{p})\cdot \alpha
\]
au sens de la topologie faible des courants, et donc~:
\[
[c_{1}\left(L_{1},\|.\|_{1}^{(n)}\right)\dotsm
c_{1}\left(L_{p},\|.\|_{p}^{(n)}\right)\cdot \alpha] \text{\quad tend vers
\quad} [c_{1}(\ov{L}_{1})\dotsm c_{1}(\ov{L}_{p})\cdot \alpha]
\]
pour la topologie naturelle sur $H_{\op{DR}}^{\ast}(X)$. On
d\'eduit le r\'esultat de l'\'egalit\'e~:
\[
[c_{1}\left(L_{1},\|.\|_{1}^{(n)}\right)\dotsm 
c_{1}\left(L_{p},\|.\|_{p}^{(n)}\right)\cdot \alpha] = c_{1}(L_{1})\dotsm
c_{1}(L_{p})\cdot [\alpha]
\]
valable pour tout $n \in \M{N}$.
\medskip

\section{Groupes de Chow arithm\'etiques g\'en\'eralis\'es}~

\subsection{Th\'eorie classique de Gillet-Soul\'e}~

\label{theorie_classique}
On suit ici (\cite{13} et \cite{3}, \S 2.1). Soit $K$ un corps de nombre, $\C{O}_{K}$ son
anneau d'entier et $S = \op{Spec}(\C{O}_{K})$ le sch\'ema associ\'e. Pour
tout plongement $\sigma: K \rightarrow \M{C}$ et pour tout sch\'ema $X$ sur
$\op{Spec}K$ ou sur $S$, on note $X_{\sigma}$ le sch\'ema sur $\M{C}$ d\'eduit
de $X$ par le changements de base $\op{Spec}\M{C} \rightarrow \op{Spec}K$. De
m\^eme, si $f: X \rightarrow Y$ est un morphisme de sch\'emas sur $K$, on note
$f_{\sigma}: X_{\sigma} \rightarrow Y_{\sigma}$ le morphisme de sch\'emas sur
$\M{C}$ induit par changement de base.

Une {\it vari\'et\'e arithm\'etique\/} est par d\'efinition un sch\'ema $\pi: X
\rightarrow S$,
plat,
projectif, int\`egre et {\em r\'egulier} sur $S$. En particulier la fibre 
g\'en\'erique $X_{K} = X \times_{S}\;
\op{Spec} K$ est lisse. On note $d = \dim_{K}X_{K}$.

Pour toute vari\'et\'e arithm\'etique $X$ et tout entier $p$ positif, on note
$Z_{p}(X)$ (resp. $Z^{p}(X)$) le groupe des cycles sur $X$ de dimension $p$
(resp. codimension $p$). Pour un tel cycle $Z$, on note $|Z| \subset X$ le
support de $Z$.

On note $X(\M{C})$ les points complexes du sch\'ema $X$; c'est
la r\'eunion disjointe $\coprod_{\sigma: K \hookrightarrow
\M{C}}X_{\sigma}(\M{C})$. Soit $F_{\infty}: X(\M{C}) \rightarrow X(\M{C})$
l'involution antiholomorphe provenant de l'action de la conjugaison complexe
sur les coordonn\'ees des points complexes de $X$. On note $A^{p,p}(X_{\R})$
(resp. $D^{p,p}(X_{\R})$) l'ensemble des formes r\'eelles $\alpha \in
A^{p,p}(X(\M{C}))$ (resp. des courants r\'eels $\alpha \in D^{p,p}(X(\M{C}))$)
tels que $F_{\infty}^{\ast}(\alpha) = (-1)^{p}\alpha$. On note
$\widetilde{A}^{p,p}(X_{\R})$ (resp. $\widetilde{D}^{p,p}(X_{\R})$) le quotient
$\left. A^{p,p}(X_{\R})\right/ (\op{Im}\partial + \op{Im}\ov{\partial})$ (resp. le
quotient $\left. D^{p,p}(X_{\R}) \right/ (\op{Im}\partial +
\op{Im}\ov{\partial})$). On note \'egalement $Z^{p,p}(X_{\R}) = \op{Ker}\{d:
A^{p,p}(X_{\R}) \rightarrow A^{2p+1}(X(\M{C}))\} \subset A^{p,p}(X_{\R})$.

Tout cycle irr\'eductible $Z \in Z^{p}(X)$ d\'efinit un courant $\delta_{Z} \in
D^{p,p}(X_{\R})$ par int\'egration sur l'ensemble de ses points complexes. Si
$Z = \sum_{i \in I}n_{i}Z_{i}$ o\`u $Z_{i}$ est irr\'eductible, on pose~:
\[
\delta_{Z} = \sum_{i \in I}n_{i}\delta_{Z_{i}(\M{C})}.
\]
Un {\it courant de Green\/} pour $Z$ est un courant $g \in D^{p-1,p-1}(X_{\R})$
tel que $\dd g + \delta_{Z}$ est une forme $C^{\infty}$.

Soit $X$ une vari\'et\'e arithm\'etique. On note $\widehat{Z}^{p}(X)$ le groupe
form\'e des couples de la forme $(Z,g)$, o\`u $Z \in Z^{p}(X)$ et $g$ est un
courant de Green pour $Z$, muni de la loi d'addition composante par
composante. Soit $\widehat{R}^{p}(X) \subset \widehat{Z}^{p}(X)$ le sous-groupe
engendr\'e par les paires de la forme $(0,\partial u + \ov{\partial} v)$ et
$\left( \op{div} (f), - \log |f|^{2}\right)$, o\`u $f \in k(Y)^{\ast}$ est une
fonction rationnelle non identiquement nulle sur un sous sch\'ema int\`egre $Y
\subset X$ de codimension $p-1$, et 
o\`u pour simplifier les notations on a \'ecrit $- \log |f|^{2}$ 
pour d\'esigner $[-\log |f|^{2}]_{Y(\M{C})}$, le courant sur
$X(\M{C})$ dont l'action sur les formes diff\'erentielles est obtenue par
restriction \`a la partie lisse de $Y(\M{C})$ puis l'int\'egration contre la
fonction $- \log |f|^{2}$.

Le {\it groupe de Chow arithm\'etique de codimension $p$ de $X$\/} est d\'efini
par~:
\[
\widehat{CH}^{p}(X) = \left. \widehat{Z}^{p}(X) \right/ \widehat{R}^{p}(X).
\]
On note $R_{\op{fin}}^{p}(X)_{\M{Q}}$ l'espace des $\M{Q}$-cycles de la 
forme $\sum_{i}q_{i}\op{div}(f_{i})$, o\`u $q_{i} \in \M{Q}$ et $f_{i} \in
k(Y_{i})^{\ast}$ est une fonction rationnelle non identiquement nulle 
sur $Y_{i}$ un sous-sch\'ema int\`egre de codimension $p-1$ contenu dans 
une fibre ferm\'ee du morphisme $\pi:X \rightarrow S$. On remarque que pour
tout $R \in R_{\op{fin}}^{p}(X)_{\M{Q}}$, la classe de $(R,0)$ dans
$\widehat{CH}^{p}(X)_{\M{Q}}$ est nulle.

Une {\it $K_{1}$-chaine\/} de codimension $p$ dans $X$ est un \'el\'ement du
groupe \\
$\bigoplus_{x \in X^{(p-1)}}k(x)^{\ast}$, o\`u $X^{(p-1)}$ d\'esigne
l'ensemble des sous-sch\'emas ferm\'es int\`egres de codimension $p-1$ dans
$X$. Si $f = \sum_{i \in I}[f_{W_{i}}]$, o\`u $W_{i} \in X^{(p-1)}$ et
$f_{W_{i}} \in k(W_{i})^{\ast}$, est une $K_{1}$-chaine de codimension $p$, on 
pose~:
\begin{align*}
\op{div}f &= \sum_{i \in I}\op{div}(f_{W_{i}}) \in Z^{p}(X), \\
- \log |f|^{2} &= \sum_{i \in I}[-\log|f_{W_{i}}|^{2}]_{W_{i}(\M{C})} \in 
D^{p-1,p-1}(X_{\R}), \\
\qquad \qquad \qquad \qquad \qquad \qquad 
\op{\widehat{\op{div}}}f &= (\op{div}f, - \log |f|^{2}) \in
\widehat{Z}^{p}(X);\qquad \qquad \qquad \qquad 
\end{align*}
et on appelle {\it support\/} de $f$ la fermeture de la r\'eunion des
$W_{i}$.
Si $\E{Z} = \{Z_{1},\dots,Z_{n}\}$ est une famille de sous-sch\'emas ferm\'es
int\`egres de $X$, on dit que la $K_{1}$-chaine $f = \sum_{i \in I}[f_{W_{i}}]$
et la famille $\E{Z}$ s'intersectent {\it presque proprement\/} si pour tout
$W_{i}$ et tout $Z_{j} \in \E{Z}$, les cycles $\op{div}(W_{i})$ et $Z_{j}$
s'intersectent proprement.
\medskip

On dispose de trois morphismes de groupes d\'efinis comme suit~:
\begin{alignat*}{3}
\zeta :\; &\widehat{CH}^{p}(X) & &\longrightarrow & &CH^{p}(X) \\
&[(Z,g)] & &\longmapsto & &[Z], \\
& & & & & \\
\omega :\; &\widehat{CH}^{p}(X) & &\longrightarrow &  &A^{p,p}(X_{\R})\\
&[(Z,g)] & & \longmapsto & &\dd g + \delta_{Z},
\end{alignat*}
et 
\begin{alignat*}{3}
a :\; &\widetilde{A}^{p-1,p-1}(X_{\R}) & &\longrightarrow &
&\widehat{CH}^{p}(X) \\
&\eta & & \longmapsto & &[(0,\eta)].
\end{alignat*}
Ces morphismes donnent lieu \`a deux suites exactes~:
\begin{align}
\label{suite_exacte1}
 &CH^{p,p-1}(X) \stackrel{\rho}{\longrightarrow}
\widetilde{A}^{p-1,p-1}(X_{\R}) \stackrel{a}{\longrightarrow}
\widehat{CH}^{p}(X) \stackrel{\zeta}{\longrightarrow} CH^{p}(X)
\longrightarrow 0, \\
\label{suite_exacte2}
 &CH^{p,p-1}(X) \stackrel{\rho}{\longrightarrow}
H^{p-1,p-1}(X_{\R}) \stackrel{a}{\longrightarrow} \widehat{CH}^{p}(X) \\
& \qquad \qquad \qquad \qquad \qquad
\stackrel{(\zeta,-\omega)}{\longrightarrow} CH^{p}(X)\oplus Z^{p,p}(X_{\R})
\stackrel{cl + [.]}{\longrightarrow} H^{p,p}(X_{\R}) \longrightarrow 0, \notag
\end{align}
o\`u $CH^{p,p-1}(X)$ est un groupe d\'efini dans (\cite{12}, \S 8),
o\`u $\rho$ est $-2$ fois 
le r\'egulateur de Beilinson (cf. \cite{13}, \S 3.3.5 et 3.5), o\`u
$[.]$ est le morphisme ``classe en cohomologie de De Rham'' et o\`u
$cl$ est l'application cycle.

Tout morphisme de sch\'emas $f: X \rightarrow Y$ entre vari\'et\'es arithm\'etiques
induit un morphisme de groupes~:
\[
f^{\ast}: \widehat{CH}^{p}(Y) \longrightarrow \widehat{CH}^{p}(X).
\]
De plus, on dispose d'un produit~:
\[
\widehat{CH}^{p}(X) \otimes \widehat{CH}^{q}(X) \longrightarrow 
\widehat{CH}^{p+q}(X)_{\M{Q}}, 
\]
que l'on peut en fait d\'efinir \`a valeur dans $\widehat{CH}^{p+q}(X)$ 
quand $X$ est lisse sur $S$. 
Les morphismes $\zeta$ et $\omega$ d\'efinis ci-dessus sont des morphismes 
d'anneaux. On a \'egalement $f^{\ast}(x\cdot y) = 
f^{\ast}(x)\cdot f^{\ast}(y)$.

Enfin la formule suivante est utile dans la pratique. Soient $\eta \in
\widetilde{A}^{\ast}(X_{\R})$ et $x \in \widehat{CH}^{\ast}(X)$, on a~:
\[
a(\eta)\cdot x = a (\eta \, \omega(x)).
\]
Pour plus de d\'etails sur ces d\'efinitions et cette th\'eorie, voir
\cite{13}.
\medskip

\subsection{Fibr\'es en droites int\'egrables sur une vari\'et\'e
arithm\'etique}~

Soit $\pi: X \rightarrow \op{Spec}\C{O}_{K}$ 
une vari\'et\'e arithm\'etique de dimension relative $d$
et $p$ un entier positif. On note
$\ov{A}^{p,p}(\XR)$ (resp. $\AAA^{p,p}(\XR)$, resp. $\AAA^{p,p}_{\log}(\XR)$)
l'intersection $\AA^{p,p}(\XC) \cap D^{p,p}(\XR)$ (resp. l'intersection
$\AAA^{p,p}(\XC) \cap D^{p,p}(\XR)$, resp. l'intersection $\AAA^{p,p}_{\log}
(\XC) \cap D^{p,p}(\XR)$).

Un {\it fibr\'e en droites hermitiens sur $X$\/} est un couple $\ov{E} =
(E,\|.\|)$ form\'e d'un fibr\'e en droites $E$ sur $X$ et d'une m\'etrique
continue $\|.\|$ sur $E$ invariante sous l'action de $F_{\infty}$. 
Un tel fibr\'e est dit {\it admissible\/} lorsque
le fibr\'e holomorphe hermitien
$(E_{\M{C}},\|.\|)$ est admissible au sens de (\ref{def_admissibilite1}).

Un fibr\'e en droites hermitien $\ov{L}$ est dit {\it int\'egrable\/} \ssi\
il existe $\ov{E}_{1}$ et $\ov{E}_{2}$ deux fibr\'es en droites hermitiens
admissibles sur $X$ tels que~:
\[
\ov{L} = \ov{E}_{1} \otimes \left(\ov{E}_{2}\right)^{-1}.
\]
\begin{rem}
Si $\ov{L} = (L,\|.\|)$ est un fibr\'e en droites hermitiens int\'egrable sur
$X$, alors le fibr\'e holomorphe hermitien $(L(\M{C}),\|.\|)$ est int\'egrable sur
$\XC$ au sens de (\ref{decomposable}).
\end{rem}
\begin{expl}~
\begin{enumerate}
\item{Soit $\ov{L} = (L,\|.\|)$ un fibr\'e en droites hermitien sur $X$; si
$\|.\|$ est $C^{\infty}$ sur $\XC$, alors $\ov{L}$ est int\'egrable sur $X$.
(\'Ecrire $\ov{L} \simeq (\ov{L}\otimes \ov{H}^{n})\otimes (\ov{H}^{n})^{\ast}$, avec 
$\ov{H}$ ample muni d'une m\'etrique
$C^{\infty}$ strictement positive sur $\XC$ et $n \gg 0$).}
\item{Soit $L$ un fibr\'e en droites ample muni d'une m\'etrique $\|.\|$ continue
et positive; le fibr\'e hermitien $(L,\|.\|)$ est admissible sur $X$ d'apr\`es
(\ref{admi_positivite}).}
\item{
Soient $\P$ une vari\'et\'e torique projective lisse et
$\ov{L}_{\infty}$ un fibr\'e en droites sur $\P$ muni de
sa m\'etrique canonique; $\ov{L}_{\infty}$ est int\'egrable sur $\P$ d'apr\`es
(\ref{exemple_decomp1}).}
\end{enumerate}
\end{expl}
On appelle {\it forme diff\'erentielle g\'en\'eralis\'ee globale de type
$(p,p)$\/} tout \'el\'ement $\alpha \in \AA^{p,p}(\XR)$ pouvant s'\'ecrire sur
$\XC$ sous la forme~:
\[
\alpha = \sum_{i = 1}^{n}\omega_{i}c_{1}(\ov{L}_{i,1})\dots c_{1}(\ov{L}_{i,q_{i}})
\]
o\`u pour tout $1 \leqslant i \leqslant n$, $\omega_{i} \in
A^{\ast}(X_{\M{R}})$ et $\ov{L}_{i,1}, \dots, \ov{L}_{i,q_{i}}$ sont des fibr\'es
en droites int\'egrables sur $X$. On note $\AA_{\op{g}}^{p,p}(X_{\R}) \subset
\AA^{p,p}(X_{\R})$ l'ensemble des formes diff\'erentielles g\'en\'eralis\'ees 
globales de type
$(p,p)$ sur $X$. C'est une $\R$-alg\`ebre.

On note $\widetilde{\AA}_{\op{g}}^{p,p}(X_{\R})$ l'image de $\AA^{p,p}_{\op{g}}(X_{\R})$
dans $\widetilde{D}^{p,p}(X_{\R})$. On note \'egalement
$\AA_{\op{g}}^{\ast}(X_{\R}) = \oplus_{p \geqslant
0}\AA_{\op{g}}^{p,p}(X_{\R})$ et $\widetilde{\AA}_{\op{g}}^{\ast}(X_{\R}) =
\oplus_{p \geqslant 0}\widetilde{\AA}_{\op{g}}^{p,p}(X_{\R})$.
\begin{rem}
L'espace $\AA_{\op{g}}^{\ast}(X_{\R})$ est stable pour le produit d\'efini au
(\ref{produit_formes}). La structure d'alg\`ebre sur $\AA^{\ast}(X(\M{C}))$
induit donc une structure d'alg\`ebre associative commutative unif\`ere et
gradu\'ee sur $\AA^{\ast}_{\op{g}}(X_{\R})$.
\end{rem}
\begin{prop}
Si $\ov{E}$ est un fibr\'e en droites int\'egrable, alors $(\ov{E})^{-1}$ 
est int\'egrable.
\end{prop}
\demo\ C'est une cons\'equence imm\'ediate des d\'efinitions.
\medskip

\begin{prop}
Soient $\ov{E}$ et $\ov{F}$ deux fibr\'es en droites admissibles (resp.
int\'egrables);
leur produit tensoriel $\ov{E}\otimes\ov{F}$ est admissible (resp.
int\'egrable).
\end{prop}
\demo\ On suppose tout d'abord que $\ov{E} = (E,\|.\|_{E})$ et $\ov{F} =
(F,\|.\|_{F})$ sont admissibles.
Le fibr\'e holomorphe hermitien $((E\otimes F)({\M{C}}),
\|.\|_{E}\otimes\|.\|_{F})$ est admissible sur $X(\M{C})$ d'apr\`es
(\ref{produit_decomposables}),
et donc $\ov{E}\otimes\ov{F}$ est admissible.

On suppose maintenant que $\ov{E}$ et $\ov{F}$ sont int\'egrables. On peut
donc trouver $\ov{E}_{1}$, $\ov{E}_{2}$, $\ov{F}_{1}$ et $\ov{F}_{2}$ des
fibr\'es admissibles sur $X$ tels que $\ov{E} \simeq
\ov{E}_{1}\otimes(\ov{E}_{2})^{-1}$ et $\ov{F} \simeq \ov{F}_{1}\otimes
(\ov{F}_{2})^{-1}$; et donc $\ov{E}\otimes\ov{F} \simeq
(\ov{E}_{1}\otimes\ov{F}_{1})\otimes(\ov{E}_{2}\otimes\ov{F}_{2})^{-1}$, et
comme $\ov{E}_{1}\otimes\ov{F}_{1}$ et $\ov{E}_{2}\otimes\ov{F}_{2}$ sont
admissibles d'apr\`es ce qui pr\'ec\`ede, on en d\'eduit que
$\ov{E}\otimes\ov{F}$ est int\'egrable.
\medskip

Les deux propositions pr\'ec\'edentes nous autorisent \`a poser la d\'efinition
suivante~:
\begin{defn}
Soit $X$ une vari\'et\'e arithm\'etique. On note
$\widehat{\op{Pic}}_{\,\op{int}}(X)$ le groupe form\'e 
des classes d'isomorphie isom\'etrique des fibr\'es hermitiens int\'egrables 
sur $X$, 
et dont la structure de groupe est donn\'ee par le produit tensoriel.
\end{defn}
\begin{prop}
Soient $X$ et $Y$ deux vari\'et\'es arithm\'etiques et $f: Y \rightarrow X$ un
morphisme. Si $\ov{L}$ est un fibr\'e en droites admissible (resp.
int\'egrable) sur $X$, alors $f^{\ast}(\ov{L})$ est admissible (resp.
int\'egrable) sur $Y$.

L'application $f^{\ast}: \ov{L} \mapsto f^{\ast}(\ov{L})$ d\'efinit 
un morphisme de groupes~:
\[
f^{\ast}: \widehat{\op{Pic}}_{\,\op{int}}(X) \longrightarrow 
\widehat{\op{Pic}}_{\,\op{int}}(Y).
\]
\end{prop}
\demo\ C'est une cons\'equence directe de (\ref{image_reciproque})
et des d\'efinitions.
\bigskip

\subsection{Groupes de Chow arithm\'etiques g\'en\'eralis\'es}~
\label{sous_section_def}

Pour tout entier $p \geqslant 0$, on pose $\ZZ^{p}(X) = Z^{p}(X)
\oplus D^{p-1,p-1}(X_{\R})$; du fait de l'inclusion 
$\widehat{R}^{p}(X) \subset \ZZ^{p}(X)$, on peut d\'efinir le
groupe quotient~:
\[
\widetilde{CH}^{p}(X) = \left. \ZZ^{p}(X)\right/
\widehat{R}^{p}(X).
\]
On a imm\'ediatement l'inclusion $\widehat{CH}^{p}(X) \subset
\widetilde{CH}^{p}(X)$. On dispose des morphismes suivants~:
\begin{alignat*}{3}
\zeta :\; &\widetilde{CH}^{p}(X) & &\longrightarrow & &CH^{p}(X) \\
&[(Z,g)] & &\longmapsto & &[Z], \\
& & & & & \\
\omega :\; &\widetilde{CH}^{p}(X) & &\longrightarrow &  &D^{p,p}(X_{\R})\\
&[(Z,g)] & & \longmapsto & &\dd g + \delta_{Z},
\end{alignat*}
et 
\begin{alignat*}{3}
a :\; &\widetilde{\ov{A}}^{p-1,p-1}(X_{\R}) & &\longrightarrow &
&\widetilde{CH}^{p}(X) \\
&\eta & & \longmapsto & &[(0,\eta)],
\end{alignat*}
dont les restrictions \`a $\widehat{CH}^{p}(X)$ et
$\widetilde{A}^{p-1,p-1}(X_{\R})$ respectivement co\"\i ncident avec ceux
d\'efinis au (\ref{theorie_classique}).
\begin{defn}
\label{classe_de_chern}
Soit $\ov{L} = (L,\|.\|)$ un fibr\'e en droites hermitien int\'egrable sur
$X$. On appelle {\it premi\`ere classe de Chern arithm\'etique de $\ov{L}$\/}
et on note $\hat{c}_{1}(\ov{L})$ la classe dans $\widetilde{CH}^{1}(X)$ de~:
\[
\op{\widehat{div}}s := (\op{div} s, - \log \|s\|^{2}), 
\]
o\`u $s$ est une section rationnelle non nulle de $L$ au-dessus de $X$.
\end{defn}
\begin{rem}
La classe $\hat{c}_{1}(\ov{L})$ est bien d\'efinie et ne d\'epend pas du choix
de $s$. En effet, si $s'$ est une section rationnelle de $L$ au-dessus de $X$,
il existe $f$ une fonction rationnelle sur $X$ telle que $s' = f\cdot s$, et
l'on a~:
\[
\op{\widehat{div}}s' = (f,-\log|f|^{2}) + \op{\widehat{div}}s.
\]
\end{rem}
\begin{prop}
Soit $\ov{L} = (L,\|.\|)$ un fibr\'e en droites int\'egrable sur $X$, on a~:
\begin{alignat*}{3}
&\omega(\hat{c}_{1}(\ov{L}))& &= c_{1}(\ov{L})& &\in \AA^{1,1}_{\op{g}}(X_{\R})
\cap C_{0}^{1,1}(X_{\R})
\\
\text{et} \quad &\zeta(\hat{c}_{1}(\ov{L}))& &= c_{1}(L) & &\in CH^{1}(X).
\end{alignat*}
\end{prop}
\demo\ La premi\`ere \'egalit\'e est une cons\'equence directe de la formule de
Poincar\'e-Lelong g\'en\'eralis\'ee (\ref{PL_generalisee}); la seconde \'egalit\'e r\'esulte
imm\'ediatement des d\'efinitions.
\medskip

Soient $p$ et $q$ deux entiers positifs (\'eventuellement nuls), et soient $Z
\in Z^{p}(X)$ un cycle de codimension $p$, $g \in D^{p-1,p-1}(X_{\R})$ un
courant de Green pour $Z$ et $\ov{L}_{1} = (L_{1},\|.\|_{1}), \dots, \ov{L}_{q}
= (L_{q}, \|.\|_{q})$ des fibr\'es en droites admissibles sur $X$. On choisit
$s_{1}, \dots, s_{q}$ des sections rationnelles non identiquement nulles sur
$X$ de $L_{1}, \dots, L_{q}$ respectivement, telles que les cycles $Y_{0} = Z$,
$Y_{1}=\op{div}s_{1}, \dots, Y_{q} = \op{div}s_{q}$ soient d'intersection
propre (i.e. tels que $\op{codim} (Y_{j_{1}} \cap \dots \cap Y_{j_{m}})
\geqslant \sum_{i = 1}^{m}\op{codim}Y_{j_{i}}$ pour tout choix d'indices $j_{1}
< \dots < j_{m}$ dans $\{0,\dots,q\}$). On pose $\tilde{\omega} = \dd g + \delta_{Z}$.
Pour tout $1 \leqslant i \leqslant q$, on note \'egalement $\delta_{i} =
\delta_{\op{div}s_{i}} \in \AAA^{1,1}(X_{\R})$, $g_{i} = - \log
\|s_{i}\|_{i}^{2}$ et $\omega_{i} = c_{1}(\ov{L}_{i}) \in \AA^{1,1}(X_{\R})$.
On d\'efinit \`a partir de ces donn\'ees un \'el\'ement de
$D^{p+q-1,p+q-1}(X_{\R})$ que l'on notera $\{g\ast (g_{1},s_{1})\ast \dots \ast
(g_{q},s_{q})\}$, par l'\'egalit\'e~:
\begin{multline*}
\{g\ast (g_{1},s_{1})\ast \dots \ast
(g_{q},s_{q})\} = g\cdot\delta_{\op{div}s_{1}\cap\dots\cap\op{div}s_{q}} + 
g_{1}\cdot\tilde{\omega}\delta_{2}\dots\delta_{q} + \dots \\
+ g_{i}\cdot\tilde{\omega}\omega_{1}\dots
\omega_{i-1}\delta_{i+1}\dots\delta_{q} 
+ \dots + g_{q-1}\cdot\tilde{\omega}\omega_{1}\dots\omega_{q-2}\delta_{q} + 
g_{q}\cdot\tilde{\omega}\omega_{1}\dots\omega_{q-1}.
\end{multline*}
On remarque que cette expression a bien un sens; en effet le premier terme $g
\cdot\delta_{\op{div}s_{1}\cap\dots\cap\op{div}s_{q}}$ est bien d\'efini (cf.
\cite{13}, 2.1.3.2), et les autres le sont gr\^ace \`a 
la proposition (\ref{exemple_produit_formes}).
\begin{prop}
\label{approximation2}
Soient $\left(\|.\|_{1}^{(n)}\right)_{n \in \N}, \dots ,
\left(\|.\|_{q}^{(n)}\right)_{n \in \N}$ des suites croissantes de m\'etriques
positives $C^{\infty}$ sur $L_{1}, \dots, L_{q}$ convergeant vers $\|.\|_{1},
\dots, \|.\|_{q}$ respectivement. Si l'on note $g_{1}^{(n)} = - \log
\left(\|s_{1}\|_{1}^{(n)}\right)^{2}, \dots, g_{q}^{(n)} = - \log 
\left(\|s_{q}\|_{q}^{(n)}\right)^{2}$, alors
$\{g \ast (g_{1}^{(n)},s_{1})\ast \dots \ast (g_{q}^{(n)},s_{q})\}$ tend vers
$\{g\ast (g_{1},s_{1}) \ast \dots \ast (g_{q},s_{q})\}$ dans
$D^{p+q-1,p+q-1}(X_{\R})$, au sens de la topologie
faible, quand $n$ tend vers $+\infty$.
\end{prop}
\demo\ C'est une cons\'equence directe de (\ref{demailly}).
\medskip

\begin{prop}
\label{naturalite2}
Supposons que,
pour un certain $k \in \{1,\dots,q\}$, 
les m\'etriques $\|.\|_{1},\dots,\|.\|_{k}$ 
soient $C^{\infty}$ sur $X(\M{C})$. 
Les courants $g_{1},\dots,g_{k}$ sont alors des courants de Green pour les cycles
$\op{div}s_{1}, \dots, \op{div}s_{k}$ et on peut former le produit ~:
\[
(g\ast g_{1}\ast \dots \ast g_{k}) := g \ast (g_{1} \ast ( g_{2} \ast ( \dotsm 
\ast (g_{k}) \dotsm))),
\]
pris au sens de Gillet-Soul\'e (cf. \cite{13}, \S
2.1), qui est un courant de Green pour le cycle intersection $Y_{0}\cap Y_{1}
\cap \dots \cap Y_{q}$.
On a alors l'\'egalit\'e dans $\widetilde{D}^{p+q-1,p+q-1}(X_{\R})$
\[
\{(g\ast g_{1} \ast \dots \ast g_{k})\ast (g_{k+1},s_{k+1}) \ast \dots \ast
(g_{q},s_{q})\} = \{g\ast (g_{1},s_{1})\ast \dots \ast (g_{q},s_{q})\}.
\]
\end{prop}
\demo\ On montre tout d'abord le r\'esultat quand $k=q$. Pour tout $1 \leqslant
i \leqslant q$, le produit $g_{i}\ast\dots\ast g_{q}$ est un courant de Green
pour $\op{div}s_{i} \cap \dots \cap \op{div}s_{q}$ et on a $\omega(
\op{div}s_{i} \cap \dots \cap \op{div}s_{q}, g_{i}\ast\dots\ast g_{q}) =
\omega_{i}\dots\omega_{q}$. On a dans $\widetilde{D}^{p+q-1,p+q-1}(X_{\R})$~:
\[
\begin{split}
g\ast g_{1}\ast\dots\ast g_{q} & = g\ast(g_{1}\ast(\dots \ast(g_{q-1}\ast
g_{q})\dots )) \\
&= g\cdot \delta_{\op{div}s_{1} \cap \dots \cap \op{div}s_{q}} + 
\tilde{\omega}(g_{1}\ast(\dots \ast(g_{q-1}\ast g_{q})\dots ))\\
= g\cdot \delta_{\op{div}s_{1} \cap \dots \cap \op{div}s_{q}} & \\
+ 
\tilde{\omega}(g_{1}\cdot &\delta_{\op{div}s_{2}\cap \dots \cap \op{div}s_{q}} + 
\omega_{1}g_{2}\cdot \delta_{\op{div}s_{3}\cap \dots\cap \op{div}s_{q}} + 
\dots + \omega_{1}\dots\omega_{q-1}g_{q}).
\end{split}
\]
Il suffit alors de prouver que pour tout $1 \leqslant i \leqslant (q-1)$, on a
dans $D^{q-i,q-i}(X_{\R})$ l'\'egalit\'e~:
\begin{equation}
\label{eq_inter_1}
g_{i}\cdot \delta_{i+1}\dots \delta_{q} = g_{i}\cdot
\delta_{\op{div}s_{i+1}\cap \dots \cap \op{div}s_{q}}, 
\end{equation}
o\`u le premier produit est pris au sens de (\ref{produit_generalise}) et le 
second au sens de
Gillet-Soul\'e (cf. \cite{13}, 2.1.3.2). 
Le probl\`eme \'etant local, on peut supposer que le diviseur $\op{div}s_{i}$
est effectif.
D'apr\`es (\cite{7}, prop. 3.4.12), on
a~:
\[
\delta_{i+1}\dots \delta_{q} = \delta_{\op{div}s_{i+1}\cap \dots \cap
\op{div}s_{q}}.
\]
Pour tout $t \in \R$, on pose $g_{i}^{(t)} = \max (g_{i}, t)$. Pour tout $t \in
\R$ fix\'e, la fonction $x \mapsto g_{i}^{(t)}(x)$ est \psh\ et born\'ee sur
$X(\M{C})$; on a donc~:
\[
g_{i}^{(t)}\cdot \delta_{i+1}\dots \delta_{q} = g_{i}^{(t)}\cdot
\delta_{\op{div}s_{i+1} \cap \dots \cap \op{div}s_{q}}.
\]
D'apr\`es (\ref{demailly}), le courant $g_{i}^{(t)}\cdot \delta_{i+1}\dots \delta_{q}$ tend
vers $g_{i}\cdot \delta_{i+1}\dots \delta_{q}$ au sens de la topologie faible
quand $t$ tend vers $-\infty$.

Pour finir de montrer (\ref{eq_inter_1}), il suffit de prouver que
$g_{i}^{(t)}\cdot\delta_{\op{div}s_{i+1} \cap \dots \cap \op{div}s_{q}}$ tend
faiblement vers $g_{i} \cdot\delta_{\op{div}s_{i+1} \cap \dots \cap \op{div}s_{q}}$
quand $t$ tend vers $-\infty$. Pour cela \'ecrivons~:
\[
(\op{div}s_{i+1})_{\M{C}} \cap \dots \cap (\op{div}s_{q})_{\M{C}}
= 
\sum_{k}n_{k}\, C_{k},
\]
o\`u les $C_{k}$ sont les composantes irr\'eductibles de
$(\op{div}s_{i+1})_{\M{C}} \cap \dots \cap
(\op{div}s_{q})_{\M{C}}$.

Nous devons v\'erifier que pour chaque $k$, 
le courant 
$g_{i}^{(t)}\cdot \delta_{C_{k}}$ tend faiblement vers $g_{i}\cdot
\delta_{C_{k}}$,
o\`u le produit $g_{i}\cdot\delta_{C_{k}}$ est d\'efini comme dans Gillet-Soul\'e
(\cite{13}, \S 2.1.3.2).
Soit $\pi : \widetilde{C}_{k} \rightarrow C_{k}$ une r\'esolution des
singularit\'es de $C_{k}$. Par d\'efinition, on a~:
\[
\delta_{C_{k}} = \op{\pi_{\ast}} 1, 
\]
et donc ~:
\[
g_{i}^{(t)}\cdot \delta_{C_{k}} = \op{\pi_{\ast}}(\op{\pi^{\ast}}g_{i}^{(t)}),
\]
car $\delta_{C_{k}}$ est un courant d'ordre $0$ et $g_{i}^{(t)}$ est born\'ee.
On a par ailleurs (cf. \cite{13}, \S 2.1.3.2)~:
\[
g_{i}\cdot \delta_{C_{k}} = \op{\pi_{\ast}}(\op{\pi^{\ast}}g_{i}).
\]
Enfin $\op{\pi^{\ast}}g_{i}$ est \psh\ et $\op{\pi^{\ast}}g_{i}^{(t)} =
\op{max}(\op{\pi^{\ast}}g_{i}, t)$ converge faiblement vers $\op{\pi^{\ast}}g_{i}$ quand 
$t$ tend vers $-\infty$, d'o\`u le
r\'esultat par continuit\'e faible de $\pi_{\ast}$.

On s'int\'eresse maintenant au cas g\'en\'eral~: Soient
$\left(\|.\|^{(n)}_{k+1}\right)_{n \in \N}, \dots,
\left(\|.\|^{(n)}_{q}\right)_{n\in \N}$ des suites croissantes de m\'etriques
positives $C^{\infty}$ convergent uniform\'ement sur $X(\M{C})$ vers
$\|.\|_{k+1}, \dots, \|.\|_{q}$ respectivement. On note, pour tout $n \in \N$,
$g_{k+1}^{(n)}, \dots, g_{q}^{(n)}$ les courants $- \log
\left(\|s_{k+1}\|^{(n)}_{k+1}\right)^{2}, \dots, - \log 
\left(\|s_{q}\|_{q}^{(n)}\right)^{2}$. D'apr\`es ce qui
pr\'ec\`ede, on a dans $\widetilde{D}^{p+q-1,p+q-1}(X_{\R})$ les \'egalit\'es~:
\[
\begin{split}
\{(g\ast g_{1}\ast \dots \ast g_{k})\ast (&g_{k+1}^{(n)},s_{k})\ast \dots \ast 
(g_{q}^{(n)},s_{q})\} \\
&= (g\ast g_{1}\ast \dots g_{k}) \ast g_{k+1}^{(n)} \ast
\dots \ast g_{q}^{(n)} \\
&= g \ast g_{1} \ast \dots g_{k} \ast g_{k+1}^{(n)}\ast \dots g_{q}^{(n)} \\
&= \{g\ast (g_{1},s_{1})\ast \dots \ast (g_{k},s_{k})\ast 
(g_{k+1}^{(n)},s_{k+1})\ast \dots \ast (g_{q}^{(n)},s_{q})\}.
\end{split}
\]
Comme $X(\M{C})$ est projective, les sous-espaces $\op{Im}\partial$ et
$\op{Im}\ov{\partial}$ sont ferm\'es
(puisque les morphismes $\partial$ et $\ov{\partial}$
sont continues et que leurs groupes de cohomologie
sont de dimension finie).
On en d\'eduit que $\widetilde{D}^{p+q-1,p+q-1}(X_{\R})$,
muni de la topologie quotient de la topologie faible sur 
${D}^{p+q-1,p+q-1}(X_{\R})$,
est s\'epar\'e. 

Pour \'etablir l'\'egalit\'e recherch\'ee~:
\[
\{(g\ast g_{1}\ast \dots \ast g_{k})\ast (g_{k+1},s_{k+1})\ast \dots \ast
(g_{q},s_{q})] = [g \ast (g_{1},s_{1}) \ast \dots \ast (g_{q},s_{q})\},
\]
dans $\widetilde{D}^{p+q-1,p+q-1}(X_{\R})$, il suffit 
d'observer que d'apr\`es 
la proposition (\ref{approximation2}) on a, pour cette topologie, les
limites~:
\begin{multline*}
\lim_{n \rightarrow +\infty}\{(g\ast g_{1}\ast \dots \ast g_{k})\ast
(g_{k+1}^{(n)},s_{k+1})\ast \dots \ast (g_{q}^{(n)},s_{q})\} \\
= 
\{(g \ast g_{1} \ast \dots \ast g_{k})\ast (g_{k+1},s_{k+1})\ast \dots \ast
(g_{q},s_{q})\}
\end{multline*}
et 
\begin{multline*}
\lim_{n \rightarrow +\infty}\{g \ast (g_{1},s_{1})\ast \dots \ast
(g_{k},s_{k})\ast (g_{k+1}^{(n)},s_{k+1})\ast \dots \ast (g_{q}^{(n)},s_{q})\}\\
= \{g \ast (g_{1},s_{1})\ast \dots \ast (g_{q},s_{q})\}.
\end{multline*}
\medskip

On en d\'eduit aussit\^ot la proposition suivante~:
\begin{prop}
\label{commutativite2}
La classe du courant $\{g\ast (g_{1},s_{1})\ast \dots \ast (g_{q},s_{q})\}$
dans\\
$\widetilde{D}^{p+q-1,p+q-1}(X_{\R})$ est ind\'ependante de l'ordre des
$\ov{L}_{1},\dots,\ov{L}_{q}$.
\end{prop}
\demo\ Soit $\sigma$ une permutation de $\{1,\dots,q\}$. Reprenant les
notations de la proposition (\ref{approximation2}), on sait que~:
\begin{multline*}
\delta_{n}:= \{g\ast (g^{(n)}_{1},s_{1})\ast \dots \ast (g^{(n)}_{q},s_{q})\} \\
- 
\{g\ast (g^{(n)}_{\sigma(1)},s_{\sigma(1)})\ast \dots \ast
(g^{(n)}_{\sigma(q)},s_{\sigma(q)})\}
\in \left(\op{Im}\partial + \op{Im}\ov{\partial}\right).
\end{multline*}
De plus, $\delta_{n}$ tend vers $\{g\ast (g_{1},s_{1})\ast \dots \ast
(g_{q},s_{q})\}
- \{g\ast (g_{\sigma(1)},s_{\sigma(1)})\ast \dots \ast
(g_{\sigma(q)},s_{\sigma(q)})\}
$ au sens de la topologie faible quand $n$ tend vers $+\infty$ d'apr\`es
(\ref{approximation2}). La proposition d\'ecoule alors imm\'ediatement du fait
que $\op{Im}\partial + \op{Im}\ov{\partial}$ est ferm\'e dans
$D^{p+q-1,p+q-1}(X_{\R})$ puisque $X(\M{C})$ est projective.
\medskip

On d\'efinit un \'el\'ement de $\widetilde{CH}^{p+q}(X)$ que l'on note
provisoirement\\
$[(Z,g)\cdot \hat{c}_{1}(\ov{L}_{1},s_{1})\dotsm
\hat{c}_{1}(\ov{L}_{q},s_{q})]$ par l'\'egalit\'e~:
\begin{multline}
\label{definition_prod2}
[(Z,g)\cdot \hat{c}_{1}(\ov{L}_{1},s_{1})\dotsm \hat{c}_{1}(\ov{L}_{q},s_{q})] \\
= [(Z \cdot \op{div}s_{1} \dotsm \op{div}s_{q},
\{g\ast(g_{1},s_{1})\ast \dots \ast(g_{q},s_{q})\})].
\end{multline}
\begin{prop}
\label{independance2}
La classe de $[(Z,g)\cdot \hat{c}_{1}(\ov{L}_{1},s_{1})\dotsm
\hat{c}_{1}(\ov{L}_{q},s_{q})]$
dans $\widetilde{CH}^{p+q}(X)$ ne d\'epend que de la classe de $(Z,g)$ dans
$\widehat{CH}^{p}(X)$ et des classes d'isomorphie isom\'etrique des fibr\'es
admissibles $\ov{L}_{1},\dots,\ov{L}_{q}$; elle ne d\'epend pas de l'ordre
des fibr\'es $\ov{L}_{1},\dots,\ov{L}_{q}$.
\end{prop}
\demo\ Le fait que la classe de $[(Z,g)\cdot \hat{c}_{1}(\ov{L}_{1},s_{1})\dotsm
\hat{c}_{1}(\ov{L}_{q},s_{q})]$ ne d\'epende pas de l'ordre des fibr\'es $
\ov{L}_{1},\dots,\ov{L}_{q}$
est une simple cons\'equence de
(\ref{commutativite2}).

On s'int\'eresse maintenant \`a la premi\`ere partie de
l'\'enonc\'e. 
Il suffit, pour prouver le r\'esultat recherch\'e, de montrer que pour toute
section $s_{1}'$ de $L_{1}$ au-dessus de $X$ et tout repr\'esentant $(Z',g')$
de $[(Z,g)]$ dans $\widehat{CH}^{p}(X)$ tels que les cycles $Z',\op{div}s_{1},
\dots, \op{div}s_{q}$ et les cycles $Z',\op{div}s_{1}',\op{div}s_{2},
\dots,\op{div}s_{q}$ soient d'intersection propre, on a l'\'egalit\'e~:
\[
[(Z',g')\cdot \hat{c}_{1}(\ov{L}_{1},s_{1}')\cdot
\hat{c}_{1}(\ov{L}_{2},s_{2})\dotsm \hat{c}_{1}(\ov{L}_{q},s_{q})]
= 
[(Z,g)\cdot \hat{c}_{1}(\ov{L}_{1},s_{1})\dotsm \hat{c}_{1}(\ov{L}_{q},s_{q})].
\]
On choisit $\left(\|.\|_{1}^{(n)}\right)_{n \in \N}, \dots,
\left(\|.\|_{q}^{(n)}\right)_{n \in \N}$, $q$-suites croissantes de m\'etriques
positives $C^{\infty}$ sur $L_{1},\dots,L_{q}$ convergeant uniform\'ement sur
$X(\M{C})$ vers $\|.\|_{1},\dots,\|.\|_{q}$ respectivement.
On note $f_{1}$ la fonction rationnelle sur $X$ telle que $s_{1}' = f_{1}\cdot
s_{1}$. On note \'egalement $g_{1}' = -\log \|s_{1}'\|_{1}^{2}$ et pour tout $1
\leqslant i \leqslant q$ et tout $n \in \N$, $g_{i}^{(n)} = - \log 
\left(\|s_{i}\|_{i}^{(n)}\right)^{2}$ et ${g_{1}'}^{(n)} = - \log
\left(\|s_{1}'\|_{1}^{(n)}\right)^{2}$.

On d\'eduit de l'\'egalit\'e $[Z'] = \zeta(Z',g') = \zeta(Z,g) = [Z]$ qu'il
existe une $K_{1}$-chaine $f$ telle que~:
\begin{equation}
\label{eq_fin1}
Z' = Z + \op{div}f.
\end{equation}
D'apr\`es le lemme de d\'eplacement pour les $K_{1}$-chaines, (cf. \cite{13},
\S 4.2.6), on peut de plus supposer que $f$ rencontre $\op{div}s_{1}, \dots,
\op{div}s_{q}$ presque proprement sur $X_{K}$.
Notons~:
\[
Y = Z'\cdot \op{div}s_{1}'\cdot \op{div}s_{2}\dotsm \op{div}s_{q}
- Z \cdot \op{div}s_{1}\dotsm \op{div}s_{q}.
\]
Il vient~:
\begin{align*}
Y &= Z'\cdot (\op{div}s_{1}' - 
\op{div}s_{1})\cdot \op{div}s_{2}\dotsm \op{div}s_{q}
+ (Z' - Z)\cdot \op{div}s_{1}\dotsm \op{div}s_{q}\\
&= 
\op{div}f_{1}\cdot Z'\cdot \op{div}s_{2} \dotsm \op{div}s_{q} + 
\op{div}f\cdot \op{div}s_{1}\dotsm \op{div}s_{q}, 
\end{align*}
ce que, d'apr\`es (\cite{13}, \S 4.2.5) on peut r\'ecrire~:
\begin{equation}
\label{eq_fin3}
Y = \op{div}(f_{1}\cdot (Z'\cdot \op{div}s_{2}\dotsm \op{div}s_{q}))
+ \op{div}(f\cdot (\op{div}s_{1}\dotsm \op{div}s_{q})) + R, 
\end{equation}
o\`u $f_{1}\cdot (Z'\cdot \op{div}s_{2}\dotsm \op{div}s_{q})$ (resp. $f\cdot
(\op{div}s_{1}\dotsm \op{div}s_{q})$) d\'esigne une $K_{1}$-chaine
repr\'esentant l'intersection de la $K_{1}$-chaine $f_{1}$ et du cycle
$Z'\cdot\op{div}s_{2}\dotsm \op{div}s_{q}$ (resp. de la $K_{1}$-chaine $f$ et
du cycle $\op{div}s_{1}\dotsm\op{div}s_{q}$) et o\`u $R \in
R_{\op{fin}}^{p+q}(X)$.
Posons~:
\[
\delta = \{g'\ast(g'_{1},s_{1}')\ast (g_{2},s_{2})\ast \dotsm \ast(g_{q},s_{q})\}
- 
\{g\ast(g_{1},s_{1})\ast\dotsm \ast(g_{q},s_{q})\},
\]
et pour tout $n \in \N$,
\[
\delta_{n} = \{g'\ast({g_{1}'}^{(n)},s_{1}')\ast(g_{2}^{(n)},s_{2})\ast \dotsm
\ast(g_{q}^{(n)},s_{q})\}
- 
\{g\ast (g_{1}^{(n)},s_{1})\ast \dotsm \ast (g_{q}^{(n)},s_{q})\}.
\]
D'apr\`es la proposition (\ref{naturalite2}), il vient~:
\begin{align*}
\delta_{n} &= g'\ast {g_{1}'}^{(n)}\ast g_{2}^{(n)}\ast \dotsm 
\ast g_{q}^{(n)}
- g \ast g_{1}^{(n)}\ast \dotsm \ast g_{q}^{(n)}\\
&= (g'\ast {g_{1}'}^{(n)} - g \ast g_{1}^{(n)})\ast g_{2}^{(n)}\ast\dotsm 
\ast g_{q}^{(n)}.
\end{align*}
Par ailleurs, on tire de l'\'egalit\'e $[(Z',g')] = [(Z,g)]$ et de la relation 
(\ref{eq_fin1}) qu'il existe des courants $u$ et $v$ tels que~:
\[
(Z',g') = (Z,g) + \op{\widehat{\op{div}}}f + (0,\partial u + 
\ov{\partial}v), 
\]
ce dont on d\'eduit que~:
\[
g' = g + (-\log |f|^{2}), 
\]
dans $\widetilde{D}^{p-1,p-1}(X_{\R})$. Ceci combin\'e \`a la relation~:
\[
{g_{1}'}^{(n)} = g_{1}^{(n)} + (- \log |f_{1}|^{2}), 
\]
et \`a (\cite{13}, \S 4.2.5.1) montre que~:
\begin{align*}
g'\ast {g_{1}'}^{(n)} - g\ast g_{1}^{(n)}
&= 
g'\ast ({g'_{1}}^{(n)} - g_{1}^{(n)}) + (g' - g)\ast g_{1}^{(n)} \\
&= (- \log |f_{1}|^{2})\ast g' + (-\log|f|^{2})\ast g_{1}^{(n)} \\
&= (- \log |f_{1}|^{2})\wedge \delta_{Z'} + (- \log|f|^{2})\wedge \op{div}s_{1}
\\
&= (- \log|f_{1}\cdot Z'|^{2}) + (- \log |f\cdot \op{div}s_{1}|^{2}); 
\end{align*}
et donc que~:
\begin{align*}
\delta_{n} &= ((-\log |f_{1}\cdot Z'|^{2}) 
+ (-\log |f \cdot \op{div}s_{1}|^{2})\ast g_{2}^{(n)} \ast \dotsm \ast
g_{q}^{(n)} \\
&= (- \log |f_{1}\cdot Z'\cdot \op{div}s_{2}\dotsm \op{div}s_{q}|^{2})
+ (- \log |f\cdot \op{div}s_{1} \dotsm \op{div}s_{q}|^{2}), 
\end{align*}
dans $\widetilde{D}^{p+q-1,p+q-1}(X_{\R})$. Comme $\delta = \lim_{n \rightarrow
+\infty}\delta_{n}$, on a de m\^eme~:
\begin{equation}
\label{eq_fin2}
\delta = (-\log|f_{1}\cdot Z'\cdot \op{div}s_{2}\dotsm \op{div}s_{q}|^{2})
+ (- \log|f\cdot \op{div}s_{1}\dotsm \op{div}s_{q}|^{2}).
\end{equation}
Finalement, on d\'eduit de la d\'efinition (\ref{definition_prod2}) et des relations
(\ref{eq_fin3}) et (\ref{eq_fin2}) que~:
\begin{multline*}
[(Z',g')\cdot \hat{c}_{1}(\ov{L}_{1},s_{1}')\cdot \hat{c}_{1}(\ov{L}_{2},s_{2})
\dotsm \hat{c}_{1}(\ov{L}_{q},s_{q})]
- [(Z,g)\cdot \hat{c}_{1}(\ov{L}_{1},s_{1})\dotsm
\hat{c}_{1}(\ov{L}_{q},s_{q})]\\
= [(Y,\delta)] 
= [\op{\widehat{\op{div}}}(f_{1}\cdot Z'\cdot \op{div}s_{2}\dotsm
\op{div}s_{q})] + 
[\op{\widehat{\op{div}}}(f\cdot \op{div}s_{1}\dotsm \op{div}s_{q})]
+ [(R,0)] = 0.
\end{multline*}
\medskip

La proposition pr\'ec\'edente nous autorise \`a noter d\'esormais $[(Z,g)]\cdot
\hat{c}_{1}(\ov{L}_{1})\dotsm \hat{c}_{1}(\ov{L}_{q})$ la classe $[(Z,g)\cdot
\hat{c}_{1}(\ov{L}_{1},s_{1})\dotsm \hat{c}_{1}(\ov{L}_{q},s_{q})]$.
\begin{prop}
\label{linearite2}
Soient $\alpha \in \widehat{CH}^{p}(X)$ et $\ov{L}_{0}$,
$\ov{L}_{1},\dots, \ov{L}_{q}$ des
fibr\'es en droites admissibles sur $X$. Pour tout $1 \leqslant i \leqslant q$, on a
dans $\widetilde{CH}^{p+q}(X)$ la relation~:
\begin{multline*}
\alpha \cdot \hat{c}_{1}(\ov{L}_{1})\dotsm \hat{c}_{1}(\ov{L}_{i-1})\hat{c}_{1}
(\ov{L}_{i}\otimes \ov{L}_{0})\hat{c}_{1}(\ov{L}_{i+1})\dotsm
\hat{c}_{1}(\ov{L}_{q}) \\
=
\alpha \cdot\hat{c}_{1}(\ov{L}_{1})\dotsm \hat{c}_{1}(\ov{L}_{q})
+ \alpha \cdot \hat{c}_{1}(\ov{L}_{1})\dotsm \hat{c}_{1}(\ov{L}_{i-1})
\hat{c}_{1}(\ov{L}_{0})\hat{c}_{1}(\ov{L}_{i+1})\dotsm \hat{c}_{1}(\ov{L}_{q}).
\end{multline*}
\end{prop}
\demo\ D'apr\`es (\ref{independance2}), il suffit de d\'emontrer le cas $i=1$.
Soit $(Z,g)$ un repr\'esentant de la classe $\alpha \in \widehat{CH}^{p}(X)$.
Soient \'egalement
$s_{1},\dots,s_{q}$ des sections rationnelles au-dessus de $X$ de $L_{1},
\dots, L_{q}$ respectivement, telles que $Z$, $\op{div}s_{1}, \dots,
\op{div}s_{q}$ s'intersectent proprement; et $s_{0}$ une section
rationnelle de $L_{0}$ au-dessus de $X$ telle que $Z$, $\op{div}s_{0}$,
$\op{div}s_{2}, \dots, \op{div}s_{q}$ s'intersectent proprement. 
Pour tout $0 \leqslant i \leqslant q$, on note $\delta_{i} =
\delta_{\op{div}s_{i}}$, $g_{i} = - \log\|s_{i}\|_{i}^{2}$ et $\omega_{i} =
c_{1}(\ov{L}_{i})$. On note \'egalement $\tilde{\omega} = \omega(Z,g)$. On a~:
\begin{multline*}
[g \ast (g_{0} + g_{1},s_{0} + s_{1})\ast (g_{2},s_{2}) \ast \dots \ast
(g_{q},s_{q})] = 
g\cdot \delta_{(\op{div}s_{0}\cup\op{div}s_{1})\cap \op{div}s_{2}\cap \dots
\cap \op{div}s_{q}}\\ + \tilde{\omega}(g_{0}+g_{1})\cdot
\delta_{2}\dots\delta_{q}
+ \tilde{\omega}(\omega_{0}+ \omega_{1})g_{2}\cdot \delta_{3}\dots \delta_{q}
+ \dots + \tilde{\omega}(\omega_{0} +
\omega_{1})\omega_{2}\dots\omega_{q-1}g_{q}
\\
= [g \ast (g_{1},s_{1})\ast \dots \ast (g_{k},s_{k})] + [g \ast
(g_{0},s_{0})\ast (g_{2},s_{2})\ast \dots \ast (g_{q},s_{q})],
\end{multline*}
et comme~:
\begin{multline*}
Z\cdot (\op{div}s_{0} + \op{div}s_{1})\cdot \op{div}s_{2}\dotsm \op{div}s_{q}
= Z \cdot \op{div}s_{1}\cdot \op{div}s_{2} \dotsm \op{div}s_{q} \\
+ Z \cdot \op{div}s_{0} \cdot \op{div}s_{2} \dotsm \op{div}s_{q}
\end{multline*}
dans $Z^{p+q}(X)$, on a bien le r\'esultat cherch\'e.
\medskip

Plus g\'en\'eralement, soient $\alpha \in \widehat{CH}^{p}(X)$ et
$\ov{L}_{1}, \dots, \ov{L}_{q}$ des fibr\'es en droites int\'egrables sur $X$. On choisit
$\ov{E}_{1},\dots,\ov{E}_{q}$ et $\ov{F}_{1},\dots,\ov{F}_{q}$ des fibr\'es
en droites admissibles sur $X$ tels que pour tout $1 \leqslant i \leqslant q$, on ait~: 
$\ov{L}_{i} = \ov{E}_{i}\otimes \left(\ov{F}_{i}\right)^{-1}$.
\begin{prop_defn}
\label{produit_general2}
La classe dans $\widetilde{CH}^{p+q}(X)$ donn\'ee par~:
\[
\sum_{\substack{S_{1},S_{2} \\ S_{1}\cup S_{2} = \{1,\dots,q\}}}
(-1)^{\#S_{2}}\;
\alpha \cdot \prod_{i \in S_{1}}\hat{c}_{1}(\ov{E}_{i})\cdot
\prod_{j \in S_{2}}\hat{c}_{1}(\ov{F}_{j})
\]
ne d\'epend que de la classe $\alpha \in \widehat{CH}^{p}(X)$ et des
classes d'isomorphie isom\'etrique des fibr\'es $\ov{L}_{1}, \dots,\ov{L}_{q}$.
On note $\alpha\cdot \hat{c}_{1}(\ov{L}_{1})\dotsm \hat{c}_{1}(\ov{L}_{q})$
cette classe. Si les fibr\'es $\ov{L}_{1}, \dots, \ov{L}_{q}$ sont admissibles,
cet \'el\'ement co\"\i ncide avec celui d\'efini en (\ref{definition_prod2}).
\end{prop_defn}
\demo\ C'est une cons\'equence directe de la proposition (\ref{linearite2}) et du fait
que toute application $q$-lin\'eaire $\varphi: \underbrace{A\times\dots\times
A}_{\text{$q$ fois}} \rightarrow B$, o\`u $A$ est un semi-groupe ab\'elien et
$B$ un groupe ab\'elien, s'\'etend de mani\`ere unique en une application
$q$-lin\'eaire $\varphi_{s}: \underbrace{A_{s}\times \dots\times
A_{s}}_{\text{$q$ fois}} \rightarrow B$, o\`u
$A_{s}$ est le groupe sym\'etris\'e de $A$.
\medskip

Diverses propri\'et\'es de cette construction sont rassembl\'ees dans le
th\'eor\`eme suivant~:
\begin{thm}
\label{th_general}
Soient $\alpha$ un \'el\'ement de
$\widehat{CH}^{p}(X)$ et $\ov{L}_{1}, \dots, \ov{L}_{q}$ 
des fibr\'es en
droites int\'egrables sur $X$. On a les propri\'et\'es suivantes~:
\begin{enumerate}
\item{La classe $\alpha\cdot 
\hat{c}_{1}(\ov{L}_{1})\dotsm
\hat{c}_{1}(\ov{L}_{q})$ ne d\'epend pas de l'ordre de
$\ov{L}_{1},\dots,\ov{L}_{q}$.}
\item{ 
L'application qui \`a $\alpha \in \widehat{CH}^{p}(X)$ et
$\ov{L}_{1},\dots,\ov{L}_{q}$ des fibr\'es en droites hermitiens int\'egrables 
sur $X$ associe $\alpha\cdot \hat{c}_{1}(\ov{L}_{1})\dotsm
\hat{c}_{1}(\ov{L}_{q})$, d\'efinit une application multilin\'eaire~:
\[ 
\widehat{CH}^{p}(X)\times \widehat{\op{Pic}}_{\,\op{int}}(X)\times \dots \times
\widehat{\op{Pic}}_{\,\op{int}}(X)
\longrightarrow  \widetilde{CH}^{p+q}(X).
\]
}
\item{Si les m\'etriques des fibr\'es $\ov{L}_{1},\dots,\ov{L}_{k}$
pour $1 \leqslant k \leqslant q$ fix\'e, sont $C^{\infty}$ sur $X(\M{C})$, 
alors on a~:
\[
\alpha \cdot \hat{c}_{1}(\ov{L}_{1})\dotsm
\hat{c}_{1}(\ov{L}_{q})
= (\alpha \cdot\hat{c}_{1}(\ov{L}_{1})\dotsm
\hat{c}_{1}(\ov{L}_{k}))
\cdot\hat{c}_{1}(\ov{L}_{k+1})\dotsm\hat{c}_{1}(\ov{L}_{q}),
\]
o\`u le produit $(\alpha\cdot\hat{c}_{1}(\ov{L}_{1})\dotsm
\hat{c}_{1}(\ov{L}_{k})) \in \widehat{CH}^{p+k}(X)$
dans le second membre est pris au sens de Gillet-Soul\'e. (Voir \cite{13},
th. 4.2.3).
En particulier si les m\'etriques sur $\ov{L}_{1},\dots,\ov{L}_{q}$ sont
$C^{\infty}$, le produit $\alpha\cdot \hat{c}_{1}(\ov{L}_{1})\dotsm
\hat{c}_{1}(\ov{L}_{q})$ d\'efini par la proposition (\ref{produit_general2}) 
co\"\i ncide avec
celui d\'efini par Gillet-Soul\'e.}
\item{Soit $\boldsymbol{1} = [(X,0)] \in \widehat{CH}^{0}(X)$. On a~:
\[
\boldsymbol{1}\cdot \hat{c}_{1}(\ov{L}_{1})
= \hat{c}_{1}(\ov{L}_{1}). 
\]
}
\item{On a les relations~:
\begin{multline*}
\qquad \quad \; \omega (\alpha\cdot \hat{c}_{1}(\ov{L}_{1})\dotsm
\hat{c}_{1}(\ov{L}_{q})) \\ =
\omega(\alpha)\cdot c_{1}(\ov{L}_{1})\dotsm c_{1}(\ov{L}_{q})
\in \AA_{\op{g}}^{p+q,p+q}(X_{\R}) \cap C_{0}^{p+q,p+q}(X_{\R})
\end{multline*}
et
\[
\zeta (\alpha\cdot \hat{c}_{1}(\ov{L}_{1})\dotsm
\hat{c}_{1}(\ov{L}_{q})) = 
\zeta(\alpha) \cdot c_{1}(L_{1})\dotsm c_{1}(L_{q}) \in CH^{p+q}(X).
\]
}
\item{Enfin, si $\alpha = (0,\varphi)$, avec $\varphi \in
A^{p-1,p-1}(X_{\R})$, alors~:
\[
\alpha \cdot \hat{c}_{1}(\ov{L}_{1})\dotsm
\hat{c}_{1}(\ov{L}_{q}) = 
[(0,\varphi\cdot c_{1}(\ov{L}_{1})\dotsm c_{1}(\ov{L}_{q}))].
\]
}
\end{enumerate}
\end{thm}
\demo\ 
\begin{enumerate}
\item{Cela d\'ecoule de la proposition (\ref{independance2}) et 
de la d\'efinition donn\'ee \`a la proposition (\ref{produit_general2}).}
\item{Cela d\'ecoule des propositions (\ref{linearite2}) 
et (\ref{produit_general2}).}
\item{On se ram\`ene gr\^ace au (2) au cas o\`u
$\ov{L}_{1},\dots,\ov{L}_{q}$ sont des fibr\'es en droites admissibles, on
applique alors (\ref{naturalite2}).}
\item{Il suffit de le faire pour $\ov{L}_{1}$ admissible et c'est alors 
\'evident.}
\item{D'apr\`es le (2), il suffit de d\'emontrer les relations pour
$\ov{L}_{1}, \dots, \ov{L}_{q}$ admissibles. 

L'\'egalit\'e
$\zeta (\alpha\cdot \hat{c}_{1}(\ov{L}_{1})\dotsm
\hat{c}_{1}(\ov{L}_{q})) = \zeta(\alpha)\cdot c_{1}(L_{1}) 
\cdots c_{1}(L_{q})$
est une cons\'equence imm\'ediate de la d\'efinition (\ref{definition_prod2}).

Soient $\left(\|.\|_{1}^{(n)}\right)_{n \in
\N},\dots,$ $\left(\|.\|_{q}^{(n)}\right)_{n \in \N}$ des suites croissantes de
m\'etriques positives $C^{\infty}$ convergeant uniform\'ement sur $X(\M{C})$
vers $\|.\|_{1},\dots,\|.\|_{q}$ respectivement. D'apr\`es le (3) 
et (\cite{13}, th. 4.2.9), on a~:
\begin{equation}
\begin{split}
\label{relation_intermediaire}
\omega (\alpha\cdot\hat{c}_{1}(L_{1},\|.\|_{1}^{(n)})\dotsm
\hat{c}_{1}(L_{q},\|.&\|_{q}^{(n)})) \\
&= \omega(\alpha)\cdot c_{1}(L_{1},\|.\|_{1}^{(n)})\dotsm
c_{1}(L_{q},\|.\|_{q}^{(n)}).
\end{split}
\end{equation}
Soient $(Z,g)$ un repr\'esentant de la classe $\alpha \in 
\widehat{CH}^{p}(X)$ et $s_{1},\dots,s_{q}$ des sections
rationnelles non identiquement nulles sur $X$ de $L_{1},\dots, L_{q}$
respectivement, telles que les cycles $Z$, $\op{div}s_{1}, \dots,
\op{div}s_{q}$ s'intersectent proprement. Pour tout $1 \leqslant i \leqslant
q$, on note $g_{i} = - \log \|s_{i}\|_{i}^{2}$ et pour tout $1 \leqslant i
\leqslant q$ et tout $n \in \N$, on note $g_{i}^{(n)} = - \log
\left(\|s_{i}\|_{i}^{(n)}\right)^{2}$. On a, pour tout $n \in \N$,
l'\'egalit\'e de courants~:
\begin{multline*}
\qquad \omega (\alpha\cdot\hat{c}_{1}(L_{1},\|.\|_{1}^{(n)})\dotsm
\hat{c}_{1}(L_{q},\|.\|_{q}^{(n)})) = \\
\dd(\{g\ast (g_{1}^{(n)},s_{1})\ast \dots \ast (g_{q}^{(n)},s_{q})\}) + 
\delta_{Z \cap \op{div}s_{1} \cap \dots \cap \op{div}s_{q}}.
\end{multline*}
Comme $\{g\ast (g_{1}^{(n)},s_{1})\ast \dots \ast (g_{q}^{(n)},s_{q})\}$ tend
vers $\{g\ast (g_{1},s_{1})\ast \dots \ast (g_{q},s_{q})\}$ quand $n$ tend vers
$+\infty$ d'apr\`es (\ref{approximation2}), on d\'eduit de la continuit\'e faible de
l'op\'erateur $\dd$ et de la relation~:
\begin{multline*}
\qquad \omega (\alpha\cdot\hat{c}_{1}(L_{1},\|.\|_{1})\dotsm
\hat{c}_{1}(L_{q},\|.\|_{q})) = \\
\dd(\{g\ast (g_{1},s_{1})\ast \dots \ast (g_{q},s_{q})\}) + 
\delta_{Z \cap \op{div}s_{1} \cap \dots \cap \op{div}s_{q}},
\end{multline*}
que la forme diff\'erentielle
$\omega (\alpha\cdot\hat{c}_{1}(L_{1},\|.\|_{1}^{(n)})\dotsm
\hat{c}_{1}(L_{q},\|.\|_{q}^{(n)}))$ tend vers \\
$\omega (\alpha\cdot\hat{c}_{1}(L_{1},\|.\|_{1})\dotsm
\hat{c}_{1}(L_{q},\|.\|_{q}))$ au sens de la convergence faible
quand $n$ tend vers $+\infty$. Enfin, d'apr\`es (\ref{demailly}), on a~:
\[
\lim_{n \rightarrow +\infty} 
\omega(\alpha)\cdot c_{1}(L_{1},\|.\|_{1}^{(n)})\dotsm
c_{1}(L_{q},\|.\|_{q}^{(n)}) = 
\omega(\alpha)\cdot c_{1}(\ov{L}_{1})\dotsm c_{1}(\ov{L}_{q}), 
\] 
ce qui joint \`a la relation (\ref{relation_intermediaire}) permet de
conclure.}
\item{D'apr\`es le (2), il suffit de d\'emontrer l'\'egalit\'e recherch\'ee
pour $\ov{L}_{1},\dots,\ov{L}_{q}$ des fibr\'es en droites admissibles sur $X$. On
reprend ici les notations de la d\'emonstration du (5).
En utilisant la d\'efinition (\ref{definition_prod2}), on obtient la relation~:
\begin{equation}
\label{int_eq_5}
[(0,\varphi)]\cdot \hat{c}_{1}(\ov{L}_{1})\dotsm \hat{c}_{1}(\ov{L}_{q})
= 
[(0,\{\varphi\ast (g_{1},s_{1})\ast \dotsm \ast (g_{q},s_{q})\})].
\end{equation}
D'apr\`es la proposition (\ref{naturalite2}) et la th\'eorie d\'evelopp\'ee dans
\cite{13}, on a l'\'egalit\'e suivante, valable pour tout $n \in \M{N}$~:
\begin{align*}
\{\varphi\ast (g_{1}^{(n)},s_{1})\ast \dotsm \ast (g_{q}^{(n)},s_{q})\}
&= 
\varphi\ast (g_{1}^{(n)}\ast ( \dotsm \ast (g_{q-1}^{(n)} \ast
g_{q}^{(n)})\dotsm )) \\
&= (\dotsm ((\varphi \ast g_{1}^{(n)})\ast g_{2}^{(n)}) \ast \dotsm ) \ast
g_{q}^{(n)}.
\end{align*}
Par ailleurs, on sait d'apr\`es (\cite{13}, \S 2.2.9) que pour toute forme
diff\'erentielle $\varphi' \in A^{\ast}(X_{\R})$ et tout $1 \leqslant i
\leqslant q$~:
\[
\varphi' \ast g_{i}^{(n)} = g_{i}^{(n)} \ast \varphi' = \varphi'\cdot 
c_{1}(L_{i},\|.\|_{i}^{(n)}).
\]
On d\'eduit alors par r\'ecurrence de ce qui pr\'ec\`ede la relation~: 
\[
\{\varphi \ast (g_{1}^{(n)},s_{1})\ast \dotsm \ast (g_{q}^{(n)},s_{q})\}
= 
\varphi\cdot c_{1}(L_{1},\|.\|_{1}^{(n)}) \dotsm c_{1}(L_{q},\|.\|_{q}^{(n)})
\]
dans $\widetilde{D}^{p+q-1,p+q-1}(X_{\R})$.

En remarquant que d'une part, $\{\varphi \ast (g_{1}^{(n)},s_{1})\ast 
\dotsm \ast (g_{q}^{(n)},s_{q})\}$ tend vers $
\{\varphi\ast (g_{1},s_{1})\ast \dotsm \ast (g_{q},s_{q})\}$
au sens de la convergence faible des courants quand $n$ tend vers $+\infty$
d'apr\`es (\ref{approximation2}), et que d'autre part~:
\[
\lim_{n \rightarrow +\infty}\varphi\cdot c_{1}(L_{1},\|.\|_{1}^{(n)}) 
\dotsm c_{1}(L_{q},\|.\|_{q}^{(n)})
= 
\varphi\cdot c_{1}(\ov{L}_{1}) 
\dotsm c_{1}(\ov{L}_{q}),
\]
\'egalement au sens de la topologie faible des courants, on conclut de ce qui
pr\'ec\`ede que~:
\[
\{\varphi \ast (g_{1},s_{1}) \ast \dotsm \ast 
(g_{q},s_{q})\} = \varphi\cdot c_{1}(\ov{L}_{1})\dotsm c_{1}(\ov{L}_{q}) 
\] 
dans $\widetilde{D}^{p+q-1,p+q-1}(X_{\R})$, ce qui joint \`a la relation
(\ref{int_eq_5}) prouve le r\'esultat \'enonc\'e.}
\end{enumerate}
\medskip

\begin{defn}
\label{def_gch_gen}
Soit $p$ un entier positif. On appelle {\it groupe de Chow arithm\'etique
g\'en\'eralis\'e de codimension $p$\/} et on note $\CH^{p}(X)$ le
$\Z$ sous-module de $\widetilde{CH}^{p}(X)$ engendr\'e par les \'el\'ements de la
forme $\alpha\cdot \hat{c}_{1}(\ov{L}_{1})\dotsm
\hat{c}_{1}(\ov{L}_{q})$, o\`u $\alpha \in \widehat{CH}^{p-q}(X)$
et o\`u $\ov{L}_{1},\dots,\ov{L}_{q}$ sont des fibr\'es en droites
int\'egrables sur $X$. On convient de noter $\CH^{\ast}(X) = 
\bigoplus_{p \geqslant 0}\CH^{p}(X)$.
\end{defn}
\begin{rem}
Le groupe $\CH^{p}(X)$ contient $\widehat{CH}^{p}(X)$.
\end{rem}
\begin{rem}
Soit $\ov{L}$ un fibr\'e int\'egrable sur $X$; la premi\`ere classe de Chern
$\hat{c}_{1}(\ov{L})$ de $\ov{L}$ appartient \`a $\CH^{1}(X)$.
\end{rem}
\begin{prop}
L'application qui \`a tout fibr\'e en droites hermitien int\'egrable $\ov{L}$
associe sa classe de Chern $\hat{c}_{1}(\ov{L}) \in \CH^{1}(X)$, d\'efinit un
isomorphisme de groupes~:
\[
\hat{c}_{1}: \widehat{\op{Pic}}_{\,\op{int}}(X) \longrightarrow \CH^{1}(X).
\]
\end{prop}
\demo\ Il d\'ecoule des d\'efinitions que $\hat{c}_{1}(\cdot)$ est bien
d\'efinie, surjective, et pr\'eserve la structure de groupe. Il nous reste \`a 
montrer que
$\hat{c}_{1}(\cdot)$ est injective. Soit $\ov{L} = (L,\|.\|) \in 
\widehat{\op{Pic}}_{\,\op{int}}(X)$ tel que $\hat{c}_{1}(\ov{L}) = 0$. 
Comme $c_{1}(L) = \zeta(\hat{c}_{1}(\ov{L})) = 0$, le fibr\'e $L$ est isomorphe
au fibr\'e trivial; on note $s$ une section constante non nulle de $L$
au-dessus de $X$. La relation~:
\[
\hat{c}_{1}(\ov{L}) = [\op{\widehat{\op{div}}}s] = [(0,-\log \|s\|^{2})] = 0, 
\]
jointe au fait que le groupe $CH^{1,0}(X)$ est toujours trivial, entra\^\i nent
que $\|s\| = 1$ identiquement sur $X(\M{C})$, ce qui nous permet de conclure.
\medskip

Soient $\alpha \in A^{r-1,r-1}(X_{\R})$ et $\ov{L}_{1},\dots,\ov{L}_{q}$ des
fibr\'es en droites int\'egrables sur $X$, avec $q$ et $r$ deux entiers
positifs tels que $q + r = p$. D'apr\`es le th\'eor\`eme (\ref{th_general}), on a~: 
\[
[(0,\alpha c_{1}(\ov{L}_{1}) \dotsm c_{1}(\ov{L}_{q}))] = [(0,\alpha)]\cdot
\hat{c}_{1}(\ov{L}_{1})\dotsm \hat{c}_{1}(\ov{L}_{q}) \in
\CH^{p}(X). 
\]
On en d\'eduit que $a\left(\AA_{\op{g}}^{p-1,p-1}(X_{\R})\right) \subset 
\CH^{p}(X)$. On dispose donc du morphisme de groupes~:
\begin{alignat*}{3}
a: \AA_{\op{g}}^{p-1,p-1}&(X_{\R})& &\longrightarrow &
&\CH^{p}(X) \\
&\beta& &\longmapsto & & [(0,\beta)].
\end{alignat*}
D'apr\`es le th\'eor\`eme (\ref{th_general}), on obtient par restriction \`a $
\CH^{p}(X)$ des morphismes $\zeta$ et $\omega$, les morphismes de groupes
suivants~:
\begin{alignat*}{3}
\zeta :\; &\CH^{p}(X) & &\longrightarrow & &CH^{p}(X) \\
&[(Z,g)] & &\longmapsto & &Z, \\
\end{alignat*}
et 
\begin{alignat*}{3}
\qquad \qquad \quad \omega :\; &\CH^{p}(X) & &\longrightarrow &  &\AA_{\op{g}}^{p,p}(X_{\R})\\
&[(Z,g)] & & \longmapsto & &\omega(Z,g) = \dd g + \delta_{Z}.
\end{alignat*}
\bigskip

\begin{prop}
\label{decomposition_classe}
Soit $\alpha \in \CH^{p}(X)$ et choisissons $(Z,g)$ un repr\'esentant de
$\alpha$ dans $\CH^{p}(X)$. Il existe $g_{Z} \in D^{p-1,p-1}(X_{\R})$ un
courant de Green pour $Z$ et $\varphi \in C_{\log,0}^{p-1,p-1}(X_{\R})$ tels que
l'on ait~:
\[
g = g_{Z} + \varphi \quad \text{dans} \quad \widetilde{D}^{p-1,p-1}(X_{\R}).
\]
\end{prop}
\demo\ 
Par lin\'earit\'e, il suffit de d\'emontrer le r\'esultat pour $\alpha =
\beta\cdot \hat{c}_{1}(\ov{L}_{1})\dotsm \hat{c}_{1}(\ov{L}_{q})$, o\`u $\beta
\in \widehat{CH}^{p-q}(X)$ et $\ov{L}_{1} = (L_{1},\|.\|_{1}), \dots,
\ov{L}_{q} = (L_{q},\|.\|_{q})$ sont des fibr\'es admissibles sur $X$.

Soient $(Z',g')$ un repr\'esentant de la classe $\beta$ et $s_{1},\dots,s_{q}$
des sections rationnelles non identiquement nulles sur $X$ de
$L_{1},\dots,L_{q}$ respectivement telles que
$Z',\op{div}s_{1},\dots,\op{div}s_{q}$ s'intersectent proprement. 
Soient \'egalement $\|.\|_{1}^{(0)}, \dots, \|.\|_{q}^{(0)}$ des m\'etriques
positives $C^{\infty}$ sur $L_{1},\dots,L_{q}$ respectivement.

On pose $\tilde{\omega}' = \dd g' + \delta_{Z'}$. Pour tout $1 \leqslant i
\leqslant q$, on note \'egalement $\delta_{i} = \delta_{\op{div}s_{i}}$, $g_{i}
= - \log \|s_{i}\|_{i}^{2}$, $\omega_{i} = c_{1}(\ov{L}_{i})$, $g_{i}^{(0)} = -
\log (\|s_{i}\|_{i}^{(0)})^{2}$ et $\omega_{i}^{(0)} =
c_{1}(L_{i},\|.\|_{i}^{(0)})$.

En appliquant les d\'efinitions, il vient~:
\begin{align*}
\{g'\ast (g_{1},s_{1})\ast \dotsm \ast (g_{q},&s_{q})\}
-
\{g'\ast (g_{1},s_{1})\ast \dotsm \ast (g_{q-1},s_{q-1})\ast
(g_{q}^{(0)},s_{q})\}\\
 &= 
(g_{q} - g_{q}^{(0)})\, \tilde{\omega}'\omega_{1}\dotsm \omega_{q-1} \\
&= - \log \left(\frac{\|.\|_{q}}{\;\, \|.\|_{q}^{(0)}}\right)^{2}
\tilde{\omega}'\omega_{1}\dots \omega_{q-1} \in C_{\log,0}^{p-1,p-1}(X_{\R}).
\end{align*}
Or, d'apr\`es la proposition (\ref{naturalite2}), on a~:
\begin{multline*}
\{g'\ast(g_{1},s_{1})\ast \dotsm \ast
(g_{q-1},s_{q-1})\ast(g_{q}^{(0)},s_{q})\}\\
= 
\{(g'\ast g_{q}^{(0)})\ast(g_{1},s_{1})\ast \dotsm \ast (g_{q-1},s_{q-1})\},
\end{multline*}
dans $\widetilde{D}^{p-1,p-1}(X_{\R})$.
En it\'erant $q$ fois ce proc\'ed\'e, on trouve $\varphi \in
C_{\log,0}^{p-1,p-1}(X)$ tel que~:
\[
\{g'\ast (g_{1},s_{1})\ast \dotsm \ast (g_{q},s_{q})\}
=
g'\ast g_{q}^{(0)}\ast \dotsm \ast g_{1}^{(0)} + \varphi,
\]
dans $\widetilde{D}^{p-1,p-1}(X_{\R})$. Il suffit alors de remarquer que $g_{Z}
= g'\ast g^{(0)}_{q} \ast \dotsm \ast g_{1}^{(0)}$ est un courant de Green pour
$Z'\cdot \op{div}s_{1}\dotsm \op{div}s_{q}$ et la proposition est
d\'emontr\'ee.
\medskip

la proposition (\ref{decomposition_classe}) justifie la d\'efinition suivante~:
\begin{defn}
Soit $p$ un entier positif. On note $\widehat{CH}_{\op{gen}}^{p}(X)$ le $\Z$
sous-module de $\widetilde{CH}^{p}(X)$ engendr\'e par les \'el\'ements de
$\widehat{CH}^{p}(X)$ et ceux de la forme $a(\varphi)$ avec $\varphi \in
C_{\log,0}^{p-1,p-1}(X_{\R})$. On convient de noter
$\widehat{CH}_{\op{gen}}^{\ast}(X) = \bigoplus_{p \geqslant
0}\widehat{CH}_{\op{gen}}^{p}(X)$.
\end{defn}
\begin{rem}
D'apr\`es la proposition (\ref{decomposition_classe}) on a $\CH^{p}(X) \subset 
\widehat{CH}_{\op{gen}}^{p}(X)$.
\end{rem}
\begin{rem}
On dispose du morphisme de groupes~:
\begin{alignat*}{3}
\omega~: \;&\widehat{CH}_{\op{gen}}^{p}(X) & & \longrightarrow &
\;&C_{0}^{p,p}(X_{\R}) \\
&[(Z,g)]& &\longmapsto & & \dd g + \delta_{Z}.
\end{alignat*}
\end{rem}
\bigskip

\subsection{L'accouplement $\CH^{p}(X) \otimes
\CH^{q}(X) \rightarrow \CH^{p+q}(X)_{\M{Q}}$}~

Soit $\pi: X \rightarrow S = \op{Spec}\C{O}_{K}$ une vari\'et\'e arithm\'etique
de dimension de Krull $d+1$.

On d\'efinit dans cette section un accouplement $
\CH^{p}(X) \otimes \CH^{q}(X) 
\rightarrow \linebreak[4] \CH^{p+q}(X)_{\M{Q}}$
qui \'etend l'accouplement $\widehat{CH}^{p}(X) \otimes \widehat{CH}^{q}(X)
\rightarrow \widehat{CH}^{p+q}(X)_{\M{Q}}$ d\'efini par Gillet-Soul\'e (voir
\cite{13}, \S 4.2; on peut \'egalement consulter \cite{3}, \S 2.2).

Soient $x \in \CH^{p}(X)$ et $y \in \CH^{q}(X)$.
On peut \'ecrire $x$ et $y$ sous la forme~:
\begin{align*}
x &= \sum_{i = 1}^{n} \alpha_{i}\cdot \hat{c}_{1}(\ov{E}_{i,1})
\dotsm \hat{c}_{1}(\ov{E}_{i,r_{i}}) \\
\intertext{et}
y &= \sum_{j =1}^{m} \beta_{j}\cdot \hat{c}_{1}(\ov{F}_{j,1})
\dotsm \hat{c}_{1}(\ov{F}_{j,s_{j}}), 
\end{align*}
o\`u pour tout $1 \leqslant i \leqslant n$ (resp. tout $1 \leqslant j \leqslant
m$), $\alpha_{i}$ est un \'el\'ement de $\widehat{CH}^{p - r_{i}}(X)$ et
$\ov{E}_{i,1}, \dots \ov{E}_{i,r_{i}}$ sont des fibr\'es en droites
int\'egrables sur $X$ (resp. $\beta_{j}$ est un \'el\'ement de
$\widehat{CH}^{q-s_{j}}(X)$ et $\ov{F}_{j,1}, \dots, \ov{F}_{j,s_{j}}$ sont des
fibr\'es en droites int\'egrables sur $X$).

On d\'efinit alors le produit $(x\cdot y) \in \CH^{p+q}(X)_{\M{Q}}$ par la formule~: 
\[
(x\cdot y) = 
\sum_{i =1}^{n}\sum_{j=1}^{m} 
(\alpha_{i}\cdot \beta_{j})\cdot \hat{c}_{1}(\ov{E}_{i,1})\dotsm 
\hat{c}_{1}(\ov{E}_{i,r_{i}})\cdot \hat{c}_{1}(\ov{F}_{j,1})\dotsm
\hat{c}_{1}(\ov{F}_{j,s_{j}}), 
\]
o\`u pour tout $1 \leqslant i \leqslant n$ et $1 \leqslant j \leqslant m$, on a
not\'e $(\alpha_{i}\cdot \beta_{j}) \in \widehat{CH}^{p+q - r_{i} -
s_{j}}(X)_{\M{Q}}$ le produit de $\alpha_{i}$ et $\beta_{j}$ au sens
de Gillet-Soul\'e.

Bien entendu, il faut montrer que cette d\'efinition ne d\'epend pas des choix
effectu\'es pour repr\'esenter $x$ et $y$. C'est l'objet du th\'eor\`eme
suivant~:
\begin{thm}
\label{produit_bien_defini}
La d\'efinition ci-dessus ne d\'epend pas des choix effectu\'es. Elle d\'efinit
un accouplement~:
\begin{equation}
\label{accouplement_1}
\CH^{p}(X) \otimes \CH^{q}(X) \longrightarrow 
\CH^{p+q}(X)_{\M{Q}}, 
\end{equation}
qui prolonge celui d\'efini par Gillet-Soul\'e.

Cet accouplement munit $\CH^{\ast}(X)_{\M{Q}}$ d'une
structure d'anneau commutatif, associatif et unif\`ere.
\end{thm}
\begin{rem}
\label{rem_acc_1}
Si $p =1$ ou $q =1$, on dispose d'un accouplement~:
\[
\CH^{p}(X) \otimes \CH^{q}(X) \longrightarrow 
\CH^{p+q}(X)
\]
qui induit l'accouplement (\ref{accouplement_1}) \`a valeur dans
$\CH^{p+q}(X)_{\M{Q}}$.
\end{rem}
\begin{rem}
\label{rem_acc_2}
Si $X$ est lisse sur $\op{Spec}\C{O}_{K}$, on peut d\'efinir les produits
$(\alpha_{i}\cdot \beta_{j})$ dans
$\widehat{CH}^{p+q-r_{i}-s_{j}}(X)$. La construction pr\'ec\'edente donne donc
un accouplement~:
\[
\CH^{p}(X) \otimes \CH^{q}(X) \longrightarrow 
\CH^{p+q}(X)
\]
qui induit par produit tensoriel avec $\M{Q}$ l'accouplement
(\ref{accouplement_1}) \`a valeur dans \linebreak[4]
$\CH^{p+q}(X)_{\M{Q}}$.
\end{rem}
\demo\ On montre tout d'abord que le produit $(x\cdot y)$ introduit
plus haut est bien d\'efini.

Soit $x = \sum_{i=1}^{n}
\alpha_{i}\cdot \hat{c}_{1}(\ov{E}_{i,1})\dotsm
\hat{c}_{1}(\ov{E}_{i,r_{i}})$ un \'el\'ement de
$\CH^{p}(X)$. Il faut montrer que si $x =0$ dans $\CH^{p}
(X)$, alors pour tout $y \in \CH^{q}(X)$, on a $(x\cdot y) = 0$.

On montre dans ce qui suit un r\'esultat plus g\'en\'eral~: Si $x = a(\varphi)$
avec $\varphi \in C_{\log,0}^{p-1,p-1}(X_{\R})$, alors $x\cdot y = a(\varphi\,
\omega(y))$.

Par lin\'earit\'e, il suffit de d\'emontrer cet \'enonc\'e pour $y = 
\beta\cdot \hat{c}_{1}(\ov{F}_{1})\dotsm \hat{c}_{1}(\ov{F}_{s})$,
o\`u $\beta$ est un \'el\'ement de $\widehat{CH}^{q-s}(X)$ et $\ov{F}_{1}
= (F_{1}, \|.\|_{1}), \dots, \ov{F}_{s} = (F_{s},\|.\|_{s})$ sont des fibr\'es
en droites admissibles sur $X$. Remarquons \'egalement que la proposition
(\ref{produit_general2}) et le th\'eor\`eme (\ref{th_general}) nous permettent de supposer que pour tout
$1 \leqslant i \leqslant n$ les fibr\'es en droites $\ov{E}_{i,1} =
(E_{i,1},\|.\|_{i,1}), \dots, \ov{E}_{i,r_{i}} = (E_{i,r_{i}},
\|.\|_{i,r_{i}})$ sont admissibles sur $X$.

Soient
$\left(\|.\|_{1}^{(k)}\right)_{k \in \N}, \dots,
\left(\|.\|_{s}^{(k)}\right)_{k \in \N}$ des suites croissantes de m\'etriques
positives $C^{\infty}$ sur $F_{1},\dots,F_{s}$ convergeant vers
$\|.\|_{1},\dots,\|.\|_{s}$ respectivement, et 
pour tout $1 \leqslant i \leqslant n$, soient 
$\left(\|.\|^{(k)}_{i,1}\right)_{k \in \N}, \dots, 
\left(\|.\|^{(k)}_{i,r_{i}}\right)_{k \in \N}$ des suites croissantes de 
m\'etriques positives
$C^{\infty}$ sur $E_{i,1},\dots,E_{i,r_{i}}$ convergeant vers
$\|.\|_{i,1},\dots,\|.\|_{i,r_{i}}$ respectivement.

On choisit $(Z',g')$, $(Z_{1},g_{1}), \dots, (Z_{n},g_{n})$ des repr\'esentants
des classes $\beta, \alpha_{1}, \dots,$ $\alpha_{n}$ respectivement, tels que
$Z_{K}'$ soit d'intersection propre avec chacun des $(Z_{1})_{K},$ $\dots, 
(Z_{n})_{K}$, et tels que $g',g_{1},\dots,g_{n}$ soient des courants de type
logarithmique.

Pour tout $1 \leqslant i \leqslant n$, on note $Z'\cdot Z_{i}$ un
repr\'esentant dans $Z^{p+q-r_{i}}(|Z'|\cap |Z_{i}|)_{\M{Q}}$ du produit
$[Z']\cdot [Z_{i}] \in CH_{|Z'|\cap|Z_{i}|}^{p+q-r_{i}}(X)_{\M{Q}}$; et on
choisit $s_{i,1},\dots,s_{i,r_{i}}$ des sections rationnelles non identiquement
nulles sur $X$ de $E_{i,1},\dots,E_{i,r_{i}}$ respectivement 
et $s_{1},\dots,s_{s}$ des sections rationnelles non identiquement nulles sur
$X$ de $F_{1},\dots,F_{s}$ respectivement
telles que les
cycles $Z_{i},\op{div}s_{i,1},\dots,\op{div}s_{i,r_{i}},
\op{div}s_{1},\dots,\op{div}s_{s}$ ainsi que les cycles
$(Z'\cdot Z_{i}), \op{div}s_{i,1},\dots,\op{div}s_{i,r_{i}},
\op{div}s_{1},\dots,\op{div}s_{s}$ et $Z_{K}, Z_{K}',
(\op{div}s_{i,1})_{K}, \dots,\linebreak[4] (\op{div}s_{i,r_{i}})_{K},
(\op{div}s_{1})_{K},\dots,(\op{div}s_{s})_{K}$ s'intersectent
proprement. \\
On convient de noter $g_{i,1} = -\log \|s_{i,1}\|_{i,1}^{2}, \dots, g_{i,r_{i}}
= - \log\|s_{i,r_{i}}\|^{2}_{i,r_{i}}$ et 
$g_{1} = - \log\|s_{1}\|_{1}^{2},\dots,g_{s} = - \log \|s_{s}\|_{s}^{2}$, 
et pour tout $k \in \N$, $g_{i,1}^{(k)}
= - \log \left(\|s_{i,1}\|_{i,1}^{(k)}\right)^{2}, \dots\linebreak[4], 
g_{i,r_{i}}^{(k)}
= - \log \left(\|s_{i,r_{i}}\|_{i,r_{i}}^{(k)}\right)^{2}$
et $g_{1}^{(k)} = - \log \left(\|s_{1}\|_{1}^{(k)}\right)^{2}, \dots, 
g_{s}^{(k)} = - \log \left( \|s_{s}\|_{s}^{(k)}\right)^{2}$.

D'apr\`es la d\'efinition (\ref{definition_prod2}) et 
(\cite{13}, \S 4.2.1 et 4.2.2), on a~:
\begin{equation}
\label{eq_prod1}
x = \sum_{i = 1}^{n}
[(\, Z_{i}\cdot \odiv s_{i,1} \dotsm \odiv s_{i,r_{i}}, 
\{g_{i}\ast (g_{i,1},s_{i,1}) \ast \dotsm \ast (g_{i,r_{i}},s_{i,r_{i}})\}
\,)]
\end{equation}
et 
\begin{multline}
\label{eq_prod2}
x\cdot y = 
{\sum_{i=1}^{n}
[(\, (Z'\cdot Z_{i})\cdot \odiv s_{i,1}\dotsm \odiv s_{i,r_{i}}
\cdot \odiv s_{1} \dotsm \odiv s_{s}, }\\
\{(g'\ast g_{i})\ast (g_{i,1},s_{i,1})\ast \dotsm \ast
(g_{i,r_{i}},s_{i,r_{i}})\ast (g_{1},s_{1})\ast \dotsm \ast (g_{s},s_{s})
\}\, )].
\end{multline}
Puisque $\zeta(x) = \sum_{i=1}^{n}
[Z_{i}\cdot \odiv s_{i,1} \dotsm \odiv s_{i,r_{i}}] = 0$, on peut trouver une
$K_{1}$-chaine $f$ telle que~:
\begin{equation}
\label{eq_prod3}
\sum_{i=1}^{n}Z_{i}\cdot \odiv s_{i,1} \dotsm \odiv s_{i,r_{i}} = \odiv f.
\end{equation}
D'apr\`es le lemme de d\'eplacement pour les $K_{1}$-chaines (cf. \cite{13}, \S
4.2.6), on peut de plus supposer que $f$ rencontre $Z'
\cdot \odiv s_{1}\dotsm \odiv s_{s}$ presque proprement dans
$X_{K}$. Cela entra\^\i ne que d'une part~:
\begin{multline}
\label{eq_prod4}
\sum_{i=1}^{n}
(Z_{i}\cdot \odiv s_{i,1} \dotsm \odiv s_{i,r_{i}}, 
\{g_{i}\ast (g_{i,1},s_{i,1})\ast \dotsm \ast(g_{i,r_{i}},s_{i,r_{i}})\}
)\\
 = \op{\widehat{\op{div}}}f + (0,\gamma_{1} + \gamma_{2}), 
\end{multline}
o\`u l'on a not\'e~:
\[
\gamma_{1} = \sum_{i=1}^{n}
\{g_{i}\ast (g_{i,1},s_{i,1}) \ast \dotsm \ast(g_{i,r_{i}},s_{i,r_{i}})\},
\qquad 
\gamma_{2} = \log |f|^{2}; 
\]
et que d'autre part, d'apr\`es (\cite{13},\S 4.2.1 et 4.2.5), 
\begin{multline}
\label{eq_prod5}
\sum_{i=1}^{n}
(\,(Z'\cdot Z_{i})\cdot \odiv s_{i,1} \dotsm \odiv s_{i,r_{i}}
\cdot \odiv s_{1}\dotsm \odiv s_{s}, \\
\shoveright{\{(g'\ast g_{i})\ast (g_{i,1},s_{i,1})\ast \dotsm \ast
(g_{i,r_{i}},s_{i,r_{i}})
\ast (g_{1},s_{1})\ast \dotsm \ast (g_{s},s_{s})\})}\\
= \op{\widehat{\op{div}}}(f\cdot (Z'
\cdot \odiv s_{1}\dotsm \odiv s_{s})) + (0,\gamma')
+ (R,0), \quad
\end{multline}
o\`u l'on a not\'e~:
\begin{multline}
\gamma' = 
\sum_{i=1}^{n}\{
(g'\ast g_{i})\ast (g_{i,1},s_{i,1}) \ast \dotsm \ast
(g_{i,r_{i}},s_{i,r_{i}})\ast (g_{1},s_{1})\ast \dotsm \ast (g_{s},s_{s})\} \\
+ \log |f\cdot (Z'\cdot \odiv s_{1}\dotsm \odiv s_{s})|^{2}, 
\end{multline}
o\`u $R \in R_{\op{fin}}^{p+q}(X)$ et
o\`u $f\cdot (Z'\cdot \odiv s_{1}\dotsm \odiv s_{s})$ d\'esigne une $K_{1}$-chaine repr\'esentant l'intersection
de la $K_{1}$-chaine $f$ et du cycle $Z'\cdot \odiv s_{1}\dotsm \odiv s_{s}$ 
(bien que $f\cdot (Z'\cdot \odiv s_{1}\dotsm \odiv s_{s})$ ne soit pas
d\'efinie de mani\`ere univoque, 
les cycles $\odiv (f\cdot (Z'\cdot \odiv s_{1}\dotsm \odiv s_{s}))$ et 
$\op{\widehat{\op{div}}}(f\cdot (Z'\cdot \odiv s_{1}\dotsm \odiv s_{s}))$ le sont,
voir \cite{13}, \S 4.2.5). 

D'apr\`es (\ref{eq_prod3}) et (\cite{13}, \S 2.2.9 et 4.2.7), on a pour tout $k
\in \N$~:
\begin{multline*}
\sum_{i=1}^{n}(g'
\ast g_{1}^{(k)} \ast \dotsm \ast g_{s}^{(k)})
\ast (g_{i}\ast g_{i,1}^{(k)}\ast \dotsm \ast
g_{i,r_{i}}^{(k)}) + \log |f\cdot (Z'\cdot \odiv s_{1}\dotsm \odiv s_{s})|^{2}  \\
\shoveleft{= (g'\ast g_{1}^{(k)} \ast \dotsm \ast g_{s}^{(k)})\cdot \delta_{\odiv
f}} \\
 + 
\omega(\beta)
\, c_{1}\big(F_{1},\|.\|_{1}^{(k)}\big)
\dotsm c_{1}\big(F_{s},\|.\|_{s}^{(k)}\big)\cdot \left(
\sum_{i=1}^{n}g_{i}\ast 
g_{i,1}^{(k)}\ast \dotsm \ast g_{i,r_{i}}^{(k)}\right) \\
\shoveright{
+ \log |f\cdot (Z'\cdot \odiv s_{1}\dotsm \odiv s_{s})|^{2}} \\
\shoveleft{= \Big((g'\ast g_{1}^{(k)} \ast \dotsm \ast g_{s}^{(k)})
\cdot \delta_{\odiv f} }
- \omega(\beta)\, c_{1}\big(F_{1},\|.\|_{1}^{(k)}\big)
\dotsm c_{1}\big(F_{s},\|.\|_{s}^{(k)}\big)\log|f|^{2}  \\
\shoveright{+ \log |f\cdot
(Z'\cdot \odiv s_{1}\dotsm \odiv s_{s})|^{2}\Big) }\\
\shoveright{+ \omega(\beta)\, c_{1}\big(F_{1},\|.\|_{1}^{(k)}\big)
\dotsm c_{1}\big(F_{s},\|.\|_{s}^{(k)}\big)
\left(
\sum_{i=1}^{n}g_{i}\ast 
g_{i,1}^{(k)}\ast \dotsm \ast g_{i,r_{i}}^{(k)} + \log|f|^{2}\right)\;\,}
\\
\shoveleft{= 
\omega(\beta)\, c_{1}\big(F_{1},\|.\|_{1}^{(k)}\big)
\dotsm c_{1}\big(F_{s},\|.\|_{s}^{(k)}\big)
\left(
\sum_{i=1}^{n}g_{i}\ast 
g_{i,1}^{(k)}\ast \dotsm \ast g_{i,r_{i}}^{(k)} + \log|f|^{2}\right),}
\end{multline*}
dans $\widetilde{D}^{p+q-1,p+q-1}(X_{\R})$, ce qui, en faisant tendre $k$ vers
$+\infty$, entra\^\i ne que~:
\[
\gamma' = \gamma_{1} \,\omega(y) +\gamma_{2} \,\omega(y) , 
\]
dans $\widetilde{D}^{p+q-1,p+q-1}(X_{\R})$
(une telle expression a bien un sens car $\omega(y) \in 
C_{0}^{q-1,q-1}(X_{\R})$).
De l'\'egalit\'e $x = [(0,\gamma_{1} + \gamma_{2})] = [(0,\varphi)]$ obtenue 
en combinant (\ref{eq_prod1})
et (\ref{eq_prod4}), on d\'eduit qu'il existe une $K_{1}$-chaine $g$, que l'on
choisit d'intersection presque propre avec $Z'
\cdot \odiv s_{1}\dotsm \odiv s_{s}$ dans $X_{K}$, telle que $\odiv
g = 0$ et $\gamma_{1} + \gamma_{2} - \varphi = - \log |g|^{2}$ dans
$\widetilde{D}^{p-1,p-1}(X_{\R})$.

On tire de ce qui pr\'ec\`ede et de (\cite{13}, \S 4.2.5 et 4.2.7)
la relation~:
\begin{multline*}
- \log |g\cdot (Z'\cdot \op{div}s_{1}\dotsm \op{div}s_{s})|^{2}\\
= (\gamma_{1} + \gamma_{2} - \varphi)\,\omega(\beta)
\, c_{1}\big(F_{1},\|.\|_{1}^{(k)}\big) 
\dotsm 
c_{1}\big(F_{s},\|.\|_{s}^{(k)}\big),
\end{multline*}
valable pour tout $k\in \N$; 
ce qui montre, en faisant tendre $k$ vers $+\infty$, que~:
\[
\gamma_{1}\,\omega(y) + \gamma_{2}\,\omega(y) - \varphi\, \omega(y)
= - \log |g\cdot (Z'\cdot \op{div}s_{1}\dotsm \op{div}s_{s})|^{2},
\]
dans $\widetilde{D}^{p+q-1,p+q-1}(X_{\R})$. Comme de plus 
$\odiv(g\cdot (Z'\cdot \op{div}s_{1}\dotsm \op{div}s_{s}))_{K} = 
(\odiv g)_{K}\cdot (Z'\cdot \op{div}s_{1}\dotsm \op{div}s_{s})_{K} = 0$, 
on peut affirmer que $[(0,\gamma')]
= [(0,\gamma_{1}\,\omega(y) + \gamma_{2}\,\omega(y))] =
[(0,\varphi\,\omega(y))]$, 
ce qui combin\'e \`a (\ref{eq_prod2}) et 
(\ref{eq_prod5}) montre que $x\cdot y = a(\varphi\,\omega(y))$.

On a donc d\'emontr\'e que l'accouplement (\ref{accouplement_1}) est bien d\'efini.
Les autres propri\'et\'es \'enonc\'ees se d\'eduisent ais\'ement du th\'eor\`eme
(\ref{th_general}) et des propri\'et\'es analogues pour
$\widehat{CH}^{\ast}(X)$.
\medskip

Au cours de la d\'emonstration pr\'ec\'edente, on a de plus prouv\'e le
r\'esultat suivant~:
\begin{prop}
\label{produit_generalise_fin}
Soient 
$x \in \CH^{\ast}(X)$ et $\varphi \in C_{\log,0}^{\ast}(X_{\R})$ 
tel que $a(\varphi) \in
\CH^{\ast}(X)$, on a~:
\[
a(\varphi)\cdot x = a(\varphi\,\omega(x)).
\]
\end{prop}
\begin{rem}
Soient $x$ et $y$ deux \'el\'ements de $\widehat{CH}_{\op{gen}}^{\ast}(X)$,
et soient $\alpha$, $\alpha' \in \widehat{CH}^{\ast}(X)$ et $\varphi$,
$\varphi' \in C_{\log,0}^{\ast}(X_{\R})$ tels que $x = \alpha + a (\varphi)$ et $y =
\alpha' + a(\varphi')$. On d\'efinit le produit de $x$ et de $y$ par la
formule~:
\[
x\cdot y = \alpha\cdot \alpha' + 
a(\varphi\,\omega(\alpha') + \varphi'\,\omega(\alpha) + \varphi\,\dd \varphi')
\in \widehat{CH}_{\op{gen}}^{\ast}(X)_{\Q}.
\]
On peut montrer que ce produit est bien d\'efini et qu'il munit
$\widehat{CH}_{\op{gen}}^{\ast}(X)_{\Q}$ d'une structure d'anneau commutatif,
associatif et unif\`ere. D'apr\`es la proposition (\ref{produit_generalise_fin}), il co\"\i ncide sur
$\CH^{\ast}(X)_{\Q}$ avec le produit d\'efini au th\'eor\`eme
(\ref{produit_bien_defini}).
\end{rem}

\begin{thm}
Soit $\pi: X \rightarrow \op{Spec}\C{O}_{K}$ une vari\'et\'e arithm\'etique.
Les morphismes~:
\begin{align*}
\zeta&: \CH^{\ast}(X)_{\M{Q}} \longrightarrow
CH^{\ast}_{\M{Q}}(X) \\
\intertext{et}
\omega&: \CH^{\ast}(X)_{\M{Q}} \longrightarrow
\AA_{\op{g}}^{\ast}(X_{\R}), 
\end{align*}
d\'efinis \`a la section (\ref{sous_section_def}) sont des morphismes
d'anneaux. De plus, le produit d\'efini \`a la remarque (\ref{rem_acc_1}) et (dans le
cas o\`u $\pi: X \rightarrow \op{Spec}\C{O}_{K}$ est lisse) celui d\'efini \`a
la remarque (\ref{rem_acc_2}) sont compatibles avec les morphismes $\zeta$ et $\omega$.
\end{thm}
\demo\ C'est une simple cons\'equence du th\'eor\`eme (\ref{th_general}),
alin\'ea (5) et des propri\'et\'es analogues pour $\widehat{CH}^{\ast}(X)$.
\medskip

\begin{rem}
Soient $\varphi \in \AA_{\op{g}}^{\ast}(X_{\R})$ et $x \in
\CH^{\ast}(X)$. On d\'eduit ais\'ement de la d\'efinition de
$\AA_{\op{g}}^{\ast}(X_{\R})$ et des alin\'eas (5) et (6) du th\'eor\`eme
(\ref{th_general}) la formule utile suivante~:
\[
a(\varphi)\cdot x = a(\varphi\, \omega(x)).
\]
\end{rem}
\medskip

\subsection{Degr\'e arithm\'etique et hauteurs}~

\subsubsection{Degr\'e arithm\'etique}~
\label{sous_section_degre}

Soit $\pi: X \rightarrow \op{Spec}\C{O}_{K} = S$ une vari\'et\'e
arithm\'etique (projective) de dimension relative $d$. 
On rappelle (voir par exemple \cite{3}, \S 2.1.3) que l'on
dispose des deux morphismes~:
\[
\op{deg}_{K} : \widehat{CH}^{0}(S) = CH^{0}(S) \longrightarrow \M{Z}
\]
et 
\[
\widehat{\op{deg}} : \widehat{CH}^{1}(S) \longrightarrow \R, 
\]
qui induisent par composition avec $\pi_{\ast}: 
\widehat{CH}^{\ast}(X) \rightarrow \widehat{CH}^{\ast -d}(S)$ les morphismes~:
\[
\op{deg}_{K} : \widehat{CH}^{d}(X) \longrightarrow \M{Z} 
\qquad \text{(degr\'e g\'eom\'etrique)} 
\]
et
\[
\widehat{\op{deg}} : \widehat{CH}^{d+1}(X) \longrightarrow \R
\qquad \text{(degr\'e arithm\'etique)}.
\]
Soit $i \in \{ 0,1\}$. Pour toute classe $\alpha \in \CH^{d+i}(X)$, choisissons
$(Z,g)$ un repr\'esentant de $\alpha$; on d\'eduit de (\cite{13}, \S 3.6)
que la classe $[(\pi_{\ast}(Z),\pi_{\ast}(g))]$ ne d\'epend que de $\alpha$.
On dispose donc d'un morphisme, encore not\'e $\pi_{\ast}$, 
qui est d\'efini comme suit~:
\begin{alignat*}{3}
\pi_{\ast}: \;&\CH^{d+i}(X)& &\longrightarrow &&\;\widehat{CH}^{i}(S) \\
&\;\;[(Z,g)]&&\longmapsto &&[(\pi_{\ast}(Z),\pi_{\ast}(g))],
\end{alignat*}
et qui prolonge \`a $\CH^{d+i}(X)$ le morphisme image directe usuel
$\pi_{\ast}: \widehat{CH}^{d+i}(X) \rightarrow \widehat{CH}^{i}(S)$.

En composant $\op{deg}_{K}$ et $\widehat{\op{deg}}$ avec 
$\pi_{\ast}$, on obtient deux nouveaux morphismes~:
\[
\op{deg}_{K} : \CH^{d}(X) \longrightarrow \M{Z} 
\]
et
\[
\widehat{\op{deg}} : \CH^{d+1}(X) \longrightarrow \R,
\]
qui prolongent ceux d\'efinis pr\'ec\'edemment.

\subsubsection{Hauteurs}~

Soient $Z \in Z_{q}(X)$ un cycle de dimension $q$ et $\ov{L}_{1} = (L_{1},
\|.\|_{1}), \dots, \ov{L}_{q} = (L_{q},\|.\|_{q})$ des fibr\'es en droites
admissibles sur $X$.
Choisissons $s_{1},\dots,s_{q}$ des sections rationnelles non identiquement
nulles sur $X$ de $L_{1},\dots,L_{q}$ respectivement, telles que les cycles
$Z,\odiv s_{1}, \dots, \odiv s_{q}$ soient d'intersection propre.
Pour tout $1 \leqslant i \leqslant q$, on note $g_{i} = - \log
\|s_{i}\|_{i}^{2}$ et $\omega_{i} = c_{1}(\ov{L}_{i})$.
On d\'efinit \`a partir de ces donn\'ees un \'el\'ement de $D^{p,p}(X_{\R})$ que
l'on note $\{(g_{1},s_{1})\ast \dotsm \ast (g_{q},s_{q})|\delta_{Z}\}$ par la
formule suivante~:
\begin{multline}
\label{int_fin_eq1}
\{(g_{1},s_{1})\ast \dotsm \ast (g_{q},s_{q})|\delta_{Z}\} = 
g_{1}\cdot \delta_{\odiv s_{2} \cap \dotsm \cap \odiv s_{q} \cap Z} 
+ \omega_{1}g_{2}\cdot \delta_{\odiv s_{3}\cap \dotsm \cap \odiv s_{q} \cap
Z}\\
+ \dots + \omega_{1}\dotsm \omega_{i-1}g_{i}\cdot \delta_{\odiv s_{i+1}\cap
\dotsm \cap \odiv s_{q}\cap Z}
+ \dots + \omega_{1} \dotsm \omega_{q-1}g_{q}\cdot \delta_{Z}.
\end{multline}
On v\'erifie que chacun des termes de la relation (\ref{int_fin_eq1})
est bien d\'efini. 
En effet $\omega_{1}\dotsm \omega_{i-1}\cdot \delta_{\odiv s_{i+1} \cap \dotsm
\cap \odiv s_{q} \cap Z}$ est bien d\'efini comme le produit de
$\omega_{1}\dotsm \omega_{i-1} \in B^{i-1,i-1}(X_{\R})$ et de $\delta_{\odiv
s_{i+1} \cap \dotsm \cap \odiv s_{q} \cap Z} \in B^{d-i+1,d-i+1}_{\odiv s_{i+1}
\cap \dotsm \cap \odiv s_{q} \cap Z}(X_{\R})$ d'apr\`es les propositions
(\ref{produit_uniforme}) et (\ref{produit_uniforme2}).
Par ailleurs le produit~:
\[
\gamma_{i} = g_{i}\,\omega_{1}\dotsm \omega_{i-1}\cdot \delta_{\odiv s_{i+1}\cap
\dotsm \cap \odiv s_{q}\cap Z},
\]
a bien un sens d'apr\`es le th\'eor\`eme (\ref{demailly}); et si l'on pose $g_{i}^{(t)}
= \op{max} (g_{i}, t)$ pour tout $t \in \R$, la limite~:
\[
\lim_{t \rightarrow - \infty}
g_{i}^{(t)}\omega_{1}\dotsm \omega_{i-1}\cdot \delta_{\odiv s_{i+1} \cap \dotsm
\cap \odiv s_{q} \cap Z} = 
g_{i} \,\omega_{1}\dotsm \omega_{i-1}\cdot \delta_{\odiv s_{i+1} \cap \dotsm
\cap \odiv s_{q} \cap Z},
\]
montre que le courant $\gamma_{i}$ ne d\'epend pas des choix effectu\'es pour
le d\'efinir. 

Ceci \'etant \'etabli, on pose~:
\[
h_{\ov{L}_{1}, \dots,\ov{L}_{q}}(Z) := \widehat{\op{deg}}
(\pi_{\ast}(Z\cdot \odiv s_{1} \dotsm \odiv s_{q}), 
\pi_{\ast}(\{(g_{1},s_{1})\ast \dotsm \ast (g_{q},s_{q})|\delta_{Z}\})) \in \R.
\]
\begin{prop_defn}
\label{int_fin_prop1}
Le nombre r\'eel $h_{\ov{L}_{1},\dots,\ov{L}_{q}}(Z)$ d\'efini ci-dessus
ne d\'epend pas du choix des sections $s_{1},\dots,s_{q}$; on
l'appelle {\rm hauteur de $Z$ relativement \`a $\ov{L}_{1},\dots,\ov{L}_{q}$\/}.
\end{prop_defn}
\begin{prop}
\label{int_fin_prop2}
Soient $Z \in Z_{q}(X)$ et $\ov{L}_{0}, \ov{L}_{1} = (L_{1}, \|.\|_{1}), \dots,
\ov{L}_{q} = (L_{q},\|.\|_{q})$ des fibr\'es admissibles sur $X$. Les
assertions suivantes sont v\'erifi\'ees~:
\begin{enumerate}
\item{La hauteur $h_{\ov{L}_{1}, \dots, \ov{L}_{q}}(Z)$ ne d\'epend que du
cycle $Z$
et des
classes d'isomorphie isom\'etrique des fibr\'es admissibles
$\ov{L}_{1},\dots,\ov{L}_{q}$; elle ne d\'epend pas de l'ordre des fibr\'es
$\ov{L}_{1},\dots,\ov{L}_{q}$.}
\item{Soient $\Big(\|.\|_{1}^{(k)}\Big)_{k \in \N}, \dots, 
\Big(\|.\|_{q}^{(k)}\Big)_{k \in \N}$ des suites croissantes de m\'etriques
positives convergeant vers $\|.\|_{1}, \dots,\|.\|_{q}$ 
sur $L_{1},\dots,L_{q}$ respectivement, on a~:
\[
\lim_{k \rightarrow + \infty}
h_{\big(L_{1},\|.\|_{1}^{(k)}\big), \dots, 
\big(L_{q},\|.\|_{q}^{(k)}\big)}
(Z) = h_{\ov{L}_{1},\dots,\ov{L}_{q}}(Z).
\]
}
\item{Si les m\'etriques $\|.\|_{1}, \dots, \|.\|_{q}$ sont $C^{\infty}$, alors
on a~:
\[
h_{\ov{L}_{1}, \dots, \ov{L}_{q}}(Z) = 
\widehat{\op{deg}}
(\hat{c}_{1}(\ov{L}_{1})\dotsm \hat{c}_{1}(\ov{L}_{q})|Z), 
\]
o\`u $(\cdot |\cdot)$ d\'esigne l'accouplement d\'efini dans (\cite{3}, \S
2.3).}
\item{Pour tout $1 \leqslant i \leqslant q$, on a~:
\[
h_{\ov{L}_{1}, \dots,\ov{L}_{i-1},\ov{L}_{i}\otimes \ov{L}_{0}, 
\ov{L}_{i+1}, \dots,\ov{L}_{q}}(Z) = 
h_{\ov{L}_{1},\dots,\ov{L}_{q}}(Z) + 
h_{\ov{L}_{1}, \dots, \ov{L}_{i-1}, \ov{L}_{0}, \ov{L}_{i+1}, \dots,
\ov{L}_{q}}(Z).
\]
}
\end{enumerate}
\end{prop}
{\bf D\'emonstration des propositions (\ref{int_fin_prop1}) et
(\ref{int_fin_prop2}).} La proposition (\ref{int_fin_prop1}) et les assertions
(1) et (4) de la proposition (\ref{int_fin_prop2}) se d\'eduisent
imm\'ediatement des assertions (2) et (3) de la proposition (\ref{int_fin_prop2})
et des propri\'et\'es analogues dans le cas classique (cf. \cite{3}, \S 2.3).
L'assertion (2) \'etant une cons\'equence du th\'eor\`eme (\ref{demailly}), il suffit
de prouver l'assertion (3).

En remarquant que $g_{1} \ast (g_{2}\ast (\dotsm \ast(g_{q})))$ est un courant
de Green pour le cycle $\odiv s_{1} \dotsm \odiv s_{q}$, on d\'eduit de
(\cite{13}, \S 1.2.4 et 1.3.5) qu'il existe $g$ un courant de Green de type
logarithmique pour $\odiv s_{1} \dotsm \odiv s_{q}$ tel que~: 
\begin{equation}
\label{int_fin_eq2}
g = g_{1} \ast (g_{2}\ast (\dotsm \ast(g_{q}))) + \partial u + \ov{\partial}v.
\end{equation}
Pour tout $\epsi \in \R^{+\ast}$, soit $\delta_{Z}^{(\epsi)}$ une
r\'egularisation de $\delta_{Z}$ comme \`a la proposition (\ref{construction_GS_lissage}). En
multipliant les deux membres de l'\'egalit\'e (\ref{int_fin_eq2}) par
$\delta_{Z}^{(\epsi)}$, il vient~:
\begin{equation}
\label{int_fin_eq3}
g \wedge \delta_{Z}^{(\epsi)} = 
g_{1} \ast (g_{2}\ast (\dotsm \ast(g_{q}))) \wedge \delta_{Z}^{(\epsi)} + 
\partial (u \delta_{Z}^{(\epsi)}) + \ov{\partial}(v \delta_{Z}^{(\epsi)}).
\end{equation}
Comme d'une part~:
\begin{multline*}
g_{1} \ast (g_{2}\ast (\dotsm \ast(g_{q}))) 
\\
= g_{1} \cdot \delta_{\odiv s_{2}
\cap \dotsm \cap \odiv s_{q}} + \omega_{1} g_{2} \cdot \delta_{\odiv s_{3} \cap
\dotsm \cap \odiv s_{q}} + \dots + \omega_{1} \dotsm \omega_{q-1}g_{q}, 
\end{multline*}
et que d'autre part, d'apr\`es (\cite{13}, \S 2.2.12), on a les limites~:
\[
\lim_{\epsi \rightarrow 0} g \wedge \delta_{Z}^{(\epsi)} = g \wedge \delta_{Z},
\]
et pour tout $1 \leqslant i \leqslant q$, 
\[
\lim_{\epsi \rightarrow 0}
(\omega_{1}\dotsm \omega_{i-1})g_{i}\cdot \delta_{\odiv s_{i+1} \cap \dotsm
\cap \odiv s_{q}}\wedge \delta_{Z}^{(\epsi)} = 
(\omega_{1}\dotsm \omega_{i-1})g_{i}\cdot \delta_{\odiv s_{i+1} \cap \dotsm
\cap \odiv s_{q}}\wedge \delta_{Z}, 
\]
au sens de la convergence faible des courants, on tire de (\ref{int_fin_eq3})
l'\'egalit\'e~:
\begin{multline*}
g \wedge \delta_{Z} = \\ 
g_{1} \cdot 
\delta_{\odiv s_{2} \cap \dotsm \cap \odiv s_{q} \cap Z}
+ \omega_{1} g_{2} \cdot 
\delta_{\odiv s_{3}\cap \dotsm \cap \odiv s_{q} \cap Z} + \dots + 
\omega_{1}\dotsm \omega_{q-1}g_{q}\delta_{Z} \\
= \{(g_{1},s_{1})\ast \dotsm \ast (g_{q},s_{q})|\delta_{Z}\},
\end{multline*}
dans $\widetilde{D}^{d,d}(X_{\R})$, ce qui termine la d\'emonstration.
\medskip

Plus g\'en\'eralement, soient $Z \in Z_{q}(X)$ et $\ov{L}_{1},
\dots,\ov{L}_{q}$ des fibr\'es int\'egrables sur $X$. Choisissons $\ov{E}_{1},
\dots,\ov{E}_{q}$ et $\ov{F}_{1}, \dots,\ov{F}_{q}$ des fibr\'es en droites
admissibles sur $X$ tels que pour tout $ 1 \leqslant i \leqslant q$, on ait~:
$\ov{L}_{i} = \ov{E}_{i}\otimes (\ov{F}_{i})^{-1}$.
\begin{prop_defn}
\label{int_fin_prop3}
On appelle {\rm hauteur de $Z$ relativement \`a $\ov{L}_{1}, \dots,
\ov{L}_{q}$} et on note $h_{\ov{L}_{1}, \dots,\ov{L}_{q}}(Z)$ le nombre
r\'eel d\'efini par la formule~:
\[
h_{\ov{L}_{1},\dots,\ov{L}_{q}}(Z)
= \sum
_{\substack{S_{1}, S_{2} \\ S_{1} \cup S_{2} = \{1,\dots,q\}}}
(-1)^{\# S_{2}}h_{ \{\ov{E}_{i}\}_{i \in S_{1}}, \{\ov{F}_{j}\}_{j \in
S_{2}}}(Z).
\]
La hauteur $h_{\ov{L}_{1},\dots,\ov{L}_{q}}(Z)$ ainsi d\'efinie ne d\'epend que
de $Z$ et des classes d'isomorphie isom\'etrique des fibr\'es
$\ov{L}_{1},\dots,\ov{L}_{q}$. Si les fibr\'es $\ov{L}_{1},\dots,\ov{L}_{q}$
sont admissibles, elle co\"\i ncide avec la hauteur d\'efinie \`a la
proposition (\ref{int_fin_prop1}).
\end{prop_defn}
\demo\ On suit {\it mutatis mutandis\/} la d\'emonstration de la proposition
(\ref{produit_general2}) en utilisant la proposition (\ref{int_fin_prop2}).
\medskip

Lorsque $\ov{L} := \ov{L}_{1} = \dots = \ov{L}_{q}$, on convient de noter
$h_{\ov{L}}(Z)$ le nombre $h_{\ov{L}_{1},\dots,\ov{L}_{q}}(Z)$ que l'on appelle
alors 
{\it hauteur de $Z$ relativement \`a $\ov{L}$}.

\begin{rem}
Si $q=0$, la hauteur $h_{\ov{L}}(Z)$ n'est autre que $h(Z) = \widehat{\op{deg}}[(Z,0)]$. 
Dans ce cas, la fonction $h: Z_{0}(X) \rightarrow \R$ est d\'efinie par $h(P) = \log \#
k(P)$.
\end{rem}
\begin{rem}
Ces d\'efinitions ont \'et\'e introduites pour la premi\`ere fois 
dans cette g\'en\'eralit\'e par Zhang (cf. \cite{21}, th. 1.4) sous une forme
diff\'erente.
\end{rem}

Le th\'eor\`eme suivant rassemble diverses propri\'et\'es de la hauteur
introduite \`a la proposition (\ref{int_fin_prop3})~:
\begin{thm}
\label{gdthm}
Soient $Z \in Z_{q}(X)$ et $\ov{L}_{1},\dots,
\ov{L}_{q-1},\ov{L}_{q} = (L_{q},\|.\|_{q})$
des fibr\'es en droites int\'egrables sur $X$. On a les propri\'et\'es
suivantes~:
\begin{enumerate}
\item{La hauteur $h_{\ov{L}_{1},\dots,\ov{L}_{q}}(Z)$ ne d\'epend pas de
l'ordre de $\ov{L}_{1},\dots,\ov{L}_{q}$.}
\item{L'application qui \`a $Z \in Z_{q}(X)$ et $\ov{L}_{1},\dots,\ov{L}_{q}$
des fibr\'es en droites int\'egrables sur $X$ associe $h_{\ov{L}_{1},\dots,
\ov{L}_{q}}(Z) \in \R$ d\'efinit une forme multilin\'eaire~:
\[
Z_{q}(X) \times \widehat{\op{Pic}}_{\op{int}}(X) \times \dotsm \times 
\widehat{\op{Pic}}_{\op{int}}(X) \longrightarrow \R.
\]
}
\item{Si les m\'etriques des fibr\'es $\ov{L}_{1},\dots,\ov{L}_{q}$ sont
$C^{\infty}$, alors on a~:
\[
h_{\ov{L}_{1},\dots,\ov{L}_{q}}(Z) = \widehat{\op{deg}}(\hat{c}_{1}(\ov{L}_{1})
\dotsm \hat{c}_{1}(\ov{L}_{q})|Z),
\]
o\`u $(\cdot | \cdot )$ d\'esigne l'accouplement d\'efini dans (\cite{3}, \S
2.3).}
\item{Si $Z$ est le diviseur d'une fonction rationnelle sur un sous-sch\'ema
int\`egre contenu dans une fibre ferm\'ee de $X$, alors
$h_{\ov{L}_{1},\dots,\ov{L}_{q}}(Z) = 0$.}
\item{Pour tout morphisme $f: X \rightarrow X'$ de vari\'et\'es
arithm\'etiques, on a~:
\[
h_{f^{\ast}(\ov{L}_{1}), \dots , f^{\ast}(\ov{L}_{q})}(Z) =
h_{\ov{L}_{1},\dots,\ov{L}_{q}}(f_{\ast}(Z)).
\]
}
\item{Soit $s_{q}$ une section rationnelle de $L_{q}$ au-dessus de $Z$ qui
n'est identiquement nulle sur aucune des composantes irr\'eductibles de $Z$. On
a la relation~:
\[
h_{\ov{L}_{1},\dots,\ov{L}_{q-1}}(Z\cdot \odiv s_{q}) = 
h_{\ov{L}_{1},\dots,\ov{L}_{q}}(Z) + 
\int_{X(\M{C})}\log \|s_{q}\|_{q}\,
c_{1}(\ov{L}_{1})\dotsm c_{1}(\ov{L}_{q-1})\cdot \delta_{Z}.
\]
}
\item{On suppose que $Z = X$, et donc que $q = d+1$. On a~:
\[
h_{\ov{L}_{1},\dots,\ov{L}_{d+1}}(X) = \widehat{\op{deg}}(
\hat{c}_{1}(\ov{L}_{1})\dotsm \hat{c}_{1}(\ov{L}_{d+1})), 
\]
o\`u $\widehat{\op{deg}}$ d\'esigne le degr\'e arithm\'etique sur
$\CH^{d+1}(X)$ introduit au \S \ref{sous_section_degre}.}
\end{enumerate}
\end{thm}
\demo\ Par (multi)lin\'earit\'e on se ram\`ene au cas o\`u
$\ov{L}_{1},\dots,\ov{L}_{q}$ sont des fibr\'es en droites admissibles.
Les assertions (1) \`a (3) sont alors une cons\'equence imm\'ediate de la
proposition (\ref{int_fin_prop2}), et les assertions (4) \`a (6) se d\'eduisent
directement de (\ref{int_fin_prop2}) alin\'ea (2) et des assertions analogues
dans le cas classique (cf. \cite{3}, prop. 2.3.1 et 3.2.1). 
Enfin l'assertion (7) est une cons\'equence des propositions
(\ref{approximation2}) et
(\ref{naturalite2}), et des d\'efinitions.
\medskip

La proposition suivante \'etend \`a notre cadre un \'enonc\'e de Faltings (cf.
\cite{32}, prop. 2.6) g\'en\'eralis\'e dans (\cite{3}, \S 3.2.3)~:

\begin{prop}{\rm \bf (Positivit\'e).} Soient $\ov{L}_{1},\dots,\ov{L}_{q}$ des
fibr\'es en droites {\em admissibles} sur $X$ tels que pour tout $1 \leqslant i
\leqslant q$, il existe une puissance tensorielle positive $\ov{L}_{i}^{n_{i}}$
de $\ov{L}_{i}$ engendr\'ee par ses sections globales de norme sup inf\'erieure
ou \'egale \`a $1$. Pour tout cycle {\em effectif} $Z \in Z_{q}(X)$, on a~: 
\begin{equation}
\label{int_eq_6}
h_{\ov{L}_{1},\dots,\ov{L}_{q}}(Z) \geqslant 0.
\end{equation}
\label{positivite_arithmetique}
\end{prop}
\demo\ On raisonne par r\'ecurrence sur la dimension de $Z$.

Si $\op{dim}Z = 0$, le calcul se fait aux places finies et il n'y a pas de
changement avec la situation classique (voir par exemple \cite{3}, prop.
3.2.4).

On se place d\'esormais dans le cas o\`u $\op{dim} Z > 0$. 

On peut supposer que $Z$ est irr\'eductible et choisir une section $s_{q}$ de
$L_{q}^{n_{q}}$ de norme sup inf\'erieure ou \'egale \`a $1$ qui ne s'annule
pas identiquement sur $Z$. On d\'eduit alors de l'alin\'ea (6) du th\'eor\`eme 
(\ref{gdthm})
appliqu\'e aux fibr\'es $\ov{L}_{1},\dots,\ov{L}_{q-1},\ov{L}_{q}^{n_{q}}$
que~:
\begin{alignat*}{2}
h_{\ov{L}_{1},\dots,\ov{L}_{q-1}}(Z\cdot \op{div}s_{q}) &\leqslant
n_{q}\,h_{\ov{L}_{1},\dots,\ov{L}_{q}}(Z)& \qquad &\text{si} \quad q \geqslant 2
\\ 
\text{et}
\qquad \qquad \qquad \qquad \qquad 
h(Z \cdot \op{div}s_{1}) &\leqslant n_{1} \, h_{\ov{L}_{1}}(Z)& \qquad
&\text{si} \quad q  = 1.\qquad \qquad \qquad \qquad ~
\end{alignat*}
Comme $Z\cdot \op{div}s_{q}$ est effectif de dimension $\op{dim}Z - 1$
et que $h$ prend des valeurs positives ou nulles sur les cycles effectifs dans
$Z_{0}(X)$, cela implique l'in\'egalit\'e (\ref{int_eq_6}) par r\'ecurrence sur
$\op{dim}Z$.
\medskip

\begin{expl}
\label{positivite_torique}
Soient $\P$ une vari\'et\'e torique projective lisse et
$\ov{E}_{1},\dots,\ov{E}_{q}$ des fibr\'es en droites sur $\P$ engendr\'es par
leurs sections globales et munis de leur m\'etrique canonique. D'apr\`es
(\ref{exemple_adm1}) les fibr\'es $\ov{E}_{1},\dots,\ov{E}_{q}$ sont admissibles et 
d'apr\`es la proposition (\ref{sections_globales}) et l'\'egalit\'e
(\ref{image_inverse})
ils
sont engendr\'es par leurs sections globales de norme sup inf\'erieure ou
\'egale \`a $1$. On d\'eduit de la proposition (\ref{positivite_arithmetique}) que pour tout cycle
effectif $Z \in Z_{q}(\P)$, on a~: 
\[
h_{\ov{E}_{1},\dots,\ov{E}_{q}}(Z) \geqslant 0.
\]
\end{expl}
\bigskip

\section{Courants de Chern canoniques sur les vari\'et\'es toriques}~

Dans cette partie et les suivantes, nous revenons \`a l'\'etude des
vari\'et\'es toriques projectives lisses sur
$\SP$.

Soit $\D$ un \'eventail r\'egulier complet admettant une fonction 
support strictement concave et
soient $\ov{L}_{1},\dots,\ov{L}_{q}$ des fibr\'es en droites sur
$\P$ munis \`a l'infini de leur m\'etrique canonique. On donne dans cette
partie une expression explicite du courant~:
\[
c_{1}(\ov{L}_{1})\dotsm c_{1}(\ov{L}_{q}).
\]
On en d\'eduit un premier r\'esultat de trivialit\'e~: Si l'on note~: 
\[
\ov{A}_{\op{f}}^{\ast}(\P_{\R}) = \op{Ker}\{d: \ov{A}_{\op{g}}^{\ast}(\P_{\R})
\longrightarrow \ov{A}_{\op{g}}^{\ast}(\P_{\R})\}
\]
et 
\[
[\cdot]: \ov{A}_{\op{f}}^{\ast}(\P_{\R}) \longrightarrow H^{2\ast}(\P(\M{C}),\R)
\]
la surjection canonique induite par l'application classe en
cohomologie des courants, alors on peut
construire de mani\`ere canonique une section (d'anneaux) du morphisme
d'anneaux $[\cdot]$.
\bigskip

\subsection{Pr\'eliminaires}~

Soit $\op{Arg}: \M{C}^{\ast} \rightarrow \R / 2\pi\M{Z}$ la fonction argument
d\'efinie pour tout $z = \rho e^{i\alpha} \in \M{C}^{\ast}$ avec $\rho \in
\R^{+\ast}$ et $\alpha \in \R$ par~:
\[
\op{Arg}(z) = [\alpha] \in \R / 2\pi \M{Z}.
\]
Pour tout $m \in M$, on note $\arg(\chi^{m})$ la fonction multiforme sur
$\TT$ d\'efinie comme l'argument du caract\`ere $\chi^{m}$. La
diff\'erentielle $d \arg(\chi^{m})$ a bien un sens sur $\TT$ et d\'efinit un
\'el\'ement de $A^{1,0}(\TT)\oplus A^{0,1}(\TT)$. On dispose donc d'un
morphisme injectif de $\M{Z}$-module~:
\begin{alignat*}{3}
\Theta:\; &M & &\longrightarrow& &A^{1,0}(\TT)\oplus A^{0,1}(\TT) \\
&m & &\longmapsto & &\frac{d \arg(\chi^{m})}{2\pi}.
\end{alignat*}
Ce morphisme s'\'etend en un morphisme d'anneaux (encore not\'e $\Theta$)~:
\[
\Theta: {\textstyle\bigwedge\nolimits^{\ast}_{\M{Z}}}M \longrightarrow
\oplus_{p,q}A^{p,q}(\TT),
\]
o\`u ${\textstyle\bigwedge\nolimits^{\ast}_{\M{Z}}}M$
d\'esigne l'alg\`ebre ext\'erieure sur $\M{Z}$ de $M$ et o\`u
la structure d'anneau sur $\oplus_{p,q}A^{p,q}(\TT)$ est celle induite par
le produit ext\'erieur. Le morphisme d'anneaux $\Theta$ est injectif.
\begin{rem}
Soit $f_{1},\dots,f_{d}$ une base de $N$ et $f_{1}^{\ast},\dots,f_{d}^{\ast}$
la base duale. Pour tout $m \in M$, on a~:
\[
\Theta(m) = \sum_{i =1}^{d}f_{i}(m) \frac{d \arg(\chi^{f_{i}^{\ast}})}{2\pi}.
\]
\end{rem}
\begin{rem}
\label{remarque_theta}
Soit $m \in M$. On a~:
\[
d^{c} \log |\chi^{m}|^{2} = \Theta(m) \qquad \text{sur $\TT$}.
\]
\end{rem}
Pour tout $\t$ \'el\'ement de $\D$, le $\M{Z}$-module
$\bigwedge^{\op{max}}_{\M{Z}}(\t^{\perp}\cap M)$ est un $\M{Z}$-module libre de rang $1$.
On choisit un g\'en\'erateur de $\bigwedge^{\op{max}}_{\M{Z}}(\t^{\perp}\cap M)$ que l'on
notera dans toute la suite $\C{M}_{\t}$. Le choix de $\C{M}_{\t}$ d\'efinit une
orientation sur la vari\'et\'e r\'eelle $C_{\t}^{\op{int}}$ de la mani\`ere
suivante~:

On note $h = \op{dim}\t$. Soient $u_{1},\dots,u_{h}$ un syst\`eme de g\'en\'erateurs du
semi-groupe $(\t \cap N)$ que l'on compl\`ete en une base $u_{1},\dots,u_{d}$
de $N$,
et $v_{h+1},\dots,v_{d}$ une base du $\M{Z}$-module
$(\t^{\perp}\cap M)$ telle que $v_{h+1}\wedge \dots \wedge v_{d} = \C{M}_{\t}$.
On note $u_{1}^{\ast},\dots,u_{d}^{\ast}$ la base duale de $u_{1},\dots,u_{d}$.
On note enfin $B(0,1) \subset \M{C}$ la boule ouverte de centre $0$ et de rayon $1$
et $S_{1} \subset \M{C}$ le cercle unit\'e. On dispose alors du
diff\'eomorphisme~:
\begin{alignat*}{3}
\eta_{\t}:\; &C_{\t}^{\op{int}}& & \longrightarrow & &B(0,1)^{h}\times
S_{1}^{d-h}\\
&x & &\longmapsto &
&(\chi^{u_{1}^{\ast}}(x),\dots,\chi^{u_{h}^{\ast}}(x),\chi^{v_{h+1}}(x),\dots,\chi^{v_{d}}(x))
\end{alignat*}
On notera $C_{\t}^{\op{int}, +}$ la vari\'et\'e r\'eelle $C_{\t}^{\op{int}}$
munie de l'orientation produit des orientations naturelles sur $B(0,1)^{h}$ et
$S_{1}^{d-h}$. Cette orientation ne d\'epend que du choix de $\C{M}_{\t}$. On
peut remarquer que pour tout $\si \in \D_{\op{max}}$, l'orientation de
$C_{\si}^{\op{int},+}$ co\"\i ncide avec celle induite par sa structure
complexe. Par ailleurs, on a pour tout $\t \in \D$ l'\'egalit\'e
$|\chi^{v_{h+1}}| = \dots = |\chi^{v_{d}}| = 1$ sur $C_{\t}^{\op{int}}$, ce qui
entra\^\i ne que la forme $\Theta(\C{M}_{\t})$ est de classe $C^{\infty}$ sur
un voisinage de $C_{\t}^{\op{int}}$; on dispose donc du courant r\'eel~:
\[
\omega \longmapsto \int_{C_{\t}^{\op{int},+}}\Theta(\C{M}_{\t})\wedge \omega.
\]
\bigskip

\subsection{Calcul de $c_{1}(\ov{L})$}~

Soit $\ov{L}$ un fibr\'e en droites sur $\P$ muni de sa m\'etrique
canonique. On donne, dans cette section, une expression du courant
$c_{1}(\ov{L})$.

On rappelle tout d'abord un r\'esultat bien connu (voir par exemple \cite{7}, \S 3.3).
\begin{prop}
\label{Green}
{\rm (Formule de Green).} Soit $X$ une vari\'et\'e complexe de
dimension $d$ et $\Omega \subset X$ un ouvert relativement compact tel que
$\ov{\Omega}$ soit une sous-vari\'et\'e r\'eelle \`a coins de $X$. Soient $f$ et
$g$ des formes diff\'erentielles de classe $C^{2}$ sur un voisinage de
$\ov{\Omega}$ et de bidegr\'es $(p,p)$ et $(q,q)$ avec $p+q = d-1$. On a~:
\[
\int_{\Omega}\left(f \wedge \dd g - \dd f \wedge g \right) 
= \int_{\partial\Omega}\left(f \wedge d^{c} g - d^{c}f \wedge g\right).
\]
\end{prop}
\medskip

Soit $D$ un diviseur (horizontal) $T$-invariant sur $\P$ tel que $\ov{L} =
\ov{\C{O}(D)}_{\infty}$.
On rappelle que pour tout $\si \in \D$, on note $m_{D,\si}$ l'\'el\'ement de $M$
donnant la restriction \`a $\si$ de la fonction support $\Psi_{D}$ de $D$. On a
alors le th\'eor\`eme suivant~:
\begin{thm}
\label{calcul_c1}
Pour toute forme test $\omega \in A^{d-1,d-1}(\P(\M{C}))$, on a~:
\[
\int_{\P(\M{C})}c_{1}(\ov{L}) \wedge \omega = - \sum_{\si \in
\D_{\op{max}}}\int_{\partial C_{\si}^{\op{int},+}}\Theta(m_{D,\si}) \wedge
\omega.
\]
\end{thm}
\demo\ On montre tout d'abord que les termes intervenant dans le second membre
sont bien d\'efinis. Pour tout $\si \in \D_{\op{max}}$, la forme
$\Theta(m_{D,\si})\wedge \omega$ est $L^{1}$ sur $\partial C_{\si}^{\op{int}}$
(cela provient du fait que la forme $d \theta = d \op{Arg}(z)$ est localement
$L^{1}$ sur $\M{C}$).
La diff\'erence $\partial C_{\si}^{\op{int}} \backslash (\partial C_{\si}^{\op{int}}
\cap \TT)$ \'etant de codimension r\'eelle sup\'erieure ou \'egale \`a $2$, elle est
n\'egligeable au sens de la th\'eorie de la mesure, ce qui entra\^\i ne que
l'int\'egrale $\int_{\partial C_{\si}^{\op{int},+}}\Theta(m_{D,\si})\wedge
\omega$ est bien d\'efinie.

Soit maintenant $\boldsymbol{1}$ la section (rationnelle) canonique de
$\C{O}(D)$. On a, d'apr\`es la formule de Poincar\'e-Lelong
g\'en\'eralis\'ee (\ref{PL_generalisee}) l'\'egalit\'e des courants~:
\[
c_{1}(\ov{L}) = c_{1}(\ov{\C{O}(D)}_{\infty}) = \delta_{D} + 
(\dd (-\log \|\boldsymbol{1}\|^{2}_{D,\infty})).
\]
On en d\'eduit que~:
\begin{equation}
\label{equa_1}
\int_{\P(\M{C})}c_{1}(\ov{L})\wedge \omega = 
\int_{\P(\M{C})} \delta_{D} \wedge \omega + \int_{\P(\M{C})}(- \log
\|\boldsymbol{1}\|^{2}_{D,\infty}) \wedge \dd \omega.
\end{equation}
Comme de plus, d'apr\`es (\ref{BT_construction}), on a pour tout $\si \in
\D_{\op{max}}$ et $x \in C_{\si}$ l'\'egalit\'e~:
\[
- \log \|\boldsymbol{1}(x)\|^{2}_{D,\infty} = \log | \chi^{m_{D,\si}}(x)|^{2}, 
\]
on tire de (\ref{equa_1}) et de (\ref{dissection}) la relation~:
\begin{equation}
\label{equa_2}
\int_{\P(\M{C})}c_{1}(\ov{L}) \wedge \omega = \int_{\P(\M{C})}\delta_{D}\wedge \omega + 
\sum_{\si \in
\D_{\op{max}}}\int_{C_{\si}^{\op{int},+}}\log|\chi^{m_{D,\si}}(x)|^{2}\wedge
\dd \omega.
\end{equation}
On fixe provisoirement $\si \in \D_{\op{max}}$.
Soit $f_{1},\dots,f_{d}$ une famille g\'en\'eratrice de $\si$ (et donc une
$\M{Z}$-base de $N$ d'apr\`es (\ref{lissite})); on note
$f_{1}^{\ast},\dots,f_{d}^{\ast}$ la base duale de $M$. On a $m_{D,\si} =
\sum_{i =1}^{d}f_{i}(m_{D,\si})f_{i}^{\ast}$, ce dont on tire~:
\[
\log | \chi^{m_{D,\si}}(x)|^{2} = \sum_{i =1}^{d}f_{i}(m_{D,\si})\log |
\chi^{f_{i}^{\ast}}(x)|^{2}.
\]
Par lin\'earit\'e, on ram\`ene ainsi le calcul de l'int\'egrale
$\int_{C_{\si}^{\op{int},+}}\log |\chi^{m_{D,\si}}(x)|^{2} \wedge \dd \omega$
\`a celui des int\'egrales $\int_{C_{\si}^{\op{int},+}}\log |
\chi^{f_{i}^{\ast}}(x)|^{2}\wedge \dd \omega$. On ne perd rien en
g\'en\'eralit\'e en posant $i = 1$. Dans la carte affine $\varphi:
U_{\si}(\M{C}) \rightarrow \M{C}^{d}$ donn\'ee par $\varphi(x) =
(\chi^{f_{1}^{\ast}}(x),\dots,\chi^{f_{d}^{\ast}}(x))$, les ensembles $C_{\si}$
et $C_{\si}^{\op{int}}$ sont d\'efinis par les conditions~:
\[
C_{\si} = \{x \in \M{C}^{d}: \quad |x_{1}| \leqslant 1, \dots, |x_{d}|
\leqslant 1\}, 
\]
et 
\[
C_{\si}^{\op{int}} = \{x \in \M{C}^{d}: \quad |x_{1}| < 1, \dots, |x_{d}| <
1\}.
\]
On remarque que~:
\[
\partial C_{\si} = \partial C_{\si}^{\op{int}} = \{x \in C_{\si}: \quad \exists
i \in \{1,\dots,d\}: \quad |x_{i}| = 1\}.
\]
Pour tout $\epsi \in \R^{+\ast}$, on pose~:
\begin{align*}
&C_{\si}^{\epsi} = \{ x \in C_{\si}^{\op{int}}: \quad \epsi < |x_{1}| < 1\}, \\
&D_{\si}^{\epsi} = \{x \in \partial C_{\si}: \quad \epsi < |x_{1}| \leqslant 1\}
\\
\intertext{et}
&E_{\si}^{\epsi} = \{x \in C_{\si}: \quad |x_{1}| = \epsi\}.
\end{align*}
\bigskip
\begin{center}
\begin{picture}(0,0)%
\includegraphics{figure5.pstex}%
\end{picture}%
\setlength{\unitlength}{0.00083300in}%
\begingroup\makeatletter\ifx\SetFigFont\undefined
\def\x#1#2#3#4#5#6#7\relax{\def\x{#1#2#3#4#5#6}}%
\expandafter\x\fmtname xxxxxx\relax \def\y{splain}%
\ifx\x\y   
\gdef\SetFigFont#1#2#3{%
  \ifnum #1<17\tiny\else \ifnum #1<20\small\else
  \ifnum #1<24\normalsize\else \ifnum #1<29\large\else
  \ifnum #1<34\Large\else \ifnum #1<41\LARGE\else
     \huge\fi\fi\fi\fi\fi\fi
  \csname #3\endcsname}%
\else
\gdef\SetFigFont#1#2#3{\begingroup
  \count@#1\relax \ifnum 25<\count@\count@25\fi
  \def\x{\endgroup\@setsize\SetFigFont{#2pt}}%
  \expandafter\x
    \csname \romannumeral\the\count@ pt\expandafter\endcsname
    \csname @\romannumeral\the\count@ pt\endcsname
  \csname #3\endcsname}%
\fi
\fi\endgroup
\begin{picture}(3687,3726)(3226,-4675)
\put(4276,-2311){\makebox(0,0)[lb]{\smash{\SetFigFont{12}{14.4}{rm}$C_{\si}^{\epsi}$}}}
\put(6076,-4636){\makebox(0,0)[lb]{\smash{\SetFigFont{12}{14.4}{rm}$E_{\si}^{\epsi}$}}}
\put(3226,-4636){\makebox(0,0)[lb]{\smash{\SetFigFont{12}{14.4}{rm}$D_{\si}^{\epsi}$}}}
\end{picture}

\end{center}
\bigskip
On a~:
\begin{multline}
\label{equaa_1}
\int_{C_{\si}^{\op{int},+}}\log | \chi^{f_{1}^{\ast}}(x)|^{2} \wedge \dd \omega
\\
= \int_{C_{\si}^{\epsi}}\log | \chi^{f_{1}^{\ast}}(x)|^{2} \wedge \dd \omega
+ \int_{C_{\si}^{\op{int}}\backslash C_{\si}^{\epsi}}
\log | \chi^{f_{1}^{\ast}}(x)|^{2} \wedge \dd \omega.
\end{multline}
Comme $\log | \chi^{f_{1}^{\ast}}(x)|^{2} = \log |x_{1}|^{2}$ est localement
$L^{1}$ sur $U_{\si}(\M{C})$, on a~:
\begin{equation}
\label{equaa_2}
\lim_{\epsi \rightarrow 0}\int_{C_{\si}^{\op{int}}\backslash C_{\si}^{\epsi}}
\log | \chi^{f_{1}^{\ast}}(x)|^{2} \wedge \dd \omega = 0.
\end{equation}
L'application $x \mapsto \log | \chi^{f_{1}^{\ast}}(x)|^{2}$ est $C^{\infty}$
sur un voisinage de $C_{\si}^{\epsi}$; il vient donc, d'apr\`es la formule de
Green (\ref{Green})~:
\begin{multline}
\label{equaa_3}
\int_{C_{\si}^{\epsi}} \log |\chi^{f_{1}^{\ast}}(x)|^{2}\wedge \dd \omega = 
\int_{C_{\si}^{\epsi}}\dd \log |\chi^{f_{1}^{\ast}}(x)|^{2}\wedge \omega \\
+
\int_{\partial C_{\si}^{\epsi}}\left(
\log | \chi^{f_{1}^{\ast}}(x)|^{2}\wedge d^{c}\omega - d^{c}\log |
\chi^{f_{1}^{\ast}}(x)|^{2}\wedge \omega\right).
\end{multline}
On a sur $C_{\si}^{\epsi}$ l'\'egalit\'e $\dd \log |\chi^{f_{1}^{\ast}}(x)|^{2}
= \dd \log |x_{1}|^{2} = 0$. De plus, d'apr\`es la remarque
(\ref{remarque_theta}), on a~:
\[
d^{c}\log |\chi^{f_{1}^{\ast}}(x)|^{2} = d^{c} \log |x_{1}|^{2} =
\Theta(f_{1}^{\ast}).
\]
On tire des d\'efinitions l'\'egalit\'e des courants~:
\[
\int_{\partial C_{\si}^{\epsi}} = \int_{D_{\si}^{\epsi}}\; + \;
\int_{E_{\si}^{\epsi}}
\]
Comme $\log |\chi^{f_{1}^{\ast}}(x)|^{2} = \log |x_{1}|^{2}$ est localement
$L^{1}$ sur $\partial C_{\si}$, il vient~:
\[
\lim_{\epsi \rightarrow 0}\int_{D_{\si}^{\epsi}}\log
|\chi^{f_{1}^{\ast}}(x)|^{2} \wedge d^{c}\omega = \int_{\partial
C_{\si}^{\op{int},+}}\log |\chi^{f_{1}^{\ast}}(x)|^{2} \wedge d^{c}\omega.
\]
De plus $\log |\chi^{f_{1}^{\ast}}(x)|^{2} = 2\log \varepsilon$ sur
$E_{\si}^{\varepsilon}$, et pour toute forme $\eta \in A^{2d-1}(\P(\M{C}))$ on a
$\int_{E_{\si}^{\varepsilon}}\eta = {O}(\varepsilon)$ quand $\varepsilon$
tend vers $0$; on en d\'eduit que~:
\[
\lim_{\epsi \rightarrow 0}\int_{E_{\si}^{\epsi}} \log 
|\chi^{f_{1}^{\ast}}(x)|^{2} \wedge d^{c}\omega = 0.
\]
Enfin, on a les limites~:
\[
\lim_{\epsi \rightarrow 0}\int_{D_{\si}^{\epsi}}d^{c} \log 
|\chi^{f_{1}^{\ast}}(x)|^{2} \wedge \omega = \lim_{\epsi \rightarrow 0}
\int_{D_{\si}^{\epsi}}\Theta(f_{1}^{\ast}) \wedge \omega = \int_{\partial
C_{\si}^{\op{int},+}}\Theta(f_{1}^{\ast}) \wedge \omega,
\]
et
\[
\lim_{\epsi \rightarrow 0} \int_{E_{\si}^{\epsi}} d^{c} 
\log |\chi^{f_{1}^{\ast}}(x)|^{2} \wedge \omega = 
\lim_{\epsi \rightarrow 0} \int_{E_{\si}^{\epsi}} 
\frac{d\arg (x_{1})}{2\pi}\wedge \omega =
\int_{C_{\si}^{\op{int},+}}\delta_{H_{1}}\wedge \omega, 
\]
o\`u $H_{1}$ d\'esigne le diviseur d'\'equation $\chi^{f_{1}^{\ast}}(x) = 0$.

En revenant au cas g\'en\'eral ($m_{D,\si}$ quelconque), on d\'eduit de
(\ref{equaa_1}), (\ref{equaa_2}), (\ref{equaa_3}) et des calculs de limites ci-dessus que~:
\begin{multline*}
\int_{C_{\si}^{\op{int},+}}\log | \chi^{m_{D,\si}}(x)|^{2} \wedge \dd \omega \\
= \int_{\partial C_{\si}^{\op{int},+}} \log | \chi^{m_{D,\si}}(x)|^{2}
\wedge d^{c}\omega - \int _{\partial C_{\si}^{\op{int},+}}
\Theta(m_{D,\si}) \wedge \omega - 
\int_{C_{\si}^{\op{int},+}}\delta_{D} \wedge \omega.
\end{multline*}
D'apr\`es (\ref{equa_2}), il vient~:
\begin{multline*}
\int_{\P(\M{C})} c_{1}(\ov{L}) \wedge \omega = 
\sum_{\si \in \D_{\op{max}}}\int_{\partial C_{\si}^{\op{int},+}}
\log | \chi^{m_{D,\si}}(x)|^{2} \wedge d^{c}\omega \\
- \sum_{\si \in \D_{\op{max}}}\int_{\partial C_{\si}^{\op{int},+}}
\Theta(m_{D,\si}) \wedge \omega
- \sum_{\si \in \D_{\op{max}}}\int_{C_{\si}^{\op{int},+}} \delta_{D}\wedge
\omega
+ \int_{\P(\M{C})}\delta_{D} \wedge \omega.
\end{multline*}
Soient $\si_{1}$ et $\si_{2}$ deux \'el\'ements de $\D_{\op{max}}$ et $\t \in
\D (d-1)$ une face commune de $\si_{1}$ et $\si_{2}$ (i.e. telle que $\t <
\si_{1}$ et $\t < \si_{2}$). Du fait de la continuit\'e de la fonction support
$\Psi_{D}$, on a~:
\begin{equation}
\label{eq_chern_1}
\log | \chi^{m_{D,\si_{1}}}(x)|^{2} = \log | \chi^{m_{D,\si_{2}}}(x)|^{2}, \qquad
(\forall x \in C_{\t}).
\end{equation}
Puisque $\D$ est complet, on a par ailleurs~:
\begin{multline*}
\sum_{\si \in \D_{\op{max}}} \int_{\partial C_{\si}^{\op{int},+}}\log
|\chi^{m_{D,\si}}(x)|^{2}\wedge d^{c} \omega = 
\sum_{\substack{\si \in \D_{\op{max}} \\ \t \in \D(d-1) \\ \t < \si}}
\int_{C_{\t}^{\op{int},+}}
\varepsilon_{\t}(\si)
\log |\chi^{m_{D,\si}}(x)|^{2}\wedge d^{c} \omega = 
\\
\sum_{\substack{\{\si_{1},\si_{2}\} \subset \D_{\op{max}} \\
\t = \si_{1} \cap \si_{2} \in \D(d-1)}}
\varepsilon_{\t}(\si_{1})
\left(
\int_{C_{\t}^{\op{int},+}}\log|\chi^{m_{D,\si_{1}}}(x)|^{2} \wedge d^{c} \omega
- 
\int_{C_{\t}^{\op{int},+}}\log|\chi^{m_{D,\si_{2}}}(x)|^{2} \wedge d^{c} \omega
\right)
\end{multline*}
o\`u pour tout $\si \in \D_{\op{max}}$ et tout $\t \in \D(d-1)$ tels que $\t <
\si$ on a pos\'e
$\varepsilon_{\t}(\si) = 1$ si les orientations de $\partial
C_{\si}^{\op{int},+}$ et de $C_{\t}^{\op{int},+}$ sont compatibles et
$\varepsilon_{\t}(\si) = -1$ sinon.

On d\'eduit de cela et de (\ref{eq_chern_1}) que~:
\[
\sum_{\si \in \D_{\op{max}}} \int_{\partial C_{\si}^{\op{int},+}}\log
|\chi^{m_{D,\si}}(x)|^{2}\wedge d^{c} \omega = 0.
\]
Enfin, comme $\op{codim}_{\R}(D\cap \partial C_{\si}) \geqslant 3$, pour tout
$\si \in \D_{\op{max}}$, on a~:
\[
\sum_{\si \in \D_{\op{max}}}\int_{C_{\si}^{\op{int},+}}
\delta_{D}\wedge \omega = \int_{\P(\M{C})}\delta_{D}\wedge \omega, 
\]
et le th\'eor\`eme est d\'emontr\'e.
\medskip

Soient $\si_{1},\dots,\si_{s}$ les \'el\'ements de $\D_{\op{max}}$. Comme $\D$
est complet, on peut associer \`a tout \'el\'ement $\t$ de $\D(d-1)$ deux
entiers $i_{\t} < j_{\t}$ dans $\{1,\dots,s\}$ tels que $\t$ soit la face
commune de $\si_{i_{\t}}$ et $\si_{j_{\t}}$. 
On munit $C_{\t}^{\op{int}}$ de l'orientation compatible avec celle de
$\partial C_{\si_{i_{\t}}}^{\op{int}}$ et on note $C_{\t}^{\op{int},++}$ la
vari\'et\'e r\'eelle $C_{\t}^{\op{int}}$ munie de cette orientation. On peut
reformuler le th\'eor\`eme (\ref{calcul_c1}) de la fa\c con suivante~:
\begin{thm}
\label{calcul2_c1}
Pour toute forme test $\omega \in A^{d-1,d-1}(\P(\M{C}))$, on a~:
\[
\int_{\P(\M{C})} c_{1}(\ov{L}) \wedge \omega = \sum_{\t \in
\D(d-1)}\int_{C_{\t}^{\op{int},++}}\Theta(m_{D,\si_{j_{\t}}} -
m_{D,\si_{i_{\t}}}) \wedge \omega.
\]
\end{thm}
\bigskip

\subsection{Calcul de $c_{1}(\ov{L}_{1})\dotsm c_{1}(\ov{L}_{q})$}~

Soient $\ov{L}_{1},\dots,\ov{L}_{q}$ des fibr\'es en droites
au-dessus de $\P$ et munis sur $\P(\M{C})$ de leur m\'etrique canonique. Le
th\'eor\`eme suivant donne une expression du produit $c_{1}(\ov{L}_{1})\dotsm 
c_{1}(\ov{L}_{q})$.
\begin{thm}
\label{calcul_prod}
Soient $\ov{L}_{1},\dots,\ov{L}_{q}$ des fibr\'es en droites
au-dessus de $\P$ munis de leur m\'etrique canonique. 
\begin{enumerate}
\item{Il existe une famille d'entiers
$(a_{\t}(\ov{L}_{1},\dots,\ov{L}_{q}))_{\t \in \D(d-q)}$ telle que, pour toute
forme test $\omega \in A^{d-q,d-q}(\P(\M{C}))$, on ait~:
\[
\int_{\P(\M{C})}c_{1}(\ov{L}_{1})\wedge \dots \wedge c_{1}(\ov{L}_{q}) \wedge \omega
= \sum_{\t \in \D(d-q)}a_{\t}(\ov{L}_{1},\dots, \ov{L}_{q})
\int_{C_{\t}^{\op{int},+}}\Theta(\C{M}_{\t})\wedge \omega.
\]
Les entiers $a_{\t}(\ov{L}_{1},\dots, \ov{L}_{q})$ sont d\'efinis de mani\`ere
unique par cette \'egalit\'e.
}
\item{Les entiers $a_{\t}(\ov{L}_{1},\dots, \ov{L}_{q})$ v\'erifient les
relations suivantes~: Pour tout $\si \in \D(d-q-1)$, on a~:
\[
\sum_{\substack{\t > \si \\ \t \in \D(d-q)}}
\epsi_{\si}(\t)a_{\t}(\ov{L}_{1},\dots, \ov{L}_{q})\C{M}_{\t} = 0;
\]
o\`u $\epsi_{\si}(\t) = 1$ si les orientations de $\partial
C_{\t}^{\op{int},+}$ et de $C_{\si}^{\op{int},+}$ sont compatibles, et $\epsi_{\si}(\t) =
-1$ sinon.
}
\item{
On suppose que $d < q$ et que $\ov{L}$ est un fibr\'e en droites
au-dessus de $\P$ muni de sa m\'etrique canonique. On note $D$ un diviseur
horizontal $T$-invariant sur $\P$ tel que $\ov{L} = \ov{\C{O}(D)}_{\infty}$. On
a~:
\[
a_{\si}(\ov{L},\ov{L}_{1},\dots,L_{q})\C{M}_{\si} = 
- \sum_{\substack{\t > \si \\ \t \in \D(d-q)}}
\epsi_{\si}(\t) a_{\t}(\ov{L}_{1},\dots,\ov{L}_{q}) \; m_{D,\t} \wedge \C{M}_{\t}.
\]
}
\end{enumerate}
\end{thm}
\demo\ On remarque tout de suite que l'unicit\'e des entiers
$a_{\t}(\ov{L}_{1},\dots,\ov{L}_{q})$ est une cons\'equence directe de la
proposition (\ref{dissection}) et du fait que le support du courant r\'eel
$\omega \mapsto \int_{C_{\t}^{\op{int},+}}\Theta(\C{M}_{\t})\wedge \omega$ est
$C_{\t}$.

D'apr\`es le th\'eor\`eme (\ref{calcul2_c1}), les assertions (1) et (2) sont
vraies pour $q=1$. On va montrer que si (1) et (2) sont vraies au rang $q$, alors
(1) est vraie au rang $q+1$,
(3) est vraie au rang $q$, et enfin (2)
est vraie au rang $q+1$.

D'apr\`es le (1), on peut trouver une famille d'entiers
$(a_{\t}(\ov{L}_{1},\dots,\ov{L}_{q}))_{\t \in \D(d-q)}$ telle que~:
\[
\int_{\P(\M{C})}c_{1}(\ov{L}_{1})\wedge \dots \wedge c_{1}(\ov{L}_{q})\wedge \omega =
\sum_{\t \in
\D(d-q)}a_{\t}(\ov{L}_{1},\dots,\ov{L}_{q})\int_{C_{\t}^{\op{int},+}}
\Theta(\C{M}_{\t})\wedge
\omega.
\]
Comme d'apr\`es (\ref{BT_construction}), on a pour tout $\si \in D$~:
\[
- \log \|\boldsymbol{1}(x)\|_{D,\infty}^{2} = \log | \chi^{m_{D,\si}}(x)|^{2},
\qquad (\forall x \in C_{\si}).
\]
On d\'eduit de la formule de Poincar\'e-Lelong g\'en\'eralis\'ee
(\ref{PL_generalisee}) et de la formule de Green
l'\'egalit\'e~:
\begin{multline}
\label{eq_chern_2}
\int_{\P(\M{C})}c_{1}(\ov{L})\wedge c_{1}(\ov{L}_{1}) \wedge \dots \wedge
c_{1}(\ov{L}_{q})\wedge \omega \\
\quad = \int_{\P(\M{C})}(\dd (- \log \|\boldsymbol{1}(x)\|_{D,\infty}^{2}) + \delta_{D})
\wedge c_{1}(\ov{L}_{1}) \wedge \dots \wedge c_{1}(\ov{L}_{q}) \wedge \omega 
\\
= \int_{\P(\M{C})}(-\log \|\boldsymbol{1}(x)\|_{D,\infty}^{2})\,
c_{1}(\ov{L}_{1}) \wedge \dots \wedge c_{1}(\ov{L}_{q}) \wedge \dd \omega 
\qquad \quad \\
\qquad \qquad + 
\int_{\P(\M{C})}
\delta_{D}\wedge c_{1}(\ov{L}_{1}) \wedge \dots \wedge c_{1}(\ov{L}_{q})\wedge \omega
\\
= \sum_{\t \in \D(d-q)} a_{\t}(\ov{L}_{1},\dots,\ov{L}_{q})
\int_{C_{\t}^{\op{int},+}}
\log | \chi^{m_{D,\t}}(x)|^{2}\; \Theta(\C{M}_{\t}) \wedge \dd \omega \\
+
\sum_{\t \in \D(d-q)}a_{\t}(\ov{L}_{1},\dots,\ov{L}_{q})
\int_{C_{\t}^{\op{int},+}}\delta_{D\cap C_{\t}^{\op{int}}}\wedge 
\Theta(\C{M}_{\t}) \wedge \omega.
\end{multline}
(Dans le dernier terme, on a utilis\'e le fait que $D$ et $C_{\t}^{\op{int}}$
s'intersectent de mani\`ere transverse).

Nous allons calculer les int\'egrales du type~:
\[
I_{\t}(m_{D,\t}) = 
\int_{C_{\t}^{\op{int},+}}
\log | \chi^{m_{D,\t}}(x)|^{2}\; \Theta(\C{M}_{\t}) \wedge \dd \omega.
\]
Fixons momentan\'ement $\t \in \D(d-q)$ et consid\'erons $\si \in \D_{\op{max}}$ tel que
$\t < \si$. Soient $f_{1},\dots,f_{d}$ un syst\`eme de g\'en\'erateurs du semi-groupe $(\si
\cap N)$ tels que $f_{1},\dots,f_{d-q}$ engendrent le semi-groupe $(\t \cap N)$
et tels que $f_{d-q+1}^{\ast}\wedge \dots \wedge f_{d}^{\ast} = \C{M}_{\t}$.
Comme $\D$ est r\'egulier, la famille
$f_{1}^{\ast},\dots,f_{d}^{\ast}$ forme une $\M{Z}$-base de $M$.

Par lin\'earit\'e, il suffit de calculer $I_{\t}(m_{D,\t})$ pour $m_{D,\t} =
f_{i}^{\ast}$ avec $i \in \{1,\dots,d\}$. 

Pour $i \in \{d-q+1,\dots,d\}$, on a
$|\chi^{f_{i}^{\ast}}(x)| = 1$ pour tout $x\in C_{\t}$, et donc
$I_{\t}(f_{i}^{\ast}) = 0$.

On choisit maintenant $i \in \{1,\dots,d-q\}$. On ne perd rien en
g\'en\'eralit\'e en supposant $i=1$.

Dans la carte affine $\varphi : U_{\si}(\M{C}) \rightarrow \M{C}^{d}$ donn\'ee
par $\varphi(x) = (\chi^{f_{1}^{\ast}}(x),\dots, \chi^{f_{d}^{\ast}}(x))$, les
ensembles $C_{\t}$ et $C_{\t}^{\op{int}}$ sont d\'efinis par les conditions~:
\begin{align*}
C_{\t} &= \{x \in \M{C}^{d}: \quad |x_{1}| \leqslant 1, \dots, |x_{d-q}|
\leqslant 1, 
|x_{d-q+1}| = 1, \dots, |x_{d}| = 1\} \\
\intertext{et}
C_{\t}^{\op{int}} &= \{x \in \M{C}^{d}: \quad |x_{1}| < 1, \dots, |x_{d-q}|
< 1, |x_{d-q+1}| =
1, \dots, |x_{d}| = 1\}.
\end{align*}
On remarque que~:
\[
\partial C_{\t} = \partial C_{\t}^{\op{int}} = \{x \in C_{\t}: \quad \exists i
\in \{1,\dots,d-q\}: \quad |x_{i}| = 1\}.
\]
Pour tout $\epsi \in \R^{+\ast}$, on pose~:
\begin{align*}
C_{\t}^{\epsi} &= \{x \in C_{\t}^{\op{int}}: \quad \epsi < |x_{1}| < 1\}, \\
D_{\t}^{\epsi} &= \{x \in \partial C_{\t}: \quad \epsi < |x_{1}| \leqslant 1\}
\\
\intertext{et}
E_{\t}^{\epsi} &= \{x \in C_{\t}: \quad |x_{1}| = \epsi\}.
\end{align*}
On a~:
\begin{multline}
\label{eqq_1}
I_{\t}(f_{1}^{\ast}) = 
\int_{C_{\t}^{\epsi}}\log | \chi^{f_{1}^{\ast}}(x)|^{2}\; \Theta(\C{M}_{\t})
\wedge \dd \omega 
\\ + 
\int_{C_{\t}^{\op{int},+}\backslash C_{\t}^{\epsi}}
\log | \chi^{f_{1}^{\ast}}(x)|^{2}\; \Theta(\C{M}_{\t})
\wedge \dd \omega.
\end{multline}
Comme $\log | \chi^{f_{1}^{\ast}}(x)|^{2}\,\Theta(\C{M}_{\t}) 
= \log |x_{1}|^{2} \,\Theta(\C{M}_{\t})$ est localement
$L^{1}$ sur $C_{\t}$, on a~:
\begin{equation}
\label{eqq_2}
\lim_{\epsi \rightarrow 0} \int_{C_{\t}^{\op{int},+}\backslash C_{\t}^{\epsi}}
\log | \chi^{f_{1}^{\ast}}(x)|^{2}\; \Theta(\C{M}_{\t})
\wedge \dd \omega = 0.
\end{equation}
Puisque $|\chi^{f_{d-q+1}^{\ast}}(x)| = \dots = |\chi^{f_{d}^{\ast}}(x)| = 1$
pour tout $x \in C_{\t}$ et que $\C{M}_{\t} \in
\bigwedge^{\op{max}}_{\M{Z}}(\t^{\perp}\cap M)$, la forme $\Theta(\C{M}_{\t})$, et donc
l'application $x \mapsto \log |
\chi^{f_{1}^{\ast}}(x)|^{2}\,\Theta(\C{M}_{\t})$,
sont $C^{\infty}$
sur un voisinage de $C_{\t}^{\epsi}$. On peut donc appliquer
la formule de Green (\ref{Green}) et en d\'eduire l'\'egalit\'e~:
\begin{multline}
\int_{C_{\t}^{\epsi}}
\log | \chi^{f_{1}^{\ast}}(x)|^{2}\; \Theta(\C{M}_{\t})
\wedge \dd \omega = \int_{C_{\t}^{\epsi}}
\dd \left(
\log | \chi^{f_{1}^{\ast}}(x)|^{2}\; \Theta(\C{M}_{\t})
\right) \wedge \omega \\
+ 
\int_{\partial C_{\t}^{\epsi}}
\left(
\log | \chi^{f_{1}^{\ast}}(x)|^{2}\; \Theta(\C{M}_{\t}) \wedge d^{c}\omega 
- d^{c}(
\log | \chi^{f_{1}^{\ast}}(x)|^{2}\; \Theta(\C{M}_{\t})
) \wedge \omega 
\right).
\end{multline}
En remarquant que $d \Theta(\C{M}_{\t}) = d^{c}\Theta(\C{M}_{\t}) = 0$, on a sur
$C_{\t}^{\varepsilon}$~:
\begin{align*}
\dd ( \log |\chi^{f_{1}^{\ast}}(x)|^{2}\wedge\Theta(\C{M}_{\t})) &= 
\dd (\log|\chi^{f_{1}^{\ast}}(x)|^{2}) \Theta(\C{M}_{\t}) = 0 \\
\intertext{et}
d^{c}(\log |\chi^{f_{1}^{\ast}}(x)|^{2}\;\Theta(\C{M}_{\t})) &= 
d^{c}(\log|\chi^{f_{1}^{\ast}}(x)|^{2})\wedge \Theta(\C{M}_{\t}) \\
&= \Theta(f_{1}^{\ast}) \wedge \Theta(\C{M}_{\t}) 
= \Theta(f_{1}^{\ast}\wedge 
\C{M}_{\t}).
\end{align*}
On tire des d\'efinitions l'\'egalit\'e des courants~:
\begin{equation}
\label{eqq_3}
\int_{\partial C_{\t}^{\epsi}} = \int_{D_{\t}^{\epsi}}\; + \;
\int_{E_{\t}^{\epsi}}
\end{equation}
Comme $\log | \chi^{f_{1}^{\ast}}(x)|^{2}\,\Theta(\C{M}_{\t})$ est localement
$L^{1}$ sur $\partial C_{\t}$, il vient~:
\[
\lim_{\epsi \rightarrow 0}
\int_{D_{\t}^{\epsi}} \log |
\chi^{f_{1}^{\ast}}(x)|^{2}\;\Theta(\C{M}_{\t})\wedge d^{c}\omega = 
\int_{\partial C_{\t}^{\op{int},+}} \log |
\chi^{f_{1}^{\ast}}(x)|^{2}\;\Theta(\C{M}_{\t})\wedge d^{c}\omega.
\]
De plus, $\log | \chi^{f_{1}^{\ast}}(x)|^{2} = 2 \log \varepsilon$ sur
$E_{\t}^{\varepsilon}$, la forme $\Theta(\C{M}_{\t})$ est $C^{\infty}$ sur un
voisinage de $C_{\t}$, et pour toute forme $\eta \in A^{2d - q -1}(\P(\M{C}))$ on a
$\int_{E_{\t}^{\varepsilon}}\eta = O(\varepsilon)$ quand $\varepsilon$ tend
vers $0$; on en d\'eduit que~:
\[
\lim_{\epsi \rightarrow 0}
\int_{E_{\t}^{\epsi}} \log |
\chi^{f_{1}^{\ast}}(x)|^{2}\;\Theta(\C{M}_{\t})\wedge d^{c}\omega = 0.
\]
Enfin, puisque $\Theta(f_{1}^{\ast}) = d\op{Arg}(x_{1})/2\pi$ et la forme
$\Theta(\C{M}_{\t})$ est $C^{\infty}$ sur un voisinage de $C_{\t}$, on a les
limites~:
\begin{align*}
\lim_{\epsi \rightarrow 0}
\int_{D_{\t}^{\epsi}} \Theta(f_{1}^{\ast}\wedge \C{M}_{\t}) \wedge \omega &= 
\int_{\partial C_{\t}^{\op{int},+}} \Theta(f_{1}^{\ast}\wedge \C{M}_{\t})
\wedge \omega \\
\intertext{et}
\lim_{\epsi \rightarrow 0}
\int_{E_{\t}^{\epsi}} \Theta(f_{1}^{\ast}\wedge \C{M}_{\t}) \wedge \omega &= 
\int_{C_{\t}^{\op{int},+}}\delta_{H_{1}\cap C_{\t}^{\op{int}}} \wedge 
\Theta(\C{M}_{\t})\wedge \omega,
\end{align*}
o\`u $H_{1}$ d\'esigne le diviseur d'\'equation $\chi^{f_{1}^{\ast}}(x) = 0$.

En revenant au cas g\'en\'eral (i.e. $m_{D,\t}$ quelconque), on d\'eduit de
(\ref{eqq_1}), (\ref{eqq_2}), (\ref{eqq_3}) et des calculs de limites ci-dessus que~:
\begin{multline*}
I_{\t}(m_{D,\t}) = \int_{\partial C_{\t}^{\op{int},+}}
\log | \chi^{m_{D,\t}}(x)|^{2}\; \Theta(\C{M}_{\t})\wedge d^{c}\omega \\
- \int_{\partial C_{\t}^{\op{int},+}}\Theta(m_{D,\t}\wedge\C{M}_{\t})\wedge
\omega
- \int_{C_{\t}^{\op{int},+}}
\delta_{D\cap C_{\t}^{\op{int}}}\wedge \Theta(\C{M}_{\t})\wedge \omega.
\end{multline*}
(On a utilis\'e ici le fait que pour tout $i \in \{d-q+1,\dots,d\}$, on a
$f_{i}^{\ast}\wedge \C{M}_{\t} = 0$). En utilisant (\ref{eq_chern_2}) il vient~:
\begin{multline*}
\int_{\P(\M{C})}c_{1}(\ov{L})\wedge c_{1}(\ov{L}_{1}) \wedge \dots \wedge
c_{1}(\ov{L}_{q})\wedge \omega \\
= \sum_{\t \in \D(d-q)}a_{\t}(\ov{L}_{1},\dots,\ov{L}_{1})\int_{\partial
C_{\t}^{\op{int},+}}
\log | \chi^{m_{D,\t}}(x)|^{2}\;\Theta(\C{M}_{\t})\wedge d^{c}\omega \\
- \sum_{\t \in \D(d-q)}a_{\t}(\ov{L}_{1},\dots,\ov{L}_{1})\int_{\partial
C_{\t}^{\op{int},+}}
\Theta(m_{D,\t}\wedge \C{M}_{\t}) \wedge \omega.
\end{multline*}
Enfin, on a~: 
\begin{multline*}
\sum_{\t \in \D(d-q)}a_{\t}(\ov{L}_{1},\dots,\ov{L}_{1})\int_{\partial
C_{\t}^{\op{int},+}}
\log | \chi^{m_{D,\t}}(x)|^{2}\;\Theta(\C{M}_{\t})\wedge d^{c}\omega  \\
= 
\sum_{\si \in \D(d-q-1)}\int_{C_{\si}^{\op{int},+}}
\Theta\left(
\sum_{\substack{\t > \si \\ \t \in \D(d-q)}}\epsi_{\si}(\t)a_{\t}(\ov{L}_{1},
\dots, \ov{L}_{q})\;\C{M}_{\t}\right)
\wedge (\log | \chi^{m_{D,\si}}(x)|^{2}d^{c}\omega) = 0
\end{multline*}
d'apr\`es l'assertion (2). On a donc~:
\begin{multline*}
\int_{\P(\M{C})}c_{1}(\ov{L})\wedge c_{1}(\ov{L}_{1}) \wedge \dots \wedge
c_{1}(\ov{L}_{q})\wedge \omega \\
= \sum_{\si \in \D(d-q-1)}\int_{C_{\si}^{\op{int},+}}
\Theta 
\left(
- \sum_{\substack{\t > \si \\ \t \in \D(d-q)}}
\epsi_{\si}(\t)
a_{\t}(\ov{L}_{1},\dots,\ov{L}_{q}) 
\; 
m_{D,\t} \wedge \C{M}_{\t}
\right)
\wedge \omega.
\end{multline*}
Pour \'etablir $(1)$ au rang $q+1$ et $(3)$ au rang $q$,
il suffit maintenant de prouver que pour tout $\si \in \D(d-q-1)$ il existe un
entier 
$a_{\si}(\ov{L},\ov{L}_{1},\dots,\ov{L}_{q})$ tel que~:
\[
a_{\si}(\ov{L},\ov{L}_{1},\dots,\ov{L}_{q})\; \C{M}_{\si} = 
- 
\sum_{\substack{\t > \si \\ \t \in \D(d-q)}}
\epsi_{\si}(\t)
a_{\t}(\ov{L}_{1},\dots,\ov{L}_{q})\;
m_{D,\t} \wedge \C{M}_{\t}.
\]
Soit $\si \in \D(d-q-1)$ et $\t_{0} \in \D(d-q)$ tel que $\si < \t_{0}$. 
On d\'eduit de l'assertion (2) au rang $q$ que l'on a~:
\[
\sum_{\substack{\t > \si \\ \t \in \D(d-q)}}
\epsi_{\si}(\t)
a_{\t}(\ov{L}_{1},\dots,\ov{L}_{q})\;
m_{D,\t_{0}} \wedge \C{M}_{\t} = 0.
\]
On a donc~:
\begin{multline*}
\sum_{\substack{\t > \si \\ \t \in \D(d-q)}}
\epsi_{\si}(\t)
a_{\t}(\ov{L}_{1},\dots,\ov{L}_{q})\;
m_{D,\t} \wedge \C{M}_{\t} \\
= \sum_{\substack{\t > \si \\ \t \in \D(d-q)}}
\epsi_{\si}(\t)
a_{\t}(\ov{L}_{1},\dots,\ov{L}_{q})\; (m_{D,\t} - m_{D,\t_{0}})\wedge
\C{M}_{\t}.
\end{multline*}
Or, comme la fonction support $\Psi_{D}$ est continue en $\si$, on a $(m_{D,\t}
- m_{D,\t_{0}}) \in \si^{\perp}\cap M$ pour tout $\t \in \D(d-q)$ tel que $\t >
\si$, et donc $(m_{D,\t}- m_{D,\t_{0}})\wedge \C{M}_{\t} \in \M{Z}\C{M}_{\si}$.

Montrons enfin que (2) est vraie au rang $q+1$.

Soit $\si' \in \D(d-q-2)$. D'apr\`es le (3), on a~:
\begin{multline*}
\sum_{\substack{\si > \si' \\ \si \in \D(d-q-1)}}
\epsi_{\si'}(\si) a_{\si}(\ov{L},\ov{L}_{1},\dots,\ov{L}_{q})\;\C{M}_{\si} \\
= - 
\sum_{\substack{\si > \si' \\ \si \in \D(d-q-1)}}
\sum_{\substack{\t > \si \\ \t \in \D(d-q)}}
\epsi_{\si'}(\si)\epsi_{\si}(\t)a_{\t}(\ov{L},\ov{L}_{1},\dots,\ov{L}_{q})\;
m_{D,\t}\wedge \C{M}_{\t} \\
= - 
\sum_{\substack{\t > \si' \\ \t \in \D(d-q)}} 
a_{\t}(\ov{L},\ov{L}_{1},\dots,\ov{L}_{q}) 
\left(
\sum_{\substack{\si \in \D(d-q-1) \\ \t > \si > \si'}}
\epsi_{\si'}(\si)\epsi_{\si}(\t)\right)
\; m_{D,\t}\wedge \C{M}_{\t}.
\end{multline*}
Or pour tout $\t \in \D(d-q)$ tel que $\t > \si'$, on a~: 
\[
\sum_{\substack{\si \in \D(d-q-1) \\ \t > \si > \si'}}
\epsi_{\si'}(\si)\epsi_{\si}(\t) = 0, 
\]
du fait de la relation $\partial \circ \partial = 0$ dans le complexe
simplicial associ\'e \`a la triangulation de $\C{S}^{d-1}$ induite par
$\D\backslash\{0\}$. 

Le th\'eor\`eme est donc d\'emontr\'e.
\medskip

\begin{rem}
\label{rem_dani}
Pour tout couple d'entiers positifs $(p,q)$, on introduit le $\M{Z}$-module~:
\[
C^{q}(\D,p) = \bigoplus_{\t \in \D(d-q)}{\textstyle\bigwedge\nolimits^{p}}(\t^{\perp}).
\]
On d\'efinit \'egalement un cobord $d: C^{q}(\D,p) \rightarrow C^{q+1}(\D,p)$
de la fa\c con suivante~:
\begin{alignat*}{3}
d:\; &C^{q}(\D,p) & &\longrightarrow & &\quad C^{q+1}(\D,p) \\
&\bigoplus_{\t \in \D(d-q)}x_{\t}& &\longmapsto & &\bigoplus_{\si \in \D(d-q-1)}
\left(
\sum_{\substack{\t > \si \\ \t \in \D(d-q)}}\epsi_{\si}(\t)x_{\t}
\right).
\end{alignat*}
Danilov a d\'emontr\'e (cf. \cite{4}, 12.4.1) que le $q^{\text{\`eme}}$ groupe de
cohomologie du complexe $(C^{\ast}(\D,p)\otimes \M{C},d)$ est isomorphe \`a
$H^{p,q}(\P(\M{C}))$. En reprenant les notations du th\'eor\`eme (\ref{calcul_prod}), on peut
associer au courant $c_{1}(\ov{L}_{1})\dotsm c_{1}(\ov{L}_{q})$ un \'el\'ement 
$C(\ov{L}_{1},\dots,\ov{L}_{q})$ de $C^{q}(\D,q)$ d\'efini par~:
\[
C(\ov{L}_{1},\dots,\ov{L}_{q}) = \bigoplus_{\t \in
\D(d-q)}a_{\t}(\ov{L}_{1},\dots,\ov{L}_{q})\,\C{M}_{\t}.
\]
D'apr\`es l'alin\'ea (2) de (\ref{calcul_prod}), on a $d C(\ov{L}_{1},\dots,\ov{L}_{q})
= 0$; et donc $C(\ov{L}_{1},\dots,\ov{L}_{q})$ d\'efinit un \'el\'ement de
$H^{q,q}(\P(\M{C}))$ dont on peut montrer qu'il co\"\i ncide avec la classe
$c_{1}(L_{1})\dotsm c_{1}(L_{q})$ par l'isomorphisme \'evoqu\'e
pr\'ec\'edemment. On retrouve ainsi gr\^ace au (3) du
th\'eor\`eme (\ref{calcul_prod}) la structure multiplicative de
$H^{2\ast}(\P(\M{C}))$.
\end{rem}
\begin{rem}
Comme l'\'eventail $\D$ est r\'egulier, la famille d'entiers \linebreak[4]
$(a_{\t}(\ov{L}_{1}, \dots,
\ov{L}_{q}))_{\t \in \D(d-q)}$ est un {\it poids de Minkowski\/} de codimension
$q$ de $\D$ au sens de \cite{23}.
\end{rem}

On d\'eduit imm\'ediatement du th\'eor\`eme (\ref{calcul_prod}) le corollaire~:
\begin{cor}
\label{support_courant}
Soient $\ov{L}_{1},\dots,\ov{L}_{q}$ des fibr\'es en droites
au-dessus de $\P$ munis \`a l'infini de leur m\'etrique canonique; on a~:
\[
\op{Supp}(c_{1}(\ov{L}_{1})\dotsm c_{1}(\ov{L}_{q})) \subset \bigcup_{\t \in
\D(d-q)}C_{\t}.
\]
\end{cor}
On note $\C{S}_{N}^{+} = C_{\{0\}}^{\op{int},+}$ le tore compact $\C{S}_{N}$
muni de l'orientation canonique induite par le choix de $\C{M}_{\{0\}}$, et 
$d \mu^{+}$ la forme volume canonique $\Theta(\C{M}_{\{0\}})$. 
On peut alors \'enoncer un second corollaire~:
\begin{cor}
\label{produit_chern_max}
Soient $\ov{L}_{1},\dots,\ov{L}_{d}$ des fibr\'es en droites
au-dessus de $\P$ munis \`a l'infini de leur m\'etrique canonique; pour
toute fonction $f$ de classe $C^{\infty}$ sur $\P(\M{C})$, on a~:
\[
\int_{\P(\M{C})}f\, c_{1}(\ov{L}_{1})\wedge \dots \wedge c_{1}(\ov{L}_{d}) =
\op{deg}(c_{1}(L_{1}) \dotsm c_{1}(L_{d})) 
\int_{\C{S}_{N}^{+}}f \, d\mu^{+}.
\]
\end{cor}
\demo\ D'apr\`es le th\'eor\`eme (\ref{calcul_prod}), il existe une constante
\linebreak[4]
$a_{\{0\}}(\ov{L}_{1},\dots,\ov{L}_{d}) \in \M{Z}$ telle que~:
\[
\int_{\P(\M{C})}f\, c_{1}(\ov{L}_{1}) \wedge \dots \wedge c_{1}(\ov{L}_{d})
= a_{\{0\}}(\ov{L}_{1},\dots,\ov{L}_{d})
\int_{\C{S}_{N}^{+}}f \, d\mu^{+}.
\]
On prend $f = 1$ et le r\'esultat d\'ecoule alors directement de
(\ref{coho_courant}).
\medskip
\begin{rem}
En prenant $\ov{L} = \ov{L}_{1} = \dots = \ov{L}_{d}$ dans le corollaire
(\ref{produit_chern_max}), on constate que la m\'etrique canonique $\|.\|_{L,\infty}$ est une
solution au premier probl\`eme de Calabi pour la forme volume singuli\`ere
$\delta_{S_{N}^{+}}\wedge d\mu^{+}$ sur $\M{P}(\D)(\M{C})$.
\end{rem}
\begin{rem}
\label{algo_efficace}
Soient $K_{1},\dots,K_{d}$ des polytopes convexes de $M_{\R}$ \`a sommets dans
$M$ tels que $K = K_{1} + \dots + K_{d}$ soit d'int\'erieur non vide. Soient
$\D'$ un raffinement r\'egulier de $\D$ l'\'eventail associ\'e \`a $K$ comme
au th\'eor\`eme
(\ref{construction_inverse}) et \`a la remarque (\ref{construction_inverse3}), et 
$E_{1}',\dots,E_{d}'$ les
diviseurs horizontaux invariants sur $\M{P}(\D')$ associ\'es \`a
$K_{1},\dots,K_{d}$ respectivement comme au
(\ref{construction_inverse3}).
En remarquant que le calcul de $c_{1}(\ov{E}_{1,\infty}')\dotsm
c_{1}(\ov{E}_{d,\infty}')$ ne n\'ecessite pas de conna\^\i tre $\D'$ mais
seulement $\D$ et les fonctions supports $\psi_{K_{1}}, \dots, \psi_{K_{d}}$, on d\'eduit
du th\'eor\`eme (\ref{calcul_prod}) un algorithme efficace 
pour le calcul du volume mixte~:
\[
V(K_{1},\dots,K_{d}) = \frac{1}{d!}\op{deg}(
c_{1}(\ov{E}_{1,\infty}')\dotsm c_{1}(\ov{E}_{d,\infty}')).
\]
\end{rem}
\bigskip

\subsection{Diviseurs \'el\'ementaires}~

Dans ce paragraphe on \'etablit un raffinement de
(\ref{support_courant}) dans le
cas o\`u les fibr\'es en droites consid\'er\'es sont des faisceaux
associ\'es \`a des diviseurs invariants \'el\'ementaires. Cela nous permet de
construire de mani\`ere canonique une section du morphisme d'anneaux~:
$[\cdot]: \ov{A}_{\op{f}}^{\ast}(\P_{\R}) \rightarrow H^{2\ast}(\P(\M{C}),\R)$.
\begin{thm}
\label{trivialite1}
Soient $D_{1} = V(\t_{1}), \dots, D_{q} = V(\t_{q})$ pour $q \leqslant d$ des
diviseurs invariants \'el\'ementaires; on a~:
\[
\op{Supp}(c_{1}(\ov{\C{O}(D_{1})}_{\infty}) \dotsm
c_{1}(\ov{\C{O}(D_{q})}_{\infty})) \subset \bigcup_{\substack{\t \in \D(d-q) \\
\t < \si \in \D_{\op{max}} \\ \si > \t_{1}, \dots,\t_{q}}}C_{\t}.
\]
\end{thm}
\demo\ On suppose dans un premier temps que $q = 1$. Par d\'efinition, la
fonction support $\Psi_{D_{1}}$ de $D_{1}$ est nulle sur tout c\^one maximal
$\si \in \D_{\op{max}}$ ne contenant pas $\t_{1}$. L'\'enonc\'e pour $q=1$ est
alors une simple cons\'equence du th\'eor\`eme (\ref{calcul2_c1}).

On suppose maintenant que $q > 1$. Le support du courant
$c_{1}(\ov{\C{O}(D_{1})}_{\infty})\dotsm \linebreak[0] c_{1}(\ov{\C{O}(D_{q})}_{\infty})$ est
inclus dans celui de $c_{1}(\ov{\C{O}(D_{1})}_{\infty})$ (on peut voir cela en
approchant la m\'etrique canonique sur $\C{O}(D_{1})$ par une m\'etrique
$C^{\infty}$). On d\'eduit de ce qui pr\'ec\`ede et du th\'eor\`eme
(\ref{calcul_prod}) que~:
\[
\op{Supp}(c_{1}(\ov{\C{O}(D_{1})}_{\infty}) \dotsm
c_{1}(\ov{\C{O}(D_{q})}_{\infty})) \subset \bigcup_{\substack{\t \in \D(d-q) \\
\t < \si \in \D_{\op{max}} \\ \si > \t_{1}}}C_{\t}.
\]
Le courant $c_{1}(\ov{\C{O}(D_{1})}_{\infty})\dotsm c_{1}(\ov{\C{O}(D_{q})}_{\infty})$
\'etant ind\'ependant de l'ordre des diviseurs $D_{1},\dots,D_{q}$, il vient~:
\[
\op{Supp}(c_{1}(\ov{\C{O}(D_{1})}_{\infty}) \dotsm
c_{1}(\ov{\C{O}(D_{q})}_{\infty})) \subset \bigcup_{\substack{\t \in \D(d-q) \\
\t < \si \in \D_{\op{max}} \\ \si > \t_{1}, \dots,\t_{q}}}C_{\t}.
\]
\medskip

On a alors le corollaire suivant~:
\begin{cor}
\label{trivialite2}
Soient $D_{1} = V(\t_{1}), \dots, D_{q} = V(\t_{q})$ pour $q \leqslant d$, des
diviseurs invariants \'el\'ementaires tels que $D_{1}\dotsm D_{q} = 0$ dans
$CH^{\ast}(\P)$ (i.e. tels que le c\^one $\t = \R^{+}\t_{1} + \dots +
\R^{+}\t_{q}$ ne soit pas un \'el\'ement de $\D$). On a~:
\[
c_{1}(\ov{\C{O}(D_{1})}_{\infty})\dotsm c_{1}(\ov{\C{O}(D_{q})}_{\infty}) = 0.
\]
\end{cor}
\demo\ Comme $\t \notin \D$, on ne peut pas trouver $\si \in \D_{\op{max}}$ tel
que $\t_{1},\dots,\t_{q}$ soient des faces de $\si$ (sinon $\t$ serait une face
de $\si$, et donc un \'el\'ement de $\D$). On d\'eduit de (\ref{trivialite1})
que $\op{Supp}(c_{1}(\ov{\C{O}(D_{1})}_{\infty}) \dotsm
c_{1}(\ov{\C{O}(D_{q})}_{\infty})) = \emptyset$, et donc que $
c_{1}(\ov{\C{O}(D_{1})}_{\infty})\dotsm c_{1}(\ov{\C{O}(D_{q})}_{\infty}) = 0$.
\medskip

On d\'eduit de (\ref{trivialite2}) et des th\'eor\`emes (\ref{anneau_chow}) et
(\ref{calcul2_c1}) qu'il existe un morphisme d'anneaux~:
\[
\varsigma: H^{2\ast}(\P(\M{C}),\R) \longrightarrow \ov{A}_{\op{f}}^{\ast}(\P_{\R})
\]
tel que pour tout $q$-uplet $(D_{1},\dots,D_{q})$ de diviseurs $T$-invariants
sur $\P$, on ait~:
\[
\varsigma(c_{1}(\C{O}(D_{1})) \dotsm c_{1}(\C{O}(D_{q}))) = 
c_{1}(\ov{\C{O}(D_{1})}_{\infty})\dotsm c_{1}(\ov{\C{O}(D_{q})}_{\infty}).
\]
D'apr\`es (\ref{coho_courant}), $\varsigma$ est une section du morphisme 
classe $[\cdot]$.
\begin{rem}
Soit $\C{A}^{q}(\D)$ le groupe des poids de Minkowski de codimension $q$ pour
$\D$ tel qu'il est d\'efini dans \cite{23}. On note $\xi: \C{A}^{q}(\D)
\rightarrow \ov{A}_{\op{f}}^{q}(\P_{\R})$ le morphisme de groupes d\'efini par
l'identit\'e~:
\[
\xi(\oplus_{\t \in \D(d-q)}a_{\t}) = \sum_{\t \in \D(d-q)}a_{\t}
\int_{C_{\t}^{\op{int},+}}\Theta(\C{M}_{\t}) \wedge \cdot
\]
On note \'egalement $\xi$ le morphisme $\C{A}^{q}(\D) \otimes \R \rightarrow 
\ov{A}_{\op{f}}^{q}(\P_{\R})$ obtenu par extension des scalaires \`a partir de
$\xi$. On peut alors factoriser l'application $\varsigma$ de la fa\c con suivante~:
\begin{center}
\mbox{
\xymatrix{
\C{A}^{q}(\D)\otimes \R \ar[r]^{\xi} & \ov{A}_{\op{f}}^{q}(\P_{\R}) \\
H^{2\ast}(\P(\M{C}),\R) \ar[u]^{\tilde{\varsigma}} \ar[ru]_{\varsigma} & }
}
\end{center}
\medskip

\noindent
D'apr\`es la remarque (\ref{rem_dani}), le morphisme $\tilde{\varsigma}$ est un
isomorphisme; de plus il induit par restriction un isomorphisme d'anneaux
gradu\'es~:
\[
\tilde{\varsigma}: H^{2\ast}(\P(\M{C}),\M{Z}) \longrightarrow \C{A}^{\ast}(\D).
\]
L'alin\'ea (3) du th\'eor\`eme (\ref{calcul_prod}) redonne la structure multiplicative de
$\C{A}^{\ast}(\D)$ induite comme dans \cite{23} par la structure d'anneau
naturelle de l'anneau de Chow op\'eratoriel de $\P$. On retrouve ainsi , sous
une forme diff\'erente, certains r\'esultats de \cite{23}.
\end{rem}
\bigskip

\section{G\'eom\'etrie d'Arakelov des vari\'et\'es toriques}~

Dans toute cette partie, $\P$ d\'esigne une vari\'et\'e torique projective
lisse de dimension absolue $d+1$. 

On montre tout d'abord que les multihauteurs ``canoniques'' de $\P$
(c'est-\`a-dire celles relatives \`a des fibr\'es en droites munis de leur
m\'etrique canonique) sont nulles. On en d\'eduit un r\'esultat remarquable~:
La hauteur d'une hypersurface dans $\P$ relativement \`a un fibr\'e en droites
muni de sa m\'etrique canonique est essentiellement donn\'ee par la {\em mesure
de Mahler} du polyn\^ome qui la d\'efinit. 

On construit enfin de mani\`ere canonique une section du morphisme d'anneaux~: 
\[
\zeta: \CH^{\ast}(X) \longrightarrow CH^{\ast}(X).
\]
L'existence d'une section canonique pour $\zeta$
\'etend les \'enonc\'es d'annulations
de nombres arithm\'etiques obtenus dans un premier temps
et conduit \`a une description de la structure de l'anneau
$\CH^{\ast}(X)$.
\bigskip

\subsection{Annulation des multihauteurs}~

L'\'enonc\'e suivant est un cas particulier de (\cite{21}, th. 2.4). 
On peut \'egalement consulter (\cite{21}, conj. 2.5, 2.6 et th. 2.9)
pour des \'enonc\'es proches de celui-ci.
\begin{prop}
\label{annulation_hauteur}
Soient $\ov{L}_{1,\infty}, \dots, \ov{L}_{d+1,\infty}$ des fibr\'es en droites
sur $\P$ munis de leur m\'etrique canonique. On a~: 
\[
h_{\ov{L}_{1,\infty}, \dots, \ov{L}_{d+1,\infty}}(\P) = 0.
\]
En particulier, pour tout fibr\'e en droites sur $\P$ muni de
sa m\'etrique canonique $\ov{L}_{\infty}$, on a~: 
\[
h_{\ov{L}_{\infty}}(\P) = 0.
\]
\end{prop}
\demo\ Soit $p$ un entier strictement sup\'erieur \`a $1$ et consid\'erons
l'endomorphisme $[p]: \P \rightarrow \P$ d\'efini au (\ref{definition_endo}). D'apr\`es la
proposition (\ref{relation_zhang}), on a pour tout $1 \leqslant i \leqslant d+1$~:
\[
[p]^{\ast}(\ov{L}_{i,\infty}) \simeq (\ov{L}_{i,\infty})^{p}.
\]
On tire de cela que pour tout $1 \leqslant i \leqslant d+1$, 
\[
\hat{c}_{1}([p]^{\ast}(\ov{L}_{i,\infty})) =
\hat{c}_{1}((\ov{L}_{i,\infty})^{p}) = p\,\hat{c}_{1}(\ov{L}_{i,\infty}).
\]
On obtient donc~:
\begin{multline}
\label{eq_canonique1}
\widehat{\op{deg}}(\hat{c}_{1}([p]^{\ast}(\ov{L}_{1,\infty})) \dotsm
\hat{c}_{1}([p]^{\ast}(\ov{L}_{d+1,\infty}))) \\
= p^{d+1}\,\widehat{\op{deg}}(\hat{c}_{1}(\ov{L}_{1,\infty}) \dotsm 
\hat{c}_{1}(\ov{L}_{d+1,\infty})).
\end{multline}
D'autre part, on a d'apr\`es les alin\'eas (5) et (7) 
du th\'eor\`eme (\ref{gdthm})~:
\begin{equation}
\begin{split}
\label{eq_canonique2}
\widehat{\op{deg}}(\hat{c}_{1}([p]^{\ast}(\ov{L}_{1,\infty})) \dotsm 
\hat{c}_{1}([p]&^{\ast}(\ov{L}_{d+1,\infty}))) \\
&= 
h_{[p]^{\ast}(\ov{L}_{1,\infty}),\dots,[p]^{\ast}(\ov{L}_{d+1,\infty})}
(\P)
\\
&= 
h_{\ov{L}_{1,\infty},\dots,\ov{L}_{d+1,\infty}}([p]_{\ast}\P) \\
&= p^{d}
h_{\ov{L}_{1,\infty},\dots,\ov{L}_{d+1,\infty}}(\P)\\
&= p^{d}\,\widehat{\op{deg}}(\hat{c}_{1}(\ov{L}_{1,\infty}) \dotsm
\hat{c}_{1}(\ov{L}_{d+1,\infty})).
\end{split}
\end{equation}
Comme $p > 1$, la conjonction de (\ref{eq_canonique1}) 
et de (\ref{eq_canonique2}) implique que~:
\[
h_{\ov{L}_{1,\infty},\dots,\ov{L}_{d+1,\infty}}(\P) =
\widehat{\op{deg}}(\hat{c}_{1}(\ov{L}_{1,\infty}) \dotsm
\hat{c}_{1}(\ov{L}_{d+1,\infty})) =0.
\]
\bigskip

\subsection{Hauteurs canoniques des hypersurfaces de $\P$}~

\begin{prop}
\label{hauteur_hypersurfaces}
Soient $\ov{L}_{1,\infty},\dots, \ov{L}_{d,\infty}$ des fibr\'es en droites sur
$\P$ munis de leur m\'etrique canonique et soit $s$ une section rationnelle non
nulle d'un fibr\'e en droites $L$ sur $\P$. Soit alors $D$ un diviseur $T$-invariant 
sur $\P$ tel que $L \simeq \C{O}(D)$ et notons $s_{D}$ la fonction 
rationnelle sur $\P$ correspondant \`a $s$ par cet isomorphisme. On a~:
\[
h_{\ov{L}_{1,\infty},\dots, \ov{L}_{d,\infty}}(\op{div}s) = \op{deg}
(c_{1}(L_{1}) \dotsm c_{1}(L_{d})) \int_{\C{S}_{N}^{+}}\log |s_{D}| \, d\mu^{+}.
\]
\end{prop}
\demo\ Soit $\|.\|_{\infty}$ la m\'etrique canonique de $L$ et notons
$\ov{L}_{\infty} = (L,\|.\|_{\infty})$. Suivant l'alin\'ea (6) du th\'eor\`eme
(\ref{gdthm}) il vient~:
\begin{multline*}
h_{\ov{L}_{1,\infty},\dots, \ov{L}_{d,\infty}}(\op{div}s) \\
= h_{\ov{L}_{1,\infty},\dots, \ov{L}_{d,\infty},\ov{L}_{\infty}}(\P) + 
\int_{\P(\M{C})}\log \|s\|_{\infty}\,c_{1}(\ov{L}_{1,\infty})\dotsm
c_{1}(\ov{L}_{d,\infty});
\end{multline*}
d'autre part on sait que $h_{\ov{L}_{1,\infty},\dots, 
\ov{L}_{d,\infty},\ov{L}_{\infty}}(\P) = 0$ du fait de la proposition
(\ref{annulation_hauteur}).
On termine alors la d\'emonstration en utilisant le corollaire
(\ref{produit_chern_max}) et en
remarquant que $\|s\|_{\infty} = |s_{D}|$ sur $\C{S}_{N}$.
\medskip

\begin{defn}
Soit $s$ une fonction rationnelle non nulle sur $\P$. On
appelle {\it mesure de Mahler\/} de $s$ et on note $M(s)$ le nombre r\'eel~: 
\[
M(s) = \int_{\C{S}_{N}^{+}}\log|s|\,d\mu^{+}.
\]
\end{defn}

La proposition suivante est une cons\'equence 
imm\'ediate de la d\'emonstration de (\cite{3}, prop. 1.5.1).
\begin{prop}
\label{mahler_remarque}
Soient $n$ un entier positif et $U$ un ouvert de $\M{C}^{n}$. Si $F$ est une
fonction m\'eromorphe sur $U \times \P(\M{C})$ dont le diviseur $\op{div}F$ est
plat sur $U$ (relativement \`a la premi\`ere projection), alors la mesure de
Mahler $M(F(u,\cdot))$ d\'epend contin\^ument de $u$.
\end{prop}
\bigskip

\subsection{Un exemple.}~

Pla\c cons nous sur $\M{P}^{n}_{\M{Z}}$ vu comme vari\'et\'e
torique de fa\c con standard, et soit $P \in
\M{Z}[X_{0},\dots,X_{n}]$ un polyn\^ome homog\`ene de degr\'e $k \in
\M{N}^{\ast}$.
Le polyn\^ome $P$ d\'efinit une section globale (encore not\'ee $P$) de $\C{O}(k)$ et la
hauteur de l'hypersurface $\op{div}P$ est donn\'ee d'apr\`es la proposition
(\ref{hauteur_hypersurfaces}) par~: 
\[
h_{\ov{\C{O}(1)}_{\infty}}(\op{div}P) = M(P) = \frac{1}{(2\pi)^{n+1}}
\int_{0}^{2\pi}\dotsi \int_{0}^{2\pi}\log |
P(e^{i\theta_{0}},\dots,e^{i\theta_{n}})|\,d\theta_{0} \dotsm d\theta_{n}.
\]
On peut trouver dans certains cas une formule explicite pour $M(P)$.

Soit par exemple la forme lin\'eaire $P(X_{0},\dots,X_{n}) = a_{0}X_{0} + \dots
+ a_{n}X_{n}$ avec $(a_{0},\dots,a_{n}) \in \M{C}^{n+1}-\{0\}$ 
et notons $I(a_{0},\dots,a_{n}) = 
M(a_{0}X_{0} + \dots + a_{n}X_{n})$ sa mesure de Mahler.

On d\'eduit de la formule de Jensen l'\'egalit\'e~: 
\[
I(a_{0},a_{1}) = \log \Sup (|a_{0}|,|a_{1}|).
\]
Le calcul de $I(a_{0},a_{1},a_{2})$ est plus d\'elicat et fait l'objet de
l'\'enonc\'e suivant, obtenu en collaboration avec J. Cassaigne~: 
\begin{prop}~
\label{polylog}
\begin{enumerate}
\item{
Si $|a_{0}|$, $|a_{1}|$ et $|a_{2}|$ sont les longueurs des c\^ot\'es d'un
triangle du plan (i.e. v\'erifient les in\'egalit\'es $|a_{i}|
\leqslant |a_{j}| + |a_{k}|$ avec $\{i,j,k\} = \{0,1,2\}$) alors~:
\[
I(a_{0},a_{1},a_{2}) = \frac{\alpha_{0}}{\pi}\log |a_{0}| + 
\frac{\alpha_{1}}{\pi}\log |a_{1}| + 
\frac{\alpha_{2}}{\pi}\log |a_{2}| + 
\frac{1}{\pi}\C{D}\left(\left|\frac{a_{1}}{a_{0}}\right|
e^{i\alpha_{2}}\right),
\]
o\`u $\C{D}(\cdot)$ d\'esigne le dilogarithme de Bloch-Wigner d\'efini par~: 
\[
\C{D}(z) = \Im (\op{li_{2}}(z) + \log|z|\log(1-z)) \qquad \text{pour}\quad z
\in \M{C}\backslash \{0,1\}, 
\]
et o\`u $\alpha_{0},\alpha_{1}$ et $\alpha_{2}$ sont les mesures principales
non orient\'ees des angles aux sommets $A_{0}$, $A_{1}$ et $A_{2}$ d'un
triangle du plan $\C{T} = (A_{0},A_{1},A_{2})$ tel que $|a_{0}| = A_{1}A_{2}$,
$|a_{1}| = A_{0}A_{2}$ et $|a_{2}| = A_{0}A_{1}$.}
\item{Sinon, 
\[
I(a_{0},a_{1},a_{2}) = \log \Sup (|a_{0}|,|a_{1}|,|a_{2}|).
\]}
\end{enumerate}
\end{prop}
\demo\ D\'emontrons tout d'abord l'assertion (1).
Les expressions consid\'er\'ees \'etant sym\'etriques
du fait de l'\'equation fonctionnelle v\'erifi\'ee par $\C{D}(\cdot)$, 
on peut supposer que
$|a_{0}| \geqslant |a_{1}| \geqslant |a_{2}|$. On a, en toute g\'en\'eralit\'e,
les relations~:
\begin{align}
I(a_{0},a_{1},a_{2}) &= \frac{1}{(2\pi)^{3}}\int_{0}^{2\pi}
\int_{0}^{2\pi}\int_{0}^{2\pi}\log \left|a_{0}e^{i\theta_{0}} + 
a_{1}e^{i\theta_{1}} + a_{2}e^{i\theta_{2}}\right|\,
d\theta_{0}d\theta_{1}d\theta_{2} \notag \\
&= 
\frac{1}{(2\pi)^{2}}\int_{0}^{2\pi}\int_{0}^{2\pi}
\log\left| |a_{0}| + |a_{1}|e^{i\theta_{1}} + |a_{2}|e^{i\theta_{2}}\right| \,
d\theta_{1}d\theta_{2} \notag \\
&= 
\frac{1}{2\pi}\int_{0}^{2\pi}\log \Sup \left(\left| |a_{0}| +
|a_{1}|e^{i\theta_{1}}\right|,
|a_{2}|\right)\, d\theta_{1}, 
\label{geo_1}
\end{align}
la derni\`ere \'egalit\'e \'etant obtenue par application de la formule de
Jensen. On pose alors $\alpha = \pi - \alpha_{2}$ et on v\'erifie que
$|\,|a_{0}| +
|a_{1}|e^{i\alpha}| = |a_{2}|$. 
On d\'eduit de cela et de (\ref{geo_1}) la relation~: 
\[
I(a_{0},a_{1},a_{2}) = \frac{\alpha_{2}}{\pi}\log |a_{2}| 
+ 
\frac{1}{2\pi}\int_{-\alpha}^{\alpha}\log 
\left||a_{0}| + |a_{1}|e^{i\theta_{1}}\right|\,d\theta_{1}, 
\]
ce dont on tire~:
\begin{equation}
\begin{split}
\label{geo_2}
I(a_{0},a_{1},a_{2}) = \frac{\alpha_{2}}{\pi}\log |a_{2}| +
\frac{\alpha_{0}}{\pi}\log |a_{0}| + \frac{\alpha_{1}}{\pi}\log |a_{0}|
\qquad \qquad \\
~\qquad \qquad \qquad \qquad + 
\frac{1}{2\pi}\int_{-\alpha}^{\alpha}\log \left| 
1 + \left| \frac{a_{1}}{a_{0}}\right|e^{i \theta_{1}}\right| \, d\theta_{1}, 
\end{split}
\end{equation}
puisque $\alpha = \alpha_{0} + \alpha_{1}$.
\bigskip
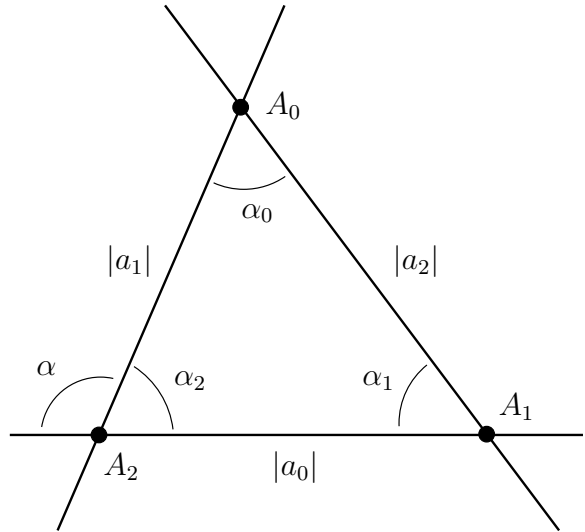
\begin{center}
\begin{picture}(0,0)%
\includegraphics{figure6.pstex}%
\end{picture}%
\setlength{\unitlength}{0.00083300in}%
\begingroup\makeatletter\ifx\SetFigFont\undefined
\def\x#1#2#3#4#5#6#7\relax{\def\x{#1#2#3#4#5#6}}%
\expandafter\x\fmtname xxxxxx\relax \def\y{splain}%
\ifx\x\y   
\gdef\SetFigFont#1#2#3{%
  \ifnum #1<17\tiny\else \ifnum #1<20\small\else
  \ifnum #1<24\normalsize\else \ifnum #1<29\large\else
  \ifnum #1<34\Large\else \ifnum #1<41\LARGE\else
     \huge\fi\fi\fi\fi\fi\fi
  \csname #3\endcsname}%
\else
\gdef\SetFigFont#1#2#3{\begingroup
  \count@#1\relax \ifnum 25<\count@\count@25\fi
  \def\x{\endgroup\@setsize\SetFigFont{#2pt}}%
  \expandafter\x
    \csname \romannumeral\the\count@ pt\expandafter\endcsname
    \csname @\romannumeral\the\count@ pt\endcsname
  \csname #3\endcsname}%
\fi
\fi\endgroup
\begin{picture}(3644,3344)(3579,-4583)
\put(6001,-2921){\makebox(0,0)[lb]{\smash{\SetFigFont{12}{14.4}{rm}$|a_{2}|$}}}
\put(4211,-2916){\makebox(0,0)[lb]{\smash{\SetFigFont{12}{14.4}{rm}$|a_{1}|$}}}
\put(4196,-4206){\makebox(0,0)[lb]{\smash{\SetFigFont{12}{14.4}{rm}$A_{2}$}}}
\put(6661,-3831){\makebox(0,0)[lb]{\smash{\SetFigFont{12}{14.4}{rm}$A_{1}$}}}
\put(4636,-3641){\makebox(0,0)[lb]{\smash{\SetFigFont{12}{14.4}{rm}$\alpha_{2}$}}}
\put(5206,-1951){\makebox(0,0)[lb]{\smash{\SetFigFont{12}{14.4}{rm}$A_{0}$}}}
\put(5056,-2596){\makebox(0,0)[lb]{\smash{\SetFigFont{12}{14.4}{rm}$\alpha_{0}$}}}
\put(5811,-3641){\makebox(0,0)[lb]{\smash{\SetFigFont{12}{14.4}{rm}$\alpha_{1}$}}}
\put(5246,-4206){\makebox(0,0)[lb]{\smash{\SetFigFont{12}{14.4}{rm}$|a_{0}|$}}}
\put(3776,-3581){\makebox(0,0)[lb]{\smash{\SetFigFont{12}{14.4}{rm}$\alpha$}}}
\end{picture}

\end{center}
\bigskip
Enfin, 
\begin{align*}
\frac{1}{2\pi}\int_{-\alpha}^{\alpha}\log \left|
1 + \left| \frac{a_{1}}{a_{0}}\right|e^{i \theta_{1}}\right| \, d\theta_{1}
&= 
\frac{1}{2\pi}\Re\left(\int_{-\alpha}^{\alpha}\log \left(1 + 
\left|\frac{a_{1}}{a_{0}}\right|e^{i\theta_{1}}\right)\,d\theta_{1}\right) \\
&= 
\frac{1}{2\pi}\Im\left(\int_{L}\frac{\log(1-z)}{z}\,dz\right), 
\end{align*}
o\`u $L$ est un chemin allant de $(- |a_{1}/a_{0}|e^{-i\alpha})$ \`a 
$(- |a_{1}/a_{0}|e^{i\alpha})$ d'indice z\'ero par rapport \`a $0$ et $1$. 

On en d\'eduit que~:
\begin{align*}
\frac{1}{2\pi}\int_{-\alpha}^{\alpha}\log 
\left|
1 + \left| \frac{a_{1}}{a_{0}}\right|e^{i \theta_{1}}\right| &\, d\theta_{1} \\
&= 
\frac{1}{2\pi}
\Im\left(\op{li_{2}}\left(- \left|\frac{a_{1}}{a_{0}}\right| e^{- i
\alpha}\right) - 
\op{li_{2}}\left(- \left|\frac{a_{1}}{a_{0}}\right| e^{i
\alpha}\right)
\right) \\
&=
\frac{1}{2\pi}\Im\left(\op{li_{2}}\left(\left|\frac{a_{1}}{a_{0}}\right|e^{i\alpha_{2}}
\right)\right)
- 
\frac{1}{2\pi}\Im\left(\op{li_{2}}\left(\left|\frac{a_{1}}{a_{0}}\right|e^{-i\alpha_{2}}
\right)\right)
\\
&= 
\frac{1}{\pi}
\Im\left(\op{li_{2}}\left(\left|\frac{a_{1}}{a_{0}}\right|e^{i\alpha_{2}}
\right)\right) \\
&=
\frac{1}{\pi}\C{D}\left(\left|\frac{a_{1}}{a_{0}}\right|e^{i\alpha_{2}}\right)
- \frac{1}{\pi}\log \left| \frac{a_{1}}{a_{0}}\right|
\Im\left(\log \left(1 -
\left|\frac{a_{1}}{a_{0}}\right|e^{i\alpha_{2}}\right)\right) \\
&= \frac{1}{\pi}\C{D}\left(\left|\frac{a_{1}}{a_{0}}\right|e^{i\alpha_{2}}\right)
+ \frac{\alpha_{1}}{\pi}\log \left|\frac{a_{1}}{a_{0}}\right|, 
\end{align*}
ce qui ajout\'e \`a (\ref{geo_2}) donne le r\'esultat annonc\'e.

Pla\c cons nous maintenant sous les hypoth\`eses de l'assertion (2). Les expressions
consid\'er\'ees \'etant sym\'etriques, on peut supposer que $|a_{0}| + 
|a_{1}| < |a_{2}|$. Le r\'esultat se d\'eduit alors directement de la relation
(\ref{geo_1}).
\medskip

\begin{rem}
La proposition (\ref{polylog}) g\'en\'eralise certains r\'esultats partiels dus
\`a Smyth (cf. \cite{26}, voir \'egalement \cite{25}). Pour un point de vue
moderne sur la mesure de Mahler, on peut consulter \cite{24}. On peut
s'inspirer de l'approche de (\cite{24}, \S 3) pour donner une interpr\'etation
de la proposition (\ref{polylog}) dans le langage de la K-th\'eorie et des
structures de Hodge-Tate mixtes.
\end{rem}
\begin{rem}
Si l'on consid\`ere le plan euclidien comme le bord du demi-espace de
Poincar\'e $\C{H}_{3}$ muni de la m\'etrique hyperbolique usuelle, le nombre
$\C{D}\left(\left|\frac{a_{1}}{a_{0}}\right|e^{i\alpha_{2}}\right)$ intervenant
dans la proposition (\ref{polylog}) est le volume du t\'etra\`edre
hyperbolique id\'eal dans $\C{H}_{3}$ de sommets $(A_{0},A_{1},A_{2},\infty)$.
\end{rem}
\bigskip

\subsection{L'anneau de Chow arithm\'etique g\'en\'eralis\'e d'une vari\'et\'e
torique.}~

L'objet de cette section est de prouver le r\'esultat suivant~:
\begin{thm}
\label{section_chow}
Soit $\P$ une vari\'et\'e torique projective lisse de dimension relative $d$
sur $\op{Spec}\M{Z}$. 
Il existe une unique section
$\hat{\varsigma}$ du morphisme $\zeta: \CH^{\ast}(\P)
\rightarrow CH^{\ast}(\P)$ telle que~:
\begin{itemize}
\item{$\hat{\varsigma}$ soit un morphisme d'anneaux;}
\item{
pour tout fibr\'e en droites $L$ sur $\P$,
on ait $\hat{\varsigma}(c_{1}(L)) = \hat{c}_{1}(\ov{L}_{\infty})$.}
\end{itemize}
\end{thm}
On montre tout d'abord le lemme suivant~: 

\begin{lem}
\label{lemme_section_chow}
Si $D_{1},\dots,D_{q}$ $(q \leqslant d)$ sont des diviseurs $T$-invariants
\'el\'ementaires tels que $D_{1}\dotsm D_{q} = 0$ dans $CH^{\ast}(\P)$, alors~:
\[
\hat{c}_{1}(\ov{\C{O}(D_{1})}_{\infty}) \dotsm \hat{c}_{1}(\ov{\C{O}(D_{q})}_{\infty})
= 0 \qquad \text{dans} \quad \CH^{q}(\P).
\]
\end{lem}
\demo\ On a~: 
\[
\zeta(\hat{c}_{1}(\ov{\C{O}(D_{1})}_{\infty}) \dotsm
\hat{c}_{1}(\ov{\C{O}(D_{q})}_{\infty}))
= c_{1}(\C{O}(D_{1})) \dots c_{1}(\C{O}(D_{q})) = 0. 
\]
De plus, 
\[
\omega(\hat{c}_{1}(\ov{\C{O}(D_{1})}_{\infty}) \dotsm
\hat{c}_{1}(\ov{\C{O}(D_{q})}_{\infty})) = 
c_{1}(\ov{\C{O}(D_{1})}_{\infty}) \dots c_{1}(\ov{\C{O}(D_{q})}_{\infty}) = 0
\]
d'apr\`es le corollaire (\ref{trivialite2}). 
Il existe donc
$\alpha \in D^{q-1,q-1}(\P_{\R})$ tel que~:
\begin{equation}
\label{equa1_sect6}
\hat{c}_{1}(\ov{\C{O}(D_{1})}_{\infty}) \dotsm 
\hat{c}_{1}(\ov{\C{O}(D_{q})}_{\infty})
= [(0,\alpha)] \qquad \text{dans} \quad \CH(\P),
\end{equation}
et comme $\dd \alpha = 0$, on peut choisir $\alpha$ dans
$Z^{q-1,q-1}(\P_{\R})$.

Soient maintenant $\ov{L}_{q+1,\infty},\dots, \ov{L}_{d+1,\infty}$ des fibr\'es
en droites sur $\P$ munis de leur m\'etrique canonique. On d\'eduit de
(\ref{equa1_sect6}) la relation~: 
\begin{multline*}
\widehat{\op{deg}}(\hat{c}_{1}(\ov{\C{O}(D_{1})}_{\infty}) \dotsm
\hat{c}_{1}(\ov{\C{O}(D_{q})}_{\infty}) 
\cdot \hat{c}_{1}(\ov{L}_{q+1,\infty})\dots \hat{c}_{1}(\ov{L}_{d+1,\infty}))
\\
= \frac{1}{2} \int_{\P(\M{C})}\alpha\, c_{1}(\ov{L}_{q+1,\infty}) \dotsm
c_{1}(\ov{L}_{d+1,\infty}),
\end{multline*}
ce qui entra\^\i ne, d'apr\`es la proposition (\ref{annulation_hauteur}) que~: 
\begin{equation}
\label{equa2_sect6}
\int_{\P(\M{C})}\alpha\, c_{1}(\ov{L}_{q+1,\infty}) \dotsm
c_{1}(\ov{L}_{d+1,\infty}) = 0,
\end{equation}
quels que soient les fibr\'es en droites $L_{q+1},\dots,L_{d+1}$ choisis. 

Comme l'anneau de cohomologie de De Rham de $\P$ est engendr\'e par les classes des
diviseurs dans $\P$ (cf. th\'eor\`eme \ref{anneau_chow}), on d\'eduit de (\ref{equa2_sect6})
et de (\ref{coho_courant}) que~: 
\[
\int_{\P(\M{C})}\alpha \wedge \beta = 0,
\]
quel que soit $\beta \in H^{d-q,d-q}(\P(\M{C})_{\R})$ et ceci implique par dualit\'e de
Poincar\'e que $\alpha = 0$ dans $\widetilde{D}^{q-1,q-1}(\P(\M{C})_{\R})$.

On passe maintenant \`a la d\'emonstration du th\'eor\`eme
(\ref{section_chow}).
En reprenant les notations du th\'eor\`eme (\ref{anneau_chow}), soit 
$\hat{\varsigma}: \C{S} \rightarrow \CH^{\ast}(\P)$ le morphisme
d'anneaux d\'efini par $\hat{\varsigma}(t_{\t}) =
\hat{c}_{1}(\ov{\C{O}(V(\t))}_{\infty})$ pour tout $\t \in \D(1)$. On d\'eduit
de la proposition (\ref{metrique_ind}) 
et du lemme (\ref{lemme_section_chow}) respectivement
les inclusions $\C{J} \subset \op{Ker}\hat{\varsigma}$
et $\C{I} \subset
\op{Ker}\hat{\varsigma}$. On conclut la d\'emonstration en utilisant le th\'eor\`eme
(\ref{anneau_chow}).
\bigskip

\section{Un th\'eor\`eme de Bernstein-Koushnirenko arithm\'etique}~

Ce chapitre est consacr\'e \`a \'etablir, comme application des r\'esultats des
pages qui pr\'ec\`edent, un analogue arithm\'etique du th\'eor\`eme de 
Bernstein-Koushnirenko. Rappelons que ce th\'eor\`eme fournit une borne du
nombre de z\'eros communs dans $(\M{C}^{\ast})^{n}$ \`a $n$ polyn\^omes de
Laurent $P_{1}, \dots, P_{n}$ dans $\M{C}[X_{1},X_{1}^{-1}, \dots, X_{n},
X_{n}^{-1}]$ en terme de volumes mixtes associ\'es aux polyh\`edres de Newton
de $P_{1},\dots,P_{n}$.

Lorsque $P_{1},\dots,P_{n}$ appartiennent \`a $\M{Z}[X_{1},X_{1}^{-1}, \dots, X_{n},
X_{n}^{-1}]$ nous d\'emontrons une majoration sur la hauteur de leurs z\'eros
communs (cf. corollaire (\ref{coro_BK})). Cette majoration fait intervenir un certain invariant r\'eel
$L(\nabla)$, associ\'e \`a un polytope convexe $\nabla$ dans $M$, que nous
d\'efinissons et \'etudions dans la section (\ref{invariant_polytope}).
\medskip

\subsection{A propos d'une constante associ\'ee \`a un polytope convexe}~
\label{invariant_polytope}

\subsubsection{D\'efinitions et propri\'et\'es}

Soit $\nabla$ un polytope convexe dans $M_{\R}$ \`a sommet dans $M$ 
que l'on suppose d'int\'erieur non vide (si
tel n'est pas le cas, on s'y ram\`ene en se pla\c cant dans $M' = M \cap
(\R\nabla+ \R(-\nabla))$.

On note $\D$ l'\'eventail dans $N$ associ\'e \`a $\nabla$ par le th\'eor\`eme
(\ref{construction_inverse}) et $\M{P}(\nabla)$ la vari\'et\'e torique $\P$.
On note \'egalement $E$ l'unique diviseur de Cartier horizontal $T$-invariant
sur $\M{P}(\nabla)$ tel que $K_{E} = \nabla$. D'apr\`es (\ref{construction_inverse}), 
on sait que le faisceau $\C{O}(E)$ est ample.

On associe alors \`a $\nabla$ la constante r\'eelle $L(\nabla)$ d\'efinie de la
mani\`ere suivante~:
\[
L(\nabla) = \Sup_{s \in \Gamma(\PPN, \C{O}(E))}\left(\Sup_{x \in
\PPN}\log \|s(x)\|_{E,\infty} - M(s) \right) .
\]
\begin{prop}
La constante $L(\Na)$ est bien d\'efinie et est positive.
\end{prop}
\demo\ 
La diff\'erence~:
\[
\Sup_{x \in \PPN}\log \|s(x)\|_{E,\infty} - \int_{S_{N}^{+}}\log |s(x)|\,
d\mu^{+}
\]
ne change pas lorsque $s$ est multipli\'ee par une constante $\lambda \in
\M{C}^{\ast}$. On a donc~:
\[
L(\nabla) = \Sup_{s \in \M{P}(\Gamma(\PPN, \C{O}(E)))}\left(\Sup_{x \in
\PPN}\log \|s(x)\|_{E,\infty} - M(s) \right) .
\]
Comme $\dim_{\M{C}}\Gamma(\PPN,\C{O}(E))$ est finie, 
$\M{P}(\Gamma(\PPN,\C{O}(E)))$ est compact et l'on d\'eduit de la proposition
(\ref{mahler_remarque}) que $|L(\Na)| < + \infty$.
Enfin pour toute section $s \in \Gamma(\PPN,\C{O}(E))$, on a~: 
\[
\int_{S_{N}^{+}}\log |s(x)| \, d\mu^{+} \leqslant \Sup_{x \in
\PPN}\log \|s(x)\|_{E,\infty} \int_{S_{N}^{+}}d\mu^{+} = \Sup_{x \in
\PPN}\log \|s(x)\|_{E,\infty}, 
\]
et donc $L(\Na) \geqslant 0$.
\medskip

\begin{prop}
\label{universalite}
Soient $\D'$ un \'eventail complet dans $N_{\R}$ et 
$E_{1}$ un diviseur de Cartier $T$-invariant sur $\M{P}(\D')$ tel que 
$\C{O}(E_{1})$ soit engendr\'e par ses sections globales. 
Si l'on note $\Na_{1} = K_{E_{1}} \subset M_{\R}$ le polytope 
convexe \`a sommets dans $M$ associ\'e \`a $E_{1}$, alors on a~: 
\[
L(\Na_{1}) = \Sup_{s \in \Gamma(\M{P}(\D')(\M{C}),\C{O}(E_{1}))}
\left(\Sup_{x \in \M{P}(\D')(\M{C})}\log \|s(x)\|_{E_{1},\infty} - 
\int_{S_{N}^{+}}\log |s(x)|\,d\mu^{+}\right) .
\]
\end{prop}
\demo\
On pose $V' = \R\Na_{1} + \R(-\Na_{1})$ et on note $M' = M \cap V'$. 
Soit $M''$ un sous-groupe de $M$ tel que l'on ait~: $M = M'\oplus M''$; on
tire~: 
\[
\TT =
\op{Spec}(\M{C}[M']) \times
\op{Spec}(\M{C}[M'']) = 
\M{G}_{\op{m}}^{\op{rg}M'}(\M{C}) \times 
\M{G}_{\op{m}}^{\op{rg}M''}(\M{C}).
\]
Notons $pr_{1}$ la premi\`ere projection et soit $m \in M'$. Par construction,
$m$ induit un caract\`ere $\chi_{\nabla_{1}}^{m}$ sur $\M{G}_{\op{m}}^{\op{rg}M'}(\M{C})
\subset \M{P}(\nabla_{1})(\M{C})$. Du fait de l'inclusion $M' \subset M$, $m$
induit \'egalement un caract\`ere $\chi_{\D'}^{m}$ sur $\TT \subset
\M{P}(\D')(\M{C})$ et les deux caract\`eres sont li\'es par la relation~: 
\[
\chi_{\D'}^{m} = \chi_{\nabla_{1}}^{m}\circ pr_{1} =
pr_{1}^{\ast}(\chi_{\nabla_{1}}^{m}).
\]
On note $E_{1}'$ le diviseur de Cartier $T$-invariant sur $\M{P}(\Na_{1})$, 
dont l'existence est assur\'ee par le th\'eor\`eme
(\ref{construction_inverse}), tel que $K_{E_{1}'} = \Na_{1}$ et $E_{1}'$ est
ample.
D'apr\`es la proposition (\ref{sections_globales}), on dispose d'un isomorphisme
canonique~:
\[
pr_{1}^{\ast}: \Gamma(\M{P}(\nabla_{1})(\M{C}),\C{O}(E_{1}')) \longrightarrow
\Gamma(\M{P}(\D')(\M{C}),\C{O}(E_{1})).
\]
Pour tout $s \in \Gamma(\M{P}(\Na_{1})(\M{C}),\C{O}(E_{1}'))$ et $x \in
\TT$, on a~:
\[
\|pr_{1}^{\ast}(s)(x)\|_{E_{1},\infty} = \|s(pr_{1}(x))\|_{E_{1}',\infty}, 
\]
et comme $\M{G}_{\op{m}}^{\op{rg}M'}(\M{C})$ (resp. $\TT$) est dense dans
$\M{P}(\Na_{1})(\M{C})$ (resp. dans $\M{P}(\D')(\M{C})$), on a~:
\[
\Sup_{x \in \M{P}(\D')(\M{C})}\|pr_{1}^{\ast}(s)(x)\|_{E_{1},\infty} = 
\Sup_{x' \in \M{P}(\Na_{1})(\M{C})}\|s(x')\|_{E_{1}',\infty}.
\]
De plus, pour tout $s \in \Gamma(\M{P}(\Na_{1})(\M{C}),\C{O}(E_{1}'))$, on a~:
\[
\int_{\C{S}_{N}^{+}} \log |pr_{1}^{\ast}(s)(x)|\,d\mu^{+} =
\int_{\C{S}_{N'}^{+}}\log |s(x')|\,d\mu^{+}.
\]
On en d\'eduit que~:
\begin{align*}
L(\Na_{1}) &= \Sup_{s \in \Gamma(\M{P}(\Na_{1})(\M{C}),\C{O}(E_{1}'))}
\left(\Sup_{x' \in \M{P}(\Na_{1})(\M{C})}\log \|s(x')\|_{E_{1}',\infty} - 
\int_{S_{N'}^{+}}\log |s(x')|\,d\mu^{+}\right) \\
 &= \Sup_{s \in \Gamma(\M{P}(\D')(\M{C}),\C{O}(E_{1}))}
\left(\Sup_{x \in \M{P}(\D')(\M{C})}\log \|s(x)\|_{E_{1},\infty} - 
\int_{S_{N}^{+}}\log |s(x)|\,d\mu^{+}\right).
\end{align*}
\medskip

\begin{prop}{\rm (Fonctorialit\'e).}
Soient $\Na_{1}$ et $\Na_{2}$ deux polytopes convexes dans $M$ et soit $\Na = \Na_{1} +
\Na_{2}$. Pour $i \in \{1,2\}$, on a~: 
\[
L(\Na_{i}) \leqslant L(\Na).
\]
\end{prop}
\demo\ Il suffit de d\'emontrer le r\'esultat pour $i =1$. 

On suppose que $\Na$ est d'int\'erieur non vide (si tel n'est pas le cas, 
on s'y ram\`ene en se pla\c cant dans $M' = M \cap (\R\Na + 
\R(-\Na))$. D'apr\`es le th\'eor\`eme (\ref{construction_inverse2}), il existe des 
diviseurs de Cartier $T$-invariants $E_{1}$, $E_{2}$ sur $\PN$ tels que
l'on ait $K_{E_{i}} = \Na_{i}$ pour $i \in \{1,2\}$.
Les faisceaux inversibles $\C{O}(E_{1})$ et $\C{O}(E_{2})$ sont engendr\'es
par leurs sections globales et l'on a $E = E_{1} + E_{2}$.

Soit $s_{1} \in \Gamma(\PPN,\C{O}(E_{1}))$ et soit $x_{0} \in \PPN$ tel que~: 
\[
\|s_{1}(x_{0})\|_{E_{1},\infty} = \Sup_{x \in \PPN}\|s_{1}(x)\|_{E_{1},\infty}.
\]
On peut trouver $\si \in \D_{\op{max}}$ tel que $x_{0} \in C_{\si}$.

On d\'eduit de l'\'egalit\'e $E = E_{1} + E_{2}$ que $m_{\Na,\si} =
m_{\Na_{1},\si} + m_{\Na_{2},\si}$. D'apr\`es le th\'eor\`eme
(\ref{polytope_et_vt}) on a
$m_{\Na_{2},\si} \in \Na_{2}$, et donc $s_{1}\otimes \chi^{m_{\Na_{2},\si}} \in
\Gamma(\PPN,\C{O}(E))$. De plus gr\^ace \`a la proposition
(\ref{BT_construction}) on peut
affirmer que~: 
\[
\|s_{1}\otimes \chi^{m_{\Na_{2},\si}}(x_{0})\|_{E,\infty} =
\|s_{1}(x_{0})\|_{E_{1},\infty}
\]
et donc que~: 
\[
\Sup_{x \in \PPN}\|s_{1}\otimes \chi^{m_{\Na_{2},\si}}(x)\|_{E,\infty}
\geqslant
\Sup_{x \in \PPN}\|s_{1}(x)\|_{E_{1},\infty}.
\]
Comme $\int_{S_{N}^{+}}\log |\chi^{m_{\Na_{2},\si}}|\,d\mu^{+} = 0$, on
d\'eduit finalement de la proposition (\ref{universalite}) la majoration~:
\begin{align*}
L(\Na_{1}) &= \Sup_{s_{1} \in \Gamma(\PPN,\C{O}(E_{1}))}\left(
\Sup_{x \in \PPN}\log \|s_{1}(x)\|_{E_{1},\infty} - 
\int_{S_{N}^{+}}\log |s_{1}(x)|\,d\mu^{+}\right) \\
&\leqslant \Sup_{s \in \Gamma(\PPN,\C{O}(E))}
\left(\Sup_{x \in \PPN}\log \|s(x)\|_{E,\infty} - 
\int_{S_{N}^{+}}\log |s(x)|\, d\mu^{+}\right) \\
&= L(\Na).
\end{align*}
\medskip

\subsubsection{Majoration de $L(\Na)$.}
Dans ce paragraphe, $\Na$ d\'esigne un polytope convexe dans $M_{\R}$ \`a
sommets dans $M$ et suppos\'e {\em d'int\'erieur non vide} (si tel n'est pas
le cas, on s'y ram\`ene en se pla\c cant dans $M' = M\cap (\R \Na + 
\R (- \Na))$). On note $\Na_{0}$ l'ensemble des sommets de $\Na$ et $\D$
l'\'eventail complet associ\'e \`a $\Na$ par la construction du th\'eor\`eme
(\ref{construction_inverse}). On note $E$ le diviseur ample sur $\P$ associ\'e \`a $\Na$ par cette
construction; on sait que $K_{E} = \Na$.

D'apr\`es la remarque (\ref{construction_inverse3}), il existe un raffinement $\D'$ de $\D$ tel que
$\M{P}(\D')$ est projective et lisse. On note $i_{\ast}: \M{P}(\D') \rightarrow
\P$ le morphisme \'equivariant induit par l'inclusion $i: \D' \hookrightarrow
\D$ et on pose $E' = (i_{\ast})^{\ast}(E)$. On sait alors que $K_{E'} = \Na$ et
le faisceau inversible $\C{O}(E')$ est engendr\'e par ses sections globales.

Suivant Lelong (cf. \cite{18}) on d\'efinit pour $n$ un entier strictement 
positif les constantes suivantes~:
\begin{align*}
C_{n} &= \frac{1}{2}\sum_{i=1}^{n-1}\frac{1}{i}\qquad \text{et} \qquad 
C_{n} = 0 \quad \text{pour} \quad n = 1, \\
C_{n}' &= \sum_{i =1}^{2n-2}\frac{1}{i}+ \sum_{i = 2n-1}^{+\infty}
\frac{1}{i2^{i}}\; .
\end{align*}
La proposition suivante est une cons\'equence directe des \'enonc\'es donn\'es
dans \cite{18}~:
\begin{prop}
\label{majoration_lelong}
Soit $P \in \M{C}[X_{1},\dots,X_{n}]$ un polyn\^ome, on a~:
\[
\Sup_{\substack{|z_{i}| \leqslant 1 \\ 1 \leqslant i \leqslant n}}
\log |P(z_{1}, \dots, z_{n})| - M(P) \leqslant (C_{n} +
C_{n}')\op{deg}P.
\]
\end{prop}
\demo\ Voir \cite{18}, th. 2, prop. 4 et \'equation (14). Voir aussi \cite{19}.
\medskip

A tout $\si \in \D_{\op{max}}'$ on attache un nombre r\'eel $L(\si)$ 
que l'on d\'efinit de la fa\c con suivante~:
\[
L(\si) = 
\Sup_{s \in \Gamma(\M{P}(\D')(\M{C}),\C{O}(E'))}
\left(
\Sup_{x \in C_{\si}} \log \|s(x)\|_{E',\infty} - \int_{S_{N}^{+}}\log
|s(x)|\, d\mu^{+} \right).
\]
L'\'eventail $\D'$ \'etant complet, on d\'eduit des propositions
(\ref{recouvrement}) et
(\ref{universalite}) l'\'egalit\'e~:
\begin{equation}
\label{eq_bern_1}
L(\Na) = \Sup_{\si \in \D'_{\op{max}}}L(\si).
\end{equation}
Pour tout $\sigma \in
\D_{\op{max}}'$, on choisit $f_{1},\dots,f_{d}$ une famille g\'en\'eratrice de
$\sigma$ (qui est donc une $\M{Z}$-base de $N$) et on note $f_{1}^{\ast},
\dots, f_{d}^{\ast}$ la base duale de $M$.

Dans la carte affine $\varphi: U_{\si}(\M{C}) \rightarrow \M{C}^{d}$ donn\'ee
par $\varphi(x) = (\chi^{f_{1}^{\ast}}(x), \dots, \chi^{f_{d}^{\ast}}(x))$,
l'ensemble $C_{\si}$ est d\'efini par les conditions~: 
\[
C_{\si} = \{z \in \M{C}^{d}: \quad |z_{1}| \leqslant 1, \dots, |z_{d}|
\leqslant 1\}.
\]
Soit $s \in \Gamma(\M{P}(\D')(\M{C}), \C{O}(E'))$. Dans la carte affine
$U_{\si}(\M{C})$, la fonction rationnelle $Q = s \cdot \chi^{-
m_{E,\si}}$ est un polyn\^ome. Comme 
\[
\int_{S_{N}^{+}}\log |Q|\, d \mu^{+} = 
\int_{S_{N}^{+}}\log |s \cdot \chi^{- m_{E,\si}}|\, d\mu^{+} = 
\int_{S_{N}^{+}}\log |s|\, d \mu^{+},
\]
et que pour tout $x \in C_{\si}$ on a~:
\[
\log \|s(x)\|_{E,\infty} = \log|Q(x)|, 
\]
on obtient d'apr\`es la proposition (\ref{majoration_lelong}) la majoration~:
\begin{align}
\label{eq_bern_2}
\Sup_{x \in C_{\si}}\log \|s(x)\|_{E,\infty} - \int_{S_{N}^{+}}\log |s|\,
d\mu^{+} &= \Sup_{\substack{|z_{i}| \leqslant 1 \\ 1 \leqslant i \leqslant d}}
\log |Q(z)| - M(Q) \notag \\
&\leqslant (C_{d} + C_{d}')\op{deg}(Q).
\end{align}
Au vu de la description des sections globales de $\C{O}(E')$ que l'on d\'eduit
de la proposition (\ref{sections_globales}), il vient~:
\[
\op{deg}Q \leqslant \Sup_{m \in (\Na - m_{E',\si})}\|m\|_{1} = 
\Sup_{m \in (\Na_{0} - m_{E',\si})}\|m\|_{1},
\]
o\`u l'on a pos\'e $\|m\|_{1} = \sum_{i = 1}^{d}|\hbox{$<f_{i},m>$}|$. 
On d\'eduit de cela et de (\ref{eq_bern_2}) la majoration~:
\begin{equation}
\label{eq_bern_3}
L(\si) \leqslant (C_{d} + C_{d}')\Sup_{m \in (\Na_{0} -
m_{E',\si})}\|m\|_{1}.
\end{equation}

Les relations (\ref{eq_bern_1}) et (\ref{eq_bern_3}) fournissent un proc\'ed\'e
th\'eorique pour l'obtention d'une majoration de l'invariant $L(\Na)$.
Malheureusement, la d\'etermination de l'\'eventail $\D'$, et {\it a
fortiori\/} de la base $f_{1},\dots,f_{d}$, repose non seulement sur la
connaissance de la g\'eom\'etrie du polytope $\Na$ mais \'egalement sur les propri\'et\'es
arithm\'etiques de l'\'eventail $\D$. Il n'est donc pas possible, en
g\'en\'eral, de trouver une majoration de $\op{deg}Q$ en fonction de la
g\'eom\'etrie du polytope $\Na$ uniquement.
\smallskip

On suppose d\'esormais que $\Na$ est un polytope absolument simple, ce qui nous
autorise \`a poser $\D' = \D$. On donne, dans ce cas, une majoration totalement
explicite de l'invariant $L(\Na)$ en terme de la combinatoire du polytope
$\Na$.

Pour tout $S \in \Na_{0}$ on
note $l_{1}(S), \dots , l_{d}(S)$ la base de $M$ associ\'ee au sommet
$S$ comme \`a la proposition (\ref{simplicite}).

Pour tout $S \in \Na_{0}$ et tout $P \in (\Na \cap M)$ on peut trouver
$a_{1}^{(S)}(P), \dots, a_{d}^{(S)}(P)$ des entiers positifs tels que~: 
\[
P-S = \sum_{i =1}^{d}a_{i}^{(S)}(P)l_{i}(S).
\]
On pose alors~: 
\[
N_{S}(P) = \sum_{i = 1}^{d} a_{i}^{(S)}(P) \in \M{N}.
\]
\begin{defn}
On appelle {\it norme\/} du polytope convexe absolument simple $\Na$ et on note
$N(\Na)$ l'entier strictement positif d\'efini par~: 
\[
N(\Na) = \Sup_{S \in \Na_{0}} \Sup_{S' \in \Na_{0}\backslash \{S\} } N_{S}(S').
\]
\end{defn}
\begin{prop}
\label{majoration_cas_AS}
Soit $\Na$ un polytope convexe absolument simple dans $M$. On a~: 
\[
L(\Na) \leqslant (C_{d} + C_{d}')N(\Na).
\]
\end{prop}
\demo\ On d\'eduit des d\'efinitions l'in\'egalit\'e~:
\[
\Sup_{m \in (\Na_{0} - m_{E,\si})}\|m\|_{1} = 
\Sup_{S \in \Na_{0}\backslash \{m_{E,\si}\}}N_{m_{E,\si}}(S) 
\leqslant N(\Na), 
\]
valable pour tout $\si \in \D_{\op{max}}$, ce qui joint aux relations
(\ref{eq_bern_1}) et (\ref{eq_bern_3}) donne le r\'esultat annonc\'e.
\bigskip

\subsection{Un th\'eor\`eme de Bernstein-Koushnirenko arithm\'etique.}~

\subsubsection{Rappels.}
Soient $X \in Z_{p}(\PN)$ et $Y \in Z_{q}(\PN)$ deux cycles effectifs tels que
$p + q \geqslant d +1$. On note $W_{1}, \dots, W_{r}$ les composantes
irr\'eductibles de l'intersection (ensembliste) $|X| \cap |Y|$, ce qui nous
permet d'\'ecrire~:
\[
|X| \cap |Y| = \bigcup_{1 \leqslant i \leqslant r} W_{i}\; .
\]
Pour tout $1 \leqslant i \leqslant r$, on a~: 
\begin{equation}
\label{int_propre}
\op{dim}W_{i} \geqslant p + q - d - 1.
\end{equation}
On dit que $W_{i}$ est {\it propre\/} si (\ref{int_propre}) est une
\'egalit\'e, et que $W_{i}$ est {\it impropre\/} sinon. On d\'efinit alors la
{\it partie propre\/} $(X\cdot Y)_{\op{pr}}$ de l'intersection 
de $X$ avec $Y$ par la formule~:
\[
(X\cdot Y)_{\op{pr}} = \sum_{\text{$W_{i}$ propre}}m_{i}W_{i} \in Z_{p+q -d
-1}(\PN),
\]
o\`u $m_{i}$ d\'esigne la multiplicit\'e d'intersection de $X$ avec $Y$ le long de
$W_{i}$ donn\'ee par la formule des ``Tor'' de Serre. En d'autres termes, $(X
\cdot Y)_{\op{pr}}$ est la somme des composantes propres de l'intersection de
$X$ avec $Y$ compt\'ees avec leur multiplicit\'e.

Si $X$, $Y$ et $Z$ sont trois cycles effectifs de $\PN$, alors on a~: 
\[
((X\cdot Y)_{\op{pr}}\cdot Z)_{\op{pr}} = (X \cdot
(Y\cdot Z)_{\op{pr}})_{\op{pr}}, 
\]
et dans toute la suite on \'ecrira $(X\cdot Y \cdot Z)_{\op{pr}}$ pour
d\'esigner l'un ou l'autre de ces cycles.

Pour plus de d\'etails sur ces questions, on peut consulter (\cite{3}, \S 5.5)
d'o\`u est extraite notre pr\'esentation. 

Ceci \'etant rappel\'e, nous \'enon\c
cons les r\'esultats que nous avons en vue~:

\subsubsection{Pr\'esentation des r\'esultats.}

On se place sur $\M{P}(\D')$ une vari\'et\'e torique
projective lisse de dimension relative $d$, associ\'ee
\`a un \'eventail $\D' \subset N_{\R}$.

On consid\`ere $E_{1}, \dots, E_{d}$ des diviseurs de Cartier horizontaux
$T$-invariants sur $\M{P}(\D')$ tels que les faisceaux inversibles
$\C{O}(E_{1}), \dots, \C{O}(E_{d})$ soient engendr\'es par leurs sections
globales. On pose $E = E_{1} + \dots + E_{d}$ et on note $\Na = K_{E}$, 
$\Na_{1} = K_{E_{1}}, \dots, \Na_{d} = K_{E_{d}}$ les polytopes convexes \`a
sommets dans $M$ associ\'es \`a $E$, $E_{1}, \dots, E_{d}$ respectivement; 
on sait d'apr\`es (\ref{additivite_1}) que $\Na = \Na_{1} + \dots + \Na_{d}$.
On supposera ici que le polytope $\Na$ est {\em d'int\'erieur non vide\/}.

\begin{thm}
\label{B_K}
Soient $s_{1},\dots,s_{d}$ des sections r\'eguli\`eres non nulles sur 
$\M{P}(\D')$ des faisceaux
$\C{O}(E_{1}), \dots, \C{O}(E_{d})$; et notons $\op{div}s_{1},
\dots, \op{div}s_{d}$ respectivement le lieu de leurs z\'eros. 
On a l'in\'egalit\'e~:
\begin{multline*}
h_{\ov{\C{O}(E)}_{\infty}}((\op{div}s_{1} \dotsm  \op{div}s_{d})_{\op{pr}}) \\
\leqslant \sum_{i=1}^{d}\frac{\op{deg}(E\cdot E_{1} \dotsm 
\widehat{E_{i}}\dotsm E_{d})}{\op{deg}(E^{d})}\, h_{\ov{\C{O}(E)}_{\infty}}
(\op{div}s_{i}) \\
+ \sum_{i =1}^{d}L(\Na_{i})\op{deg}(E\cdot E_{1} \dotsm \widehat{E_{i}} 
\dotsm E_{d}), 
\end{multline*}
o\`u le symbole $\widehat{E_{i}}$ signifie que l'on omet ce terme dans le
produit d'intersection consid\'er\'e.
\end{thm}
En utilisant le fait que~:
\[
\op{deg}(E\cdot E_{1}\dotsm \widehat{E_{i}} \dotsm E_{d}) = d!\, V(\Na,\Na_{1},
\dots, \widehat{\Na_{i}}, \dots, \Na_{d}), 
\]
o\`u $V(\Na,\Na_{1},\dots, \widehat{\Na_{i}}, \dots, \Na_{d})$ d\'esigne le
volume mixte des polytopes $\Na,\Na_{1}, \dots, \widehat{\Na_{i}},
\linebreak[4] \dots,
\Na_{d}$ (voir par exemple \cite{20}, \S A.4 et p. 78-79, et aussi \cite{11},
\S 5.4) et l'\'egalit\'e~:
\[
h_{\ov{\C{O}(E)}_{\infty}}(\op{div}s_{i}) = \op{deg}(E^{d})\, M(s_{i}), 
\]
valable pour tout $1 \leqslant i \leqslant d$ d'apr\`es la proposition
(\ref{hauteur_hypersurfaces}), on peut r\'e\'ecrire la 
majoration du th\'eor\`eme (\ref{B_K}) sous la
forme~:
\begin{multline}
\label{BK_intermediaire}
h_{\ov{\C{O}(E)}_{\infty}}((\op{div}s_{1} \dotsm  
\op{div}s_{d})_{\op{pr}}) \\
\leqslant d!\, \sum_{i
=1}^{d}V(\Na,\Na_{1},\dots,\widehat{\Na_{i}},\dots,\Na_{d})\, M(s_{i}) \\
 + 
d!\, \sum_{i=1}^{d}L(\Na_{i})V(\Na,\Na_{1},\dots,\widehat{\Na_{i}},\dots,
\Na_{d}).
\end{multline}
\medskip

\begin{defn}
Soit $P \in \M{Z}[X_{1},1/X_{1}, \dots, X_{d},1/X_{d}]$ un polyn\^ome de 
Laurent, et notons $(a_{m})_{m \in \M{Z}^{d}}$ la
famille presque nulle de ses coefficients (i.e. d\'efinie par l'\'egalit\'e $P(X)
= \sum_{m \in \M{Z}^{d}}a_{m}X^{m}$). On appelle {\it support\/} du polyn\^ome
$P$ et on note $\op{Supp}P$ le sous-ensemble fini de $\M{Z}^{d}$ d\'efini par~:
\[
\op{Supp}P = \{ m \in \M{Z}^{d}: \;\; a_{m} \neq 0\}.
\]
Le {\it polyh\`edre de Newton\/} du polyn\^ome $P$ est l'enveloppe
convexe de son support $\op{Supp}P$. C'est un polytope convexe \`a sommets dans
$\M{Z}^{d}$.
\end{defn}

Pour tout corps de nombres $K \subset \ov{\M{Q}}$, notons $S_{K}$ l'ensemble
canonique des places de $K$ (voir par exemple \cite{31}, II \S 1). 
Pour tout $\nu \in S_{K}$, on note
$|\cdot|_{\nu}$ la valeur absolue {\em normalis\'ee} sur $K$ en la place $\nu$,
celle-ci \'etant
d\'efinie par~: $|\cdot|_{\nu} = |\sigma(\cdot)|$ si $\nu$ est associ\'ee \`a
un plongement r\'eel $\sigma: K \hookrightarrow \M{R}$, $|\cdot|_{\nu} =
|\sigma(\cdot)|^{2}$ si $\nu$ est associ\'ee \`a un plongement complexe
$\sigma: K \hookrightarrow \M{C}$, et $|\cdot|_{\nu} = (N\mathfrak{p})^{-
v_{\mathfrak{p}}(\cdot)}$
si $\nu$ est une place non-archim\'edienne associ\'ee \`a un id\'eal premier
$\mathfrak{p}$ de $\C{O}_{K}$.

On peut alors \'enoncer le corollaire suivant, 
de forme plus \'el\'ementaire~:
\begin{cor}
\label{coro_BK}
Soient $P_{1},\dots,P_{d} \in \M{Z}[X_{1},1/X_{1}, \dots, X_{d},1/X_{d}]$ des 
polyn\^omes de Laurent \`a coefficients entiers, et notons
$\Na_{1}, \dots ,\Na_{d}$ respectivement leur polyh\`edre de Newton.
Notons $\Na = \Na_{1} +
\dots + \Na_{d}$ leur somme de Minkowski, que nous supposons d'int\'erieur non vide. 
Notons $Z_{1}, \dots, Z_{d}$ le lieu
des z\'eros de $P_{1},\dots,P_{d}$ respectivement dans
$(\ov{\M{Q}}^{\ast})^{d}$ et $(Z_{1} \cap \dots \cap Z_{d})_{\op{pr}}$
l'ensemble des points isol\'es du sch\'ema $Z_{1}\cap \dots \cap Z_{d}$, et
pour chaque point $x$ de cet ensemble notons $l(x)$ sa multiplicit\'e et
$h_{\Na}(x)$ sa hauteur d\'efinie par l'\'egalit\'e~:
\[
h_{\Na}(x) = \sum_{\nu \in S_{K}}\log \left(\max_{m \in \Na\cap
M}|\chi^{m}(x)|_{\nu}\right),
\]
o\`u $K$ d\'esigne un corps de nombres contenant
les coordonn\'ees de $x$. On a~:
\[
\frac{1}{d!}\sum_{x \in (Z_{1}\cap \dots \cap Z_{d})_{\op{pr}}}l(x)h_{\Na}(x) 
\leqslant \sum_{i =1}^{d}V(\Na,\Na_{1}, \dots, \widehat{\Na_{i}}, \dots,
\Na_{d})
(M(P_{i}) + L(\Na_{i})).
\]
\end{cor}
\demo\ 
Soit $\D$ l'\'eventail dans $N$ que l'on associe \`a $\Na$ gr\^ace au
th\'eor\`eme (\ref{construction_inverse}) et soit $E$ l'unique diviseur de Cartier horizontal
$T$-invariant sur $\P$ tel que $K_{E} = \Na$ et $E$ est ample. D'apr\`es le
th\'eor\`eme (\ref{construction_inverse2}) il existe des diviseurs de Cartier horizontaux
$T$-invariants $E_{1},\dots,E_{d}$ sur $\M{P}(\D)$ tels que pour tout $1 \leqslant i \leqslant
d$, $K_{E_{i}} = \Na_{i}$ et le faisceau inversible $\C{O}(E_{i})$ est
engendr\'e par ses sections globales. De plus on a $E = E_{1} + \dots + E_{d}$.

D'apr\`es la remarque (\ref{construction_inverse3}), il existe un raffinement
$\D'$ de $\D$ tel que $\M{P}(\D')$ est {\em projective} et {\em lisse}. On note 
$i_{\ast}: \M{P}(\D') \rightarrow \M{P}(\D)$ le morphisme \'equivariant 
induit par l'inclusion $i: \D' \hookrightarrow \D$ et on pose $E' =
(i_{\ast})^{\ast}(E)$, $E_{1}' = (i_{\ast})^{\ast}(E_{1}), \dots, E_{d}' =
(i_{\ast})^{\ast}(E_{d})$. On sait alors que $K_{E'} = \Na$, $K_{E_{1}'} =
\Na_{1}, \dots, K_{E_{d}'} = \Na_{d}$ et que les faisceaux inversibles
$\C{O}(E')$, $\C{O}(E_{1}'), \dots, \C{O}(E_{d}')$ sont engendr\'es par leurs
sections globales.

Pour tout $1 \leqslant i \leqslant d$, on d\'eduit de l'inclusion
$\op{Supp}P_{i} \subset \Na_{i}$ que $P_{i}$ s'\'etend en une section
r\'eguli\`ere non nulle
$s'_{i}$ de $E'_{i}$ sur $\M{P}(\D')$.
On peut donc appliquer l'in\'egalit\'e
(\ref{BK_intermediaire}) sur $\M{P}(\D')$ et on trouve~:
\begin{multline*}
h_{\ov{\C{O}(E')}_{\infty}}((\op{div}s'_{1} \dotsm  
\op{div}s'_{d})_{\op{pr}}) \\
\leqslant d!\, \sum_{i
=1}^{d}V(\Na,\Na_{1},\dots,\widehat{\Na_{i}},\dots,\Na_{d})\, M(s'_{i}) \\
 + 
d!\, \sum_{i=1}^{d}L(\Na_{i})V(\Na,\Na_{1},\dots,\widehat{\Na_{i}},\dots,
\Na_{d}).
\end{multline*}
On se ram\`ene de $\M{P}(\D')(\ov{\M{Q}})$ \`a $T(\ov{\M{Q}}) = 
(\ov{\M{Q}}^{\ast})^{d}$
en utilisant la positivit\'e de $h_{\ov{\C{O}(E')}_{\infty}}$ (cf. exemple
\ref{positivite_torique}). Enfin l'\'egalit\'e $h_{\ov{\C{O}(E')}_{\infty}}(x) = 
h_{\Na}(x)$ pour
tout $x \in (\ov{\M{Q}}^{\ast})^{d}$ est une cons\'equence directe de la
construction par image inverse de $\|.\|_{E',\infty}$ (cf.
\S\ref{chap_image_inverse}).
\medskip

\subsubsection{D\'emonstration du th\'eor\`eme.}~

On d\'emontre tout d'abord le lemme suivant~:
\begin{lem}
\label{lemme_BK}
Soit $k$ un entier compris entre $1$ et $d$, et soit $Z \in Z_{k}(\M{P}(\D'))$ un
cycle effectif tel que la section $s_{k}$ de $\C{O}(E_{k})$ ne soit identiquement nulle
sur aucune des composantes irr\'eductibles de $Z$. On a~: 
\begin{multline*}
h_{\ov{\C{O}(E)}_{\infty}, \ov{\C{O}(E_{1})}_{\infty}, \dots, 
\ov{\C{O}(E_{k-1})}_{\infty}}(Z \cdot \op{div}s_{k}) \\
\leqslant 
h_{\ov{\C{O}(E)}_{\infty}, \ov{\C{O}(E_{1})}_{\infty}, \dots, 
\ov{\C{O}(E_{k})}_{\infty}}(Z) \qquad \qquad \qquad \qquad\\
+ \op{deg}(E\cdot E_{1}\dotsm E_{k-1}\cdot Z)\left(
L(\Na_{k}) + \frac{h_{\ov{\C{O}(E)}_{\infty}}(\op{div}s_{k})}{\op{deg}(E^{d})}
\right).
\end{multline*}
\end{lem}
\demo\ D'apr\`es le th\'eor\`eme (\ref{gdthm}) alin\'ea (6), on a~: 
\begin{multline*}
h_{\ov{\C{O}(E)}_{\infty}, \ov{\C{O}(E_{1})}_{\infty}, \dots, 
\ov{\C{O}(E_{k-1})}_{\infty}}(Z \cdot \op{div}s_{k}) 
= 
h_{\ov{\C{O}(E)}_{\infty}, \ov{\C{O}(E_{1})}_{\infty}, \dots, 
\ov{\C{O}(E_{k})}_{\infty}}(Z) \\
+
\int_{Z(\M{C})}\log \|s_{k}\|_{E_{k},\infty}
c_{1}(\ov{\C{O}(E)}_{\infty})
c_{1}(\ov{\C{O}(E_{1})}_{\infty})
\dotsm 
c_{1}(\ov{\C{O}(E_{k-1})}_{\infty}).
\end{multline*}
On d\'eduit de la proposition (\ref{universalite}) que~: 
\[
\Sup_{x \in \M{P}(\D')(\M{C})}\log \|s_{k}(x)\|_{E_{k},\infty} \leqslant M(s_{k}) +
L(\Na_{k}).
\]
De l'in\'egalit\'e pr\'ec\'edente et de la positivit\'e des courants 
$c_{1}(\ov{\C{O}(E)}_{\infty}),c_{1}(\ov{\C{O}(E_{1})}_{\infty}), 
\linebreak[4]
\dots, c_{1}(\ov{\C{O}(E_{k-1})}_{\infty})$ (cf. exemple 
\ref{exemple_adm1}) on tire que~:
\begin{multline*}
\int_{Z(\M{C})}\log \|s_{k}\|_{E_{k},\infty}\,
c_{1}(\ov{\C{O}(E)}_{\infty})
c_{1}(\ov{\C{O}(E_{1})}_{\infty}) \dotsm 
c_{1}(\ov{\C{O}(E_{k-1})}_{\infty}) \\
\leqslant
\left(
\int_{Z(\M{C})}
c_{1}(\ov{\C{O}(E)}_{\infty})
c_{1}(\ov{\C{O}(E_{1})}_{\infty}) \dotsm 
c_{1}(\ov{\C{O}(E_{k-1})}_{\infty})
\right)
(M(s_{k}) + L(\Na_{k})) \\
= \op{deg}(E\cdot E_{1} \dotsm E_{k-1})(M(s_{k}) + L(\Na_{k})).
\end{multline*}
Enfin, on a d'apr\`es la proposition (\ref{hauteur_hypersurfaces}) l'\'egalit\'e~:
\[
M(s_{k}) = \frac{h_{\ov{\C{O}(E)}_{\infty}}(\op{div}s_{k})}{\op{deg}(E^{d})}, 
\]
ce qui suffit \`a \'etablir le r\'esultat.
\medskip

On passe maintenant \`a la d\'emonstration du th\'eor\`eme~:

On construit par r\'ecurrence une suite finie $Z_{0}, Z_{1}, \dots, Z_{d}$ de
cycles effectifs dans $\M{P}(\D')$ tels que $Z_{i} \in Z^{i}(\M{P}(\D'))$ pour tout $0
\leqslant i \leqslant d$ de la fa\c con suivante~:
\begin{itemize}
\item{On pose $Z_{0} = \M{P}(\D')$.}
\item{Pour $i \geqslant 1$,
le cycle $Z_{i}$ est d\'efini comme la somme
avec multiplicit\'es
des composantes de $Z_{i-1}\cdot \op{div}s_{d+1-i}$ dont l'intersection avec
$\op{div}s_{d-i}$ est propre.}
\end{itemize}
Comme dans $\M{P}(\D')$ l'intersection d'un cycle de codimension $k$ par une
hypersurface est soit vide, soit un cycle de codimension $k+1$, toute
composante non vide de~:
\[
((|Z_{i-1}|\cap |\op{div}s_{d+1-i}|) - |Z_{i}|)\cap |\op{div}s_{d-i}|
\cap \dotsm \cap |\op{div}s_{1}|
\]
est au moins de dimension $1$. On est donc assur\'e de l'\'egalit\'e~:
\[
Z_{d} = (\op{div}s_{1} \dotsm \op{div}s_{d})_{\op{pr}}.
\]
Pour tout $1 \leqslant k \leqslant d$, le cycle $(Z_{d-k}\cdot 
\op{div}s_{k}) - Z_{d -k +1}$ est effectif. On d\'eduit de cela et de
l'exemple (\ref{positivite_torique}) que~:
\[
h_{\ov{\C{O}(E)}_{\infty}, \ov{\C{O}(E_{1})}_{\infty}, \dots, 
\ov{\C{O}(E_{k-1})}_{\infty}}(Z_{d-k+1}) 
\leqslant 
h_{\ov{\C{O}(E)}_{\infty}, \ov{\C{O}(E_{1})}_{\infty}, \dots, 
\ov{\C{O}(E_{k-1})}_{\infty}}(Z_{d-k}\cdot \op{div}s_{k}).
\]
Du fait de la positivit\'e des fibr\'es $\C{O}(E), \C{O}(E_{1}), \dots, 
\C{O}(E_{k-1})$ on a de m\^eme~: 
\[
\op{deg}(E\cdot E_{1} \dotsm E_{k-1}\cdot Z_{d-k}) \leqslant 
\op{deg}(E\cdot E_{1} \dotsm \widehat{E_{k}} \dotsm E_{d}).
\]
En appliquant le lemme (\ref{lemme_BK}), on obtient pour tout $1 \leqslant k \leqslant
d$ l'in\'egalit\'e~:
\begin{multline*}
h_{\ov{\C{O}(E)}_{\infty}, \ov{\C{O}(E_{1})}_{\infty}, \dots, 
\ov{\C{O}(E_{k-1})}_{\infty}}(Z_{d-k+1}) 
\\
\leqslant 
h_{\ov{\C{O}(E)}_{\infty}, \ov{\C{O}(E_{1})}_{\infty}, \dots, 
\ov{\C{O}(E_{k})}_{\infty}}(Z_{d-k})
\\
+ \op{deg}(E\cdot E_{1}\dotsm \widehat{E_{k}} \dotsm E_{d})
\left(
L(\Na_{k}) + 
\frac{h_{\ov{\C{O}(E)}_{\infty}}(\op{div}s_{k})}{\op{deg}(E^{d})}
\right).
\end{multline*}
Une r\'ecurrence finie permet alors de conclure.

\bigskip
\begin{flushleft}
{
D\'epartement de Math\'ematiques et d'Informatique,
\'Ecole Normale Sup\'erieure, \\
45 rue d'Ulm, 75230 Paris cedex 05, FRANCE. \\
e-mail: vmaillot@ens.fr
}
\end{flushleft}

\end{document}